\documentclass[11pt]{amsart}

\usepackage{amsmath,amsthm, amscd, amssymb, amsfonts, mathtools}
\usepackage[all]{xy}
\usepackage{enumitem}

\usepackage{t1enc}
\usepackage[latin1]{inputenc}

\usepackage{fontsmpl}
\usepackage{srcltx}

\renewcommand{\_}[1]{_{\left( #1 \right)}}

\allowdisplaybreaks

\newcommand{\supera}[2]{{\mathbf A}(#1|#2)}
\newcommand{\superb}[2]{{\mathbf B}(#1|#2)}

\newcommand{\superd}[2]{{\mathbf D}(#1|#2)}
\newcommand{\superda}[1]{{\mathbf D}(2,1;#1)}
\newcommand{\superf}{{\mathbf F}(4)}
\newcommand{\superg}{{\mathbf G}(3)}

\newcommand{\bq}{{\mathfrak q}}
\def\qb{\mathfrak{q}}
\newcommand{\toba}{\widetilde{\mathcal E}_{\bq}}

\newcommand{\cO}{{\mathcal O}}
\newcommand{\vsp}{\vspace*{-0.7cm}}

\newcommand{\Ss}{{\mathcal S}}
\newcommand{\ztu}{\overline{\zeta}}

\usepackage{color}

\numberwithin{equation}{section}\theoremstyle{plain}

\newtheorem{theorem}{Theorem}[section]

\newtheorem{lem}[theorem]{Lemma}
\newtheorem{cor}[theorem]{Corollary}

\newtheorem{pro}[theorem]{Proposition}

\theoremstyle{definition}

\theoremstyle{remark}
\newtheorem{rem}[theorem]{Remark}

\newcommand{\ydh}{{}^{H}_{H}\mathcal{YD}}

\def\G{\mathbb{G}}

\def\Jb{\mathbb{J}}
\newcommand\id{\operatorname{id}}

\newcommand\ord{\operatorname{ord}}

\newcommand\Alg{\operatorname{Alg}}
\newcommand\Cleft{\operatorname{Cleft}}
\newcommand\Res{\operatorname{Res}}
\newcommand\End{\operatorname{End}}

\newcommand\Isom{\operatorname{Isom}}

\newcommand\gr{\operatorname{gr}}
\newcommand\co{\operatorname{co}}

\newcommand\ad{\operatorname{ad}}
\newcommand\GL{\operatorname{GL}}

\def\hm{{}_{H}\mathcal{M}}

\def\q{\mathbf{q}}

\def\k{\Bbbk}
\def\ku{\Bbbk}

\def\ot{\otimes}

\def\Z{\mathbb{Z}}
\def\N{\mathbb{N}}

\def\B{\mathfrak{B}}

\def\Lb{\mathbf{L}}
\def\eps{\epsilon}
\def\mT{\mathcal{T}}

\def\mL{\mathcal{L}}
\def\mP{\mathcal{P}}

\def\mR{\mathcal{R}}
\def\mJ{\mathcal{J}}
\def\mF{\mathcal{F}}

\def\mA{\mathcal{A}}

\def\mE{\mathcal{E}}
\def\mtE{\widetilde{\mE}_{\bq}}

\def\mH{\mathcal{H}}

\def\Sym{\mathbb{S}}

\def\br{\mathfrak{br}}
\def\brj{\mathfrak{brj}}
\def\bgl{\mathfrak{wk}}
\def\el{\mathfrak{el}}
\def\g{\mathfrak{g}}
\def\ufo{\mathfrak{ufo}}
\def\u{\mathfrak{u}}
\def\U{\mathfrak{U}}

\def\by{\mathbf{y}}

\newcommand{\J}{{\mathcal J}}
\newcommand{\Gc}{{\mathcal G}}
\renewcommand{\Mc}{{\mathcal M}}
\newcommand{\D}{{\mathcal D}}
\newcommand\I{\mathbb I}

\def\bs{\boldsymbol}

\def\lg{\langle}
\def\rg{\rangle}

\def\pf{\begin{proof}}
\def\epf{\end{proof}}

\newcommand{\Dchaintwo}[3]{\xymatrix@C-4pt{\overset{#1}{\underset{\ }{\circ}}\ar
@{-}[r]^{#2}
& \overset{#3}{\underset{\ }{\circ}}}}
\newcommand{\Dchainthree}[5]{\xymatrix@C-6pt{
\overset{#1}{\underset{\ }{\circ}}\ar  @{-}[r]^{#2}  & \overset{#3}{\underset{\
}{\circ}}\ar  @{-}[r]^{#4}
& \overset{#5}{\underset{\ }{\circ}} }}

\newcommand{\Dchainfour}[7]{\xymatrix@C-6pt{\overset{#1}{\underset{\ }{\circ}}\ar
@{-}[r]^{#2}
& \overset{#3}{\underset{\ }{\circ}}\ar  @{-}[r]^{#4}  & \overset{#5}{\underset{\
}{\circ}} \ar  @{-}[r]^{#6}
& \overset{#7}{\underset{\ }{\circ}}}}

\newcommand{\Dchainfive}[9]{\xymatrix@C-6pt{\overset{#1}{\underset{\ }{\circ}}\ar
@{-}[r]^{#2}  & \overset{#3}{\underset{\ }{\circ}}\ar  @{-}[r]^{#4}  &
\overset{#5}{\underset{\ }{\circ}}
\ar  @{-}[r]^{#6}  & \overset{#7}{\underset{\ }{\circ}}\ar  @{-}[r]^{#8}  &
\overset{#9}{\underset{\ }{\circ}}}}

\newcommand{\Dtriangle}[6]{
\xymatrix@R-12pt{  &    \overset{#2}{\circ} \ar  @{-}[dl]_{#4}\ar  @{-}[dr]^{#5} & \\
\overset{#1}{\circ} \ar  @{-}[rr]^{#6}  &  &\overset{#3}{\circ} }}

\newcommand{\Dthreefork}[8]{
\rule[-9\unitlength]{0pt}{12\unitlength}
\begin{picture}(28,12)(0,9)
\put(2,10){\ifthenelse{\equal{#1}{l}}{\circle*{2}}{\circle{2}}}
\put(3,10){\line(1,0){10}}
\put(14,10){\ifthenelse{\equal{#1}{m}}{\circle*{2}}{\circle{2}}}
\put(15,10){\line(1,1){7}}
\put(15,10){\line(1,-1){7}}
\put(22,18){\ifthenelse{\equal{#1}{t}}{\circle*{2}}{\circle{2}}}
\put(22,2){\ifthenelse{\equal{#1}{b}}{\circle*{2}}{\circle{2}}}
\put(2,12){\makebox[0pt]{\scriptsize #2}}
\put(8,11){\makebox[0pt]{\scriptsize #3}}
\put(14,12){\makebox[0pt]{\scriptsize #4}}
\put(19,16){\makebox[0pt][r]{\scriptsize #5}}
\put(19,4){\makebox[0pt][r]{\scriptsize #6}}
\put(24,17){\makebox[0pt][l]{\scriptsize #7}}
\put(24,2){\makebox[0pt][l]{\scriptsize #8}}
\end{picture}}

\newcommand{\Drightofway}[8]{\xymatrix@R-6pt{  &    \overset{#6}{\circ} \ar
@{-}[d]_{#4}\ar  @{-}[dr]^{#7} & \\
\overset{#1}{\circ} \ar  @{-}[r]^{#2}  &\overset{#3}{\circ} \ar  @{-}[r]^{#5}
&\overset{#8}{\circ} }}

\begin{document}


 \title[All liftings are cocycle deformations]{\small  Liftings of Nichols algebras of diagonal type \\ II. All liftings are cocycle deformations}
\author[Angiono; Garc\'ia Iglesias]{Iv\'an Angiono; Agust\'in Garc\'ia Iglesias}

\address{FaMAF-CIEM (CONICET), Universidad Nacional de C\'ordoba,
Medina A\-llen\-de s/n, Ciudad Universitaria, 5000 C\' ordoba, Rep\'
ublica Argentina.} \email{angiono@famaf.unc.edu.ar, aigarcia@famaf.unc.edu.ar}

\thanks{\noindent 2010 \emph{Mathematics Subject Classification.}
16T05. \newline The work was partially supported by CONICET,
FONCyT-ANPCyT, Secyt (UNC), the MathAmSud project
GR2HOPF}

\begin{abstract}
We classify finite-dimensional pointed Hopf algebras with abelian coradical, up to isomorphism, and show that they are cocycle deformations of the associated graded Hopf algebra. 

More generally, for any braided vector space of diagonal type $V$ with a principal realization in the category of Yetter-Drinfeld modules of a cosemisimple Hopf algebra $H$ and such that the Nichols algebra $\B(V)$ is finite-dimensional, thus presented by a finite set $\Gc$ of relations, we define a family of Hopf algebras $\u(\bs\lambda)$, $\bs\lambda\in \k^{\Gc}$, which are cocycle deformations of $\B(V)\# H$ and such that $\gr\u(\bs\lambda)\simeq \B(V)\# H$.
\end{abstract}

\maketitle

\section{Introduction}\label{sec:intro}

The program for the classification of finite-dimensional pointed Hopf algebra with abelian coradical $\k\Gamma$ was initiated by N.~Andruskiewitsch and H.J.~Schneider in \cite{AS-lift-meth}, where they proposed a {\it lifting method} to deal with this problem. The method proved to be very fruitful as they completed the classification in the case when $|\Gamma|$ was not divisible by small primes in \cite{AS-annals}.

In \cite{AAG}, we proposed a strategy to compute all liftings of braided vector spaces of diagonal type $V\in\ydh$, that is all Hopf algebras $\u$ such that $\gr \u=\B(V)\# H$. This is a large class of Hopf algebras, which includes f.d.~pointed Hopf algebras with abelian coradical. In this article we show that this algorithm is exhaustive. Our method recovers the examples in \cite{AS-annals} and completes the classification, as it becomes {\it independent} of $\Gamma$ and holds more generally for any coradical $H$. Moreover, it also applies for Nichols algebras of infinite dimension which are finitely presented.

A remarkable feature of this method is that all liftings are constructed as cocycle deformations of $\u_0=\B(V)\# H$. 
This has important implications on the study of tensor categories and more general Hopf algebras. On the one hand, it is equivalent to the fact that the categories of comodules of $\u$ and $\u_0$ are tensor equivalent; in the finite-dimensional case this equivalence extends to the categories of Yetter-Drinfeld modules ${}^{\u}_{\u}\mathcal{YD}$ and ${}^{\u_0}_{\u_0}\mathcal{YD}$. This permits the construction of new examples of finite-dimensional Nichols algebras over the basic Hopf algebra $\u^\ast$ \cite{AA-new}, and \emph{a fortiori} of new Hopf algebras by bosonization with $\u^\ast$. 


We fix
\begin{itemize} [leftmargin=*]\renewcommand{\labelitemi}{$\circ$}
\item A cosemisimple Hopf algebra $H$.
\item A f.d.~braided vector space  of diagonal type $(V,c)$ with a principal realization in $\ydh$ and such that $\dim \B(V)<\infty$.
\end{itemize}
Let $\Gc\subset T(V)$ be a (minimal) set of relations defining $\B(V)$. We produce, for each family $\bs\lambda=(\lambda_r)_{r\in\Gc}$ of scalars in a subset $\mR\subset\k^{\Gc}$, 
\begin{itemize} [leftmargin=*]\renewcommand{\labelitemi}{$\bullet$}
\item a Hopf algebra $\u(\bs\lambda)$ such that
  \[
  \gr\u(\bs\lambda)\simeq \B(V)\#H.
  \]
  that is, a {\it lifting} of $V$. Moreover, $\u(\bs\lambda)$ is a cocycle deformation of $\B(V)\#H$.
\end{itemize}

Our main result is the following:
\begin{theorem}\label{thm:cocycle-defs}
Let $H$ be a cosemisimple Hopf algebra and let $(V,c)$ be a braided vector space  of diagonal type such that $\dim\B(V)<\infty$.
Let $L$ be a Hopf algebra with coradical $L_0\simeq H$ and infinitesimal braiding given by a principal realization
$V\in \ydh$. Then
\begin{enumerate}
\item There is $\bs\lambda\in\mR$ such that $L\simeq \u(\bs\lambda)$.
\item $L$ is a cocycle deformation of $\B(V)\#H$.
\end{enumerate}
\end{theorem}
See \S\ref{sec:def-u} below for the definition of the algebras $\u(\bs\lambda)$ and \S\ref{sec:proof-u} for comments on the proof of this result. See also Theorem \ref{thm:isos} for the description of the isomorphism classes.
\pf
Let $\pi\colon \gr L\twoheadrightarrow L_0\simeq H$ be the canonical projection. Then $R=(\gr L)^{\co\pi}$ is a coradically graded Hopf algebra in $\ydh$ and $\gr L=R\# H$. 
By \cite[Theorem 2]{A-nichols}, $R=\B(V)$. The result now follows by Theorem \ref{teo:lift}.
\epf

When $H$ is the group algebra of a finite abelian group $\Gamma$, Theorem \ref{thm:cocycle-defs} is equivalent to the following.
\begin{cor}\label{cor:pointed}
Let $A$ be a finite-dimensional pointed Hopf algebra with abelian coradical. Then $A$ is a cocycle deformation of $\gr A$ and there is $\bs\lambda\in\mR$ such that $A\simeq \u(\bs\lambda)$.\qed
\end{cor}

\subsection{Definition of $\u(\bs\lambda)$}\label{sec:def-u}

Set $\theta=\dim V$. In this setting there exist 
 \begin{align*}
g_i&\in G(H), & \chi_i&\in \Alg(H,\k), \qquad 1\leq i\leq\theta,
 \end{align*}
 such that $\Gamma \coloneqq \lg g_1,\dots,g_\theta\rg < Z(G(H))$ and there is a linear basis $\{x_1,\dots,x_\theta\}$ of $V$ so that $x_i\in V_{g_i} ^{\chi_i}=V_{g_i}\cap V^{\chi_i}$, $1\leq i\leq \theta$, cf.~\eqref{eqn:components}. Moreover, we may assume that $\Gc$ is composed of $\Z^\theta$-homogeneous elements, so that if $r\in\Gc$ is of degree $(a_1,\dots,a_\theta)\in \Z^\theta$, then $r\in T(V)_{g_r}^{\chi_r}$, where
 \begin{align}\label{eqn:gr}
g_r&:=g_1^{a_1}\dots g_\theta^{a_\theta}\in\Gamma, & \chi_r&:=\chi_1^{a_1}\dots \chi_\theta^{a_\theta}\in \Alg(H,\k).
 \end{align}
For each $\bs\lambda=(\lambda_r)_{r\in\Gc}$ in the set
 \begin{align}\label{eqn:R}
 \mR\coloneqq\{\bs\lambda=(\lambda_r)_{r\in \Gc}\in\k^{\Gc}:\lambda_r=0\text{ if } \chi_r\neq \eps\text{ or } g_r=1 \},
 \end{align} we define a Hopf algebra $\u(\bs\lambda)$ 
 recursively as follows: 
\begin{itemize} [leftmargin=*] 
\item We fix a stratification $\Gc=\Gc_0\sqcup\dots \sqcup G_\ell$ of the set of relations of $\B(V)$ in such a way that $\Gc_{k+1}$ is a set of primitive elements of 
\[\B_{k+1}\coloneqq T(V)/\lg \Gc_0\cup\dots\cup\Gc_k\rg.\]
\item We set $\mE_0(\bs\lambda)=T(V)$ and proceed stepwise by defining an algebra 
\[\mE_{k+1}(\bs\lambda)\coloneqq \mE_k(\bs\lambda)/\lg \gamma_k(r)-\lambda_r|r\in\Gc_{k+1}\rg.\]
In this setting, $\mA_k(\bs\lambda)=\mE_k(\bs\lambda)\# H$ is a right $\B_{k}\#H$-cleft object and $\gamma_k:\B_{k}\#H\to \mA_k(\bs\lambda)$ is a section.
\item We set $\u_0(\bs\lambda)=T(V)\#H$ and proceed stepwise by defining an algebra 
\[\u_{k+1}(\bs\lambda)\coloneqq \u_k(\bs\lambda)/\lg \tilde{r}-\lambda_r(1-g_r)|r\in\Gc_{k+1}\rg\]
in such a way that $\mA_k(\bs\lambda)$ becomes a left $\u_{k}(\bs\lambda)$-cleft object and $\tilde{r}\in\u_k(\bs\lambda)$ is defined so that the coaction satisfies  $\gamma_k(r)_{(-1)}\ot \gamma_k(r)_{(0)}=\tilde{r}\ot 1+g_r\ot \gamma_k(r)$. Thus,
\begin{align}\label{eqn:u}
\u(\bs\lambda)\coloneqq T(V)\# H/\lg \tilde{r}-\lambda_r(1-g_r)|r\in\Gc_{k}\rg_{0\leq k\leq\ell}.
\end{align}
\end{itemize}

The explicit presentation of the algebras $\u(\bs\lambda)$ is a hard computational problem. When the braiding is of Cartan type $A$, a complete description was achieved in \cite{AAG} and the Cartan type $G_2$ was completed in \cite{JG}.

\subsection{On the proof}\label{sec:proof-u}
The proof of our main theorem relies on the fact that we are able to show that, for each braiding matrix $\qb$, a given algebra $\mtE(\bs\lambda)$ is nonzero. This is a technical condition pointed in our previous work \cite{AAG} as an obstruction to prove that the list of algebras $\u(\bs\lambda)$ was exhaustive. 

This is achieved in Theorem \ref{thm:E}, whose proof takes over more than 50 pages in the Appendix, as it is a case-by-case analysis of a series of algebras associated to each generalized Dynkin diagram in the classification of matrices $\qb$, in terms of Dynkin diagrams, from \cite{H-classif}. Each example is preceded by the corresponding Dynkin diagram.

Although being tedious and of considerable length (relative to the main body of the paper), we decided to include the full proof here, as there is no general argument underlying it, and each example is analyzed with a subtle combination of different ad-hoc techniques, which are developed in the present article.

\subsection{Organization} In Section \ref{sec:prelim} we review some general Hopf-related concepts and the generic setup for the strategy in \cite{AAGMV}. In \S\ref{sec:cocycles} we present a series of lemmas which set the ground to prove our key Theorem \ref{thm:E}, from which Theorem \ref{thm:cocycle-defs} follows. As stated, the proof of Theorem \ref{thm:E} is completed in the Appendix. In \S\ref{sec:non-connected} we show that the non-connected diagrams can be easily understood from the connected case. Finally, in Section \ref{sec:isos} we study the isomorphism classes of the algebras $\u(\bs\lambda)$ from the classification.

\section*{Acknowledgments} We thank Nicol\'as Andruskiewitsch for his constant support and council. 
We also thank Cristian Vay for pointing to us a mistake in a previous version of this article. We thank the referee for his/her comments, that we believe have helped to improve the presentation of the article, as well as the scope of potential readers.

\section{Preliminaries}\label{sec:prelim}
We work over an algebraically closed field of characteristic zero $\k$. If $A$ is a $\k$-algebra and $S\subseteq A$ is a set, then we denote by $\lg S\rg$ the ideal generated by $S$ in $A$. We denote by $\Alg(A,\k)$ the set of algebra maps $A\to\k$. We write $\G_N$ for the group of $N$th roots of 1 in $\k$, and $\G'_N\subset \G_N$ for the subset of primitive roots.

If $\Gamma$ is an abelian group, then we denote by $\widehat{\Gamma}$ the group of characters of $\Gamma$.  If $n\in \N$, then we set $\I_n\coloneqq\{1,\dots,n\}$, $\I_n^{\circ}\coloneqq\{0,1,\dots,n\}$. If $G$ is a group and $S\subseteq G$ is a subset, then $\lg S\rg$ shall denote the subgroup generated by $S$.

If $H$ is a Hopf algebra, then we denote by $(H_n)_{n\geq 0}$ the coradical filtration of $H$ and set $\gr H=\oplus_{n\geq 0}H_n/H_{n-1}$ the associated graded coalgebra; here $H_{-1}\coloneqq0$.  We denote by $\mP(H)$ the subspace of primitive elements in $H$. 
If $H'$ is another Hopf algebra, then $\Isom(H,H')$ denotes the set of Hopf algebra isomorphisms $H\to H'$. 

We say that $H'$ {\it is a cocycle deformation of} $H$ if there a is a Hopf cocycle $\sigma:H\ot H\to \k$ in such a way that $H'\simeq H_{\sigma}$, the Hopf algebra obtained by deforming the multiplication of $H$. We recall that this is equivalent to the existence of a $(H,H')$-bicleft object $A$. This is the approach we follow in the series. See \cite{S} for details. 

\medbreak

We write $\hm$ for the category of $H$-modules; an $H$-module algebra is thus an algebra in $\hm$.
In turn, $\ydh$ stands for the category of Yetter-Drinfeld modules over $H$; this is a braided tensor category. Given $V\in\ydh$, $g\in G(H)$ and $\chi\in\Alg(H,\k)$, then we set
\begin{align}\label{eqn:components}
V_g&\coloneqq\{x\in V| \delta(x)= g\ot x\}, & V^\chi&\coloneqq\{x\in V| h\cdot x=\chi(h)x,\, h\in H\}.
\end{align}

If $V\in\ydh$, then we denote by $\B(V)$ the {\it Nichols algebra} of $V$: This is the quotient $\B(V)=T(V)/\mJ(V)$ where $\mJ(V)$ is the maximal graded Hopf ideal generated by homogeneous elements of degree $\geq 2$. In particular, it inherits a graduation $\B(V)=\oplus_{n\geq 0}\B^n(V)$ with $\B^0(V)=\k$ and $\B^1(V)=V$. We recall that it only depends on the underlying braided vector space $(V,c_{V,V})$; see \cite{AS-annals} and references therein for details.

\medbreak 

A {\it lifting} of $V\in\ydh$ is a Hopf algebra $L$ such that $\gr L\simeq \B(V)\#H$. A {\it lifting map} is a Hopf algebra projection $\phi:T(V)\#H\twoheadrightarrow L$ such that
\begin{align}\label{eqn:lifting-map}
\phi_{|H}&=\id_H & \text{and} && \phi_{|V\#H}&\colon V\#H\stackrel{\simeq}{\longrightarrow} P_1 &\text{as  Hopf bimodules over }H,
\end{align}
where $P_1=L_1\cap\ker\Pi$  and $\Pi:L\twoheadrightarrow H$ is a $H$-bimodule coalgebra projection such that $\Pi_{|H}=\id_H$, see \cite[\S2.5]{AV} for details.

\subsection{Diagonal braidings}\label{sec:diagonal}
Let $(V,c)$ be a braided vector space. We set $\theta=\dim V$, $\I=\I_\theta$. Assume that $V$ is of diagonal type, with braiding matrix $\qb=(q_{ij})_{i,j\in\I}$; i.e.~there is a basis $\{x_1,\dots,x_\theta\}$ such that $c=c^\qb$ satisfies:
\[
c^\qb(x_i\ot x_j)=q_{ij}x_j\ot x_i, \quad i,j\in\I.
\]
We shall denote by $\B_{\qb}$ the Nichols algebra associated to $(V,c^\qb)$.

\subsubsection{} There is a (generalized) Dynkin diagram $\D$ attached to the braiding matrix $\qb$ \cite{H-classif}. This is a graph whose set of vertices is $\I$ and there is an edge connecting $i\neq j\in \I$ if and only if $\widetilde{q_{ij}}\coloneqq q_{ij}q_{ji}\neq 1$. It is decorated with a label $q_{ii}$ on each vertex $i\in \I$ and a label $\widetilde{q_{ij}}$ on the edge between $i$ and $j$.

Let $\D=\bigsqcup_{i=1}^m \D_i$ be the decomposition of $\D$ into $m$ connected components, $m\in\N$. Given two vertices $i,j\in\I$ of $\D$, then we write $i\sim j$ if they belong to the same connected component and $i\not\sim j$ otherwise. If $V=\bigoplus_{i=1}^m V_i$ is the corresponding decomposition of $V$ into braided subspaces $V_i$, associated to the connected diagram $\D_i$, $i\in\I_m$, then we refer to each $V_i$, $i\in \I_m$, as a {\it (connected) component} of $V$. 

We shall say that the matrix ${\qb}$ is connected when the corresponding Dynkin diagram is so.

In \cite{H-classif}, Heckenberger classified all connected matrices ${\qb}$ such that the Nichols algebra $\B_{\qb}$ is finite-dimensional\footnote{More generally, \cite{H-classif} contains the classification of all matrices ${\qb}$ with finite root system.}. A list of the corresponding Dynkin diagrams can be found in loc.cit. We shall make explicit use of many of the diagrams of this list, and hence they will appear in the case-by-case analysis in the Appendix. We choose not to reproduce the complete list here, as our computations rely on connected subdiagrams of small rank.

\subsubsection{}\label{sec:presentation}

We fix $\bq$ a (connected) braiding matrix of the lists in \cite{H-classif}, $(V,c)$ the associated braided vector space of diagonal type and $(x_i)_{i\in \I}$ a basis of $V$ such that $c(x_i\ot x_j)=q_{ij} x_j\ot x_i$.

Next we give a presentation for the Nichols algebra $\B_{\bq}$. We recall that $\B_{\bq}$ is a quotient of $T(V)$ by a $\Z^\I$-homogeneous ideal $\J(V)$, where the degree on $T(V)$ is determined by $\deg x_i=\alpha_i$, $(\alpha_i)_{i\in\I}$ the canonical basis of $\Z^\I$.
Following \cite{AA-new} we introduce more notation needed to describe the relations:

\begin{itemize}[leftmargin=*]
\item By abuse of notation $\bq$ also denotes the $\Z$-bilinear form $\Z^\I\times\Z^\I\to \Bbbk^\times$ such that $\bq_{\alpha_i,\alpha_j}=q_{ij}$, $i,j\in\I$.
\item Given $x,y\in T(V)$ homogeneous of degrees $\alpha,\beta\in\Z^\I$, the braided commutator is defined by
$[x,y]_c = xy- \bq_{\alpha,\beta} yx$.
\item In particular, given $i\in \I$, $y\in T(V)$, we denote $(\ad_c x_i) y := [x_i,y]_c$.
\item To short the notation, we define recursively
\begin{align*}
x_{ij} &= (\ad_c x_i) x_j, & x_{i \, i_1 \dots x_k} &= (\ad_c x_i) x_{i_1 \dots x_k}, & &i,j,i_1,\dots i_k \in\I.
\end{align*}
\item The Nichols algebra $\B_{\bq}$ admits a PBW basis whose generators are $\Z^\I$-homogeneous. The set $\Delta^{\bq}_+$ of degrees of the PBW generators is unique (under some mild conditions, cf. \cite{AA-new} and the references therein), and is known as the set of positive roots of $\B_{\bq}$.
\item Let $a_{ij}^{\bq}:=-\max \{n\in \N_0: \alpha_j+n \,\alpha_i\in\Delta_+^{\bq} \}$, $i\neq j\in\I$. We say that $i\in\I$ is a Cartan vertex if $\widetilde{q_{ij}}=q_{ii}^{a_{ij}^{\bq}}$ for all $j\neq i$.
\item There exists a groupoid (called the Weyl groupoid, see the precise definition in \cite{AA-new}) acting on the set of roots $\Delta^{\bq} = \Delta^{\bq}_+ \cup (-\Delta^{\bq})$. We say that $\alpha\in\Delta^{\bq}$ is a Cartan root if it is the image of $\alpha_i$ for some Cartan root $i\in\I$. The set of all (positive) Cartan roots is denoted by $\cO^{\bq}$.
\end{itemize}

\begin{theorem}\label{thm:presentacion minima}\cite{A-nichols}
$\B_{\qb}$ is the quotient of $T(V)$ by the relations:
\begin{align*}
&x_\alpha^{N_\alpha}, &\alpha \in \cO^{\bq};
\\ &(\ad_cx_i)^{1-a_{ij}^{\bq}}x_j, & q_{ii}^{1-a_{ij}^{\bq}} \neq 1;
\\ &x_i^{N_i}, & i \mbox{ is not a Cartan vertex};
\end{align*}

\begin{itemize} [leftmargin=*]\renewcommand{\labelitemi}{$\circ$}

\item for $i,j \in \I$ such that
$q_{ii}=\widetilde{q_{ij}}=q_{jj}=-1$, and there exists $k\neq i,j$
such that $\widetilde{q_{ik}}^2\neq1$ or
$\widetilde{q_{jk}}^2\neq1$,
$$ x_{ij}^2; $$

\item for $i,j,k \in \I$ such that
$q_{jj}=-1$,
$\widetilde{q_{ik}}=\widetilde{q_{ij}}\widetilde{q_{kj}}=1$,
$\widetilde{q_{ij}}\neq -1$,
$$
\left[ x_{ijk} , x_j \right]_c;
$$

\item for $i,j \in \I$ such that
$q_{jj}=-1$, $q_{ii}\widetilde{q_{ij}}\in \G_6$,
$\widetilde{q_{ij}}\neq -1$, and either $q_{ii}\in \G_3$ or else 
$a_{ij}^{\bq}\leq -3$,
$$
\left[ x_{iij}, x_{ij} \right]_c;
$$

\item for $i,j,k \in \I$ such that $q_{ii}=\pm \widetilde{q_{ij}}\in\G_3$, $\widetilde{q_{ik}}=1$, and either
$-q_{jj}=\widetilde{q_{ij}}\widetilde{q_{jk}}=1$ or else  $q_{jj}^{-1}=\widetilde{q_{ij}}=\widetilde{q_{jk}}\neq -1$,
$$
\left[ x_{iijk} , x_{ij} \right]_c;
$$

\item for $i,j,k \in \I$ such that $\widetilde{q_{ik}}, \widetilde{q_{ij}}, \widetilde{q_{jk}} \neq 1$,
$$
x_{ijk}-\frac{1-\widetilde{q_{jk}}}{q_{kj}(1-\widetilde{q_{ik}})}\left[x_{ik},x_j\right]_c-q_{ij}(1-\widetilde{q_{jk}}) \ x_jx_{ik};
$$

\item for $i,j,k\in\I$ such that one of the following situations hold
\begin{enumerate}
\item $q_{ii}=q_{jj}=-1$, $\widetilde{q_{ij}}^2= \widetilde{q_{jk}}^{-1}$, $\widetilde{q_{ik}}=1$, 
\item $\widetilde{q_{ij}}=q_{jj}=-1$, $q_{ii}= -\widetilde{q_{jk}}^2\in\G_3$, $\widetilde{q_{ik}}=1$, 
\item $q_{kk}=\widetilde{q_{jk}}=q_{jj}=-1$, $q_{ii}= -\widetilde{q_{ij}}\in\G_3$, $\widetilde{q_{ik}}=1$, 
\item $q_{jj}=-1$, $\widetilde{q_{ij}}=q_{ii}^{-2}$, $\widetilde{q_{jk}}=-q_{ii}^{-3}$, $\widetilde{q_{ik}}=1$,
\item $q_{ii}=q_{jj}=q_{kk}=-1$, $\pm\widetilde{q_{ij}}=\widetilde{q_{jk}}\in\G_3$,
$\widetilde{q_{ik}}=1$,
\end{enumerate}
$$
\left[ \left[x_{ij}, x_{ijk} \right]_c, x_j \right]_c;
$$

\item for $i,j,k\in\I$ such that $q_{ii}=q_{jj}=-1$, $\widetilde{q_{ij}}^3=\widetilde{q_{jk}}^{-1}$, $\widetilde{q_{ik}}=1$,
$$
\left[ \left[x_{ij}, \left[x_{ij}, x_{ijk} \right]_c \right]_c, x_j \right]_c;
$$

\item for $i,j,k\in\I$ such that
$q_{jj}=\widetilde{q_{ij}}^2=\widetilde{q_{jk}}\in \G_3$,
$\widetilde{q_{ik}}=1$,
$$
\left[ \left[ x_{ijk} , x_j \right]_c x_j \right]_c;
$$

\item for $i,j,k\in\I$ such that
$q_{kk}=q_{jj}=\widetilde{q_{ij}}^{-1}=\widetilde{q_{jk}}^{-1}\in
\G_9$, $\widetilde{q_{ik}}=1$, $q_{ii}=q_{kk}^6$
$$
\left[ \left[ x_{iij} , x_{iijk} \right]_c, x_{ij} \right]_c;
$$

\item for $i,j,k\in\I$ such that
$q_{ii}=\widetilde{q_{ij}}^{-1}\in \G_9$,
$q_{jj}=\widetilde{q_{jk}}^{-1}=q_{ii}^5$, $\widetilde{q_{ik}}=1$,
$q_{kk}=q_{ii}^6$
$$
[\left[x_{ijk}, x_{j} \right]_c, x_k]_c -(1 +
\widetilde{q_{jk}})^{-1}q_{jk} \left[ \left[x_{ijk}, x_{k} \right]_c
, x_{j} \right]_c; $$

\item for $i,j,k\in\I$ such that
$q_{jj}=\widetilde{q_{ij}}^3=\widetilde{q_{jk}}\in \G_4$,
$\widetilde{q_{ik}}=1$,
$$
\left[ \left[ \left[ x_{ijk} , x_j \right]_c, x_j \right]_c, x_j
\right]_c;
$$

\item for $i,j,k\in\I$ such that
$q_{ii} = \widetilde{q_{ij}} =-1$, $q_{jj}= \widetilde{q_{jk}}^{-1}
\neq-1$, $\widetilde{q_{ik}}=1$,
$$
\left[x_{ij}, x_{ijk} \right]_c;
$$

\item for $i,j,k \in\I$ such that
$q_{ii}= q_{kk} =-1$, $\widetilde{q_{ik}}=1$, $\widetilde{q_{ij}}
\in \G_3$, $q_{jj}= -\widetilde{q_{jk}} = \pm \widetilde{q_{ij}}$,
$$
[x_i, x_{jjk}]_c -(1 + q_{jj}^2)q_{kj}^{-1} \left[x_{ijk}, x_{j}
\right]_c - (1 + q_{jj}^2)(1 + q_{jj}) q_{ij} x_j x_{ijk};
$$

\item for $i,j,k,l\in\{1,\ldots,\theta \}$ such that
$q_{jj}\widetilde{q_{ij}}= q_{jj}\widetilde{q_{jk}}=1$, $q_{kk}=-1$,
$\widetilde{q_{ik}}=\widetilde{q_{il}}=\widetilde{q_{jl}}=1$,
$\widetilde{q_{jk}}^2= \widetilde{q_{lk}}^{-1}= q_{ll}$,
$$
\left[\left[\left[x_{ijkl},x_k\right]_c, x_j \right]_c, x_k \right]_c;
$$

\item for $i,j,k,l\in\{1,\ldots,\theta \}$ such that
$\widetilde{q_{jk}}= \widetilde{q_{ij}}= q_{jj}^{-1}\in
\G_4'\cup\G_6'$, $q_{ii}=q_{kk}=-1$,
$\widetilde{q_{ik}}=\widetilde{q_{il}}=\widetilde{q_{jl}}=1$,
$\widetilde{q_{jk}}^3= \widetilde{q_{lk}}$,
$$
\left[\left[x_{ijk},\left[x_{ijkl}, x_k \right]_c\right]_c, x_{jk}
\right]_c;
$$

\item for $i,j,k,l\in\{1,\ldots,\theta \}$ such that
$q_{ll}=\widetilde{q_{lk}}^{-1}=
q_{kk}=\widetilde{q_{jk}}^{-1}=q^2$, $\widetilde{q_{ij}}=
q_{ii}^{-1}=q^3$ for some $q\in \ku^\times$, $q_{jj}=-1$,
$\widetilde{q_{ik}}=\widetilde{q_{il}}=\widetilde{q_{jl}}=1$,
$$
\left[\left[\left[x_{ijk},x_j\right]_c, \left[x_{ijkl},x_j\right]_c
\right]_c, x_{jk} \right]_c;
$$

\item for $i,j,k,l\in\{1,\ldots,\theta \}$ such that one
of the following situations
\begin{enumerate}
\item $q_{kk}=-1$, $q_{ii}=\widetilde{q_{ij}}^{-1}= q_{jj}^2$,
$\widetilde{q_{kl}}= q_{ll}^{-1}= q_{jj}^3$, $\widetilde{q_{jk}}=
q_{jj}^{-1}$,
$\widetilde{q_{ik}}=\widetilde{q_{il}}=\widetilde{q_{jl}}=1$,
\item $q_{ii}=\widetilde{q_{ij}}^{-1}= -q_{ll}^{-1}=-\widetilde{q_{kl}}$,
$q_{jj}=\widetilde{q_{jk}}=q_{kk}=-1$,
$\widetilde{q_{ik}}=\widetilde{q_{il}}=\widetilde{q_{jl}}=1$,
\end{enumerate}
$$
\left[\left[x_{ijkl}, x_j \right]_c, x_k \right]_c-
q_{jk}(\widetilde{q_{ij}}^{-1}-q_{jj})
\left[\left[x_{ijkl},x_k\right]_c, x_j \right]_c;
$$

\item for $i,j,k\in\I$ such that
$\widetilde{q_{jk}}=1$,
$q_{ii}=\widetilde{q_{ij}}=-\widetilde{q_{ik}}\in \G_3$,
$$
\left[x_i, \left[ x_{ij},x_{ik} \right]_c
\right]_c+q_{jk}q_{ik}q_{ji} \left[ x_{iik} ,x_{ij} \right]_c+q_{ij}
\, x_{ij} x_{iik};
$$

\item for $i,j,k\in\I$ such that
$q_{jj}=q_{kk}=\widetilde{q_{jk}}=-1$,
$q_{ii}=-\widetilde{q_{ij}}\in\G_3$, $\widetilde{q_{ik}}=1$,
$$
\left[x_{iijk}, x_{ijk} \right]_c;
$$

\item for $i,j\in\I$ such that $-q_{ii}, -q_{jj}, q_{ii}\widetilde{q_{ij}}, q_{jj}\widetilde{q_{ij}} \neq 1$,
$$
(1-\widetilde{q_{ij}})q_{jj}q_{ji}\left[x_i, \left[ x_{ij}, x_j \right]_c \right]_c - (1+q_{jj})(1-q_{jj}\widetilde{q_{ij}})x_{ij}^2;
$$

\item for $i,j\in\I$ such that either $a_{ij}^{\bq}\le -4$, or else  $a_{ij}^{\bq}=-3$, $q_{jj}=-1$, $q_{ii} \in \G_4$,
$$
\left[x_i,x_{3\alpha_i+2\alpha_j}\right]_c-\frac{1-q_{ii}\widetilde{q_{ij}}-q_{ii}^2\widetilde{q_{ij}}^2q_{jj}}{(1-q_{ii}\widetilde{q_{ij}})q_{ji}}
x_{iij}^2;
$$

\item for $i,j\in\I$ such that $3\alpha_i+2\alpha_j \in \Delta_+^\bq$, $4\alpha_i+3\alpha_j\notin \Delta_+^\bq$, either $q_{jj}=-1$ or else $a_{ji}^{\bq}\leq-2$, and either $a_{ij}^{\bq}\leq -3$,
or else $a_{ij}^{\bq}=-2$, $q_{ii}\in\G_3$,
$$
x_{4\alpha_i+3\alpha_j}=[x_{3\alpha_i+2\alpha_j}, x_{ij} ]_c;
$$

\item for $i,j\in\I$ such that $3\alpha_i+2\alpha_j\in\Delta_+^\bq$, $5\alpha_i+3\alpha_j\notin\Delta_+^\bq$, and
$q_{ii}^3\widetilde{q_{ij}}, q_{ii}^4\widetilde{q_{ij}}\neq 1$,
$$
[x_{iij}, x_{3\alpha_i+2\alpha_j}]_c;
$$

\item for $i,j\in\I$ such that $4\alpha_i+3\alpha_j\in\Delta_+^\bq$, $5\alpha_i+4\alpha_j\notin\Delta_+^\bq$,
$$
x_{5\alpha_i+4\alpha_j}=[x_{4\alpha_i+3\alpha_j}, x_{ij} ]_c;
$$

\item for $i,j\in\I$ such that $5\alpha_i+2\alpha_j\in\Delta_+^\bq$, $7\alpha_i+3\alpha_j \notin \Delta_+^\bq$,
$$
[[x_{iiij}, x_{iij}], x_{iij} ]_c;
$$

\item for $i,j\in\I$ such that $q_{jj}=-1$, $5\alpha_i+4\alpha_j\in\Delta_+^\bq$,

\begin{align*}
& [x_{iij},x_{4\alpha_i+3\alpha_j}]_c- \frac{b-(1+q_{ii})(1-q_{ii}\widetilde{q_{ij}})(1+\widetilde{q_{ij}}+q_{ii}\zeta^2)q_{ii}^6\widetilde{q_{ij}}^4} {a\ q_{ii}^3q_{ij}^2q_{ji}^3} x_{3\alpha_i+2\alpha_j}^2,
\\
& \begin{aligned}
&\text{where }&  a&=(1-\widetilde{q_{ij}})(1-q_{ii}^4\widetilde{q_{ij}}^3)-(1-q_{ii}\widetilde{q_{ij}})(1+q_{ii})q_{ii}\widetilde{q_{ij}}, 
\\ 
&& b&=(1-\widetilde{q_{ij}})(1-q_{ii}^6\widetilde{q_{ij}}^5)-a\ q_{ii}\widetilde{q_{ij}}.
\end{aligned}
\end{align*}

\end{itemize}
\qed
\end{theorem}

Now we assume that $\bq$ is not necessarily connected but every connected component belongs to the list of \cite{H-classif}. We take $\Gc$ as the set of defining relations of $\B_{\bq}$ as computed in Theorem \ref{thm:presentacion minima} for each component, union $\Gc(0)$, where
\begin{equation}\label{eqn:G(0)}
\Gc(0)\coloneqq\{x_ix_j-q_{ij}x_jx_i| i<j \text{ and }i\not\sim j\}
\end{equation}
denotes the set of $q$-commutators of vertices in different components of $\bq$.

\begin{rem}\label{rem:distinguished}
Let $\widetilde \Gc := \Gc -\{x_\alpha^{N_\alpha}:\alpha \in \cO^{\bq}\}$. Then $\widetilde{\B}_{\bq} = T(V)/\langle \widetilde \Gc \rangle$ is an intermediate dostinguished quotient between $T(V)$ and $\B_{\bq}$ called the distinguished pre-Nichols algebra of $\bq$ \cite{A-distinguished}.
\end{rem}

We stress that \cite{AA-new} presents a different approach for the defining relations. The authors list the relations of each Nichols algebra $\B_{\qb}$, for every matrix $\qb$ in \cite{H-classif}. We will mainly follow this second approach. Also, a partial presentation of $\B_{\qb}$ is given on the Appendix.

\subsubsection{}
Let $H$ be a Hopf algebra. A principal YD-realization of $V$ over $H$ is a family $(\chi_i,g_i)_{i\in\I}\in \Alg(H,\k)\times G(H)$ such that $\chi_i(g_j)=q_{ji}, i,j\in \I$, and $\chi_i(h)g_i=\chi_i(h\_{2})h\_{1}g_i\Ss(h\_{3})$ for all $h\in H$, $i\in\I$. In particular, it follows that $\Gamma=\lg g_1,\dots,g_\theta\rg<Z(G(H))$. These data realize $(V,c^\qb)$ as an object in $\ydh$ in such a way that $c^\qb$ coincides with the categorical braiding of $\ydh$ and $x_i\in V_{g_i}^{\chi_i}=V_{g_i}\cap V^{\chi_i}$, $i\in\I$,  cf.~\eqref{eqn:components}.

\subsection{The strategy}\label{sec:strategy}

We recall the setting of a generic step in the algorithm behind our strategy to compute the liftings as cocycle deformations. We refer the reader to \cite{AAGMV,AAG} for details. We fix
\begin{itemize} [leftmargin=*]\renewcommand{\labelitemi}{$\circ$}
\item A cosemisimple Hopf algebra $H$.
\item  $V\in\ydh$ such that the ideal $\J(V)$ defining the Nichols algebra $\B(V)$ is generated by a finite set $\Gc$.
\end{itemize}

Let $\Gc=\Gc_0\sqcup\dots \sqcup\Gc_\ell$ be an {\it adapted stratification} of $\Gc$ cf.~\cite[\S 5.1]{AAGMV}; in particular $\Gc_k$ is a basis of a Yetter-Drinfeld submodule $\Mc_{k}\subset \mP(\B_k)$, for each\footnote{$\Gc_\ell$ can be chosen so that it is not properly contained in  $\mP(\B_\ell)$, see loc.cit.} $k$.
We set $\B_0=T(V)$ and denote by $\B_k=\B_{k-1}/\lg \Gc_{k-1}\rg$, $k\in \I_{\ell+1}$ the corresponding pre-Nichols algebras; also $\mH_k=\B_k\# H$. In loc.~cit.~we introduced an algorithm to produce at each step $k=0,\dots,\ell$ a collection $\Cleft'(\mH_k)$ of cleft objects $\mA_k$ in such a way that the Schauenburg left Hopf algebra $\mL_k$ associated to each pair $(\mH_k,\mA_k)$, see \cite{S}, satisfies $\gr \mL_k\simeq \mH_k$; moreover this algebra is a cocycle deformation of $\mH_k$ by construction.

In this article we show that a technical condition we introduced in \cite{AAG} for the diagonal case is always satisfied, which shows that the family of Hopf algebras resulting from the $\ell+1$th step is a complete family of liftings of $V$.

Set $\mT\coloneqq T(V)\#H$. The algorithm starts with $\B_0=T(V)$ and $\Cleft'(\mH_0)$ is the singleton $\{\mA_0=\mT\}$, with coaction $\rho_0=\Delta_{\mT}$; so $\mL_0=\mT$. At a given step $k$ we compute, for each $\mA_{k-1}\in \Cleft'(\mH_{k-1})$ a collection of quotients $\mA_{k-1}\twoheadrightarrow \mA_k$.
This setup is  depicted in the following {\it snapshot}
from \cite[p.696]{AAGMV}:
\begin{equation}\label{eqn:snapshot}
\xymatrix{
\mT\ar@{->}[rr]^{\gamma_0=\id}\ar@{->>}[dd]^{\pi_{k-1}}\ar@/^-1.5pc/@{->>}[ddd]_{\pi_{k}}
^{\,\equiv} &&
\mT\ar@{->>}[dd]_{\tau_{k-1}}\ar@/^1.5pc/@{->>}[ddd]^{\tau_{k}}_{\equiv\,}\ar@{~>}[rr]
&&L(\mT,\mT)\simeq \mT\ar@{->>}[dd]\ar@/^3.5pc/@{->>}[ddd]_{\equiv\quad}\\
&&&&&\\
\mH_{k-1}\ar@{->}[rr]^{\gamma_{k-1}}\ar@{->>}[d]^{\pi'_{k}} &&
\mA_{k-1}\ar@{->>}[d]_{\tau'_{k}}\ar@{~>}[rr]&&L(\mA_{k-1},\mH_{k-1})\ar@{->>}[d]\\
\mH_{k}\ar@{->}[rr]^{\gamma_{k}} && \mA_{k}\ar@{~>}[rr]&&L(\mA_{k},\mH_{k})
}
\end{equation}
The maps $\gamma_{k-1}$, $\gamma_k$ denote the corresponding sections, and the vertical arrows represent
the (Hopf) algebra projections. The coaction $\rho_k\colon\mA_k\to \mA_k\ot \mH_k$ is given by
$\rho_k=(\tau_k'\ot\pi_k')\rho_{k-1}=(\tau_k \ot\pi_k)\Delta_{\mT}$ and the section satisfies $\gamma_k(xh)=\gamma_k(x)h$, $x\in\B_k$, $h\in H$.

In the final step $k=\ell+1$, $\mH_{\ell+1}=\B(V)\#H$ and $\mL=L(\mA,\B(V)\#H)$ is a cocycle deformation of $\B(V)\#H$ satisfying $\gr \mL\simeq \B(V)\# H$, that is a lifting of $V$ for every $\mA\in\Cleft'(\B(V)\#H)$.

\subsection{On the computation of $\mA_{k}$}
By \cite[Proposition 5.8]{AAGMV}, each algebra $\mA_{k-1}$ can be decomposed as $\mA_{k-1}=\mE_{k-1}\# H$, where $\mE_0=T(V)$ and each $\mE_{k-1}$ is an $H$-module algebra. Hence the problem of computing all quotients $\mA_{k-1}\twoheadrightarrow\mA_k$ translates into that of finding all suitable quotients $\mE_{k-1}\twoheadrightarrow \mE_k$.

For each $\bs\lambda=(\lambda_r)_{r\in\Gc}\in\k^{\Gc}$, $\mE_k(\bs\lambda)$ is defined, as a $\k$-linear space, via
\begin{equation}\label{eqn:E(lambda)}
\mE_k=\mE_{k}(\bs\lambda)=\mE_{k-1}/\lg \gamma_{k-1}(r)-\lambda_r : r\in \Gc_{k-1}\rg,
\end{equation}
If $\mE_{k}(\bs\lambda)\neq 0$, then it is a $\k$-algebra. If, additionally, $\mE_{k}(\bs\lambda)$ is an $H$-module 
algebra, then $\mA_k(\bs\lambda)\coloneqq\mE_{k}(\bs\lambda)\#H\in\Cleft(\mH_k)$, by \cite[Proposition 5.8]{AAGMV}.

The next lemma indicates when $\mE_k(\bs\lambda)$ is an $H$-module algebra. In \S \ref{sec:cocycles} we study when $\mE_{k}(\bs\lambda)\neq 0$ for the case of diagonal braidings.

We introduce some terminology first. For each $r\in\Gc_{k}$ and each $h\in H$, we write the $H$-action as $h\cdot r=\sum_{s\in\Gc_{k}}\mu_{rs}(h)s$, some $(\mu_{rs}(h))_{s\in\Gc_{k}}\in\k^{\Gc_{k}}$.

\begin{lem}\label{lem:ideal}
Assume that $\mE_k(\bs\lambda)$ is an $H$-module algebra and that the section $\gamma_{k}\colon \mH_{k}\to \mA_{k}(\bs\lambda)$ is $H$-linear. 
If $\mE_{k+1}(\bs\lambda)\neq 0$, then it is an $H$-module algebra if and only if 
\begin{align}\label{eqn:condition-ideal}
\eps(h)\lambda_r=\sum\limits_{s\in\Gc_{k-1}}\mu_{rs}(h)\lambda_s, \qquad \forall\,r\in\Gc_{k},h\in H.
\end{align}
\end{lem}
\pf If \eqref{eqn:condition-ideal} holds, then the ideal
$
\lg \gamma_{k}(r)-\lambda_r : r\in \Gc_{k}\rg\subset \mE_{k}
$
is an $H$-submodule of $\mE_{k}$ and thus $\mE_{k+1}(\bs\lambda)$ inherits an $H$-module structure. Conversely,  set $\kappa(h,r)=\sum\limits_{s\in\Gc_{k-1}}\mu_{rs}(h)\lambda_s-\eps(h)\lambda_r\in\k$, for $h\in H$, $r\in \Gc_k$. Then
\[
h\cdot (\gamma_k(r)-\lambda_r)=\sum_{s\in\Gc_{k-1}}\mu_{rs}(h)\left(\gamma(s)-\lambda_s\right) + \kappa(h,r).
\]
Since $\mE_{k+1}(\bs\lambda)\neq 0$, $\kappa(h,r)=0$ for all $h,r$ and  \eqref{eqn:condition-ideal} holds.
\epf

\begin{rem}
If $k=0$, then $\gamma_k=\id$ thus it is $H$-linear. On the other hand, if $H$ is semisimple, then $\gamma_k$ is $H$-linear for every $k$, \cite[Proposition 5.8(c)]{AAGMV}. This is also the case when the braiding is diagonal, see \S \ref{sec:cocycles}.
\end{rem}

\section{All liftings are cocycle deformations}\label{sec:cocycles}

In this section we further assume that
\begin{itemize} [leftmargin=*]\renewcommand{\labelitemi}{$\circ$}
\item
$V\in\ydh$ comes from a principal realization $(g_i,\chi_i)_{i\in\I_{\theta}}$ of a braided vector space of diagonal type $(V,c)$ with matrix $\bq$.
\end{itemize}  We show that every lifting of $V$ is a cocycle deformation of $\B_{\qb}\# H$, by showing that the algebras $\mE_k(\bs\lambda)\in\hm$ in \eqref{eqn:E(lambda)} are nonzero.

We consider $\Gc$ as in \S \ref{sec:diagonal}. We set $\Gamma=\lg g_1,\dots,g_\theta\rg < G(H)$, so each $r\in\Gc$ belongs to a unique component $r\in T(V)_{g_r}^{\chi_r}$, for some  $g_r\in\Gamma$ and $\chi_r\in\Alg(H,\k)$, cf.~\eqref{eqn:gr}. We fix a stratification $\Gc_0\sqcup\dots \sqcup\Gc_\ell$ in such a way that $\Gc_0=\Gc(0)$ as in \eqref{eqn:G(0)} and such that $\k\lg\Gc_\ell\rg$ is the normal braided Hopf subalgebra of $\B_\ell$ generated by the powers of the Cartan root vectors, cf.~\cite[Theorem 31]{A-distinguished}. When $\mE_k(\bs\lambda)\neq 0$ for every $k$, we shall denote
\begin{align}\label{eqn:Es}
\mtE(\bs\lambda)&\coloneqq \mE_{\ell}(\bs\lambda),& \mE_{\qb}(\bs\lambda)&\coloneqq \mE_{\ell+1}(\bs\lambda).
\end{align}

\smallbreak

In this context, Lemma \ref{lem:ideal} is equivalent to the following. Recall the notation from \eqref{eqn:snapshot}.
\begin{lem}\label{lem:decomposition}
Assume that $\mE_k(\bs\lambda)\neq 0$. Then $\gamma:\mH_k\to \mA_k(\bs\lambda)$ is $H$-linear and
$\mE_k(\bs\lambda)$ is an $H$-module algebra if and only if
\begin{equation}\label{eqn:condition-lambda-k}
\lambda_r=0 \quad\text{ if }\chi_r\neq \eps, \quad \text{for all }r\in\Gc_{k-1}.
\end{equation}
Let $\bs\lambda=(\lambda_r)_{r\in\Gc}\in\k$ be a family of scalars satisfying \eqref{eqn:condition-lambda-k}. For each $k\geq 0$,
$\mE=\mE_k(\bs\lambda)$ decomposes as a sum of $\Gamma$-eigenspaces
$\mE=\bigoplus_{\chi\in\widehat{\Gamma}}\mE^\chi$.
\end{lem}
\pf
We proceed by induction on $k$. As $V=\bigoplus_{i\in\I_\theta} V^{\chi_i}$, the $H$-module algebra $\mE_0=T(V)$ decomposes as a sum of $\Gamma$-eigenspaces
 $T(V)=\bigoplus_{\chi\in\widehat{\Gamma}}T(V)^\chi$ and $\gamma=\id:\mH_0\to\mA_0$ is $H$-linear.
 
Fix $k\geq 0$ for which the statement holds. By definition, the pre-Nichols algebra $\B_k$ also decomposes as a sum $\B_k=\bigoplus_{\chi\in\widehat{\Gamma}}\B_k^\chi$. Let $\tau_k:T(V)\# H\to \mA_k$ be the algebra projection; by assumption, this is an $H$-linear map --which restricts to an $H$-linear map $\tau_k:T(V)\to \mE_k=\bigoplus_{\chi\in\widehat{\Gamma}}\mE^\chi$. Let $x\#t\in\mH_k$, $x\in \B_k^\chi$, some $\chi\in \widehat{\Gamma}$, $t\in H$. As $\gamma(xt)=\gamma(x)t$ cf.~\eqref{eqn:snapshot}, we may assume $t=1$. Now, let $y\in T(V)^\chi$ be such that $\tau_k(y)=\gamma(x)$. Then 
\[
h\cdot \gamma(x)=h\cdot \tau_k(y)=\tau_k(h\cdot y)=\chi(h)\tau_k(y)=\chi(h)\gamma(x), \ \forall\,h\in H.
\]
Hence $\gamma:\mH_k\to \mA_k$ is $H$-linear and the first assertion follows by Lemma \ref{eqn:condition-ideal}, as  \eqref{eqn:condition-ideal} becomes \eqref{eqn:condition-lambda-k} in this setting. This also shows that the ideal defining $\mE_{k+1}$, see \eqref{eqn:E(lambda)}, is generated by $\widehat{\Gamma}$-homogeneous elements, and thus the $\widehat{\Gamma}$-graduation descends to $\mE_{k+1}$. The lemma follows.
\epf

Next we show that condition \eqref{eqn:condition-lambda-k} is sufficient. Namely, recall the definition of the set $\mR\subset\k^{\Gc}$ from \eqref{eqn:R}, then we have:

\begin{theorem}\label{thm:E}
$\mE_{\qb}(\bs\lambda)\neq 0$ if and only if $\bs\lambda\in\mR$.
\end{theorem}

To prove this theorem, we need a series of results, see \S\ref{sec:proof} for the proof. 

\begin{cor}\label{cor:A}
For every $\bs\lambda\in\mR$, $\mA_k(\bs\lambda)\coloneqq\mE_k(\bs\lambda)\# H$ is a right $\B_k\#H$-cleft object. In particular, 
\[
\mA(\bs\lambda)\coloneqq\mE_{\qb}(\bs\lambda)\# H\in\Cleft(\B(V)\# H), \quad \forall\,\bs\lambda\in\mR.
\]
\end{cor}
\pf
Follows by \cite[Proposition 5.8]{AAGMV}.
\epf

The following result will be useful for the recursive proof in Theorem \ref{thm:E}. 

\begin{lem}\label{lem:nonzero-k<l}
Let $k\leq \ell$, $\bs\lambda\in\mR$. Assume that the following conditions hold:
\begin{align*}
\mE_{k-1}/\lg \gamma_{k-1}(r) : r\in \Gc_{k-1}\rg&\neq 0, & \B_{k-1}/\lg r-\lambda_r : r\in \Gc_{k-1}\rg&\neq 0.
\end{align*}
Then  $\mE_k(\bs\lambda)\neq 0$.
\end{lem}

In other words, it says that if for the $k-1$-th step we have proved that all the algebras $\mE_{k-1}$ are non-zero and we are able to extend to the $k$-th step in the following two cases:
\begin{itemize}
	\item $\lambda_r$ arbitrary for $r\in \Gc_{t}$, $t<k-1$, and $\lambda_r=0$ for $r\in \Gc_{k-1}$,
	\item $\lambda_r=0$ for $r\in \Gc_{t}$, $t<k-1$, and $\lambda_r$ arbitrary for $r\in \Gc_{k-1}$,
\end{itemize}
then $\mE_k(\bs\lambda)\neq 0$ for all $\bs\lambda$.

\pf
We consider the algebra map $\varrho:\mE_{k-1}\to \mE_{k-1}\ot \B_{k-1}$ given by the braided coaction as in \cite[Lemma 4.1]{AG}. As $\Gc_{k-1}\subset \mP(\B_{k-1})$, 
\[\varrho(\gamma_{k-1}(r))=\gamma_{k-1}(r)\ot 1+1\ot r, \quad r\in\Gc_{k-1}.\]
Now $\varrho$ induces an algebra map 
\[
\varrho':\mE_{k-1}\to \mE_{k-1}/\lg \gamma_{k-1}(r) : r\in \Gc_{k-1}\rg \ot \B_{k-1}/\lg r-\lambda_r : r\in \Gc_{k-1}\rg
\]
such that $\varrho'(\gamma_{k-1}(r)-\lambda_r)=0$, $r\in\Gc_{k-1}$. Hence $\mE_k(\bs\lambda)\neq 0$.
\epf

\begin{lem}\label{lem:nonzero-k=l+1}
If $\mtE(\bs\lambda)\coloneqq \mE_\ell(\bs\lambda)\neq 0$, then $\mE_{\ell+1}(\bs\lambda)\neq 0$.
\end{lem}
\pf
We recall that the normal Hopf subalgebra $\k\lg\Gc_{\ell}\rg\subset \B_\ell$ is a $q$-poly\-nomial algebra \cite[Proposition 21]{A-distinguished}. Hence, by \cite[Theorem 3.1]{AAGMV}, also \cite[Theorem 4]{Gu}, $\mE_{\ell+1}(\bs\lambda)\neq 0$.
\epf

We fix $\D=\bigsqcup_{i=1}^m \D_i$ and respectively $V=\bigoplus_{i=1}^m V_i$ the decompositions of $\D$  and $V$ into connected components.
We denote by $\Gc(i)$ the set of generators of the ideal defining the Nichols algebra $\B(V_i)$, $i\in \I_m$. As a result, the set $\Gc$ of generators of $\J(V)$, see \S \ref{sec:diagonal}, decomposes as a disjoint union, cf.~\eqref{eqn:G(0)}:
\begin{equation}\label{eqn:Gc=union}
\Gc=\Gc(0)\sqcup\Gc(1)\sqcup \dots \sqcup \Gc(m).
\end{equation}

Let $\mR\subseteq\k^{\Gc}$, resp.~$\mR(i)\subseteq \k^{\Gc(i)}$, $i\in \I_m$,  be the set of deformation parameters for $V$, resp.~$V_i$, cf. \eqref{eqn:R}. Then we have a restriction map:
\begin{equation}\label{eqn:restriction}
\mR\to \mR(i), \qquad \bs\lambda\mapsto \bs\lambda_{|i}\coloneqq (\lambda_r)_{r\in\Gc(i)}\in\mR(i).
\end{equation}
In particular, $\bs\lambda$ splits into a product of subfamilies
\begin{align}\label{eqn:lambda-decomp}
\bs\lambda=\bs\lambda_{|0}\times \bs\lambda_{|1}\times\dots\times\bs\lambda_{|m}
\end{align}
in such a way that
\begin{itemize}
\item[(i)] $\bs\lambda_{|0}\in\k^{\Gc(0)}$ is the subset of deformation parameters linking vertices in different connected components;
\item[(ii)] $\bs\lambda_{|i}\in\k^ {\Gc(i)}$ is the subset of deformation parameters for the $i$th connected component of $V$, $i\in\I_m$.
\end{itemize}

\begin{lem}\label{lem:lambda'}
Let $\bs\lambda=(\lambda_r)_{r\in\Gc}\in\k$ be a family of scalars satisfying \eqref{eqn:condition-lambda-k}. 
If $\mE(\bs\lambda_{|i})\neq 0$ for every $i\in\I_m$ and $\bs\lambda_{|0}={\bf 0}$, then $\mtE(\bs\lambda)\neq 0$.
\end{lem}
\pf
Let $\rho:T(V)\to \End \otimes_{t=1}^m \mE(\bs\lambda_{|t})$ be the algebra map defined by
\begin{align*}
\rho(y_k)(\by_1\ot \dots \ot \by_m):= g_k\cdot\by_1\ot \dots \ot 
g_k\cdot\by_{i-1} \ot y_k\by_i \ot \by_{i+1}\ot \dots \ot \by_m,
\end{align*}
where $\by_t\in \mE(\bs\lambda_{|t})$, $k\in\D_i$, $i\in\I_m$.
We claim that $\rho(y_{k\ell})=0$ for all $k\not\sim \ell$, $k<\ell$. Indeed, if $k\in\D_i$, $\ell\in\D_j$,  $i<j\in\I_m$, then
\begin{align*}
\rho(y_k y_\ell)(\by_1\ot & \dots \ot \by_m):= g_k g_\ell \cdot\by_1\ot \dots \ot 
g_k g_\ell \cdot\by_{i-1} \ot y_k(g_\ell \cdot \by_i) \\ 
& \quad \ot g_\ell \cdot \by_{i+1}\ot \cdots \ot
g_\ell \cdot \by_{j-1} \ot y_\ell \by_j \ot \by_{j+1}\ot \dots \ot \by_m, \\
\rho(y_\ell y_k )(\by_1\ot & \dots \ot \by_m):= g_k g_\ell \cdot\by_1\ot \dots \ot 
g_k g_\ell \cdot\by_{i-1} \ot g_\ell \cdot( y_k \by_i) \\ 
& \quad \ot g_\ell \cdot \by_{i+1}\ot \cdots \ot
g_\ell \cdot \by_{j-1} \ot y_\ell \by_j \ot \by_{j+1}\ot \dots \ot \by_m.
\end{align*}
As $k\not\sim \ell$, $q_{\ell k}q_{k\ell}=1$ and hence
\begin{align*}
g_\ell \cdot( y_k \by_i) & =  ( g_\ell \cdot y_k) (g_\ell \cdot \by_i)
= q_{\ell k} y_k (g_\ell \cdot \by_i) = q_{k\ell}^{-1} y_k (g_\ell \cdot \by_i),
\end{align*}
so $\rho(y_{k\ell})(\by_1\ot \dots \ot \by_m)=0$ for all $\by_t\in \mE(\bs\lambda_{|t})$.

Let $r\in\Gc(p)$. We claim that $\rho(\gamma_{p-1}(r))(\by_1\ot \dots \ot \by_m)=0$.
It suffices to consider $\by_t\in \mE(\bs\lambda_{|t})$ such that
$\by_t \in\pi_t(T(V_t)_{{\bf g}_t}^{\bs\chi_t} )$ for some ${\bf g}_t\in\Gamma$, $\bs\chi_t\in\widehat{\Gamma}$, where $\pi_t:T(V_t)\to \mE(\bs\lambda_{|t})$ is the canonical projection. For each $k\in\D_i$,
\begin{multline*}
\rho(y_k)(\by_1\ot \dots \ot \by_m)= \bs\chi_1\dots \bs\chi_{i-1} (g_k) 
\by_1\ot \dots \ot \by_{i-1} \ot y_k\by_i \ot \by_{i+1}\ot \dots \ot \by_m \\
= \chi_k^{-1} ({\bf g}_1\dots {\bf g}_{i-1}) 
\by_1\ot \dots \ot \by_{i-1} \ot y_k\by_i \ot \by_{i+1}\ot \dots \ot \by_m.
\end{multline*}
Since $\gamma_{p-1}(r)\in \mE(\bs\lambda_{|p})^{\chi_r}$, we have 
\begin{multline*}
\rho(\gamma_{p-1}(r))(\by_1\ot \dots \ot \by_m)= \\
\chi_r^{-1} ({\bf g}_1\dots {\bf g}_{i-1}) 
\by_1\ot \dots \ot \by_{i-1} \ot \gamma_{p-1}(r) \by_i \ot \by_{i+1}\ot \dots \ot \by_m=0.
\end{multline*}
Hence $\rho$ induces an algebra map $\widetilde{\rho}:\mtE(\bs\lambda)\to \End \otimes_{t=1}^m \mE(\bs\lambda_{|t})$ and therefore $\mtE(\bs\lambda)\neq 0$.
\epf

\begin{lem}\label{lem:lambda''}
If $\bs\lambda_{|i}={\bf 0}$ for all $i>0$, then $\mtE(\bs\lambda)\neq 0$.
\end{lem}
\pf
Set $N_j=\ord \chi_j(g_j)$, $y_{jk}\coloneqq y_jy_k-q_{jk}y_ky_j$, $j<k\in\I$. Let $\mF(\bs\lambda_{|0})$ be the quotient of $T(V)$ by the ideal generated by
\begin{align*}
& y_{jk}-\lambda_{jk}, \ j\not\sim k, & & y_{jk}, \ j\sim k, j<k, && y_j^{N_j}, \ j\in \I.
\end{align*}
An easy application of the Diamond Lemma shows that $\mF(\bs\lambda_{|0})\neq 0$.

For each connected component $V_i$, the relations in $\Gc(i)$ hold in $\mtE$, since $\bs\lambda_{|i}={\bf 0}$, cf.~\eqref{eqn:E(lambda)}. Recall that the set $\Gc(i)$ consists of relations $y_j^{N_j}$, for certain $j\in \I_{\theta_i}$, and relations involving at least two letters $y_j, y_k$ contained in the ideal $\lg y_{jk} | j<k\in \I_{\theta_i}\rg$. Hence $\mtE(\bs\lambda)$ projects onto $\mF(\bs\lambda_{|0})$.
\epf

\begin{lem}\label{lem:cocycle-defs-nc}
Let $\bs\lambda=(\lambda_r)_{r\in\Gc}\in\k$ be a family of scalars satisfying \eqref{eqn:condition-lambda-k}. 
If $\mE(\bs\lambda_{|i})\neq 0$ for every $i\in\I_m$, then $\mE(\bs\lambda)\neq 0$.
\end{lem}
\pf
The case  $\bs\lambda_{|0}={\bf 0}$ follows by Lemma \ref{lem:lambda'}, while the case $\bs\lambda_{|i}={\bf 0}$ for all $i>0$ is Lemma \ref{lem:lambda''}. Hence we  make use of Lemma \ref{lem:nonzero-k<l} to combine both cases and conclude that $\mE(\bs\lambda)\neq 0$ for any $\bs\lambda$.
\epf

\subsection{Proof of Theorem \ref{thm:E}}\label{sec:proof}

We show that $\mE_k(\bs\lambda)\neq 0$, for every $k\geq 0$.
We may restrict to the connected case by Lemma \ref{lem:cocycle-defs-nc}. Moreover, we just need to check that $\mtE(\bs\lambda)\neq 0$, by Lemma \ref{lem:nonzero-k=l+1}. This is done in the Appendix, by a case-by-case analysis of all Nichols algebras following the explicit list given in \cite{AA}. We explain how to define $\mtE(\bs\lambda)$ recursively and show that $\mtE(\bs\lambda)\neq 0$ for each $\qb$.
\qed

\subsection{The liftings $\u(\bs\lambda)$}\label{sec:u(lambda)}
We review \cite[\S 3.1]{AAG} where the liftings are defined. 

We fix $\bs\lambda\in\mR$ as in \eqref{eqn:R}: this defines a right $\B(V)\#H$-cleft object $\mA(\bs\lambda)$, by Corollary \ref{cor:A}. We consider the left Schauenburg Hopf algebra $\mL(\bs\lambda)\coloneqq L(\mA(\bs\lambda),\B(V)\#H)$. This is a cocycle deformation of $\B(V)\#H$ for which $\mA(\bs\lambda)$ becomes a $(\B(V)\#H,\mL(\bs\lambda)$-bicleft object. 

We provide a recursive process to define a Hopf algebra $\u(\bs\lambda)$ such that $\u(\bs\lambda)\simeq \mL(\bs\lambda)$. Our aim is to get a quotient $T(V)\#H\twoheadrightarrow\u(\bs\lambda)$, that is, a definition of the liftings by generators and relations.

We will proceed inductively by considering intermediate quotients $\u_k(\bs\lambda)\simeq \mL_k(\bs\lambda)\coloneqq L(\mH_k,\mA_k(\bs\lambda))$, recall that $\mH_k=\B_k\#H$.

Recall that $\mE_0(\bs\lambda)=T(V)$; hence $\mA_0(\bs\lambda)=T(V)\# H$ and thus we set $\u_0(\bs\lambda)=T(V)\#H\simeq \mL_0(\bs\lambda)$.

Fix $k>0$ and assume that $\u_{k-1}(\bs\lambda)$ is already defined. In particular, there is a left coaction $\mA_{k-1}(\bs\lambda)\to \u_{k-1}(\bs\lambda)\ot \mA_{k-1}(\bs\lambda)$. Let $\gamma:\mH_{k-1}\to\mA_{k-1}(\bs\lambda)$ be the section corresponding to the right coaction $\mA_{k-1}(\bs\lambda)\to \mA_{k-1}(\bs\lambda)\ot \mH_{k-1}$ and set, for each $r \in \Gc_{k-1}$,
\begin{align*}
\nabla(r) &:=\gamma(r)_{(-1)}\otimes \gamma(r)_{(0)}
 - g_r \otimes \gamma(r).
\end{align*}
\begin{cor}\cite[Corollary 5.12]{AAGMV}
$\nabla(r)\in \u_{k-1}(\bs\lambda) \otimes 1\subset \u_{k-1}(\bs\lambda) \otimes \mA(\bs\lambda)$.
\end{cor}
Hence this defines an element $ \tilde r\in \u_{k-1}(\bs\lambda)$ so that  $\nabla (r) = \tilde r \otimes 1$. We set
\begin{align}\label{eqn:Lk}
\u_k(\bs\lambda):= \u_{k-1}/\langle \tilde
r-\lambda_r(1-g_r): r\in \Gc_{k-1}\rangle.
\end{align}
By \cite[Proposition 3.3]{AAG} and Corollary \ref{cor:A}, $\u_k(\bs\lambda)\simeq L(\mA_k(\bs\lambda),\mH_k)$. Thus, $\u_k(\bs\lambda)$ is a quotient of $T(V)\#H=\mL_0(\bs\lambda)$, for each $k\geq 0$. We set, cf.~\eqref{eqn:u}:
\[
\u(\bs\lambda)\coloneqq \u_{\ell+1}(\bs\lambda).
\]
Our Theorem \ref{thm:E} readily implies the main result of the article, which in this context reads as follows.
\begin{theorem}\label{teo:lift}
Let $L$ be a lifting of $V$. Then there is $\bs\lambda\in\mR$ such that $L\simeq \u(\bs\lambda)$.
\end{theorem}
\pf
\cite[Theorem 3.5]{AAG} applies, by Theorem \ref{thm:E}.
\epf

\section{Diagonal case: non-connected diagrams}\label{sec:non-connected}

We fix a pair $(H,V)$ as in \S \ref{sec:diagonal}. We focus on the case in which the generalized Dynkin diagram $\D$  associated to $V$ is not connected. Our main result is Proposition \ref{pro:joint}, which shows that we can restrict to the connected case as the liftings for non-connected diagrams are defined via linking relations \eqref{eqn:linking} for the liftings of the connected components.

We fix $\{x_1,\dots,x_\theta\}$ a basis of $V$ with $x_i\in V_{g_i}^{\chi_i}$, $i\in\I$. In particular, recall that every $r\in\Gc$ belongs to a component $ T(V)_{g_r}^{\chi_r}$, for some $g_r\in\Gamma$, $\chi_r\in\Alg(H,\k)$.

We fix $\D=\bigsqcup_{i=1}^m \D_i$ and respectively $V=\bigoplus_{i=1}^m V_i$ the decompositions of $\D$  and $V$ into connected components.
We keep the notation  as in \S \ref{sec:cocycles}. Notice that \eqref{eqn:restriction} associates a lifting $\u(\bs\lambda_{|i})$ of $V_i$ to each lifting $\u(\bs\lambda)$ of $V$.
For each $i\in \I_m$, we set $\J(\bs\lambda_{|i})\subset T(V_i)$ as the ideal such that
\[\u(\bs\lambda_{|i})=T(V_i)\# H/\J(\bs\lambda_{|i}).\]
Let $\Gc(i)=\Gc(i)_0\sqcup \dots \sqcup \Gc(i)_{\ell(i)}$ be a stratification of the set of defining relations $\Gc(i)$ of $\B(V_i)$. Then
\begin{align}\label{eqn:strat-union}
\Gc=\Gc(1)_0\sqcup \dots \sqcup \Gc(1)_{\ell(1)}\sqcup\cdots \sqcup\Gc(m)_0\sqcup \dots \sqcup \Gc(m)_{\ell(m)}\sqcup\Gc(0)
\end{align}
is a stratification of $\Gc$ with $N+1$ steps, $N=r+\ell(1)+\dots+\ell(m)$.

We denote by $\u_{|i}(\bs\lambda)$ the subalgebra of $\u(\bs\lambda)$ generated by $H\oplus V_i\#H$; this is indeed a Hopf subalgebra.
\begin{lem}\label{lem:subalgebra}
$\u_{|i}(\bs\lambda)\simeq \u(\bs\lambda_{|i})$ as Hopf algebras.
\end{lem}
\pf
Let $\phi:T(V)\# H\to \u(\bs\lambda)$ be a lifting map, cf.~\S \ref{sec:prelim}. We assume that $i=1$, for short. We consider the map
$\phi_1:T(V_1)\# H\to \u(\bs\lambda)$ given by the composition of $\phi$ with the canonical inclusion $T(V_1)\# H\hookrightarrow T(V)\# H$. Using the stratification as in \eqref{eqn:strat-union}, it follows that $\phi_1(\J(\bs\lambda_{|1}))=0$ and thus $\phi_1$ factors through a map $\varphi:\u(\bs\lambda_{|1})\to \u(\bs\lambda)$.

By construction, the first term $\u(\bs\lambda_{|1})_1$ of the coradical filtration of $\u(\bs\lambda_{|1})$ is (isomorphic to) $H\oplus V_1\#H$; similarly $\u(\bs\lambda)_1\simeq H\oplus V\#H$. Now, the restriction of $\varphi$ to $\u(\bs\lambda_{|1})_1$ is the canonical inclusion  $H\oplus V_1\#H\hookrightarrow H\oplus V\#H$ since $\phi$ is a lifting map. Hence $\varphi$ is injective by \cite[Theorem 5.3.1]{Mo}.
\epf

Now we proceed in the opposite direction, and show in Proposition \ref{pro:joint} that we can combine the liftings of each of the $V_i$'s to find (all) liftings of $V$.

We denote by $\J(\bs\lambda_{>0})$ the ideal of $T(V)\#H $ generated by $\bigcup_{i\in\I_m}\J(\bs\lambda_{|i})$; here we consider the natural identification $T(V_i)\#H\hookrightarrow T(V)\#H$. We set
\[
\U(\bs\lambda)=T(V)\#H/\J(\bs\lambda_{>0}).
\]

The subfamily $\bs\lambda_{|0}$ consists of scalars $(\lambda_{i,j})_{i\not\sim j,i<j}$ subject to
\begin{align*}
\lambda_{i,j}=0 \qquad  \text{ if } \qquad \chi_{i}\chi_{j}\neq\eps \text{ or } g_ig_j=1,
\end{align*}
cf.~\eqref{eqn:condition-lambda-k}. We consider the Hopf ideal $\J(\bs\lambda_{|0})\subset \U(\bs\lambda)$ generated by
\begin{align}\label{eqn:linking}
 \ad(x_i)(x_j)-\lambda_{i,j}(1-g_ig_j),\qquad i<j,i\not\sim j
\end{align}
and set $\u'(\bs\lambda)=\U(\bs\lambda)/\J(\bs\lambda_{|0})$.

\begin{pro}\label{pro:joint}
$\u'(\bs\lambda)$ is a lifting of $V$; moreover, $\u'(\bs\lambda)\simeq \u(\bs\lambda)$.
\end{pro}
\pf
We fix a stratification as in \eqref{eqn:strat-union}. Then, for each step $k=0,\dots,N$, there is a surjective Hopf algebra map $\mL_k\twoheadrightarrow \U(\bs\lambda)$, cf.~\eqref{eqn:Lk}; hence there is a Hopf algebra projection $\tau':\mL_N\twoheadrightarrow \u'(\bs\lambda)$. Now, the final step $\mL_{N+1}=\u(\bs\lambda)$ is given precisely by the quotient of $\mL_N$ by the relations \eqref{eqn:linking}. Therefore, $\tau'$ factors through a map $\tau:\u(\bs\lambda)\twoheadrightarrow \u'(\bs\lambda)$; we shall denote by
$\tau_1$ the  restriction of $\tau$ to the first term $\u(\bs\lambda)_1\simeq H\oplus V\# H$ of the coradical filtration of $\u(\bs\lambda)$.

Conversely, we fix a lifting map $\phi:T(V)\# H\twoheadrightarrow \u(\bs\lambda)$, cf.~\S \ref{sec:prelim}. Now, by Lemma \ref{lem:subalgebra}, $\phi(\J(\bs\lambda_{|i}))=0$ for each $i\in\I_m$, and thus $\phi$ factors through a Hopf algebra epimorphism $\phi':\U(\bs\lambda)\twoheadrightarrow \u(\bs\lambda)$. We observe that $\Gc(0)$ is composed of primitive elements in $T(V)$, hence relations \eqref{eqn:linking} hold in $\u(\bs\lambda)$ by definition, see \eqref{eqn:Lk}. Hence $\phi'(\J(\bs\lambda_{|0}))=0$ and this induces a Hopf algebra projection $\varphi:\u'(\bs\lambda)\twoheadrightarrow \u(\bs\lambda)$ such that $\varphi_1\coloneqq\varphi_{|H\oplus V\#H}$ is injective, by \eqref{eqn:lifting-map}.

It follows that $\tau_1\circ\varphi_1=\varphi_1\circ \tau_1=\id_{H\oplus V\#H}$, hence $\tau:\u(\bs\lambda)\to\u'(\bs\lambda)$ defines a Hopf algebra isomorphism, with inverse $\varphi$.
\epf

\section{On isomorphism classes}\label{sec:isos}

Let $(V,H)$ be as in \S \ref{sec:diagonal}.  If $\bs\lambda\in\mR$, cf.~\eqref{eqn:R}, then we say that $(H,V,\bs\lambda)$ is a {\it lifting datum}.
Let $(V',H')$ a new pair as in \S \ref{sec:diagonal}; we denote by $((g'_i,\chi'_i))_{i\in\I_{\theta'}}$, $\theta'=\dim V'$, $\mR'$, the corresponding input information.

When $V,V'$ are connected, \cite[Theorem 3.9]{AAG} describes the set of isomorphisms classes $\Isom(\u(\bs\lambda),\u(\bs\lambda'))$ in terms of certain symmetries of the lifting data.  This is in turn inspired by \cite[Theorem 7.2]{AS-annals}, where this is done for certain families of liftings of Cartan type with minor restrictions on the parameters; see also \cite[Theorem 7.1]{AS} for some families of  type $A$.

The claim of \cite[Theorem 3.9]{AAG} extends verbatim to the non-connected case, see Theorem \ref{thm:isos}; also Remark \ref{rem:hypothesis}.
In addition, we review this case in the spirit of Proposition \ref{pro:joint}, that is we study the isomorphism classes of the algebras $\u(\bs\lambda)$ in terms of the components of $V$, see Proposition \ref{pro:iso}.

Consider the  decompositions $V=\bigoplus_{i=1}^m V_i$ and $V'=\bigoplus_{i=1}^{m'} V'_i$ of $V$ and $V'$
as sums of connected braided vector subspaces; set $\theta_i=\dim V_i$, $\theta'_j=\dim V_j'$, $i\in\I_m, j\in \I_{m'}$. We fix $\I=\I_\theta$ and set $\I(i)=\{j\in\I | g_j=g_i\text{ and }\chi_j=\chi_i\}$, for each $i\in \I$ and consider the groups
\begin{align*}
\Sym_{\qb}&=\{\sigma\in\Sym_\theta| q_{ij} = q_{\sigma(i) \sigma(j)} \,\forall i, j \in \I \};\\
\Lb&=\{s\in \GL_\theta(\k)|s_{ij}=0 \text{ if } j\notin\I(i)\}.
\end{align*}
These groups act on $\mR$ \cite[Lemmas 3.10 \& 3.11]{AAG}; we write each action as
\begin{align}\label{eqn:actions}
\bs\lambda&\mapsto \bs\lambda^\sigma, & \bs\lambda&\mapsto s\cdot \bs\lambda; & \text{for }\bs\lambda\in\mR, \ \sigma\in \Sym_{\qb}, \ s\in \Lb.
\end{align}

\begin{rem}\label{rem:action}
Let $\bs\lambda\in\mR$. Then the actions \eqref{eqn:actions} restrict to each component $\bs\lambda_{|i}$, $i\geq 0$, as in \eqref{eqn:lambda-decomp}, cf.~\cite[Lemmas 3.10 \& 3.11]{AAG}, also \cite[Example 3.12]{AAG}. In particular,
\begin{align}\label{eqn:actions-0}
(s\cdot \bs\lambda^\sigma)_{|0}=(s_{ij}^{-1}\lambda_{\sigma^{-1}(i),\sigma^{-1}(j)})_{i,j}, \qquad \sigma\in \Sym_{\qb}, s\in \Lb.
\end{align}
Here, if $\bs\lambda_{|0}=(\lambda_r)_{r\in\Gc(0)}$, then we set $\lambda_{i,j}\coloneqq\lambda_r$ for $i<j\in\I$ defining the q-commutator $r$ cf.~\eqref{eqn:G(0)}.
\end{rem}

\begin{rem}\label{rem:hypothesis}
When $V$ connected, the sets $\I(i)$ and $\Lb$ are explicitly described  in \cite[\S 3.2.1]{AAG}. This is indeed the only occurrence of the connectedness hypothesis in loc.cit.
\end{rem}

We recall from \cite[Lemma 3.7]{AAG} that if $\psi\in\Isom(\u(\bs\lambda), \u(\bs\lambda'))$, then $\varphi\coloneqq\psi_{|H}\in\Isom(H,H')$ and $T\coloneqq \psi_{|V}\colon V\to V'$ is an isomorphism of braided vector spaces; so $\theta=\theta'$, $m=m'$ and there is $s=(s_{ij})\in  \Lb$ with
\begin{align*}
T(a_i)=\sum\limits_{j\in\I(\sigma(i))} s_{ij}a'_j, \ i\in\I,
\end{align*}
where we denote by $\{a_i\}_{i\in\I}$, resp. $\{a'_i\}_{i\in\I}$, the images of the elements $x_i\in V$, resp. $x_i'\in V'$, $i\in\I$, via the natural projections $T(V)\#H\twoheadrightarrow \u(\bs\lambda)$, resp. $T(V')\#H'\twoheadrightarrow \u(\bs\lambda')$.
In turn, $\varphi$ defines an element $\sigma\in\Sym_{\qb}$ by the identities
\begin{align}\label{eqn:sigma}
\varphi(g_{j})&=g'_{\sigma(j)} \quad \text{ and }\quad \chi'_{\sigma(j)}\circ \varphi=\chi_j, \ j\in\I.
\end{align}
These data actually determine the isomorphism $\psi$, see Theorem \ref{thm:isos}. We let $\Isom(\bs\lambda,\bs\lambda')$ be the set of all triples $(\varphi,\sigma,s)\in\Isom(H,H')\times \Sym_{\qb}\times \Lb$ satisfying \eqref{eqn:sigma} and such that $\bs\lambda'=s\cdot\bs\lambda^\sigma$.

\begin{theorem}\cite[Theorem 3.9]{AAG}\label{thm:isos}
Let $(H,V,\bs\lambda)$ and $(H',V',\bs\lambda')$ be two lifting data. Then
$\Isom(\u(\bs\lambda),\u(\bs\lambda'))\simeq \Isom (\bs\lambda,\bs\lambda').$\qed
\end{theorem}

\begin{rem}\label{rem:not-connected}
If $s\in\Lb$ and $i\not\sim j\in\I$, then one of the following holds:
\begin{itemize}
\item[(a)]$s_{ij}=0$;
\item[(b)] $\k\{x_i\}\simeq\k\{x_j\}\simeq V_p$ for some $p\in\I_m$ and $q_{ii}=q_{ij}=q_{ji}=q_{jj}=-1$.
\end{itemize}
Indeed, given such $i,j$, if there is $k\neq i,j$ connected with one of them, say $k\sim i$, and $s_{ij}\neq 0$, then $j\in\I(i)$, that is $\chi_j=\chi_i$ and $g_j=g_i$. Now, as $k\not\sim j$:
\[
1=\chi_j(g_k)\chi_k(g_j)=\chi_i(g_k)\chi_k(g_i),
\]
which contradicts the fact that $k\sim i$. Hence $s_{ij}=0$. If, on the contrary, both $x_i$ and $x_j$ determine a one-dimensional connected component and $s_{ij}\neq 0$, then $j\in\I(i)$ implies $q_{ii}=q_{ij}=q_{ji}=q_{jj}=-1$.
\qed
\end{rem}

We recall from p.~\pageref{lem:subalgebra} that $\u_{|i}(\bs\lambda)\subseteq \u(\bs\lambda)$ stands for the subalgebra generated by $V_i$ and $H$. The following is straightforward.

\begin{lem}\label{lem:iso}
Let $\psi\in\Isom(\u(\bs\lambda), \u(\bs\lambda'))$. Then  there is $\tau\in\Sym_m$ such that $T_i\coloneqq T_{|V_i}:V_i\to V'_{\tau(i)}$ is an isomorphism of braided vector spaces for each $i\in \I_m$. Moreover, \[\psi_{|i}\coloneqq \psi_{|\u_{|i}(\bs\lambda)}:\u_{|i}(\bs\lambda)\to \u_{|\tau(i)}(\bs\lambda')\] is a Hopf algebra isomorphism.
\qed
\end{lem}
By Lemma \ref{lem:subalgebra}, $\u_{|i}(\bs\lambda)\simeq \u(\bs\lambda_{|i})$ and hence
each $\psi\in\Isom(\u(\bs\lambda), \u(\bs\lambda'))$ induces isomorphisms $\u(\bs\lambda_{|i})\to \u(\bs\lambda_{|\tau(i)}')$ for all $i\in \I_m$. We denote by
\[\Isom(\bs\lambda_{|1},\bs\lambda_{|\tau(1)}')\times_H \dots\times_H \Isom(\bs\lambda_{|m},\bs\lambda_{|\tau(m)}')\subset \bigtimes_{i=1}^m\Isom(\bs\lambda_{|i},\bs\lambda_{|\tau(i)}')
\]
the subset of all $m$-tuples $\left((\varphi_i,\sigma_i,s_i)\right)_{i\in\I_m}$ with $\varphi_i=\varphi_j$, $i,j\in\I_m$. Therefore, Lemma \ref{lem:iso} defines a map
\begin{align}\label{eqn:times-H}
\Res\colon\Isom(\bs\lambda,\bs\lambda')&\longrightarrow \Isom(\bs\lambda_{|1},\bs\lambda_{|\tau(1)}')\times_H \dots\times_H \Isom(\bs\lambda_{|m},\bs\lambda_{|\tau(m)}')
\end{align}

\begin{lem}\label{lem:iso-inj}
The map $\Res$ is injective.
\end{lem}
\pf
Let $(\varphi,\sigma,s)\in \Isom(\bs\lambda,\bs\lambda')$ and denote by $\psi:\u(\bs\lambda)\to \u(\bs\lambda')$ the corresponding Hopf algebra isomorphism from Theorem \ref{thm:isos}.
The collection of maps $(\psi_{|i})_{i\in\I_m}$ controls the behavior of $\psi$ when restricted to the first term of the coradical filtration of $\u(\bs\lambda)$ and hence determines the morphism $\psi$.
\epf
To describe the image of the map $\Res$, we set, for each $i\in\I_m$,
\[
\Isom_0(\bs\lambda_{|i},\bs\lambda_{|\tau(i)}')=\{(\varphi,\sigma,s)\in \Isom(\bs\lambda_{|i},\bs\lambda_{|\tau(i)}')|\lambda_0'=s\cdot \lambda_0^{\sigma_i}\}
\]
and define $\Isom_0(\bs\lambda_{|1},\bs\lambda_{|\tau(1)}')\times_H \dots\times_H \Isom_0(\bs\lambda_{|m},\bs\lambda_{|\tau(m)}')$ accordingly.
By definition, the image of $\Res$ lies inside this subset.
\begin{pro}\label{pro:iso}
The map $\Res$ defines a bijective correspondence
\begin{align}\label{eqn:times-H-0}
\Isom(\bs\lambda,\bs\lambda')\overset{\sim}{\longrightarrow}\Isom_0(\bs\lambda_{|1},\bs\lambda_{|\tau(1)}')\times_H \dots\times_H \Isom_0(\bs\lambda_{|m},\bs\lambda_{|\tau(m)}').
\end{align}
\end{pro}
\pf
First, $\Res$ is injective by Lemma \ref{lem:iso-inj}. Now, let $\left((\varphi,\sigma_i,s_i)\right)_{i\in\I_m}$ belong to the right hand side of \eqref{eqn:times-H-0}, for some  $\varphi\in\Isom(H,H')$. Then, via the natural identification $\I_{\theta_1}\times \dots \times \I_{\theta_m}=\I_\theta$:
\begin{itemize}
\item[(a)] $\sigma=\sigma_1\times\dots\times\sigma_m\in\Sym_{\qb}$.
\item[(b)] The block matrix
$
s=\left[\begin{smallmatrix}
s_1 & 0 & \cdots &0 \\
0& s_2 & & 0  \\
0 & 0 & \ddots & \vdots \\
0 & 0 & \cdots  & s_m
\end{smallmatrix}\right]$ defines an element in $\Lb$.

\end{itemize}
Indeed, (a) is clear from the definition and (b) follows by Remark \ref{rem:not-connected}. Hence, $(\varphi,\sigma,s)\in\Isom (\bs\lambda,\bs\lambda')$ and it is mapped onto $\left((\varphi,\sigma_i,s_i)\right)_{i\in\I_m}$ by $\Res$.
\epf

\appendix

\section{Nonzero cleft objects}

In this section we show that each algebra $\toba=\toba(\bs\lambda)$ is nonzero for each connected matrix $\bq$ and each $\bs\lambda\in\mR$, see  \eqref{eqn:R}. This completes the proof of Theorem \ref{thm:E}.
We refer to \S\ref{sec:presentation} for the notation $x_{ij}$, etc. We use the corresponding notation $y_{ij}$, etc., for the analogous elements in terms of the generators $y_i$, $i\in\I$. See also \cite[\S 2.2]{AA}.

\subsection{The algebras $\mtE$}

\subsubsection{} To define $\widetilde{\mathcal{E}}_{\bq}(\bs\lambda)$, we follow \eqref{eqn:E(lambda)}, by deforming the corresponding defining relations of $\B_{\qb}$, as computed in Theorem \ref{thm:presentacion minima} from \cite{A-nichols}, for each diagram in the classification of \cite{H-classif}. This is a case-by-case analysis of all Nichols algebras following the explicit list given in \cite{AA}. 

Namely, given the presentation of a Nichols algebra $\B_{\qb}$ in \cite{AA} by a set of relations $\Gc$, we keep only those relations which are not powers of Cartan root vectors (which is, in a generic way, a presentation of the distinguished pre-Nichols algebra $\widetilde{\B}_{\qb}$, see Remark \ref{rem:distinguished}), we stratify this set $\widetilde{\Gc}$ and inductively deform each relation according to \eqref{eqn:E(lambda)}.

\begin{rem}
As a result, we reproduce here the (truncated, for our pourposes) presentation of each algebra $\B_{\qb}$. Moreover, each one of these presentations is preceded by the Dynkin diagram attached to $\qb$, and hence we also reproduce the list in \cite{H-classif}. 
\end{rem}

\subsection{Tools and techniques}
For each $\bq$ we show that $\toba \neq0$ using a combination of the following results and resources. For each one of these tools, we include a reference to a specific example where it is used. In this example, the technique is explained in detail, while subsequent uses of the same techinque are simply mentioned, without further explanation.

\subsubsection{}\label{subsubsec:gap} For diagrams of small rank we explicitly compute the Gr\"obner basis of the defining ideal of $\toba$. This is achieved using the computer program \cite{GAP}, together with the package \texttt{GBNP} by \cite{CG}. The \texttt{log} files are stored in the authors' webpages: 
\texttt{http://www.famaf.unc.edu.ar/$\sim$(angiono|aigarcia)}.

This is the approach in the last line of Example \ref{subsubsec:type-B-N>4}, where a Gr\"obner basis associated to Cartan type $B_2$, $N=5$, is presented.

\subsubsection{}\label{subsubsec:oldresults} We invoke some results from previous articles. The case in which only the powers of simple roots are deformed follows from \cite[Lemma 5.16]{AAGMV}. For instance, the second case in \ref{subsec:type-G-st} follows directly from that result. 

On the other hand, in \cite[Proposition 3.2]{AGI} we determined the matrices $\qb$ for which the  quantum Serre relations admit  a nontrivial deformation. In many cases, this restricts our computations to a single matrix $\qb$; see for instance \ref{subsubsec:type-B-N>4}.

\subsubsection{}\label{subsubsec:cut} 
Let $i<j\in \I$. 
Assume that $\widetilde{q_{ij}}\neq 1$ (i.e. $i$ and $j$ are adjacent vertices). By looking at the lists in \cite{H-classif}, 
we have two cases: 
\begin{enumerate}[leftmargin=*]
\item either $(\widetilde{q_{ik}}-1) (\widetilde{q_{jk}}-1)=0$ (that is, $k$ is not adjacent to both $i$ and $j$) for all $k\neq i,j$ , or else
\item there exists a unique $k\neq i,j$ such that $(\widetilde{q_{ik}}-1) (\widetilde{q_{jk}}-1)\neq 0$ (that is, $k$ is the unique vertex adjacent to both $i$ and $j$).
\end{enumerate}

Let $\Gc[ij]\subset \Gc$ be the subset of generalized quantum Serre relations containing both letters $x_i$ and $x_j$.

We analize the two possible cases separately in the next two lemmas. 
In both cases we conclude that $\toba(\bs\lambda)\neq 0$ since it projects onto $\widetilde{\mathcal{E}}_{\bq'}(\bs\lambda')$, where $\bq'$ is a new matrix constructed from $\bq$ such that $i$ and $j$ belong to different connected components, and $\bs\lambda'$ is a new set of parameters for $\bq'$ constructed from $\bs\lambda$; as $\bq'$ has exactly two connected components, we may argue recursively using Lemma \ref{lem:lambda'}.

\begin{lem}\label{lem:cut1}
Assume that there is no $k\neq i,j\in \I$ adjacent to both $i$ and $j$. Let $\bq'=(q_{st}')_{s,t\in\I}$ be the matrix 
\begin{align*}
q_{st}'=\begin{cases}
q_{st}, & (s,t)\neq(j,i),\\
q_{ij}^{-1}, & (s,t)=(j,i).
\end{cases}
\end{align*}
Then
\begin{enumerate}[leftmargin=*]
\item\label{item:cut1-1} (The diagram associated to) $\bq'$ has two connected components, which arise by erasing the edge between $i$ and $j$. 
\item\label{item:cut1-2} The set of defining relations of $\B_{\bq'}$ is $\Gc'=\left(\Gc\setminus \Gc[ij]\right)\cup \{x_{ij}\}$.
\item\label{item:cut1-3} 
Assume that all relations $r\in \Gc[ij]$ are not deformed 
(in particular, $\lambda_r=0$ for all $r\in \Gc[ij]$).
Then $\toba(\bs\lambda)$ projects onto $\widetilde{\mathcal{E}}_{\bq'}(\bs\lambda')$ for $\bs\lambda'=(\lambda'_r)_{r\in \Gc'}$ with $\lambda'_r=\lambda_r $ if $r\in\Gc \setminus\Gc[ij]$ and $\lambda'_{ij}=0$. 
\end{enumerate}
\end{lem}

\pf
For \eqref{item:cut1-1}, we notice that all the vertices of the Dyinkin diagram of $\bq'$ and all the edges different from the one between $i$ and $j$ coincide with those of $\bq$ (with the same labels), while $\bq'$ has no edge between $i$ and $j$ and thus, by inspection, the corresponding diagram has two connected components: one containing $i$ and one containing $j$.

Now \eqref{item:cut1-2} and \eqref{item:cut1-3} follow from \eqref{item:cut1-1} and Theorem \ref{thm:presentacion minima}: Each relation $r\in\Gc$ depends on the subdiagram attached to the vertices of the letters appearing in $r$. If $r$ does not contain $x_i$ and $x_j$ simultaneously (that is, $r\notin\Gc[ij]$), then the subdiagram of letters of $r$ is contained in the Dynkin diagram of $\bq'$ and then $r\in\Gc'$; moreover, $r$ belongs to the ideal generated by $x_{ij}$. On the other hand, if $x_i$, $x_j$ appear in $r$ (that is, $r\in\Gc[ij]$), then the associated subdiagram is not contained in the Dynkin diagram of $\bq'$ and thus $r\notin \Gc'$; the relation $r$ is written as an iteration of braided brackets, all of them containing $x_{ij}$ as we can check in Theorem \ref{thm:presentacion minima}, and thus $r$ belongs to the ideal generated by $x_{ij}$. Finally $x_{ij} \in\Gc'$ since $\widetilde{q'_{ij}}=1$ by definition.
\epf

As an example, we apply this result for type $\superf$, see \S \ref{subsec:type-F-super}. Let $N=\ord q$. 
For the diagram $\xymatrix{ \overset{\,\, q^2}{\underset{\ }{\circ}}\ar  @{-}[r]^{q^{-2}}  & \overset{\,\,
q^2}{\underset{\ }{\circ}} \ar  @{-}[r]^{q^{-2}}  & \overset{q}{\underset{\ }{\circ}} & \overset{-1}{\underset{\ }{\circ}} \ar  @{-}[l]_{q^{-1}}}$, there are 2 cases: $N>4$ and $N=4$. When $N>4$, $\toba$ is generated by $y_i$, $i\in\I_4$, with 
relations
\begin{align*}
y_{13}&=0; & y_{14}&=0; & y_{24}&=0; \\
y_{112}&=0; & y_{221}&=0; & y_{223}&=\lambda_{223};\\
y_{334}&=0; & y_{3332}&=\lambda_{3332}; & y_{4}^2&=\mu_4.
\end{align*}
Here, $\Gc[34]=\{x_{334}\}$ and this relation is not deformed. Hence we consider the corresponding matrix $\bq'$ as in the Lemma, which has diagram $\xymatrix{ \overset{\,\, q^2}{\underset{\ }{\circ}}\ar  @{-}[r]^{q^{-2}}  & \overset{\,\,q^2}{\underset{\ }{\circ}} \ar  @{-}[r]^{q^{-2}}  & \overset{q}{\underset{\ }{\circ}} & \overset{-1}{\underset{\ }{\circ}}}$. Thus $\bq'$ has two components, one of type $B_{3}$ and another of type $A_1$, and the defining relations are the same as above but we put $y_{334}$ in place of $y_{34}$. By a recursive argument on the rank of connected components, we may assume that $\widetilde{\mathcal E}_{\bq'}\neq 0$. As $\toba$ projects onto $\widetilde{\mathcal E}_{\bq'}$ by the previous result, we conclude that $\toba\neq 0$.

The case $N=4$ is solved analogously.

\begin{lem}\label{lem:cut2}
Assume that there is a vertex $k\neq i,j\in \I$ adjacent to both $i$ and $j$. Let $\bq'=(q_{st}')_{s,t\in\I}$ be the matrix 
\begin{align*}
q_{st}'=\begin{cases}
q_{st}, & (s,t)\neq(j,i), (k,i),\\
q_{it}^{-1}, & s=i, \, t=j,k.
\end{cases}
\end{align*}
Then
\begin{enumerate}[leftmargin=*]
\item\label{item:cut2-1} (The diagram associated to) $\bq'$ has two connected components, which arise by erasing the edges between $i$ and $j$, and between $i$ and $k$. 
\item\label{item:cut2-2} The set of defining relations of $\B_{\bq'}$ is $\Gc'=\left(\Gc\setminus (\Gc[ij] \cup \Gc[ik])\right)\cup \{x_{ij}\}$.
\item\label{item:cut2-3} 
Assume that all relations $r\in \Gc[ij]\cup \Gc[ik]$ are not deformed.
Then $\toba(\bs\lambda)$ projects onto $\widetilde{\mathcal{E}}_{\bq'}(\bs\lambda')$ for $\bs\lambda'=(\lambda'_r)_{r\in \Gc'}$ with $\lambda'_r=\lambda_r $ if $r\in\Gc \setminus (\Gc[ij] \cup \Gc[ik])$ and $\lambda'_{ij}=\lambda'_{ik}=0$. 
\end{enumerate}
\end{lem}

\pf
The proof is similar to the previous Lemma.
For \eqref{item:cut1-1}, we notice that all the vertices of the Dyinkin diagram of $\bq'$ and all the edges different from the ones between $i$ and $j$ and between $i$ and $k$ coincide with those of $\bq$ (with the same labels), while $\bq'$ has no edges between $i$ and $j$ and between $i$ and $k$. The diagram (without labels) of $\bq$ contains a triangle with vertices $i,j,k$ and ramifications in one or two vertices, which are \emph{lines}; hence, 
after removing the edges between $i$ and $j$  and between $i$ and $k$ we obtain exactly two connected components: one containing $i$ and one containing $j$ and $k$.
Now \eqref{item:cut1-2} and \eqref{item:cut1-3} follow from \eqref{item:cut1-1} and Theorem \ref{thm:presentacion minima}.
\epf

As an example, we apply this result for some diagrams in type $\g(2,6)$, see \S \ref{subsec:type-g(2,6)}. More precisely, for the diagram in p. \pageref{page:g26-example}
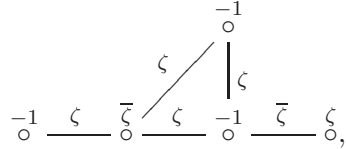
$$ \xymatrix@R-8pt{  &  &  \overset{-1}{\circ} \ar  @{-}[dl]_{\zeta}\ar  @{-}[d]^{\zeta} &
	\\
	\overset{-1}{\underset{\ }{\circ}} \ar  @{-}[r]^{\zeta}  & \overset{\ztu}{\underset{\
		}{\circ}} \ar  @{-}[r]^{\zeta}  & \overset{-1}{\underset{\ }{\circ}} \ar  @{-}[r]^{\ztu}
	& \overset{\zeta}{\underset{\ }{\circ}},}$$
$\toba$ 
is generated by $y_i$, $i\in\I_5$, with defining relations
$y_{ij}=\lambda_{ij}, \, i<j, \, \widetilde{q}_{ij}=1;$
\begin{align*}
y_{1}^2&=\mu_1; & [y_{(24)},&y_{3}]_c=\nu_2; & y_{3}^2&=\mu_3; & y_{221}&=0;\\
y_{223}&=0; & y_{5}^2&=\mu_5; & y_{443}&=0; & [y_{435},&y_{3}]_c=0;
\end{align*}
\vsp
\begin{align*}
y_{225}&=0; & y_{235}-q_{35}\ztu[y_{25},y_3]_c-q_{23}(1-\zeta)y_3y_{25}=0.
\end{align*}
Here, $\Gc[25] \cup \Gc[35]=\{x_{225} , [x_{435},x_{3}]_c , x_{235}-q_{35}\ztu[x_{25},x_3]_c-q_{23}(1-\zeta)x_3x_{25} \}$, 
and these relations are not deformed. Hence we consider the corresponding matrix $\bq'$ as in the Lemma, which has diagram 

$$ \xymatrix@R-8pt{  &  &  \overset{-1}{\circ} &
\\
\overset{-1}{\underset{\ }{\circ}} \ar  @{-}[r]^{\zeta}  & \overset{\ztu}{\underset{\ }{\circ}} \ar  @{-}[r]^{\zeta}  & \overset{-1}{\underset{\ }{\circ}} \ar  @{-}[r]^{\ztu}
& \overset{\zeta}{\underset{\ }{\circ}} \, .}$$ 
Thus $\bq'$ has two components, one of super type $A_4$ and another of type $A_1$, and the defining relations are the same as above but we remove $\Gc[25]\cup \Gc[35]$ and add $x_{25}$, $x_{35}$. By a recursive argument on the rank of connected components, we may assume that $\widetilde{\mathcal E}_{\bq'}\neq 0$. As $\toba$ projects onto $\widetilde{\mathcal E}_{\bq'}$ by the previous result, we conclude that $\toba\neq 0$.

\subsubsection{}\label{subsubsec:projection} Assume that all the relations involving a letter $y_i$, $i\in \I$, are not deformed. Let $\Gc[i]\subset \Gc$ be the subset containing these relations. 
If all the relations in $\Gc[i]$ are not deformed, then we may deal again recursively using the next result, see for instance \ref{subsubsec:type-F-N>4}.

\begin{lem}\label{lem:cut3}
Let $\bq'$ be the submatrix of $\bq$ obtained by removing the $i$-th row and column. Then
\begin{enumerate}[leftmargin=*]
\item\label{item:cut3-1} The diagram associated to $\bq'$ arises by erasing the vertex $i$ and all the edges adjacent to $i$. 
\item\label{item:cut3-2} The set of defining relations of $\B_{\bq'}$ is $\Gc'=\Gc\setminus \Gc[i]$.
\item\label{item:cut3-3} 
Assume that all relations $r\in \Gc[i]$ are not deformed.
Then $\toba(\bs\lambda)$ projects onto $\widetilde{\mathcal{E}}_{\bq'}(\bs\lambda')$ for $\bs\lambda'=(\lambda'_r)_{r\in \Gc'}$. 
\end{enumerate}
\end{lem}

\pf
Analogous to Lemmas \ref{lem:cut1} and \ref{lem:cut2}.
\epf

\subsubsection{}\label{subsubsec:casebycase} In many cases, we can restrict to the case in which a single relation $r\in\Gc$ is deformed. Then an iterative application of Lemma \ref{lem:nonzero-k<l} shows that $\toba\neq 0$ by deforming one relation at a time. In each case, we can restrict to the subdiagram involving just the letters in the relation, using the technique described in \S \ref{subsubsec:projection}.

\subsection{Cartan type}\label{sec:by-diagram-cartan}

In what follows $q\in \k^{\times}$ has order $N > 1$. All generalized Dynkin diagrams are
numbered from the left to the right and from bottom to top.
If a numbered display contains several equalities (or diagrams), they will be referred to
with roman letters from the left to the right.

\subsubsection{Type $A_{\theta}$, $\theta \ge 1$}\label{subsec:type-A}
The Dynkin diagram is
\begin{align*}
\xymatrix{ \overset{q}{\underset{\ }{\circ}}\ar  @{-}[r]^{q^{-1}}  &
	\overset{q}{\underset{\ }{\circ}}\ar  @{-}[r]^{q^{-1}} &  \overset{q}{\underset{\
		}{\circ}}\ar@{.}[r] & \overset{q}{\underset{\ }{\circ}} \ar  @{-}[r]^{q^{-1}}  &
	\overset{q}{\underset{\ }{\circ}}}
\end{align*}
For $N > 2$, the algebra $\toba$ is generated by $y_i$, $i\in\I_{\theta}$, with defining
relations
\begin{align*}
y_{ij} &= 0, \quad i < j - 1; & y_{iij} &= \lambda_{iij}, \quad \vert j - i\vert = 1.
\end{align*}
This algebra is nonzero, see \cite[Proposition 5.2]{AAG}.

If $N = 2$, then $\toba$ is generated by $y_i$, $i\in\I_{\theta}$, with defining
relations
\begin{align*}
y_{ij} &= \lambda_{ij}, \quad i < j - 1; & [y_{(i\, i+2)}, y_{i+1}]_c &= \nu_i.
\end{align*}
This algebra is nonzero, see \cite[Proposition 4.2]{AAG}.

\subsubsection{Type $B_{\theta}$, $\theta \ge 2$, $N> 2$}\label{subsec:type-B}
The Dynkin diagram is
\begin{align*}
\xymatrix{ \overset{\,\,q^2}{\underset{\ }{\circ}}\ar  @{-}[r]^{q^{-2}}  &
	\overset{\,\,q^2}{\underset{\ }{\circ}}\ar  @{-}[r]^{q^{-2}} &
	\overset{\,\,q^2}{\underset{\ }{\circ}}\ar@{.}[r] & \overset{\,\,q^2}{\underset{\
		}{\circ}} \ar  @{-}[r]^{q^{-2}}  & \overset{q}{\underset{\ }{\circ}}}
\end{align*}
Now we consider different cases according with $N$.

\noindent
$\diamond${$N>4$}.\label{subsubsec:type-B-N>4}
The algebra $\toba$ is generated by $y_i$, $i\in\I_{\theta}$, with defining
relations
\begin{align*}
y_{ij} &= 0, \quad i < j - 1; & y_{iij}&=\lambda_{iij}, \quad i < \theta, \, \vert j - i\vert = 1;  &
y_{\theta\theta\theta\theta-1}&=\lambda_{\theta\theta\theta\theta-1}.
\end{align*}
If $N>5$, then $\lambda_{\theta\theta\theta\theta-1}=\lambda_{iij}=0$ for all $i,j$, so $\toba=\widetilde{\B}_{\bq}$.
For $N=5$ all these scalars are 0 again, unless $\bq=\left(\begin{smallmatrix} q^2&q\\q^2&q\end{smallmatrix}\right)$, in which case $\lambda_{\theta-1\theta-1\theta}$, $\lambda_{\theta\theta\theta\theta-1}$ can be non-zero. This algebra projects onto the corresponding algebra of type $B_2$, which is nonzero, see \texttt{b2.log}. 

Explicitly, the corresponding Gr\"obner basis is given by:
\begin{align*}
&y_{221} -\lambda_{221}, \qquad  y_{1112}-\lambda_{1112},\\
& y_1^2y_2y_1y_2 + (1+q)y_1y_2y_1^2y_2 + (1+q^2)y_1y_2y_1y_2y_1 - q(2+q+2q^2)y_2y_1y_2y_1^2\\
&- q(1+q)y_2y_1^2y_2y_1  + 2y_2^2y_1^3 -q^2(1-q+q^2)\lambda_{221} y_1^2 + q(1+q)\lambda_{1112} y_2.
\end{align*}
\noindent
$\diamond${$N=4$}.\label{subsubsec:type-B-N=4}
The algebra $\toba$ is generated by $y_i$, $i\in\I_{\theta}$, with defining
relations
\begin{align*}
y_{ij} &= \lambda_{ij}, \quad i < j - 1; & [y_{(i\, i+2)}, y_{i+1}]_c&=\nu_i, \quad i < \theta;  &
y_{\theta\theta\theta\theta-1}&=0.
\end{align*}
Notice that $\nu_{\theta-1}=0$ since $\chi_{\theta-2}\chi_{\theta-1}^2 \chi_{\theta} (g_{\theta-2}g_{\theta-1}^2 g_{\theta})=-q$. Also, $\lambda_{i\theta}=0$ for all $i\leq \theta-2$, so Lemma \ref{lem:cut3} applies for $i=\theta$ and then $\toba$ projects over $\widetilde{\mathcal E}_{\bq'}$, $\bq'$ of type $A_{\theta-1}$.

\noindent
$\diamond${$N=3$}.\label{subsubsec:type-B-N=3}
The algebra $\toba$ is generated by $y_i$, $i\in\I_{\theta}$, with defining
relations
\begin{align*}
y_{ij} &= 0, \quad i+1 < j < \theta; & y_{iij}&=\lambda_{iij}, \quad i < \theta, \, \vert j - i\vert = 1; \\
y_{i\theta} &= \lambda_{i\theta}, \quad i+1 < \theta;   &
[y_{\theta\theta\theta-1\theta-2}&,y_{\theta\theta-1}]_c=0.
\end{align*}
Notice that $\lambda_{\theta-1 \, \theta-1 \, \theta}=0$. If $\lambda_{\theta-1 \, \theta-1 \, \theta-2}=\lambda_{\theta-2 \, \theta-2 \, \theta-1}=0$, then $\toba$ projects onto $\widetilde{\mathcal E}_{\bq'}$, $\bq'$ of type $A_{\theta-2} \times A_1$ by Lemma \ref{lem:cut1}. If either $\lambda_{\theta-1 \, \theta-1 \, \theta-2}\neq0$ or $\lambda_{\theta-2 \, \theta-2 \, \theta-1}\neq 0$, then
$\lambda_{\theta-2 \, \theta-2 \, \theta-3}=\lambda_{\theta-3 \, \theta-3 \, \theta-2}=0$; by Lemma \ref{lem:cut1}, $\toba$ projects onto $\widetilde{\mathcal E}_{\bq'}$, $\bq'$ of type $A_{\theta-3} \times B_3$, which is nonzero, see \texttt{b3.log}.

\subsubsection{Type $C_{\theta}$, $\theta \ge 3$, $N> 2$}\label{subsec:type-C}

The associated Dynkin diagram is
\begin{align*}
\xymatrix{ \overset{q}{\underset{\ }{\circ}}\ar  @{-}[r]^{q^{-1}}  &
	\overset{q}{\underset{\ }{\circ}}\ar  @{-}[r]^{q^{-1}} &  \overset{q}{\underset{\
		}{\circ}}\ar@{.}[r] & \overset{q}{\underset{\ }{\circ}} \ar  @{-}[r]^{q^{-2}}  &
	\overset{\,\,q^2}{\underset{\ }{\circ}}}
\end{align*}

\noindent
$\diamond${$N>3$}.\label{subsubsec:type-C-N>3}
The algebra $\toba$ is generated by $y_i$, $i\in\I_{\theta}$, with defining
relations
\begin{align*}
y_{ij} &= 0, \quad i < j - 1; &
y_{iij}&=0, \quad i \leq \theta-2, \, \vert j - i\vert = 1;
\\
y_{\theta-1\theta-1\theta-1\theta}&=\lambda_{\theta-1\theta-1\theta-1\theta}; &
y_{\theta\theta\theta-1}&=\lambda_{\theta\theta\theta-1}.
\end{align*}
If $N\neq 5$, then $\toba=\widetilde{\B}_{\bq}$.
For $N=5$, $\bs\lambda=0$ again, unless $\bq=\left(\begin{smallmatrix} q & q^2 \\ q & q^2 \end{smallmatrix}\right)$, in which case $\lambda_{\theta-1\theta-1\theta-1\theta}$, $\lambda_{\theta\theta\theta-1}$ can be non-zero. Thus $\toba$ projects onto the corresponding algebra of type $B_2$, which is nonzero, see \texttt{b2.log}.

\noindent
$\diamond${$N=3$}.\label{subsubsec:type-C-N=3}
The algebra $\toba$ is generated by $y_i$, $i\in\I_{\theta}$, with defining
relations
\begin{align*}
y_{ij} &= 0, \quad i+1 < j <\theta; &
y_{iij}&=\lambda_{iij}, \quad j=i\pm 1, i<\theta-1, & y_{\theta\theta\theta-1}&=0;
\\
y_{i\theta} &= \lambda_{i\theta}, \quad i < \theta - 1; &
&[[ y_{(\theta-2\theta)}, y_{\theta-1}]_c, y_{\theta-1}]_c
=\lambda_{(\theta-2\theta)}.
\end{align*}
If $\lambda_{(\theta-2\theta)}=0$, then $\toba$ projects onto $\widetilde{\mathcal E}_{\bq'}$, $\bq'$ of type $A_{\theta-1}\times A_1$, see Lemma \ref{lem:cut1}. If $\lambda_{(\theta-2\theta)}\neq 0$, then we may restrict to rank 3 following \S \ref{subsubsec:casebycase} and $\bq=\left(\begin{smallmatrix}
\zeta & b\zeta^2 & b^{-3}\zeta^2 \\
 b^{-1} & \zeta & b \\
b^3\zeta & b^{-1}\zeta & \zeta^2
\end{smallmatrix}\right)$; this algebra is not zero by \texttt{c3.log}.

\subsubsection{Type $D_{\theta}$, $\theta \ge 4$}\label{subsec:type-D}
The associated Dynkin diagram is
\begin{align*}
\xymatrix{ & & & &  \overset{q}{\circ} &\\
	\overset{q}{\underset{\ }{\circ}}\ar  @{-}[r]^{q^{-1}}  & \overset{q}{\underset{\
		}{\circ}}\ar  @{-}[r]^{q^{-1}} &  \overset{q}{\underset{\ }{\circ}}\ar@{.}[r] &
	\overset{q}{\underset{\ }{\circ}} \ar  @{-}[r]^{q^{-1}}  & \overset{q}{\underset{\
		}{\circ}} \ar @<0.7ex> @{-}[u]_{q^{-1}}^{\qquad} \ar  @{-}[r]^{q^{-1}} &
	\overset{q}{\underset{\ }{\circ}}}
\end{align*}

\noindent
$\diamond${$N > 2$}.\label{subsubsec:type-D-N>2}
The algebra $\toba$ is generated by $y_i$, $i\in\I_{\theta}$, with defining
relations
\begin{align*}
y_{\theta-1\theta} &=0; & y_{ij} &= 0, \quad i < j - 1, (i,j)\neq (\theta-2,\theta);  &
y_{\theta\theta\theta-2}&=\lambda_{\theta\theta\theta-2}; \\
&& y_{iij} &= \lambda_{iij}, \quad \vert j - i\vert = 1, i,j\neq \theta; &
y_{\theta-2\theta-2\theta}&=\lambda_{\theta-2\theta-2\theta}.
\end{align*}
If $N\neq 3$, then $\bs\lambda=0$, so $\toba=\widetilde{\B}_{\bq}$.
If $N=3$, then either $\lambda_{\theta-2\theta-2\theta-1}= \lambda_{\theta-1\theta-1\theta-2}=0$ or else $\lambda_{\theta-2\theta-2\theta}= \lambda_{\theta\theta\theta-2}=0$ by \cite[Lemma 5.1 (2)]{AAG}; in both cases we apply Lemma \ref{lem:cut3}: $\toba$ projects onto $\widetilde{\mathcal E}_{\bq'}$, $\bq'$ of type $A_{\theta-1}$.

\noindent
$\diamond${$N = 2$}.\label{subsubsec:type-D-N=2}
The algebra $\toba$ is generated by $y_i$, $i\in\I_{\theta}$, with defining
relations
\begin{align*}
y_{ij} &= \lambda_{ij}, \quad i < j - 1, (i,j)\neq (\theta-2,\theta);  &  y_{\theta-1\theta} &=\lambda_{\theta-1\theta}; \\
  [y_{(ii+2)},& y_{i+1}]_c = \nu_i, \quad i\leq \theta-3; &
[y_{\theta-3\theta-2\theta},y_{\theta-2}]_c&=\nu_{\theta-2}'.
\end{align*}
By Diamond Lemma \cite{B - diamond} $\toba\neq 0$: the defining relations are as for type $A$, so the proof of \cite[Proposition 4.2]{AAG} applies to solve the ambiguities.

\subsubsection{Type $E_{\theta}$, $6 \le \theta \le 8$}\label{subsec:type-E}
The associated Dynkin diagram is
\begin{align*}
\xymatrix{ &  &   \overset{q}{\circ} &\\
	\overset{q}{\underset{\ }{\circ}}\ar  @{-}[r]^{q^{-1}}  &  \overset{q}{\underset{\
		}{\circ}}\ar@{.}[r]  & \overset{q}{\underset{\ }{\circ}} \ar @<0.7ex> @{-}[u]_{q^{-1}} \ar
	@{-}[r]^{q^{-1}} &  \overset{q}{\underset{\ }{\circ}}  \ar  @{-}[r]^{q^{-1}} &
	\overset{q}{\underset{\ }{\circ}}}
\end{align*}
The proof is analogous to the case $D_{\theta}$.



\subsubsection{Type $F_{4}$, $N> 2$}\label{subsec:type-F}
The associated Dynkin diagram is 
\begin{align*}
\xymatrix{ \overset{\,\,q^2}{\underset{\ }{\circ}}\ar  @{-}[r]^{q^{-2}}  &
	\overset{\,\,q^2}{\underset{\ }{\circ}}\ar  @{-}[r]^{q^{-2}} &   \overset{q}{\underset{\
		}{\circ}} \ar  @{-}[r]^{q^{-1}}  &  \overset{q}{\underset{\ }{\circ}}. }
\end{align*}

\noindent
$\diamond${$N>4$}.\label{subsubsec:type-F-N>4}
The algebra $\toba$ is generated by $y_i$, $i\in\I_4$, with defining
relations
\begin{align*}
y_{ij} &= 0, \quad i < j - 1; & y_{112}&=0, & y_{221}&=0, & y_{334}&=0, \\ y_{443}&=0, & y_{223}&=\lambda_{223}, &  y_{3332}&=\lambda_{3332}.
\end{align*}
Here we may apply Lemma \ref{lem:cut3} twice (to vertices $i=4$ and $i=1$ respectively) to see that $\toba$ projects onto $\widetilde{\mathcal E}_{\bq'}$, $\bq'$ of type $B_2$. 

\noindent
$\diamond${$N=4$}.\label{subsubsec:type-F-N=4}
The algebra $\toba$ is generated by $y_i$, $i\in\I_4$, with defining
relations
\begin{align*}
y_{ij} &= 0, \quad i < j - 1; & [y_{(13)}, y_{2}]_c&=0; & y_{3332}&=0  \\
y_{334}&=0; & y_{443}&=0.
\end{align*}
Hence, $\toba=\widetilde{\B}_{\bq} \neq 0$.

\noindent
$\diamond${$N=3$}.\label{subsubsec:type-F-N=3}
The algebra $\toba$ is generated by $y_i$, $i\in\I_4$, with defining
relations
\begin{align*}
y_{ij} &= \lambda_{ij}, \quad i < j - 1; & [[y_{(24)}, y_3]_c, y_3]_c&=\lambda_{(24)};\\
y_{iij}&=\lambda_{iij}, \quad  j=i\pm 1, (i,j)\neq (3,2);   &
[y_{3321},y_{32}]_c&=\lambda_{321}.
\end{align*}
Notice that $\lambda_{223}=0$.
If $\lambda_{321}=\lambda_{(24)}=0$, then $\toba$ projects onto $\widetilde{\mathcal E}_{\bq'}$, $\bq'$ of type $A_2\times A_2$, see Lemma \ref{lem:cut1}. If $\lambda_{321}\neq0$ but all the other $\lambda$'s are 0, then $\toba$ projects onto $\widetilde{\mathcal E}_{\bq'}$, $\bq'$ of type $B_3\times A_1$ by Lemma \ref{lem:cut1} again. If $\lambda_{(24)}\neq0$ but all the other $\lambda$'s are 0, then $\toba$ projects onto $\widetilde{\mathcal E}_{\bq'}$, $\bq'$ of type $C_3\times A_1$, again by Lemma \ref{lem:cut1}.  The general case follows by applying Lemma
\ref{lem:nonzero-k<l} twice.

\subsubsection{Type  $G_{2}$, $N > 3$}\label{subsec:type-G}
The associated Dynkin diagram is $\xymatrix{  \overset{\,\,q}{\underset{\ }{\circ}} \ar  @{-}[r]^{q^{-3}} & \overset{q^3}{\underset{\ }{\circ}}}$.
The algebra $\toba$ is generated by $y_1, y_2$ with defining
relations
\begin{align*}
y_{112} &= \lambda_{1}; &y_{22221}&=\lambda_{2}
\end{align*}
with $\lambda_1=\lambda_2=0$ unless $N=7$ and $\bq=\left(\begin{smallmatrix}
q&q^3\\q&q^3
\end{smallmatrix}\right)$. This algebra is nonzero, see \texttt{g2.log}.

\subsection{Standard type}\label{sec:by-diagram-standard}

We start by recalling some notation from \cite[\S 5]{AA}.
Given $q\in\Bbbk^\times - \G_2$ and $\Jb\subset\I_{\theta}$, 
${\bf A}_{\theta}(q;\Jb)$\label{page:diagram-super-A} is the generalized Dynkin diagram 
\begin{align*}
\xymatrix{ \overset{q_{11}}{\underset{\ }{\circ}}\ar  @{-}[rr]^{\widetilde{q}_{12}}  &&
	\overset{q_{22}}{\underset{\ }{\circ}}\ar@{.}[r] &  \overset{\quad q_{\theta-1 \theta-1}}{\underset{\ }{\circ}} \ar  @{-}[rr]^{\widetilde{q}_{\theta-1 \theta}}  &&
	\overset{\quad q_{\theta\theta}}{\underset{\ }{\circ}}}
\end{align*}
subject to the following requirements:
\begin{enumerate} [leftmargin=*]
	\item\label{it:aqj-1} $q = q_{\theta\theta}^2\widetilde{q}_{\theta-1 \theta}$;
	
	\smallbreak\item\label{it:aqj-2} if $i\in\Jb$, then $q_{ii}=-1$ and $\widetilde{q}_{i-1 i}=\widetilde{q}_{i i+1}^{-1}$;
	
	\smallbreak\item\label{it:aqj-4} if $i\notin\Jb$, then $\widetilde{q}_{i-1 i}= q_{ii}^{-1} = \widetilde{q}_{i i+1}$ (only the second equality if $i=1$, only the first if $i= \theta$).
\end{enumerate}

\subsubsection{Standard type $B_{\theta}$, $\theta \ge 2$, $\zeta \in
\G'_3$}\label{subsec:type-B-standard}
The associated Dynkin diagram is
$$ \xymatrix{ {\bf A}_{\theta-1}(-\ztu;\Jb) \ar @{-}[r]^(.7){-\zeta }  & \overset{\zeta}{\underset{\ }{\circ}}.}$$
The algebra $\toba$ is generated by $y_i$, $i\in\I_{\theta}$, with defining
relations
\begin{align*}
y_{ij} &= \lambda_{ij}, \quad i < j - 1; &  [y_{(i-1i+1)},y_i]_c&=\nu_i, \quad
q_{ii}=-1;  \\
 y_{iii\pm1} &= 0, \quad q_{ii}=-\zeta^{\pm
1}; & [y_{\theta\theta\theta-1\theta-2},y_{\theta\theta-1}]_c&=0;  \\
 y_\theta^3&=\mu_\theta;  & [y_{\theta\theta\theta-1},y_{\theta\theta-1}]_c
 &=\lambda_{\theta\theta-1},
\quad q_{\theta-1\theta-1}=-1; \\
  y_i^2&=\mu_i, \quad
q_{ii}=-1.
\end{align*}
As $\chi_{\theta-1}\chi_\theta^2\chi_{\theta+1} (g_{\theta-1}g_\theta^2 g_{\theta+1})\neq 1$, $\nu_{\theta-1}=0$. Hence, if either $q_{\theta-1\theta-1}\neq -1$ or $q_{\theta-1\theta-1}=-1$ but $\lambda_{\theta\theta-1}=0$, then by Lemma \ref{lem:cut1}, $\toba$ projects onto $\widetilde{\mathcal E}_{\bq'}$, $\bq'$ with two components, one of super type $A$ and another of type $A_1$; thus $\toba$ is nonzero.

If $q_{\theta-1\theta-1}=-1$ and $\lambda_{\theta\theta-1}\neq 0$, then either $q_{\theta-2\theta-2}\neq 1$ or else $q_{\theta-2\theta-2}\neq 1$, $\nu_{\theta-2}=-1$. In both cases, all the relations in $\Gc[\theta-2 \, \theta-1]$ are not deformed so Lemma \ref{lem:cut1} applies: $\toba$ projects onto $\widetilde{\mathcal E}_{\bq'}$, $\bq'$ with two components, one of super type $A$ and another of type $B_2$ with matrix
$\left(\begin{smallmatrix} -1 & \zeta \\ -1 & \zeta
\end{smallmatrix}\right)$; this algebra is nonzero, see \texttt{b2st.log}.

\subsubsection{Standard type $G_{2}$, $\zeta \in \G'_8$}\label{subsec:type-G-st}
The Weyl groupoid has three objects.

For $\xymatrix{ \overset{\,\, \zeta^2}{\underset{\ }{\circ}} \ar  @{-}[r]^{\zeta}  &
\overset{\ztu}{\underset{\ }{\circ}}}$, $\toba$ is generated by $y_1, y_2$ with defining relations
\begin{align*}
y_1^4&=\lambda_1; & y_{221}&=\lambda_2; & [y_{1112},y_{112}]&=\lambda_3+2\lambda_2(1+q)^2y_2y_{12},
\end{align*}
with $\lambda_1=\lambda_2=\lambda_3=0$ unless $\bq=\left(\begin{smallmatrix}
\zeta^2&\zeta^{-1}\\\zeta^2&\zeta^{-1}
\end{smallmatrix}\right)$; see \texttt{g2-st-a1.log} for the deformation of the generalized quantum Serre relation.
This algebra is nonzero, see \texttt{g2-st-a.log}.

For $\xymatrix{ \overset{\,\, \zeta^2}{\underset{\ }{\circ}} \ar  @{-}[r]^{\zeta^3}  &
\overset{-1}{\underset{\ }{\circ}}}$, $\toba$ is generated by $y_1, y_2$ with defining relations
\begin{align*}
y_1^4&=\lambda_1; & y_2^2&=\lambda_2; & [y_1,y_{3\alpha_1+2\alpha_2}]_c+\frac{q_{12}}{1-\zeta}y_{112}^2&=0,
\end{align*}
with $\lambda_1\lambda_2=0$. This algebra is nonzero, using \cite[Lemma 5.16]{AAGMV}. The last relation is not primitive: it remains undeformed by direct computation.

For $\xymatrix{ \overset{\zeta}{\underset{\ }{\circ}} \ar
 @{-}[r]^{\zeta^5}  & \overset{-1}{\underset{\ }{\circ}}}$, $\toba$ is generated by $y_1, y_2$ with defining relations
\begin{align*}
y_{11112}&=\lambda_1; & y_2^2&=\lambda_2, & [y_{112},y_{12}]_c,y_{12}]_c&=\lambda_3;
\end{align*}
with $\lambda_1=\lambda_2=\lambda_3=0$ unless $\bq=\left(\begin{smallmatrix}
\zeta&-1\\\zeta&-1
\end{smallmatrix}\right)$.
This algebra is nonzero, see \texttt{g2-st-c.log}.

\subsection{Super type}\label{sec:by-diagram-super}

\subsubsection{Type $\supera{j-1}{\theta - j}$, $\theta \ge 1$, $1\le j\le \theta$, $N >
2$}\label{subsec:type-A-super}
The Dynkin diagram is of the shape ${\bf A}_{\theta}(q;\Jb)$, see p.~\pageref{page:diagram-super-A}. 
The algebra $\toba$ is generated by $y_i$, $i\in\I_{\theta}$, with defining
relations
\begin{align*}
y_{ij} &= \lambda_{ij}, \quad i < j - 1; & y_{iii\pm1} &= \lambda_{iii\pm1}, \quad
q_{ii}\neq-1;  \\
 y_i^2&=0, \quad q_{ii}=-1;  & [y_{(i-1i+1)},y_i]_c&=\nu_i, \quad
q_{ii}=-1.
\end{align*}
We argue by induction on $\theta$. If $\theta=2$, then there are 3 diagrams.
\begin{description}[leftmargin=*]
\item[$\xymatrix{ \overset{q}{\underset{\ }{\circ}} \ar@{-}[r]^{q^{-1}}  & \overset{q}{\underset{\ }{\circ}}}$] It is of Cartan type, already solved in \S \ref{subsec:type-A}.
\item[$\xymatrix{\overset{-1}{\underset{\ }{\circ}} \ar@{-}[r]^{q}  & \overset{-1}{\underset{\ }{\circ}}}$] $y_i^2=\mu_i$, $i=1,2$. Hence \cite[Lemma 5.16]{AAGMV} applies.
\item[$\xymatrix{ \overset{q}{\underset{\ }{\circ}} \ar@{-}[r]^{q^{-1}}  & \overset{-1}{\underset{\ }{\circ}}}$] $y_{112}=\lambda_{112}$, $y_2^2=\mu_2$. If $\lambda_{112}=0$, then \cite[Lemma 5.16]{AAGMV} applies. Otherwise, $\bq=\left(\begin{smallmatrix} q & -1 \\ q & -1 \end{smallmatrix}\right)$, $q\in\G_4'$; $\toba\neq 0$ by \texttt{a2-super.log}.
\end{description}

For $\theta=3$, we have six cases:
\begin{description}[leftmargin=*]
	\item[$\xymatrix{ \overset{q}{\underset{\ }{\circ}} \ar@{-}[r]^{q^{-1}}  & \overset{q}{\underset{\ }{\circ}} \ar@{-}[r]^{q^{-1}}  & \overset{q}{\underset{\ }{\circ}}}$] It is of Cartan type, already solved.
	\item[$\xymatrix{ \overset{q}{\underset{\ }{\circ}} \ar@{-}[r]^{q^{-1}}  & \overset{q}{\underset{\ }{\circ}} \ar@{-}[r]^{q^{-1}}  & \overset{-1}{\underset{\ }{\circ}} }$] If $\lambda_{223}=0$, then we quotient by $y_{23}$. If $\lambda_{223}\neq 0$, then $q\in\G_4'$, so $\lambda_{112}=\lambda_{221}=0$ and we quotient by $y_{12}$.
	\item[$\xymatrix{ \overset{q}{\underset{\ }{\circ}} \ar@{-}[r]^{q^{-1}}  & \overset{-1}{\underset{\ }{\circ}} \ar@{-}[r]^{q}  & \overset{q^{-1}}{\underset{\ }{\circ}}}$] First assume $\nu_2=0$. If $\lambda_{112}\lambda_{332}=0$, then it projects over a rank 2 case. Otherwise, $\bq=\left(\begin{smallmatrix} q & -1 & b \\ q & -1 & -q \\ -b & -1 & -q
	\end{smallmatrix}\right)$, $b=\pm q\in\G_4'$; this algebra is nonzero,
	see \texttt{a3-super-1.log}. If $\nu_2\neq 0$, then $\bq=\left(\begin{smallmatrix} q & a & q^{-1}a^{-2} \\ q^{-1}a^{-1} & -1 & qa \\ qa^2 & a^{-1} & q^{-1} \end{smallmatrix}\right)$, $a\neq 0$; this algebra is nonzero, see \texttt{a3-super-2.log}. The general case follows by Lemma \ref{lem:nonzero-k<l}.
	\item[$\xymatrix{\overset{-1}{\underset{\ }{\circ}} \ar@{-}[r]^{q}  & \overset{-1}{\underset{\ }{\circ}} \ar@{-}[r]^{q^{-1}}  & \overset{q}{\underset{\ }{\circ}}}$] as $\nu_2=0$, we quotient by $y_{12}$.
	\item[$\xymatrix{\overset{-1}{\underset{\ }{\circ}} \ar@{-}[r]^{q^{-1}}  & \overset{q}{\underset{\ }{\circ}} \ar@{-}[r]^{q^{-1}}  & \overset{-1}{\underset{\ }{\circ}}}$] If $\lambda_{221}\lambda_{223}=0$, then we apply Lemma \ref{lem:cut1} either to $(i,j)=(1,2)$ or $(i,j) =(2,3)$. Otherwise, $\bq=\left(\begin{smallmatrix} -1 & q & -1 \\ -1 & q & -1 \\ -1 & q & -1
	\end{smallmatrix}\right)$, $q\in\G_4'$; this algebra is nonzero, see \texttt{a3-super-3.log}.
	\item[$\xymatrix{\overset{-1}{\underset{\ }{\circ}} \ar@{-}[r]^{q^{-1}}  & \overset{-1}{\underset{\ }{\circ}} \ar@{-}[r]^{q}  & \overset{-1}{\underset{\ }{\circ}}}$] If $\nu_2=0$, then we quotient by $y_{12}$. If $\nu_2\neq 0$, then $\bq=\left(\begin{smallmatrix} -1 & a & -a^{-2} \\ q^{-1}a^{-1} & -1 & qa \\ -a^2 & a^{-1} & -1 \end{smallmatrix}\right)$, $a\neq 0$; this algebra is nonzero, see \texttt{a3-super-4.log}.
\end{description}

For the inductive step, if the relations in $\Gc[\theta \, \theta+1]$ are not deformed, then we apply Lemma \ref{lem:cut1}: $\toba$ projects onto $\widetilde{\mathcal E}_{\bq'}$, $\bq'$ with two components, of super type $A$ and $A_1$. If the relations in $\Gc[\theta \, \theta+1]$ are the unique deformed relations, then it projects over the rank 2 or rank 3 cases by Lemma \ref{lem:cut3}. The general case follows by Lemma \ref{lem:nonzero-k<l}.

\subsubsection{Type $\superb{m}{n}$, $q\notin\G_4$}\label{subsec:type-B-super}
The Dynkin diagram has the shape 
$$ \xymatrix@R-6pt{ {\bf A}_{\theta-1}(q^2;\Jb) \ar @{-}[rr]^{\hspace*{1.2cm} q^{-2}} & & \overset{q}{\underset{\	}{\circ}}} .$$

\noindent
$\diamond${$N>4$}.\label{subsubsec:type-B-super-N>4}
The algebra $\toba$ is generated by $y_i$, $i\in\I_{\theta}$, with defining
relations
\begin{align*}
y_{ij} &= \lambda_{ij}, \, i < j - 1; \qquad y_i^2=\mu_i, \, q_{ii}=-1;   &  y_{\theta\theta\theta\theta-1}&=\lambda_{\theta\theta\theta\theta-1}; \\
y_{iij}&=\lambda_{iij}, \, i < \theta, |j-i|=1, q_{ii}\neq -1; & [y_{(i-1\, i+1)},y_i]_c&=\nu_i, \, q_{ii}= -1.
\end{align*}
First assume $q_{\theta-1\theta-1}\neq -1$. We apply Lemma \ref{lem:cut1} in all the cases. If $\lambda_{\theta\theta\theta\theta-1} =\lambda_{\theta-1\theta-1\theta}=0$, then $\toba$ projects onto $\widetilde{\mathcal E}_{\bq'}$, $\bq'$ with two components, one of super type $A$ and another of type $A_1$. If either $\lambda_{\theta\theta\theta\theta-1}\neq 0$ or $\lambda_{\theta-1\theta-1\theta}\neq 0$, then $q\in\G_5'$ and $\nu_{\theta-2}=0$; hence
$\toba$ projects onto $\widetilde{\mathcal E}_{\bq'}$, $\bq'$ with two components, one of super type $A$ and another of Cartan type $B_2$.

When $q_{\theta-1\theta-1}=-1$, $\nu_{\theta-1}=0$ since $\chi_{\theta-2}\chi_{\theta-1}^2\chi_\theta\neq \eps$. If  $\lambda_{\theta\theta\theta\theta-1}=0$, then we apply Lemma \ref{lem:cut1} to $(i,j)=(theta-1 \, \theta)$.
If $\lambda_{\theta\theta\theta\theta-1}\neq 0$, then $q\in\G_6'$ and $\nu_{\theta-2}=0$; by Lemma \ref{lem:cut1} again,
$\toba$ projects onto $\widetilde{\mathcal E}_{\bq'}$, $\bq'$ with two components, one of super type $A$ and another of Cartan type $B_2$ with matrix $\left(\begin{smallmatrix} -1 & q \\ -1 & q \end{smallmatrix}\right)$; this algebra is nonzero, see \texttt{b2supera.log}.

\noindent
$\diamond${$N=3$}.\label{subsubsec:type-B-super-N=3}
The algebra $\toba$ is generated by $y_i$, $i\in\I_{\theta}$, with defining
relations
\begin{align*}
y_{ij} &= \lambda_{ij}, \quad i < j - 1; &
[y_{(i-1i+1)},y_i]_c&=\nu_i, \quad q_{ii}=-1;  \\
y_{iii\pm1} &= \lambda_{iii\pm1}, \quad q_{ii}\neq -1; &
[y_{\theta\theta\theta-1},y_{\theta\theta-1}]_c&=\lambda_{\theta\theta-1}, \quad q_{\theta-1\theta-1}=-1;\\
[y_{\theta\theta\theta-1\theta-2},&y_{\theta\theta-1}]_c
=\lambda_{\theta\theta-1\theta-2};  &
y_i^2&=\mu_i, \quad q_{ii}=-1.
\end{align*}
First assume $\lambda_{\theta\theta-1\theta-2}\neq 0$. Then $q_{ii}=q^2$ for $i=\theta-2,\theta-1$. Moreover, $\lambda_{\theta-2\theta-2\theta-3}=0$, $\nu_{\theta-3}=0$ if $q_{\theta-3\theta-3}=-1$, $\lambda_{\theta-3\theta-3\theta-2}=0$ if $q_{\theta-3\theta-3}=q^2$; hence Lemma \ref{lem:cut1} applies for $i=\theta-3$, $j=\theta-2$, and
$\toba$ projects onto $\widetilde{E}_{\bq'}$, $\bq'$ with two components, one of super type $A$ and of Cartan type $B_3$.

Now assume $\lambda_{\theta\theta-1\theta-2}=0$. If either $q_{\theta-1\theta-1}=q^2$ or $q_{\theta-1\theta-1}=-1$, $\lambda_{\theta\theta-1}=0$, then Lemma \ref{lem:cut1} applies for $i=\theta-1$, $j=\theta$ and $\toba$ projects onto $\widetilde{\mathcal E}_{\bq'}$, $\bq'$ with two components, one of super type $A$ and another of type $A_1$. If $q_{\theta-1\theta-1}=-1$, $\lambda_{\theta\theta-1}\neq 0$, then $\nu_{\theta-1}=0$, $\nu_{\theta-2}=0$ if $q_{\theta-2\theta-2}=-1$, $\lambda_{\theta-2\theta-2\theta-1}=0$ if $q_{\theta-2\theta-2}\neq -1$;
again by Lemma \ref{lem:cut1}, $\toba$ projects onto $\widetilde{\mathcal E}_{\bq'}$, $\bq'$ with two components, one of super type $A$ and another of Cartan type $B_2$ with matrix $\left(\begin{smallmatrix} -1 & \zeta \\ 1 & \zeta \end{smallmatrix}\right)$; this algebra is nonzero, see \texttt{b2superb.log}.

\subsubsection{Type $\superd{m}{n}$, $N > 2$}\label{subsec:type-D-super}
The Weyl groupoid has objects with generalized Dynkin diagrams of six possible shapes.

For the diagram $\xymatrix@R-6pt{ {\bf A}_{\theta-2}(q;i_1,\dots, i_j) \ar @{-}[r]^{\hspace*{1.5cm}
 q^{-1}}  & 	\overset{q}{\underset{\ }{\circ}} \ar  @{-}[r]^{q^{-2}}  & \overset{q^2}{\underset{\ }{\circ}}}$, there are 3 cases:
	
\noindent $\diamondsuit N>4$.
Then $\toba$ is generated by $y_i$, $i\in\I_{\theta}$, with defining
relations
\begin{align*}
& & y_{iii\pm1}&=0, \quad i < \theta-1, q_{ii}\neq -1;  \\
	y_{\theta-1\theta-1\theta-2}&=0; &
	y_{ij} &= \lambda_{ij}, \quad i < j - 1;
	\\
	y_{\theta\theta\theta-1}&=\lambda_{\theta\theta\theta-1}; & [y_{(i-1i+1)},y_i]_c&=\nu_i, \quad q_{ii}=-1, 2\leq i\leq \theta-2;\\
	y_{\theta-1\theta-1\theta-1\theta}&=\lambda_{\theta-1\theta-1\theta-1\theta}; &
	y_i^2&=\mu_i, \quad q_{ii}=-1.
\end{align*}
If $\lambda_{\theta-1\theta-1\theta-1\theta} =\lambda_{\theta\theta\theta-1} =0$, then we apply Lemma \ref{lem:cut1} for $i=\theta-1$, $j=\theta$ and $\toba$ projects onto $\widetilde{\mathcal E}_{\bq'}$, $\bq'$ with two components, one of super type $A$ and another of Cartan type $A_1$. If either $\lambda_{\theta-1\theta-1\theta-1\theta}\neq 0$ or $\lambda_{\theta\theta\theta-1}\neq 0$, then $\nu_{\theta-2}=0$ if $q_{\theta-2\theta-2}=-1$;
hence  we apply Lemma \ref{lem:cut1} for $i=\theta-2$, $j=\theta-1$ and $\toba$ projects onto $\widetilde{\mathcal E}_{\bq'}$, $\bq'$ with two components, one of super type $A$ and another of Cartan type $B_2$.
	
\noindent $\diamondsuit N=4$.
Then $\toba$ is generated by $y_i$, $i\in\I_{\theta}$, with defining
relations
\begin{align*}
& & y_{iii\pm1}&=\lambda_{iii\pm1}, \,  i < \theta-1, q_{ii}\neq -1;  \\
y_{\theta-1\theta-1\theta-2}&=\lambda_{\theta-1\theta-1\theta-2}; &
y_{ij} &= \lambda_{ij}, \,  i < j - 1; \\
y_{\theta-1\theta-1\theta-1\theta}&=0; &
[y_{(i-1i+1)},y_i]_c&=\nu_i, \,  q_{ii}=-1,  i\neq \theta-1;\\
[y_{(\theta-2\theta)},y_{\theta-1\theta}]_c&=0; &
y_i^2&=\mu_i, \,  q_{ii}=-1.
\end{align*}
Here  we apply Lemma \ref{lem:cut1} for $i=\theta-1$, $j=\theta$ and $\toba$ projects onto $\widetilde{\mathcal E}_{\bq'}$, $\bq'$ with two components, one of super type $A$ and another of Cartan type $A_1$.

\noindent $\diamondsuit N=3$.
Then $\toba$ is generated by $y_i$, $i\in\I_{\theta}$, with defining
relations
\begin{align*}
& & y_{iii\pm1}&=\lambda_{iii\pm1}, \, i < \theta-1, q_{ii}\neq -1;  \\
y_{\theta-1\theta-1\theta-2}&=\lambda_{\theta-1\theta-1\theta-2}; &
y_{ij} &= \lambda_{ij}, \, i < j - 1; \\
y_{\theta\theta\theta-1}&=0; &
[y_{(i-1i+1)},y_i]_c&=\nu_i, \, q_{ii}=-1, i\neq \theta-1;\\
y_i^2&=\mu_i, \, q_{ii}=-1; &
[[y_{(\theta-2\theta)}&,y_{\theta-1}]_c,y_{\theta-1}]_c=\lambda_{(\theta-2\theta)}.
\end{align*}
Here, $q_{\theta-2\theta-2}=q$ if $\lambda_{(\theta-2\theta)}\neq 0$; hence we work as for Cartan type $C_\theta$.

For the diagram $\xymatrix@R-6pt{ {\bf A}_{\theta-2}(q;i_1,\dots, i_j) \ar @{-}[r]^{\hspace*{1.5cm} q} &
\overset{-1}{\underset{\ }{\circ}} \ar  @{-}[r]^{q^{-2}}  & \overset{q^2}{\underset{\ }{\circ}}}$, there are 2 cases:

\noindent $\diamondsuit N \neq 4$. Then $\toba$ is generated by $y_i$, $i\in\I_{\theta}$, with defining
relations
\begin{align*}
 y_i^2=\mu_i,& \, q_{ii}=-1; &
 y_{iii\pm1}&=0, \, i < \theta-1, q_{ii}\neq -1;  \\
 y_{ij} = \lambda_{ij}, &\, i < j - 1; &
 [[y_{\theta-2\theta-1},y_{(\theta-2\theta)}&]_c,y_{\theta-1}]_c =0, \, q_{\theta-2\theta-2}=-1; \\
 & & [[[y_{(\theta-3\theta)},y_{\theta-1}]_c, y_{\theta-2} &]_c,y_{\theta-1}]_c =0, \, q_{\theta-2\theta-2}\neq-1;\\
 y_{\theta\theta\theta-1}&=\lambda_{\theta\theta\theta-1}; &
[y_{(i-1i+1)},y_i]_c&=\nu_i, \, q_{ii}=-1,  i\neq \theta-1.
\end{align*}
If $\lambda_{\theta\theta\theta-1}=0$, then we apply Lemma \ref{lem:cut1} for $i=\theta-1$, $j=\theta$ and $\toba$ projects onto $\widetilde{\mathcal E}_{\bq'}$, $\bq'$ with two components, one of super type $A$ and another of type $A_1$. If $\lambda_{\theta\theta\theta-1}\neq 0$, then either $q_{\theta-2\theta-2} \neq -1$, or else $q_{\theta-2\theta-2}=-1$, $\nu_{\theta-2}=0$; in both cases we apply Lemma \ref{lem:cut1} for $i=\theta-2$, $j=\theta-1$ and $\toba$ projects onto $\widetilde{\mathcal E}_{\bq'}$, $\bq'$ with two components, of super types $A_{\theta-2}$ and $A_2$.

\noindent $\diamondsuit N=4$.
Then $\toba$ is generated by $y_i$, $i\in\I_{\theta}$, with defining
relations
\begin{align*}
& y_{ij} = \lambda_{ij}, \, i < j - 1; \qquad y_{iii\pm1}=\lambda_{iii\pm1}, \, i < \theta-1, q_{ii}\neq -1;  \\
\quad &
[[y_{\theta-2\theta-1},y_{(\theta-2\theta)}]_c,y_{\theta-1}]_c =\lambda_{(\theta-2\theta)}, \, q_{\theta-2\theta-2}=-1; \\
& y_{\theta-1\theta}^2=\mu_{\theta-1\theta}; \qquad
[[[y_{(\theta-3\theta)},y_{\theta-1}]_c,y_{\theta-2}]_c,y_{\theta-1}]_c=0, \, q_{\theta-2\theta-2}\neq-1;\\
& y_i^2=\mu_i, \, q_{ii}=-1; \qquad
[y_{(i-1i+1)},y_i]_c=\nu_i, \, q_{ii}=-1,  i\neq \theta-1.
\end{align*}
First consider $q_{\theta-2\theta-2}\neq-1$. If $\lambda_{iii\pm1}\neq 0$, $i=\theta-2$, $\mu_{\theta-1\theta}\neq 0$,  then we apply Lemma \ref{lem:cut1} for $i=\theta-4$, $j=\theta-3$ and $\toba$ projects onto $\widetilde{\mathcal E}_{\bq'}$, $\bq'$ with two components, one of super type $A_{\theta-4}$ and another with matrix $\left( \begin{smallmatrix}
-1 & q & -1 & 1 \\ -1 & q & -1 & \pm 1 \\ -1 & q & -1 & -1 \\ 1 & \pm 1 & 1 & -1  \end{smallmatrix} \right)$; this
 algebra is not zero, see \texttt{CDrk4a.log}. If $\mu_{\theta-1\theta}=0$, then we apply Lemma \ref{lem:cut1} for $i=\theta-1$, $j=\theta$ and $\toba$ projects onto $\widetilde{\mathcal E}_{\bq'}$, $\bq'$ with two components, one of super type $A_{\theta-1}$ and another of type $A_1$. If $\lambda_{\theta-2\theta-2\theta-1}=0$, then we apply Lemma \ref{lem:cut1} for $i=\theta-2$, $j=\theta-1$ and $\toba$ projects onto $\widetilde{\mathcal E}_{\bq'}$, $\bq'$ with two components, one of super type $A_{\theta-2}$ and another of Cartan type $A_2$. If $\lambda_{\theta-2\theta-2\theta-3}=0$, then we apply Lemma \ref{lem:cut1} for $i=\theta-3$, $j=\theta-2$ and $\toba$ projects onto $\widetilde{\mathcal E}_{\bq'}$, $\bq'$ with two components, one of super type $A_{\theta-3}$ and the other is a subalgebra of \texttt{CDrk4a.log}.

Now fix $q_{\theta-2\theta-2}=-1$. A similar analysis reduces to the case $\nu_{\theta-2}\neq 0$, $\lambda_{(\theta-2\theta)}\neq 0$; here, we apply Lemma \ref{lem:cut1} for $i=\theta-4$, $j=\theta-3$ and $\toba$ projects onto $\widetilde{\mathcal E}_{\bq'}$, $\bq'$ with two components, one of super type $A_{\theta-4}$ and another with matrix $\left( \begin{smallmatrix}
-1 & \overline{q}a^{-1} & a^2 & -a^{-4} \\ a & -1 & a^{-1} & a^3 \\ a^{-2} & qa & -1 & a^{-2} \\ -a^4 & a^{-3} & -a^2 & -1  \end{smallmatrix} \right)$, $a\neq 0$; this
algebra is not zero, see \texttt{CDrk4b.log}.

For the diagram $\xymatrix@R-6pt{ 
&   \overset{-1}{\circ} \ar  @{-}[d]_{q^{-1}}\ar  @{-}[dr]^{q^2} &
\\
{\bf A}_{\theta-3}(q;i_1,\dots, i_j) \ar @{-}[r]_{\hspace*{1.5cm} q^{-1}}  &  \overset{q}{\underset{\ }{\circ}} \ar  @{-}[r]_{q^{-1}}  & \overset{-1}{\underset{\ }{\circ}} } $ there are 2 cases:

\noindent $\diamondsuit N\neq 4$.
Then $\toba$ is generated by $y_i$, $i\in\I_{\theta}$, with defining relations
\begin{align*}
& & y_{iii\pm1}=0, \quad  i <& \theta-1, q_{ii}\neq -1; \\
y_{\theta-2\theta-2\theta}&=0;  &
[y_{(i-1i+1)},y_i]_c&=\nu_i, \quad q_{ii}=-1, i\leq \theta-3;\\
y_i^2&=\mu_i, \quad q_{ii}=-1; &
y_{ij}&=\lambda_{ij}, \quad i < j - 1, i\neq \theta-2;
\end{align*}
\vspace*{-0.7cm}
\begin{align*}
&y_{(\theta-2\theta)}+q_{\theta-1\theta}(1+q^{-1})[y_{\theta-2\theta},y_{\theta-1}]_c-
q_{\theta-2\theta-1}(1-q^2)y_{\theta-1}y_{\theta-2\theta}=0.
\end{align*}
Here we apply Lemma \ref{lem:cut2} for $i=\theta-2$, $j=\theta-1$, $k=\theta$, and $\toba$ projects onto $ \widetilde{\mathcal E}_{\bq'}$, $\bq'$ with two components, of super type $A_{\theta-1}$ and of type $A_1$.

\noindent $\diamondsuit N=4$.
Then $\toba$ is generated by $y_i$, $i\in\I_{\theta}$, with defining
relations
\begin{align*}
y_i^2&=\mu_i, \, q_{ii}=-1; &
y_{ij} &= \lambda_{ij}, \,  i < j - 1, i\neq \theta-2;  \\
y_{\theta-2\theta-2\theta}&=\lambda_{\theta-2\theta-2\theta};&
y_{iii\pm1}&=\lambda_{iii\pm1}, \,   i < \theta-1, q_{ii}\neq -1; \\
y_{\theta-1\theta}^2&=\mu_{\theta-1\theta}; &
[y_{(i-1i+1)},y_i]_c&=\nu_i, \,  q_{ii}=-1, i\leq \theta-3;
\end{align*}
\vspace*{-0.7cm}
\begin{align*}
&y_{(\theta-2\theta)}+q_{\theta-1\theta}(1+q^{-1})[y_{\theta-2\theta},y_{\theta-1}]_c-
q_{\theta-2\theta-1}(1-q^2)y_{\theta-1}y_{\theta-2\theta}=0.
\end{align*}
Note that $\lambda_{\theta-2\theta-2\theta} \lambda_{\theta-2\theta-2\theta-1} =0$, so we may assume $\lambda_{\theta-2\theta-2\theta}=0$. If $\mu_{\theta-1\theta}=0$ (respectively, $\mu_{\theta-1\theta} \neq 0$, $\lambda_{\theta-2\theta-2\theta-1}=0$), then we apply Lemma \ref{lem:cut2} for $i=\theta-2$, $j=\theta-1$, $k=\theta$ (respectively $i=\theta-1$, $j=\theta-2$, $k=\theta$) and $\toba$ projects onto $\widetilde{\mathcal E}_{\bq'}$, $\bq'$ with two components, one of super type $A_{\theta-1}$ and another of type $A_1$. 
If $\mu_{\theta-1\theta}\lambda_{\theta-2\theta-2\theta-1} \neq 0$, $\lambda_{\theta-2\theta-2\theta-3}=0$, then we apply Lemma \ref{lem:cut1} for $i=\theta-1$, $j=\theta-3$ and $\toba-2$ projects onto $\widetilde{\mathcal E}_{\bq'}$, $\bq'$ with two components, one of super type $A_{\theta-2}$ and a subalgebra of that in \texttt{CDrk4c.log}.
If $\mu_{\theta-1\theta} \lambda_{\theta-2\theta-2\theta-1} \lambda_{\theta-2\theta-2\theta-3} \neq 0$, then  $\nu_{\theta-3}=0$ and we apply Lemma \ref{lem:cut1} for $i=\theta-4$, $j=\theta-3$: $\toba$ projects onto $\widetilde{\mathcal E}_{\bq'}$, $\bq'$ with two components, one of super type $A_{\theta-4}$ and another with matrix $\left( \begin{smallmatrix}
-1 & q & -1 & -1 \\ -1 & q & -1 & \pm 1 \\ -1 & q & -1 & 1 \\ -1 & \pm q^{-1} & -1 & -1  \end{smallmatrix} \right)$; this
algebra is not zero, see \texttt{CDrk4c.log}.

For the diagram $\xymatrix@R-6pt{
&   \overset{-1}{\circ} \ar  @{-}[d]_{q^{-1}}\ar  @{-}[dr]^{q^2} & 
\\	
	{\bf A}_{\theta-3}(q;i_1,\dots, i_j) \ar @{-}[r]_{\hspace*{1.5cm} q}  &
	\overset{-1}{\underset{\
		}{\circ}} \ar  @{-}[r]_{q^{-1}}  & \overset{-1}{\underset{\ }{\circ}}
}$ there are 2 cases:

\noindent $\diamondsuit N\neq 4$.
Then $\toba$ is generated by $y_i$, $i\in\I_{\theta}$, with defining
relations
\begin{align*}
&[y_{\theta-3\theta-2\theta},y_{\theta-2}]_c=\nu_{\theta-2}';&
& y_{iii\pm1}=\lambda_{iii\pm1}, \quad  i <\theta-2, q_{ii}\neq -1; \\
&& &[y_{(i-1i+1)},y_i]_c=\nu_i, \quad q_{ii}=-1, i\leq \theta-3;\\
&y_i^2=\mu_i, \quad q_{ii}=-1 &
&y_{ij} = \lambda_{ij}, \quad i < j - 1, i\neq \theta-2;
\end{align*}
\vspace*{-0.7cm}
\begin{align*}
&y_{(\theta-2\theta)}+q_{\theta-1\theta}(1+q^{-1})[y_{\theta-2\theta},y_{\theta-1}]_c-
q_{\theta-2\theta-1}(1-q^2)y_{\theta-1}y_{\theta-2\theta}=0.
\end{align*}
Here, $\nu_{\theta-2}\nu_{\theta-2}'=0$, we assume $\nu_{\theta-2}'=0$. Thus we apply Lemma \ref{lem:cut2} for $i=\theta-2$, $j=\theta-1$, $k=\theta$ and $\toba$ projects onto $\widetilde{\mathcal E}_{\bq'}$, $\bq'$ with two components, of super type $A_{\theta-1}$ and of type $A_1$.

\noindent $\diamondsuit N=4$.
Then $\toba$ is generated by $y_i$, $i\in\I_{\theta}$, with defining
relations
\begin{align*}
& y_{ij} = \lambda_{ij}, \quad i < j - 1, i\neq \theta-2; &
& y_{\theta-1 \theta}^2=\mu_{\theta-1 \theta} \\
& y_{iii\pm1}=\lambda_{iii\pm1}, \quad  i<\theta-2, q_{ii}\neq -1; &
& [y_{\theta-3\theta-2\theta},y_{\theta-2}]_c=\nu_{\theta-2}'; \\
&[y_{(i-1i+1)},y_i]_c=\nu_i, \quad q_{ii}=-1, i\leq \theta-3; &
& y_i^2=\nu_i, \quad q_{ii}=-1;
\end{align*}
\vspace*{-0.7cm}
\begin{align*}
&y_{(\theta-2\theta)}+q_{\theta-1\theta}(1+q^{-1})[y_{\theta-2\theta},y_{\theta-1}]_c-
2 q_{\theta-2\theta-1}y_{\theta-1}y_{\theta-2\theta}=0.
\end{align*}
Note that $\nu_{\theta-2} \lambda_{\theta-2}'=0$, so we may assume $\nu_{\theta-2}'=0$. If either $\mu_{\theta-1\theta}=0$ or else $\mu_{\theta-1\theta} \neq 0$, $\nu_{\theta-2}=0$, then we apply Lemma \ref{lem:cut2} and $\toba$ projects onto $\widetilde{\mathcal E}_{\bq'}$, $\bq'$ with two components, one of super type $A_{\theta-1}$ and another of type $A_1$. 
If $\mu_{\theta-1\theta}\nu_{\theta-2} \neq 0$, then $q_{\theta-3\theta-3}=-1$ and $\nu_{\theta-3}=0$; hence we apply Lemma \ref{lem:cut1} for $i=\theta-4$, $j=\theta-3$ and $\toba$ projects onto $\widetilde{\mathcal E}_{\bq'}$, $\bq'$ with two components, one of super type $A_{\theta-4}$ and another of Cartan type $C_4$.

For the diagram $\xymatrix@R-6pt{
&   \overset{q}{\circ} \ar  @{-}[d]^{q^{-1}} & 
\\
	{\bf A}_{\theta-3}(q;i_1,\dots, i_j) \ar @{-}[r]_{\hspace*{1.5cm} q^{-1}}  &
	\overset{q}{\underset{\
		}{\circ}} \ar  @{-}[r]_{q^{-1}}  & \overset{q}{\underset{\ }{\circ}}
}$ the algebra $\toba$ is generated by $y_i$, $i\in\I_{\theta}$, with defining
relations
\begin{align*}
y_i^2=\nu_i, \quad & q_{ii}=-1 &
y_{iii\pm1}&=\lambda_{iii\pm1}, \quad  i < \theta-1, q_{ii}\neq -1; \\
y_{\theta-2\theta-2\theta}&=\lambda_{\theta-2\theta-2\theta}; &
y_{\theta-1\theta}&=0;\\
y_{\theta-1\theta-1\theta-2}&=\lambda_{\theta-1\theta-1\theta-2}; &
y_{ij} &= \lambda_{ij}, \quad i < j - 1, i\neq \theta-2; \\
y_{\theta\theta\theta-2}&=\lambda_{\theta\theta\theta-2}; &
[y_{(i-1i+1)},y_i]_c&=\nu_i, \,  q_{ii}=-1, i\leq \theta-3.
\end{align*}
If $N\neq 3$, then $\lambda_{\theta-2\theta-2\theta}= \lambda_{\theta\theta\theta-2}= \lambda_{\theta-1\theta-1\theta-2}= \lambda_{\theta-2\theta-2\theta-1}=0$; if $N=3$, then either $\lambda_{\theta-2\theta-2\theta}= \lambda_{\theta\theta\theta-2}=0$, or else $\lambda_{\theta-1\theta-1\theta-2}= \lambda_{\theta-2\theta-2\theta-1}=0$. Thus we may assume $\lambda_{\theta-2\theta-2\theta}= \lambda_{\theta\theta\theta-2}=0$: We apply Lemma \ref{lem:cut1} for $i=\theta-2$, $j=\theta$ and $\toba$ projects onto $\widetilde{\mathcal E}_{\bq'}$, $\bq'$ with two components, one of super type $A_{\theta-1}$ and another of type $A_1$.

 For the diagram $\xymatrix@R-6pt{
&   \overset{q}{\circ} \ar  @{-}[d]^{q^{-1}} &
\\
{\bf A}_{\theta-3}(q;i_1,\dots, i_j) \ar @{-}[r]_{\hspace*{1.5cm} q}  &
\overset{-1}{\underset{\
}{\circ}} \ar  @{-}[r]_{q^{-1}}  & \overset{q}{\underset{\ }{\circ}} }$
the algebra $\toba$ is generated by $y_i$, $i\in\I_{\theta}$, with defining
relations
\begin{align*}
y_i^2=\mu_i, \,  & q_{ii}=-1;  &
y_{iii\pm1}&=\lambda_{iii\pm1}, \,   i < \theta-2, q_{ii}\neq -1; \\
[y_{\theta-3\theta-2\theta},y_{\theta-2}]_c&=\nu_{\theta-2}'; &
y_{\theta-1\theta}&=0;\\
y_{\theta-1\theta-1\theta-2}&=\lambda_{\theta-1\theta-1\theta-2}; &
y_{ij} &= \lambda_{ij}, \,  i < j - 1, i\neq \theta-2; \\
y_{\theta\theta\theta-2}&=\lambda_{\theta\theta\theta-2}; &
[y_{(i-1i+1)}&,y_i]_c=\nu_i, \,  q_{ii}=-1, i\leq \theta-2.
\end{align*}
As above, we may restrict to $\theta=4$ and $N=4$. As $\nu_{2}\nu_{2}'=0$, hence we may assume $\nu_{2}'=0$. Assume that $q_{11}=-1$, then $\nu_2=0$ and we apply Lemma \ref{lem:cut1} for $(i,j)=(1,2)$. Otherwise, $q_{11}=q^{-1}$. If only $\nu_2\neq0$, we can project onto smaller ranks by Lemma \ref{lem:cut3}. Hence, we can restrict to the case $\nu_2=0$ as described in \S \ref{subsubsec:casebycase}. We can further project to smaller rank unless $\lambda_{112}\lambda_{332}\lambda_{442}\neq 0$, when $\bq=\left( \begin{smallmatrix}
q^{-1} & -1 & a^{-1}q & b^{-1} \\ q^{-1} & -1 & q & q \\ a & -1 & q & c^{-1} \\ b & -1 & c & q \end{smallmatrix} \right)$, for $a,b,c\in\{\pm q\}$. In any case, $\toba\neq 0$, see \texttt{CDrk4d.log}.

\subsubsection{Type $\superda{\alpha}$, $q,r,s\neq 1$, $qrs=1$}\label{subsec:type-D2-1-alpha}

The Weyl groupoid has four objects.

For $\xymatrix@C-4pt{\overset{-1}{\underset{\ }{\circ}} \ar  @{-}[r]^{s} \ar@/^2pc/
@{-}[rr]^{q}
& \overset{-1}{\underset{\ }{\circ}} \ar  @{-}[r]^{r}& \overset{-1}{\underset{\
}{\circ}}}$, the algebra $\toba$ is generated by $y_i$, $i\in\I_3$, with defining
relations
\begin{align*}
y_1^2&=\lambda_1; & y_2^2&=\lambda_2; & y_3^2&=\lambda_3;
\end{align*}
\vspace*{-0.7cm}
\begin{align*}
y_{(13)}-\frac{1-s}{q_{23}(1-r)}[y_{13},y_2]_c-q_{12}(1-s)y_2y_{13}=0.
\end{align*}
This algebra is nonzero by \cite[Lemma 5.16]{AAGMV}.

For $\xymatrix{ \overset{q}{\underset{\ }{\circ}}\ar  @{-}[r]^{q^{-1}}  &
\overset{-1}{\underset{\ }{\circ}} \ar  @{-}[r]^{r^{-1}}  & \overset{r}{\underset{\
}{\circ}}}$, $q,r,s\neq -1$, the algebra $\toba$ is generated by $y_i$, $i\in\I_3$, with defining
relations
\begin{align*}
y_2^2&=\lambda_2; & y_{112}&=\lambda_1; & y_{332}&=\lambda_3; & y_{13}&=0,
\end{align*}
with $\lambda_1\lambda_3=0$. This algebra is nonzero by Lemma \ref{lem:cut3}applied to either $i=1$ or $i=3$, according to $\lambda_1=0$ or $\lambda_3=0$.

For the same diagram and $q=-1,r,s\neq -1$, $\toba$ is generated by $y_i$, $i\in\I_3$, with defining
relations
\begin{align*}
y_1^2&=\lambda_1; & y_2^2&=\lambda_2; &
y_{12}^2&=\lambda_4; & y_{332}&=\lambda_3; & y_{13}&=0.
\end{align*}
Here $\lambda_3\neq0$ only if $r\in\G_4$. If $\lambda_3=0$, then this algebra is nonzero as it projects onto an algebra of Cartan type $A_2$.
Similarly, a projection to the rank-two case shows that $\toba\neq0$ if $\lambda_1=\lambda_2=\lambda_4=0$.

Assume $\lambda_3\neq0$ and pick $\eps\in\k$ such that $\eps^2=1$, i.e.~$\eps=\pm r^2$. Then:
 $\lambda_1\lambda_4\neq 0$ only if $\bq=\left(\begin{smallmatrix}
-1&1& \eps \\-1&-1&r\\\eps &r^2&r
\end{smallmatrix}\right)$.
 $\lambda_2\neq 0$ only if $\bq=\left(\begin{smallmatrix}
-1&\eps & x \\-\eps &-1&r\\ x^{-1} &r^2&r
\end{smallmatrix}\right)$, with $x^2=\eps$.
In any case, $\toba\neq0$, see \texttt{d21a.log}.

\medbreak

The remaining two generalized Dynkin diagrams are of the shape of of the last one, but with $s$
instead of $q$, respectively with $s$ instead of $r$.

\subsubsection{Type $\superf$, $N > 3$}\label{subsec:type-F-super}
The Weyl groupoid has objects with generalized Dynkin diagrams of six possible shapes.

	 For the diagram $\xymatrix{ \overset{\,\, q^2}{\underset{\ }{\circ}}\ar  @{-}[r]^{q^{-2}}  & \overset{\,\,
			q^2}{\underset{\ }{\circ}} \ar  @{-}[r]^{q^{-2}}  & \overset{q}{\underset{\ }{\circ}} & \overset{-1}{\underset{\ }{\circ}} \ar  @{-}[l]_{q^{-1}}}$, there are 2 cases:
	
	\noindent $\diamondsuit N>4$.
	Then $\toba$ is generated by $y_i$, $i\in\I_4$, with defining
	relations
	\begin{align*}
	y_{13}&=0; & y_{14}&=0; & y_{24}&=0; \\
	y_{112}&=0; & y_{221}&=0; & y_{223}&=\lambda_{223};\\
	y_{334}&=0; & y_{3332}&=\lambda_{3332}; & y_{4}^2&=\mu_4.
	\end{align*}
	Here we apply Lemma \ref{lem:cut1} for $(i,j)=(3,4)$ and $\toba$ projects onto $\widetilde{\mathcal E}_{\bq'}$, $\bq'$ with two components, one of type $B_{3}$ and another of type $A_1$.

	\noindent $\diamondsuit N=4$.
	Then $\toba$ is generated by $y_i$, $i\in\I_4$, with defining
	relations
	\begin{align*}
	y_{13}&=0; & y_{14}&=\lambda_{14}; & y_{24}&=\lambda_{24}; \\
	y_{334}&=\lambda_{334}; & [y_{(13)},y_2]_c&=0; & [y_{23},y_{(24)}]_c&=\lambda_{(24)}; \\
	y_{3332}&=0; & y_{4}^2&=\mu_4; & y_{12}^2&=\mu_{12}.
	\end{align*}
If $\lambda_{(24)}=0$, then we apply Lemma \ref{lem:cut1} for $(i,j)=(2,3)$ and $\toba$ projects onto $ \widetilde{\mathcal E}_{\bq'}$, $\bq'$ with two components, one of type $A_2$ and another of super type $A_2$.
If $\lambda_{(24)}\neq0$ but all the other $\lambda$'s are 0, then we apply Lemma \ref{lem:cut3} for $i=1$ and $\toba$ projects onto $\widetilde{\mathcal E}_{\bq'}$, $\bq'$ of type super $C_3$. The general case follows by applying Lemma
\ref{lem:nonzero-k<l}.

	 For the diagram $\xymatrix{ \overset{\,\, q^2}{\underset{\ }{\circ}}\ar  @{-}[r]^{q^{-2}}  & \overset{\,\,
			q^2}{\underset{\ }{\circ}} \ar  @{-}[r]^{q^{-2}}  & \overset{-1}{\underset{\ }{\circ}} &
		\overset{-1}{\underset{\ }{\circ}} \ar  @{-}[l]_{q}}$, there are 2 cases:
	
	\noindent $\diamondsuit N>4$.
	Then $\toba$ is generated by $y_i$, $i\in\I_4$, with defining
	relations
	\begin{align*}
	y_{13}&=0; & y_{112}&=\lambda_{112}; & y_{24}&=0;  & y_{3}^2&=\mu_3; \\
	y_{14}&=0; & y_{221}&=\lambda_{221}; & y_{223}&=0; & y_{4}^2&=\mu_4; &
	[[y_{43},y_{432}]_c,y_3]_c&=0.
	\end{align*}
	Here we apply Lemma \ref{lem:cut1} for $i=2$, $j=\theta-3$ and $\toba$ projects onto $\widetilde{\mathcal E}_{\bq'}$, $\bq'$ with two components, one of type $A_2$ and another of super type $A_2$.

	\noindent $\diamondsuit N=4$.
	Then $\toba$ is generated by $y_i$, $i\in\I_4$, with defining
	relations
	\begin{align*}
	y_{13}&=\lambda_{13}; & y_{14}&=\lambda_{14}; &
	[y_{(13)},y_2]_c&=\nu_2; &  y_{3}^2&=\mu_3;  \\
	y_{23}^2&=\mu_{23}; & y_{24}&=\lambda_{24}; &
	[[y_{43},y_{432}]_c,y_3]_c&=\lambda_{432}; & y_{4}^2&=\mu_4.
	\end{align*}
	If $\nu_{2}=0$, then we apply Lemma \ref{lem:cut1} for $(i,j)=(1,2)$ and $\toba$ projects onto $ \widetilde{\mathcal E}_{\bq'}$, $\bq'$ with two components, one of type $A_1$ and another of super type $CD$. If $\nu_2\neq 0$, then $\lambda_{432}=0$ and we apply Lemma \ref{lem:cut1} for $(i,j)=(3,4)$: $\toba$ projects onto $\widetilde{\mathcal E}_{\bq'}$, $\bq'$ with two components, of types $A_3$ and $A_1$.

	 For the diagram $\xymatrix{ \overset{\,\, q^2}{\underset{\ }{\circ}}\ar  @{-}[r]^{q^{-2}}  & \overset{\,\,
			q^2}{\underset{\ }{\circ}} \ar  @{-}[r]^{q^{-2}}  & \overset{-1}{\underset{\ }{\circ}}\ar @{-}[r]^{q^3} & \overset{q^{-3}}{\underset{\ }{\circ}}}$, there are 3 cases:
	
	\noindent $\diamondsuit N\neq 4,6$.
	Then $\toba$ is generated by $y_i$, $i\in\I_4$, with defining
	relations
	\begin{align*}
	y_{443}&=\lambda_{443}; & y_{3}^2&=\mu_3; & y_{13}&=\lambda_{13}; & y_{14}&=0; \,  y_{24}=0; \\
	y_{112}&=0; & y_{221}&=0; & y_{223}&=\lambda_{223}; &
	[[[y_{432},y_3]_c,&[y_{4321},y_3]_c]_c,y_{32}]_c=0.
	\end{align*}
	Herewe apply Lemma \ref{lem:cut1} for $(i,j)=(1,2)$ and $\toba$ projects onto $\widetilde{\mathcal E}_{\bq'}$, $\bq'$ with two components, one of type $A_1$ and another of super type $D(2,1;\alpha)$.
	
	\noindent $\diamondsuit N=6$.
	Then $\toba$ is generated by $y_i$, $i\in\I_4$, with defining
	relations
	\begin{align*}
	y_{3}^2&=\mu_2; & y_{34}^2&=0;  & y_{13}&=0; & y_{14}&=0; \quad y_{24}=0; \\
	y_{112}&=\lambda_{112}; & y_{221}&=\lambda_{221}; & y_{223}&=0; &
	[[[ y_{432} & ,y_3]_c,[y_{4321},y_3]_c]_c,y_{32}]_c=0.
	\end{align*}
	Here, $\toba/ \langle y_{23} \rangle \simeq \widetilde{\mathcal E}_{\bq'}$, $\bq'$ with two components of type $A_2$.
	
	\noindent $\diamondsuit N=4$.
	Then $\toba$ is generated by $y_i$, $i\in\I_4$, with defining
	relations
	\begin{align*}
	y_{3}^2&=\mu_3;  & y_{13}&=\lambda_{13}; & y_{14}&=0; & [y_{(13)},y_2]_c&=\nu_2; \\
	y_{443}&=\lambda_{443}; &  y_{23}^2&=\mu_{23}; & y_{24}&=0; &
	[[[y_{432},y_3]_c,&[y_{4321},y_3]_c]_c,y_{32}]_c=0.
	\end{align*}
	If $\nu_{2}=0$, then we apply Lemma \ref{lem:cut1} for $(i,j)=(1,2)$ and $\toba$ projects onto $ \widetilde{\mathcal E}_{\bq'}$, $\bq'$ with two components, one of type $A_1$ and another of super type $CD$. If $\nu_2\neq 0$, then $\lambda_{443}=0$; now we apply Lemma \ref{lem:cut1} for $(i,j)=(3,4)$ and $\toba$ projects onto $ \widetilde{\mathcal E}_{\bq'}$, $\bq'$ with two components, of types $A_3$ and
	$A_1$.
	
	 For the diagram $\xymatrix{ \overset{\,\, q^2}{\underset{\ }{\circ}}\ar  @{-}[r]^{q^{-2}}  &
		\overset{q}{\underset{\ }{\circ}} \ar  @{-}[r]^{q^{-1}}  & \overset{-1}{\underset{\ }{\circ}}\ar  @{-}[r]^{q^3} & \overset{q^{-3}}{\underset{\ }{\circ}}}$, there are 3 cases:
	
	\noindent $\diamondsuit N\neq 4,6$.
	Then $\toba$ is generated by $y_i$, $i\in\I_4$, with defining
	relations
	\begin{align*}
	y_{3}^2&=\mu_3; &  y_{13}&=0; & y_{14}&=0; & y_{24}&=0; \\
	y_{443}&=\lambda_{443}; &  y_{112}&=\lambda_{221}; & y_{2221}&=\lambda_{2221}; & y_{223}&=0;
	\end{align*}
	\vspace*{-0.7cm}
	\begin{align*}
	[[y_{(14)},y_2]_c,y_3]_c-q_{23}(q^2-q)[[y_{(14)},y_3]_c,y_2]_c=0.
	\end{align*}
	Here we apply Lemma \ref{lem:cut1} for $(i,j)=(2,3)$ and $\toba$ projects onto $\widetilde{\mathcal E}_{\bq'}$, $\bq'$ with two components, one of type $B_2$ and another of super type $A_2$.
	
	\noindent $\diamondsuit N=6$.
	Then $\toba$ is generated by $y_i$, $i\in\I_4$, with defining
	relations
	\begin{align*}
	y_{3}^2&=\mu_3; & y_{13}&=0; & y_{14}&=0; & y_{24}&=0; \\
	y_{34}^2&=\mu_{34}; &  y_{112}&=0; & y_{2221}&=0; & y_{223}&=0;
	\end{align*}
	\vspace*{-0.7cm}
	\begin{align*}
	[[y_{(14)},y_2]_c,y_3]_c-q_{23}(q^2-q)[[y_{(14)},y_3]_c,y_2]_c=0.
	\end{align*}
	The proof is as for $N\neq 4,6$.
	
	\noindent $\diamondsuit N=4$.
	Then $\toba$ is generated by $y_i$, $i\in\I_4$, with defining
	relations
	\begin{align*}
	y_{3}^2&=\mu_3; &  y_{13}&=\lambda_{13}; & y_{14}&=0; & y_{24}&=0; \\
	y_{443}&=\lambda_{443}; &  [y_{12},y_{(13)}]_c&=\lambda_{(13)}; & y_{2221}&=0; & y_{223}&=\lambda_{223};
	\end{align*}
	\vspace*{-0.7cm}
	\begin{align*}
	[[y_{(14)},y_2]_c,y_3]_c-q_{23}(q^2-q)[[y_{(14)},y_3]_c,y_2]_c=0.
	\end{align*}
If $\lambda_{(13)}=0$, then we apply Lemma \ref{lem:cut1} for $(i,j)=(1,2)$ and $\toba$ projects onto $ \widetilde{\mathcal E}_{\bq'}$, $\bq'$ with two components, one of type $A_1$ and another of type $D(2,1;\alpha)$.
If $\lambda_{(13)}\neq0$ but all the other $\lambda$'s are 0, then we apply Lemma \ref{lem:cut1} for $(i,j)=(3,4)$ and $\toba$ projects onto $\widetilde{\mathcal E}_{\bq'}$, $\bq'$ with two components of types $B_3$ and $A_1$. The general case follows by applying Lemma \ref{lem:nonzero-k<l}.

	 For the diagram $\xymatrix@R-6pt@C-4pt{ & & \overset{q}{\circ}\ar  @{-}[ld]_{q^{-1}} \ar
		@{-}[d]^{q^{-1}}
		\\ \overset{q^2}{\underset{\ }{\circ}} \ar  @{-}[r]^{q^{-2}}  &  \overset{-1}{\underset{\ }{\circ}} \ar  @{-}[r]^{q^{2}}  & \overset{-1}{\underset{\ }{\circ}}}$ there are 2 cases:
	
	\noindent $\diamondsuit N>4$.
	Then $\toba$ is generated by $y_i$, $i\in\I_4$, with defining
	relations
	\begin{align*}
	y_{2}^2&=\mu_2; &  y_{13}&=0; & y_{14}&=0; & y_{112}&=\lambda_{112}; \\
	y_{3}^2&=\mu_3; &  y_{442}&=0; & y_{443}&=0; & [y_{(13)},y_2]_c&=0;
	\end{align*}
	\vspace*{-0.7cm}
	\begin{align*}
	y_{(24)}- q_{34}q[y_{24},y_3]_c-q_{23}(1-q^{-1})y_3y_{24}=0.
	\end{align*}
	Here we apply Lemma \ref{lem:cut2} for $(i,j,k)=(2,3,4)$ and $\toba$ projects onto $\widetilde{\mathcal E}_{\bq'}$, $\bq'$ with two components of super type $A_2$.

\noindent $\diamondsuit N=4$.
Then $\toba$ is generated by $y_i$, $i\in\I_4$, with defining relations
\begin{align*}
y_{2}^2&=\mu_2; &  y_{13}&=\lambda_{13}; & y_{14}&=0; & y_{12}^2&=\mu_{12}; \\
y_{3}^2&=\mu_3; &  y_{442}&=\lambda_{442}; & y_{443}&=\lambda_{443}; & [y_{(13)},y_2]_c&=\nu_2;
\end{align*}
\vspace*{-0.7cm}
\begin{align*}
y_{(24)}- q_{34}q[y_{24},y_3]_c-q_{23}(1-q^{-1})y_3y_{24}=0.
\end{align*}
If $\lambda_{442}\neq 0$, then $\nu_{2}=\lambda_{443}=0$, so we apply Lemma \ref{lem:cut2} for $(i,j,k)=(3,2,4)$ and $\toba$ projects onto $\widetilde{\mathcal E}_{\bq'}$, $\bq'$ with two components, of types $A_1$ and $C_3$.
If $\lambda_{443}\neq 0$, then $\nu_{2}=\lambda_{442}=0$, so we apply Lemma \ref{lem:cut2} for $(i,j,k)=(2,3,4)$ and $\toba$ projects onto $ \widetilde{\mathcal E}_{\bq'}$, $\bq'$ with two components, of types $A_2$ and super $A_2$.
	If $\nu_{2}\neq 0$, then $\lambda_{442}=\lambda_{443}=0$, so we apply Lemma \ref{lem:cut1} for $(i,j,k)=(4,3,2)$ and $\toba$ projects onto $\widetilde{\mathcal E}_{\bq'}$, $\bq'$ with two components, of types $A_1$ and $A_3$.
	
For the diagram $\xymatrix@R-6pt@C-4pt{ & & \overset{-1}{\circ}\ar  @{-}[ld]_{q^{2}} \ar  @{-}[d]^{q^{-3}}
		\\ \overset{q^2}{\underset{\ }{\circ}} \ar  @{-}[r]^{q^{-2}}  &  \overset{-1}{\underset{\
			}{\circ}} \ar  @{-}[r]^{q}  & \overset{-1}{\underset{\ }{\circ}}}$ there are 2 cases:
	
	\noindent $\diamondsuit N>4$.
	Then $\toba$ is generated by $y_i$, $i\in\I_4$, with defining
	relations
	\begin{align*}
	y_{2}^2&=\mu_2; & y_{13}&=0; & y_{14}&=0; &   [y_{124},y_2]_c&=0; \\
	y_{3}^2&=\mu_3; &   y_{112}&=\lambda_{112}; & y_{4}^2&=\mu_4; &  [[y_{32},y_{321}]_c,y_2]_c&=0;
	\end{align*}
	\vspace*{-0.7cm}
	\begin{align*}
	y_{(24)}+ q_{34}\frac{1-q^3}{1-q^2}[y_{24},y_3]_c-q_{23}(1-q^{-3})y_3y_{24}=0.
	\end{align*}
	Here we apply Lemma \ref{lem:cut2} for $(i,j,k)=(2,3,4)$ and $\toba$ projects onto $\widetilde{\mathcal E}_{\bq'}$, $\bq'$ with two components of super type $A_2$.
	
	\noindent $\diamondsuit N=4$.
	Then $\toba$ is generated by $y_i$, $i\in\I_4$, with defining
	relations
	\begin{align*}
	y_{2}^2&=\mu_2; & y_{13}&=\lambda_{13}; & y_{14}&=\lambda_{14}; &   [y_{124},y_2]_c&=\nu_2; \\
	y_{3}^2&=\mu_3; &   y_{12}^2&=\mu_{12}; & y_{4}^2&=\mu_4; &  [[y_{32},y_{321}]_c,y_2]_c&=\lambda_{(321)};
	\end{align*}
	\vspace*{-0.7cm}
	\begin{align*}
	y_{(24)}+ q_{34}\frac{1-q^3}{1-q^2}[y_{24},y_3]_c-q_{23}(1-q^{-3})y_3y_{24}=0.
	\end{align*}
	If $\lambda_{(321)}=0$, then we apply Lemma \ref{lem:cut2} for $(i,j,k)=(3,2,4)$ and $\toba$ projects onto $\widetilde{\mathcal E}_{\bq'}$, $\bq'$ with two components, of types $A_1$ and $A_3$.
	If $\lambda_{(321)}\neq 0$, then $\nu_{2}=0$, so we apply Lemma \ref{lem:cut2} for $(i,j,k)=(4,2,3)$ and $\toba$ projects onto $\widetilde{\mathcal E}_{\bq'}$, $\bq'$ with two components, of types $A_1$ and $A_3$.

\subsubsection{Type  $\superg$, $N > 3$}\label{subsec:type-G-super}
The Weyl groupoid has four objects.

For $\xymatrix{ \overset{-1}{\underset{\ }{\circ}}\ar  @{-}[r]^{q^{-1}}  &
\overset{q}{\underset{\ }{\circ}} \ar  @{-}[r]^{q^{-3}}  & \overset{\,\, q^3}{\underset{\
}{\circ}}}$, $N\neq 4,6$, the algebra $\toba$ is generated by $y_i$, $i\in\I_3$, with defining
relations
\begin{align*}
y_{13}&=0; & y_{221}&=0; &  y_{332}&=0; &
 y_{22223}&=\lambda_2; & y_{1}^2&=\lambda_1,
\end{align*}
with $\lambda_1\lambda_2=0$. If $\lambda_1\neq 0$, then $\toba$ projects onto $\k\langle y_1|y_1^2=\lambda_1\rangle$ by Lemma \ref{lem:cut3} applied twice, to $i=2,3$. If $\lambda_2\neq0$, then $\toba$ projects onto a algebra related to Cartan case $G_2$ by the same Lemma applied to $i=1$.

When $N =6$, $\toba$ is generated by $y_i$, $i\in\I_3$, with defining
relations
\begin{align*}
y_{13}&=\lambda_{3}; & y_{221}&=0; &  [y_{12},y_{(13)}]_c&=0; &
 y_{22223}&=0; & y_{1}^2&=\lambda_1.
\end{align*}
Now $\toba\neq 0$ by applying Lemma \ref{lem:cut3} for $i=2$.

For $N =4$, $\toba$ is generated by $y_i$, $i\in\I_3$, with defining
relations
\begin{align*}
y_{13}&=0; & y_{221}&=\lambda_2; &  [[[y_{(13)},&y_2]_c,y_2]_c,y_2]_c=0; &
 y_{332}&=0; & y_{1}^2&=\lambda_1.
\end{align*}
The longest relation is undeformed by \texttt{g3aN4.log}. Here $\toba$ projects onto a nonzero algebra of Cartan case $A_2$, by Lemma \ref{lem:cut3} applied to $i=3$.

\medskip

For $\xymatrix{ \overset{-1}{\underset{\ }{\circ}} \ar  @{-}[r]^{q}  &
\overset{-1}{\underset{\ }{\circ}} \ar  @{-}[r]^{q^{-3}}  & \overset{\,\, q^3}{\underset{\
}{\circ}}}$, $N\neq 6$, the algebra $\toba$ is generated by $y_i$, $i\in\I_3$, with defining
relations
\begin{align*}
y_{13}&=0; & y_{332}&=\lambda_3; &  [[y_{12},&[y_{12},y_{(13)}]_c]_c,y_2]_c=0; &
 y_{1}^2&=\lambda_1; & y_{2}^2&=\lambda_2.
\end{align*}
The third relation is undeformed by \texttt{g3b.log}. Here $\lambda_3\neq 0$ only when $N=12$. In any case $\toba\neq0$ by Lemma \ref{lem:cut1} applied to $(i,j)=(1,2)$.

When $N=6$, $\toba$ is generated by $y_i$, $i\in\I_3$, with defining
relations
\begin{align*}
y_{13}&=\lambda_3; & y_{23}^2&=\lambda_5; &  [[y_{12},&[y_{12},y_{(13)}]_c]_c,y_2]_c=\lambda_4; &
 y_{1}^2&=\lambda_1; & y_{2}^2&=\lambda_2.
\end{align*}
If $\lambda_4=0$, then we apply again Lemma \ref{lem:cut1} for $(i,j)=(1,2)$.
If $\lambda_4\neq 0$, then $\bq=\left(\begin{smallmatrix}
-1      & a  & -a^{-4}\\ 
a^{-1}q &-1  & -a^3\\ 
-a^4&  a^{-3}  &-1
\end{smallmatrix}\right)$. Here, $\toba\neq 0$ by \texttt{g3b2.log}.

For $\xymatrix{ \overset{-q^{-1}}{\underset{\ }{\circ}} \ar  @{-}[r]^{q^2}  &
\overset{-1}{\underset{\ }{\circ}} \ar  @{-}[r]^{q^{-3}}  & \overset{\,\, q^3}{\underset{\
}{\circ}}}$, $N\neq 6$, the algebra $\toba$ is generated by $y_i$, $i\in\I_3$, with defining
relations
\begin{align*}
y_{13}&=\lambda_3; \qquad y_{332}=\lambda_4; \qquad y_{1112}=0; \qquad y_{2}^2=\lambda_2; \\
[y_1, & [y_{123},y_2]_c]_c =\frac{q_{12}q_{32}}{1+q}[y_{12},y_{123}]_c-q^{-1}(1-q^{-1})q_{12}q_{13} y_{123}y_{12}.
\end{align*}
Here $\lambda_3\neq 0$ only if $N=4$ and $\lambda_4\neq0$ only if $N=12$; also $\lambda_2\lambda_3=0$. Hence, we can project onto a case of smaller rank and  $\toba\neq0$. The longest relation is undeformed by \texttt{g3c.log}.

When $N =6$, $\toba$ is generated by $y_i$, $i\in\I_3$, with defining
relations
\begin{align*}
y_{13}&=0; \qquad y_{23}^2=\lambda_3; \qquad [y_{112},y_{12}]_c=\lambda_1; \qquad y_{2}^2=\lambda_2;  \\
[y_1, & [y_{123},y_2]_c]_c = \frac{q_{12}q_{32}}{1+q}[y_{12},y_{123}]_c-q^{-1}(1-q^{-1})q_{12}q_{13} y_{123}y_{12}.
\end{align*}
The longest relation is undeformed by \texttt{g3c.log}. If $\lambda_1=0$, then we apply Lemma \ref{lem:cut1} for $(i,j)=(1,2)$ and $\toba\neq 0$. If $\lambda_1\neq 0$ but $\lambda_2=\lambda_{3}=0$, then we apply Lemma \ref{lem:cut3} for $i=3$. The general case follows from Lemma \ref{lem:nonzero-k<l}.

\medskip

For $\xymatrix@C-4pt{\overset{q}{\underset{\ }{\circ}} \ar  @{-}[r]^{q^{-2}} \ar@/^2pc/
@{-}[rr]^{q^{-1}}
& \overset{-1}{\underset{\ }{\circ}} \ar  @{-}[r]^{q^{3}}& \overset{-1}{\underset{\
}{\circ}}}$, the algebra $\toba$ is generated by $y_i$, $i\in\I_3$, with defining
relations
\begin{align*}
y_{1112}&=\lambda_1; & y_{2}^2&=\lambda_2; & y_3^2&=\lambda_3; &
 y_{113}&=\lambda_4;
\end{align*}
\vspace*{-0.6cm}
\begin{align*}
y_{(13)}+q^{-2}q_{23}\frac{1-q^3}{1-q}[y_{13},y_2]_c-q_{12}(1-q^3)y_2y_{13}=0.
\end{align*}
If $\lambda_1=\lambda_4=0$, then $\toba\neq 0$ by \cite[Lemma 5.16]{AAGMV}. Now, $\lambda_1\neq0$ only if $N=6$ and $\lambda_4\neq0$ only if $N=4$.
If $\lambda_1\neq 0$ and $\lambda_3=0$ or if $\lambda_4\neq0$ and $\lambda_2=0$ we may project onto a case of smaller rank. Indeed, $\lambda_2\lambda_4\neq 0$ and we can have
$\lambda_1\lambda_3\neq 0$ only if
$\bq
=\left(\begin{smallmatrix}
q&-1&\pm1\\ q&-1&\pm1\\ \pm q^{-1}& \pm q^{-1} & -1
\end{smallmatrix}\right)$. In this case, $\toba\neq0$ by \texttt{g3d.log}.

\subsection{Modular type, characteristic
3}\label{sec:by-diagram-modular-char3}

\subsubsection{Type $\br(2, a)$, $\zeta \in \G_3$, $q\notin \G_3$}\label{subsec:type-br(2,a)}
The Weyl groupoid has two objects: $\xymatrix{ \overset{\zeta}{\underset{\ }{\circ}} \ar  @{-}[r]^{q^{-1}}  &
\overset{q}{\underset{\ }{\circ}}}$ and $\xymatrix{ \overset{\zeta}{\underset{\ }{\circ}} \ar  @{-}[r]^{\zeta^2q}  &
\overset{\zeta q^{-1}}{\underset{\ }{\circ}}}$. Notice that the right hand diagram is of the shape of the left hand one, with $\zeta q^{-1}$
instead of $q$. Therefore, we only have to analyze one of them: we concentrate on the later.

If $q\neq -1$, then the algebra  $\toba$ is generated by $y_1, y_{2}$ with defining relations
\begin{align*}
y_{1}^3 &= \lambda_1; &  y_{221} &=\lambda_2,
\end{align*}
with $\lambda_1\lambda_2=0$. Also, $\lambda_2=0$ unless $q=-\zeta$ and $\bq=\left(\begin{smallmatrix}
\zeta& -\zeta \\\zeta& -\zeta
\end{smallmatrix}\right)$.
This algebra is nonzero: If $\lambda_1\neq 0$, this is \cite[Lemma 5.16]{AAGMV}. When $\lambda_2\neq0$, see \texttt{br2a-q.log}.

If $q=-1$, then the algebra $\toba$ is generated by $y_1, y_{2}$
with defining relations
\begin{align*}
y_{1}^3 &= \lambda_1; & y_{2}^{2} &=\lambda_2; & [y_{112}, y_{12}]_c&=\lambda_3,
\end{align*}
with, $\lambda_3=0$ unless $\bq=\left(\begin{smallmatrix}
\zeta& -1 \\ 1& -1
\end{smallmatrix}\right)$. 
This algebra is nonzero, see \texttt{br2a-1.log}.

\subsubsection{Type $\br(3)$, $\zeta \in \G'_{9}$}\label{subsec:type-br(3)}
The Weyl groupoid has two objects.

For $\xymatrix{ \overset{\zeta}{\underset{\ }{\circ}}\ar  @{-}[r]^{\ztu}  &
\overset{\zeta}{\underset{\ }{\circ}} \ar  @{-}[r]^{\ztu}  & \overset{\ztu^{\,
3}}{\underset{\ }{\circ}}}$, the algebra $\toba$ is generated by $y_i$, $i\in\I_3$, with defining
relations
\begin{align*}
y_{13}&=0; & y_{112}&=0; & y_3^3&=\lambda_1; & [[y_{332},& y_{3321}]_c,y_{32}]_c=0; &
 y_{221}&=0; & y_{223}&=0.
\end{align*}
This algebra is nonzero by \cite[Lemma 5.16]{AAGMV}.

For $\xymatrix{ \overset{\zeta}{\underset{\ }{\circ}}\ar  @{-}[r]^{\ztu}  &
\overset{\ztu^{\, 4}}{\underset{\ }{\circ}} \ar  @{-}[r]^{\zeta^{4}}  & \overset{\ztu^{\,
3}}{\underset{\ }{\circ}}}$, the algebra $\toba$ is generated by $y_i$, $i\in\I_3$, with defining
relations
\begin{align*}
y_{13}&=0; & y_{112}&=0; & y_3^3&=\lambda_1; &
 y_{2221}&=0; & y_{223}&=0;
\end{align*}
\vspace*{-0.7cm}
\begin{align*}
 [[y_{(13)}, y_{2}]_c,y_{3}]_c-(1+\zeta^4)^{-1}q_{23}[[y_{(13)}, y_{3}]_c,y_{2}]_c=0.
\end{align*}
This algebra is nonzero by \cite[Lemma 5.16]{AAGMV}.

\subsection{Super modular type, characteristic
3}\label{sec:by-diagram-super-modular-char3}

\subsubsection{\ Type $\brj(2; 3)$, $\zeta \in \G'_9$}\label{subsec:type-brj(2,3)}
The Weyl groupoid has three objects.

For $\xymatrix{ \overset{-\zeta}{\underset{\ }{\circ}} \ar  @{-}[r]^{\ztu^{\, 2}}  &
\overset{\zeta^3}{\underset{\ }{\circ}}}$, $\toba$ is generated by $y_1, y_2$ with defining relations
\begin{align*}
y_2^3&=\lambda_2; &  [y_1,[y_{12},y_2]_c]_c+\frac{\zeta^5(1-\zeta)q_{12}}{1-\zeta^7}
y_{12}^2&=0;&
 y_{1112}&=0.
\end{align*}
This algebra is nonzero as it projects onto $\k\langle y_2|y_2^3=\lambda_2\rangle$.

For $\xymatrix{ \overset{\zeta^3}{\underset{\ }{\circ}} \ar  @{-}[r]^{\ztu}  &
\overset{-1}{\underset{\ }{\circ}}}$, $\toba$ is generated by $y_1, y_2$ with defining relations
\begin{align*}
y_1^{3}&=\lambda_1; & y_2^2&=\lambda_2; &  [y_{112},[y_{112},y_{12}]_c]_c&=0,
\end{align*}
with $\lambda_1\lambda_2=0$. The last relation is undeformed, see \texttt{brj23.log}.
Hence  this algebra is nonzero by \cite[Lemma 5.16]{AAGMV}.

For $\xymatrix{ \overset{-\zeta^2}{\underset{\ }{\circ}} \ar  @{-}[r]^{\zeta}  &
\overset{-1}{\underset{\ }{\circ}}}$, $\toba$ is generated by $y_1, y_2$ with defining relations
\begin{align*}
y_2^2&=\lambda_2; &  [y_{112},y_{12}]_c&=0; &
 y_{111112}&=0.
\end{align*}
This algebra is nonzero by \cite[Lemma 5.16]{AAGMV}.

\subsubsection{\ Type $\g(1,6)$, $\zeta \in \G'_{3} \cup \G'_{6}$}\label{subsec:type-g(1,6)}
The Weyl groupoid has two objects.

	 For $\xymatrix{ \overset{\zeta}{\underset{\ }{\circ}}\ar  @{-}[r]^{\ztu}  &
		\overset{\zeta}{\underset{\ }{\circ}} \ar  @{-}[r]^{\ztu^{\, 2}}  &
		\overset{-1}{\underset{\ }{\circ}}}$, $\zeta \in \G'_{6}$, the algebra $\toba$ is generated by $y_i$, $i\in\I_3$, with defining relations
	\begin{align*}
	y_{13}&=0; &y_{112}&=0; &  y_{221}&=0; &
	y_{2223}&=\lambda_{2223}; &y_{3}^2&=\mu_3.
	\end{align*}
	Here we apply Lemma \ref{lem:cut1} for $(i,j)=(1,2)$ and $\toba$ projects onto $\widetilde{\mathcal E}_{\bq'}$, $\bq'$ with two components of super type $A_2$ and $A_1$.
	
	 For $\xymatrix{ \overset{\zeta}{\underset{\ }{\circ}}\ar  @{-}[r]^{\ztu}  &
		\overset{-\ztu}{\underset{\ }{\circ}} \ar  @{-}[r]^{\zeta^{2}}  & \overset{-1}{\underset{\
			}{\circ}}}$, $\zeta \in \G'_{6}$, the algebra $\toba$ is generated by $y_i$, $i\in\I_3$, with defining relations
	\begin{align*}
	y_{13}&=0; & y_{112}&=0;  & y_{3}^2&=\mu_3; & [y_{223},y_{23}]_c&=\lambda_{23};
	\end{align*}
	\vspace*{-0.7cm}
	\begin{align*}
	[y_2, [y_{21},y_{23}]_c]_c+q_{13}q_{23}q_{21}[y_{223},y_{21}]_c+q_{21}y_{21}y_{223}=0.
	\end{align*}
	Here we apply Lemma \ref{lem:cut1} for $(i,j)=(1,2)$ and $\toba$ projects onto $ \widetilde{\mathcal E}_{\bq'}$, $\bq'$ with two components of super type $B_2$ and $A_1$.

	 For $\xymatrix{ \overset{\zeta}{\underset{\ }{\circ}}\ar  @{-}[r]^{\ztu}  &
		\overset{\zeta}{\underset{\ }{\circ}} \ar  @{-}[r]^{\ztu^{\, 2}}  &
		\overset{-1}{\underset{\ }{\circ}}}$, $\zeta \in \G'_{3}$, the algebra $\toba$ is generated by $y_i$, $i\in\I_3$, with defining relations
	\begin{align*}
	y_{13}&=0; & y_{112}&=\lambda_{112}; &  [[y_{(13)},&y_2]_c,y_2]_c=0;\\
	y_{3}^2&=\mu_3; & y_{221}&=\lambda_{221}; & [y_{223},y_{23}]_c&=\lambda_{23}.
	\end{align*}
	If $\lambda_{23}=0$, then we apply Lemma \ref{lem:cut1} for $(i,j)=(2,3)$ and $\toba$ projects onto $ \widetilde{\mathcal E}_{\bq'}$, $\bq'$ with two components, of types $A_2$ and $A_1$.
	If $\lambda_{23}\neq 0$, then $\lambda_{112}=\lambda_{221}=0$, so we apply Lemma \ref{lem:cut1} for $(i,j)=(1,2)$ and $\toba$ projects onto $\widetilde{\mathcal E}_{\bq'}$, $\bq'$ with two components, of types $A_1$ and super $B_2$.

	 For $\xymatrix{ \overset{\zeta}{\underset{\ }{\circ}}\ar  @{-}[r]^{\ztu}  &
		\overset{-\ztu}{\underset{\ }{\circ}} \ar  @{-}[r]^{\zeta^{2}}  & \overset{-1}{\underset{\
			}{\circ}}}$, $\zeta \in \G'_{3}$, the algebra $\toba$ is generated by $y_i$, $i\in\I_3$, with defining relations
	\begin{align*}
	y_{13}&=0; &y_{112}&=0; &  y_{2221}&=0; &
	y_{2223}&=\lambda_{2223}; &y_{3}^2&=\mu_3;
	\end{align*}
	\vsp
	\begin{align*}
	[y_{1},y_{223}]_c+q_{23}[y_{(13)},y_{2}]_c+(\zeta^2-\zeta)q_{12}y_{2}y_{(13)}=0.
	\end{align*}
	Here we apply Lemma \ref{lem:cut1} for $(i,j)=(1,2)$ and $\toba$ projects onto $\widetilde{\mathcal E}_{\bq'}$, $\bq'$ with two components of super type $B_2$ and $A_1$.

\subsubsection{\ Type $\g(2,3)$, $\zeta \in \G'_{3}$}\label{subsec:type-g(2,3)}
The Weyl groupoid has four objects.

	 For $\xymatrix{ \overset{-1}{\underset{\ }{\circ}}\ar  @{-}[r]^{\ztu}  &
		\overset{\zeta}{\underset{\ }{\circ}} \ar  @{-}[r]^{\zeta}  & \overset{-1}{\underset{\
			}{\circ}}}$, the algebra $\toba$ is generated by $y_i$, $i\in\I_3$, with defining relations
	\begin{align*}
	y_{13}&=\lambda_{13}; & y_{221}&=0; & y_1^2&=\mu_1; &  y_3^2&=\mu_3;
	\end{align*}
	\vspace*{-0.7cm}
	\begin{align*}
	[y_{223},y_{23}]_c&=\lambda_{23}; &
	[[y_{(13)},y_2]_c,y_2]_c&=\lambda_{(13)}.
	\end{align*}
	If $\lambda_{(13)}=0$, then we apply Lemma \ref{lem:cut1} for $(i,j)=(1,2)$ and $\toba$ projects onto $ \widetilde{\mathcal E}_{\bq'}$, $\bq'$ with two components, of types $A_1$ and super $B_2$.
	If $\lambda_{(13)}\neq 0$, then $\bq=\left( \begin{smallmatrix}
	-1 & b & -b^{-3} \\ \zeta^2b^{-1} & \zeta & \zeta b \\ -b^3 & b^{-1} & -1 \end{smallmatrix} \right)$, $b\neq 0$; this
	algebra is not zero, see \texttt{g23-1.log}.

	 For $\xymatrix{ \overset{-1}{\underset{\ }{\circ}}\ar  @{-}[r]^{\zeta}  &
		\overset{-1}{\underset{\ }{\circ}} \ar  @{-}[r]^{\zeta}  & \overset{-1}{\underset{\
			}{\circ}}}$, the algebra $\toba$ is generated by $y_i$, $i\in\I_3$, with defining relations
	\begin{align*}
	y_{13}&=\lambda_{13}; & y_{1}^2&=\mu_1; & y_2^2&=\mu_2; &  y_3^2&=\mu_3;
	\end{align*}
	\vspace*{-0.7cm}
	\begin{align*}
	[[y_{12},y_{(13)}]_c,y_2]_c&=\lambda_{123}; & [[y_{32},y_{321}]_c,y_2]_c&=\lambda_{321}.
	\end{align*}
	If $\lambda_{123}=\lambda_{321}=0$, then we apply Lemma \ref{lem:cut1} for $(i,j)=(1,2)$ and $\toba$ projects onto $\widetilde{\mathcal E}_{\bq'}$, $\bq'$ with two components, of types $A_1$ and super $A_2$.
	If $\lambda_{123}\neq 0$ (similar for $\lambda_{321}\neq 0$), then $\bq=\left( \begin{smallmatrix}
	-1 & b & b^{-3} \\ \zeta b^{-1} & -1 & -\zeta b^2 \\ b^3 & -b^{-2} & -1 \end{smallmatrix} \right)$, $b\neq 0$; this
	algebra is not zero, see \texttt{g23-2.log}.  Observe that $\lambda_{123}\lambda_{321}\neq 0$ only if $b^3=-1$.

	 For $\xymatrix{ \overset{-1}{\underset{\ }{\circ}}\ar  @{-}[r]^{\ztu}  &
		\overset{-\ztu}{\underset{\ }{\circ}} \ar  @{-}[r]^{\ztu}  & \overset{-1}{\underset{\
			}{\circ}}}$, the algebra $\toba$ is generated by $y_i$, $i\in\I_3$, with defining relations
	\begin{align*}
	y_{13}&=\lambda_{13}; & y_{2221}&=\lambda_{2221}; & y_{2223}&=\lambda_{2223}; &  y_1^2&=\mu_1; & y_3^2&=\mu_3;
	\end{align*}
	\vsp
	\begin{align*}
	[y_{1},y_{223}]_c+q_{23}[y_{(13)},y_2]_c-(1-\zeta)q_{12}y_2y_{(13)}=0.
	\end{align*}
	As $\lambda_{2221}\lambda_{2223}=0$, we may assume $\lambda_{2221}=0$. Hence we apply Lemma \ref{lem:cut1} for $(i,j)=(1,2)$ and $\toba$ projects onto $\widetilde{\mathcal E}_{\bq'}$, $\bq'$ with two components, of types super $B_2$ and $A_1$.
	
	 For $\xymatrix{ \overset{\zeta}{\underset{\ }{\circ}} \ar  @{-}[r]^{\ztu} \ar@/^2pc/
		@{-}[rr]^{\ztu}  & \overset{-1}{\underset{\ }{\circ}} \ar  @{-}[r]^{\ztu}
		& \overset{\zeta}{\underset{\ }{\circ}}  }$, the algebra $\toba$ is generated by $y_i$, $i\in\I_3$, with defining relations
	\begin{align*}
	y_{112}&=0; & y_{113}&=0; & y_{331}&=0; &  y_{332}&=0; & y_2^2&=\mu_2;
	\end{align*}
	\vsp
	\begin{align*}
	y_{(13)} -
	q_{23}\zeta [y_{13},y_{2}]_c-q_{12}(1-\ztu)y_2y_{13}=0.
	\end{align*}
	Here we apply Lemma \ref{lem:cut2} for $(i,j,k)=(2,1,3)$ and $\toba$ projects onto $\widetilde{\mathcal E}_{\bq'}$, $\bq'$ with two components, of types $A_2$ and $A_1$.

\subsubsection{\ Type $\g(3, 3)$, $\zeta \in \G'_3$}\label{subsec:type-g(3,3)}
The Weyl groupoid has six objects.

	 For $\xymatrix{ \overset{\ztu}{\underset{\ }{\circ}}\ar  @{-}[r]^{\zeta}  &
		\overset{\ztu}{\underset{\ }{\circ}} \ar  @{-}[r]^{\zeta}  & \overset{\zeta}{\underset{\
			}{\circ}}
		\ar  @{-}[r]^{\ztu}  & \overset{-1}{\underset{\ }{\circ}}}$
	the algebra $\toba$ is generated by $y_i$, $i\in\I_4$, with defining relations
	\begin{align*}
	y_{223}&=0; & y_{13}&=\lambda_{13}; &y_{14}&=0; & y_{112}&=\lambda_{112}; & [y_{3321},&y_{32}]_c=\lambda_{321};\\
	y_{334}&=0; &  y_{4}^2&=\mu_4 & y_{24}&=0; & y_{221}&=\lambda_{221}; & [[y_{(24)},&y_{3}]_c,y_3]_c=0.
	\end{align*}
	Here we apply Lemma \ref{lem:cut1} for $(i,j)=(3,4)$ and $\toba$ projects onto $\widetilde{\mathcal E}_{\bq'}$, $\bq'$ with two components, of types $B_3$ and $A_1$.

	 For $ \xymatrix{ \overset{\ztu}{\underset{\ }{\circ}}\ar  @{-}[r]^{\zeta}  &
		\overset{\ztu}{\underset{\ }{\circ}} \ar  @{-}[r]^{\zeta}  & \overset{-1}{\underset{\
			}{\circ}}
		\ar  @{-}[r]^{\zeta}  & \overset{-1}{\underset{\ }{\circ}}}$
	the algebra $\toba$ is generated by $y_i$, $i\in\I_4$, with defining relations
	\begin{align*}
	y_{13}&=0; & y_{3}^2&=\mu_3; &y_{14}&=0; & y_{112}&=\lambda_{112}; & [[y_{43},&y_{432}]_c,y_3]_c=0;\\
	& & y_{4}^2&=0; & y_{24}&=0; & y_{221}&=\lambda_{221}; & y_{223}&=0.
	\end{align*}
	Here we apply Lemma \ref{lem:cut1} for $(i,j)=(3,4)$ and $\toba$ projects onto $\widetilde{\mathcal E}_{\bq'}$, $\bq'$ with two components, of types super $A_3$ and $A_1$.

	 For $\xymatrix@R-6pt{  &    \overset{-1}{\circ} \ar  @{-}[d]^{\ztu} & \\
		\overset{\ztu}{\underset{\ }{\circ}} \ar  @{-}[r]^{\zeta}  & \overset{-1}{\underset{\
			}{\circ}} \ar  @{-}[r]^{\zeta}  & \overset{-1}{\underset{\ }{\circ}}}$
	the algebra $\toba$ is generated by $y_i$, $i\in\I_4$, with defining relations
	\begin{align*}
	y_{13}&=0; & y_{3}^2&=\mu_3; &y_{14}&=0; & y_{34}&=\lambda_{34}; & [y_{124},&y_{2}]_c=0;\\
	&  & y_{4}^2&=\mu_4; & y_{112}&=0; & y_{2}^2&=\mu_2; & [y_{324},&y_{2}]_c=\nu_2.
	\end{align*}
	Here we apply Lemma \ref{lem:cut1} for $(i,j)=(1,2)$ and $\toba$ projects onto $\widetilde{\mathcal E}_{\bq'}$, $\bq'$ with two components, of types super $A_3$ and $A_1$.
	
	 For $\xymatrix@R-6pt{  &    \overset{-1}{\circ} \ar  @{-}[d]^{\zeta} & \\
		\overset{\ztu}{\underset{\ }{\circ}} \ar  @{-}[r]^{\zeta}  & \overset{\ztu}{\underset{\
			}{\circ}} \ar  @{-}[r]^{\zeta}  & \overset{-1}{\underset{\ }{\circ}}}$
	the algebra $\toba$ is generated by $y_i$, $i\in\I_4$, with defining relations
	\begin{align*}
	y_{13}&=0; & y_{3}^2&=\mu_3; &y_{14}&=0; & y_{34}&=\lambda_{34}; & y_{112}&=\lambda_{112};\\
	&& y_{4}^2&=\mu_4; & y_{221}&=\lambda_{221}; & y_{223}&=0; & y_{224}&=0.
	\end{align*}
	Here we apply Lemma \ref{lem:cut1} for $(i,j)=(3,4)$ and $\toba$ projects onto $\widetilde{\mathcal E}_{\bq'}$, $\bq'$ with two components, of types super $A_3$ and $A_1$.
	
	 For $\xymatrix@R-6pt{  &    \overset{-1}{\circ} \ar  @{-}[d]^{\ztu} & \\
		\overset{\ztu}{\underset{\ }{\circ}} \ar  @{-}[r]^{\zeta}  & \overset{\zeta}{\underset{\
			}{\circ}} \ar  @{-}[r]^{\ztu}  & \overset{-1}{\underset{\ }{\circ}}}$
	the algebra $\toba$ is generated by $y_i$, $i\in\I_4$, with defining relations
	\begin{align*}
	y_{13}&=0; &  y_{3}^2&=\mu_3; & y_{14}&=0; & y_{34}&=\lambda_{34}; & [[y_{(13)},&y_2]_c,y_2]_c=0;\\
	y_{112}&=0; &  y_{4}^2&=\mu_4; & y_{223}&=0; & y_{224}&=0; & [[y_{124},&y_2]_c,y_2]_c=0.
	\end{align*}
	Here we apply Lemma \ref{lem:cut1} for $(i,j)=(1,2)$ and $\toba$ projects onto $\widetilde{\mathcal E}_{\bq'}$, $\bq'$ with two components, of types super $A_3$ and $A_1$.
	
	 For $\xymatrix@R-6pt{  &    \overset{-1}{\circ} \ar  @{-}[d]_{\ztu}\ar  @{-}[dr]^{\ztu} & \\
		\overset{\ztu}{\underset{\ }{\circ}} \ar  @{-}[r]^{\zeta}  & \overset{-1}{\underset{\
			}{\circ}} \ar  @{-}[r]^{\ztu}  & \overset{\zeta}{\underset{\ }{\circ}}}$
	the algebra $\toba$ is generated by $y_i$, $i\in\I_4$, with defining relations
	\begin{align*}
	y_{13}&=\lambda_{13}; & y_{2}^2&=\mu_2; & y_{14}&=0; & y_{112}&=0; & [y_{(13)},&y_{2}]_c=\nu_2;\\
	&  & y_{4}^2&=\mu_4; & y_{332}&=0; & y_{334}&=0; & [y_{124},&y_{2}]_c=0;
	\end{align*}
	\vsp
	\begin{align*}
	y_{(24)}-\zeta q_{34}[y_{24},y_3]_c-q_{23}(1-\ztu)y_3y_{24}=0.
	\end{align*}
	Here we apply Lemma \ref{lem:cut2} for $(i,j,k)=(4,2,3)$ and $\toba$ projects onto $\widetilde{\mathcal E}_{\bq'}$, $\bq'$ with two components, of types super $A_3$ and $A_1$.

\subsubsection{\ Type $\g(4,3)$, $\zeta \in \G'_3$}\label{subsec:type-g(4,3)}
The Weyl groupoid has ten objects.

	 For $\xymatrix{ \overset{\ztu}{\underset{\ }{\circ}}\ar  @{-}[r]^{\zeta}  &
		\overset{-1}{\underset{\ }{\circ}} \ar  @{-}[r]^{\ztu}  & \overset{\zeta}{\underset{\
			}{\circ}}
		\ar  @{-}[r]^{\zeta}  & \overset{-1}{\underset{\ }{\circ}}}$,
	the algebra $\toba$ is generated by $y_i$, $i\in\I_4$, with defining relations
	\begin{align*}
	y_{13}&=\lambda_{13};  & y_{2}^2&=\mu_2; &y_{14}&=0; & y_{24}&=\lambda_{14}; & [[y_{(24)},&y_3]_c,y_3]_c=\lambda_{(24)};\\
	[y_{(13)}&,y_2]_c=\nu_2;  & y_{4}^2&=\mu_4; & y_{112}&=0; & y_{332}&=0; & [y_{334},&y_{34}]_c=\lambda_{34}.
	\end{align*}
	If $\nu_2=0$, then we apply Lemma \ref{lem:cut1} for $(i,j)=(1,2)$ and $\toba$ projects onto $\widetilde{\mathcal E}_{\bq'}$, $\bq'$ with two components, of types super $\g(2,3)$ and $A_1$. If $\nu_2\neq 0$, then $\lambda_{(24)}=\lambda_{34}=0$, so we apply Lemma \ref{lem:cut1} for $(i,j)=(3,4)$ and $\toba$ projects onto $\widetilde{\mathcal E}_{\bq'}$, $\bq'$ with two components, of types super $A_3$ and $A_1$.

	 For $\xymatrix{ \overset{\ztu}{\underset{\ }{\circ}}\ar  @{-}[r]^{\zeta}  &
		\overset{-1}{\underset{\ }{\circ}} \ar  @{-}[r]^{\ztu}  & \overset{-\ztu}{\underset{\
			}{\circ}}
		\ar  @{-}[r]^{\ztu}  & \overset{-1}{\underset{\ }{\circ}}}$,
	the algebra $\toba$ is generated by $y_i$, $i\in\I_4$, with defining relations
	\begin{align*}
	y_{13}&=0; & y_{4}^2&=\mu_4; &y_{14}&=0; & y_{24}&=\lambda_{24}; & [y_{(13)},y_2]_c&=0;\\
	&  & y_{2}^2&=\mu_2;  & y_{112}&=0; & y_{3332}&=\lambda_{3332}; & y_{3334}&=\lambda_{3334};
	\end{align*}
	\vsp
	\begin{align*}
	[y_2,y_{334}]_c-q_{34}[y_{(24)},y_3]_c+(\zeta^2-\zeta)q_{23}y_3y_{(24)}=0.
	\end{align*}
	Here we apply Lemma \ref{lem:cut1} for $(i,j)=(1,2)$ and $\toba$ projects onto $\widetilde{\mathcal E}_{\bq'}$, $\bq'$ with two components, of types super $\g(2,3)$ and $A_1$.

	 For $\xymatrix{ \overset{-1}{\underset{\ }{\circ}}\ar  @{-}[r]^{\ztu}  &
		\overset{-1}{\underset{\ }{\circ}} \ar  @{-}[r]^{\zeta}  & \overset{-1}{\underset{\
			}{\circ}}
		\ar  @{-}[r]^{\zeta}  & \overset{-1}{\underset{\ }{\circ}}}$,
	the algebra $\toba$ is generated by $y_i$, $i\in\I_4$, with defining relations
\begin{align*}
y_{1}^2&=\mu_1; & y_{2}^2&=\mu_2; & y_{13}&=\lambda_{13};  & [y_{(13)}&,y_2]_c=\nu_2; \\
& & y_{3}^2&=\mu_3; &y_{14}&=\lambda_{14}; & [[y_{23},&y_{(24)}]_c,y_3]_c=\lambda_{234};\\
& & y_{4}^2&=\mu_4; & y_{24}&=\lambda_{24}; &  [[y_{43},&y_{432}]_c,y_3]_c=\lambda_{432}.
\end{align*}
If $\nu_2=0$, then we apply Lemma \ref{lem:cut1} for $(i,j)=(1,2)$ and $\toba$ projects onto $ \widetilde{\mathcal E}_{\bq'}$, $\bq'$ with two components, of types super $\g(2,3)$ and $A_1$. If $\nu_2\neq 0$, then $\lambda_{234}=\lambda_{432}=0$, so we apply Lemma \ref{lem:cut1} for $(i,j)=(3,4)$ and $\toba$ projects onto $\widetilde{\mathcal E}_{\bq'}$, $\bq'$ with two components, of types super $A_3$ and $A_1$.

	 For $\xymatrix{ \overset{-1}{\underset{\ }{\circ}}\ar  @{-}[r]^{\zeta}  &
		\overset{\ztu}{\underset{\ }{\circ}} \ar  @{-}[r]^{\zeta}  & \overset{-1}{\underset{\
			}{\circ}}
		\ar  @{-}[r]^{\zeta}  & \overset{-1}{\underset{\ }{\circ}}}$,
	the algebra $\toba$ is generated by $y_i$, $i\in\I_4$, with defining relations
	\begin{align*}
	y_{13}&=\lambda_{13}; & y_{3}^2&=\mu_3; &y_{14}&=\lambda_{14}; & y_{24}&=0; & [[y_{43},&y_{432}]_c,y_3]_c=0;\\
	&  & y_{4}^2&=\mu_4; & y_{221}&=0; & y_{223}&=0; & y_1^2&=\mu_2.
	\end{align*}
	This algebra is nonzero by \cite[Lemma 5.16]{AAGMV}.

 For $\xymatrix{ \overset{-1}{\underset{\ }{\circ}}\ar  @{-}[r]^{\ztu}  &
		\overset{-1}{\underset{\ }{\circ}} \ar  @{-}[r]^{\zeta}  & \overset{\zeta}{\underset{\
			}{\circ}}
		\ar  @{-}[r]^{\ztu}  & \overset{-1}{\underset{\ }{\circ}}}$,
the algebra $\toba$ is generated by $y_i$, $i\in\I_4$, with defining relations
\begin{align*}
[y_{(13)},y_2]_c&=0; & y_{1}^2&=\mu_1; & y_{13}&=0; &  [[y_{(24)},y_{3}]_c,y_3]_c&=\lambda_{(24)}; \\
y_{334}&=0; & y_{2}^2&=\mu_2; & y_{14}&=\lambda_{14}; & [y_{3321},y_{32}]_c&=0; \\
& & y_{4}^2&=\mu_4; & y_{24}&=\lambda_{24}; & [y_{332},y_{32}]_c&=\lambda_{32}.
\end{align*}
Here we apply Lemma \ref{lem:cut1} for $(i,j)=(1,2)$ and $\toba$ projects onto $\widetilde{\mathcal E}_{\bq'}$, $\bq'$ with two components, of types super $\g(2,3)$ and $A_1$.
	
	 For $\xymatrix{ \overset{-1}{\underset{\ }{\circ}}\ar  @{-}[r]^{\zeta}  &
		\overset{\ztu}{\underset{\ }{\circ}} \ar  @{-}[r]^{\zeta}  & \overset{\zeta}{\underset{\
			}{\circ}}
		\ar  @{-}[r]^{\ztu}  & \overset{-1}{\underset{\ }{\circ}}}$,
	the algebra $\toba$ is generated by $y_i$, $i\in\I_4$, with defining relations
	\begin{align*}
	y_{13}&=0;  & y_{1}^2&=\mu_1; &y_{14}&=\lambda_{14}; & y_{24}&=0; & [[y_{(24)},&y_{3}]_c,y_3]_c=0;\\
	y_{334}&=0;  & y_{4}^2&=\mu_4; & y_{221}&=0; & y_{223}&=0; & [y_{3321},&y_{32}]_c=0.
	\end{align*}
	This algebra is nonzero by \cite[Lemma 5.16]{AAGMV}.
	
	 For $\xymatrix@R-6pt{  &    \overset{\ztu}{\circ} \ar  @{-}[d]^{\zeta} & \\
		\overset{\zeta}{\underset{\ }{\circ}} \ar  @{-}[r]^{\ztu}  & \overset{-1}{\underset{\
			}{\circ}} \ar  @{-}[r]^{\zeta}  & \overset{-1}{\underset{\ }{\circ}}}$,
	the algebra $\toba$ is generated by $y_i$, $i\in\I_4$, with defining relations
	\begin{align*}
	y_{13}&=0;  & y_2^2&=\mu_2; &y_{14}&=\lambda_{14}; & y_{34}&=0; & [y_{124},&y_2]_c=\nu_2';\\
	[y_{(13)},y_2]_c&=0;  & y_3^2&=\mu_3; &  y_{112}&=0; & y_{442}&=0; & [[y_{32},&y_{324}]_c,y_2]_c=0.
	\end{align*}
	Here we apply Lemma \ref{lem:cut1} for $(i,j)=(1,2)$ and $\toba$ projects onto $\widetilde{\mathcal E}_{\bq'}$, $\bq'$ with two components, of types super $CD$ and $A_1$.

	 For $\xymatrix@R-6pt{  &    \overset{-1}{\circ} \ar  @{-}[d]_{\ztu}\ar  @{-}[dr]^{\ztu} &
		\\
		\overset{-1}{\underset{\ }{\circ}} \ar  @{-}[r]^{\ztu}  & \overset{\zeta}{\underset{\
			}{\circ}} \ar  @{-}[r]^{\ztu}  & \overset{\zeta}{\underset{\ }{\circ}}}$,
	the algebra $\toba$ is generated by $y_i$, $i\in\I_4$, with defining relations
	\begin{align*}
	y_{13}&=0; &y_{14}&=\lambda_{14}; & y_{221}&=0; & y_{223}&=\lambda_{223}; & y_{224}&=0;\\
	y_{332}&=\lambda_{332}; & y_{334}&=0;  & y_1^2&=\mu_1; & y_4^2&=\mu_4;
	\end{align*}
	\vsp
	\begin{align*}
	y_{(24)}-\zeta q_{34}[y_{24},y_3]_c-(1-\ztu)q_{23}y_3y_{24}=0.
	\end{align*}
	Here we apply Lemma \ref{lem:cut2} for $(i,j,k)=(4,2,3)$ and $\toba$ projects onto $ \widetilde{\mathcal E}_{\bq'}$, $\bq'$ with two components, of types super $A_3$ and $A_1$.
	
	 For $\xymatrix@R-6pt{  &    \overset{\ztu}{\circ} \ar  @{-}[d]^{\zeta} & \\
		\overset{\zeta}{\underset{\ }{\circ}} \ar  @{-}[r]^{\ztu}  & \overset{\zeta}{\underset{\
			}{\circ}} \ar  @{-}[r]^{\ztu}  & \overset{-1}{\underset{\ }{\circ}}}$,
	the algebra $\toba$ is generated by $y_i$, $i\in\I_4$, with defining relations
	\begin{align*}
	y_{13}&=0;  & y_{442}&=0; &y_{14}&=\lambda_{14}; & y_{34}&=0; & [[y_{124},&y_2]_c,y_2]_c=\lambda_{124}; \\
	y_3^2&=0; & y_{223}&=0; & y_{112}&=\lambda_{112}; & y_{221}&=\lambda_{221}; & [[y_{324},&y_2]_c,y_2]_c=0.
	\end{align*}
	Here we apply Lemma \ref{lem:cut1} for $(i,j)=(2,3)$ and $\toba$ projects onto $\widetilde{\mathcal E}_{\bq'}$, $\bq'$ with two components, of types $C_3$ and $A_1$.

	 For $\xymatrix@R-6pt{  &    \overset{-1}{\circ} \ar  @{-}[d]_{\ztu}\ar  @{-}[dr]^{\ztu} & \\
		\overset{-1}{\underset{\ }{\circ}} \ar  @{-}[r]^{\zeta}  & \overset{-1}{\underset{\
			}{\circ}} \ar  @{-}[r]^{\ztu}  & \overset{\zeta}{\underset{\ }{\circ}}}$,
	the algebra $\toba$ is generated by $y_i$, $i\in\I_4$, with defining relations
	\begin{align*}
	y_{13}&=0; & y_2^2&=\mu_2; & y_{14}&=\lambda_{14}; & y_{332}&=0; & y_{334}&=0; \\
	&  & y_4^2&=\mu_4; & [y_{(13)},&y_2]_c=0; & y_{1}^2&=\mu_1; & [y_{124},&y_2]_c=\nu_2';
	\end{align*}
	\vsp
	\begin{align*}
	y_{(24)}-\zeta q_{34}[y_{24},y_3]_c-(1-\ztu)q_{23}y_3y_{24}=0.
	\end{align*}
	Here we apply Lemma \ref{lem:cut2} for $(i,j,k)=(3,2,4)$ and $\toba$ projects onto $ \widetilde{\mathcal E}_{\bq'}$, $\bq'$ with two components, of types super $A_3$ and $A_1$.

\subsubsection{\ Type $\g(3, 6)$, $\zeta \in \G'_3$}\label{subsec:type-g(3,6)}
The Weyl groupoid has seven objects.

	 For $\xymatrix{ \overset{-1}{\underset{\ }{\circ}}\ar  @{-}[r]^{\ztu}  &
		\overset{\zeta}{\underset{\ }{\circ}} \ar  @{-}[r]^{\ztu}  & \overset{\zeta}{\underset{\
			}{\circ}}\ar  @{-}[r]^{\zeta}  & \overset{-1}{\underset{\ }{\circ}}}$,
	the algebra $\toba$ is generated by $y_i$, $i\in\I_4$, with defining relations
	\begin{align*}
	y_{13}&=0; & y_{4}^2&=\mu_4; &y_{14}&=\lambda_{14}; & y_{24}&=0; & [[y_{(24)},&y_{3}]_c,y_{3}]_c=0;\\
	y_{332}&=\lambda_{332}; &y_{1}^2&=0;  & y_{221}&=0; & y_{223}&=\lambda_{223}; & [y_{334},&y_{34}]_c=\lambda_{34}.
	\end{align*}
	Here we apply Lemma \ref{lem:cut1} for $(i,j)=(1,2)$ and $\toba$ projects onto $\widetilde{\mathcal E}_{\bq'}$, $\bq'$ with two components, of types super $\g(1,6)$ and $A_1$.

	 For $\xymatrix{ \overset{-1}{\underset{\ }{\circ}}\ar  @{-}[r]^{\ztu}  &
		\overset{\zeta}{\underset{\ }{\circ}} \ar  @{-}[r]^{\ztu}  & \overset{-\ztu}{\underset{\
			}{\circ}}
		\ar  @{-}[r]^{\ztu}  & \overset{-1}{\underset{\ }{\circ}}}$,
	the algebra $\toba$ is generated by $y_i$, $i\in\I_4$, with defining relations
	\begin{align*}
	y_{13}&=0;  & y_{1}^2&=\mu_1;  &y_{14}&=\lambda_{14}; & y_{24}&=0; & y_{221}&=0;\\
	& & y_{4}^2&=\mu_4; & y_{223}&=0; & y_{3332}&=0; & y_{3334}&=\lambda_{3334};
	\end{align*}
	\vsp
	\begin{align*}
	[y_2,y_{334}]_c+q_{34}[y_{(24)},y_3]_c+(\zeta^2-\zeta)q_{23}y_3y_{(24)}=0.
	\end{align*}
	Here we apply Lemma \ref{lem:cut1} for $(i,j)=(1,2)$ and $\toba$ projects onto $\widetilde{\mathcal E}_{\bq'}$, $\bq'$ with two components, of types super $\g(1,6)$ and $A_1$.
	
For $\xymatrix{ \overset{-1}{\underset{\ }{\circ}}\ar  @{-}[r]^{\zeta}  & \overset{-1}{\underset{\ }{\circ}} \ar  @{-}[r]^{\ztu}  & \overset{\zeta}{\underset{\ }{\circ}} \ar  @{-}[r]^{\zeta}  & \overset{-1}{\underset{\ }{\circ}}}$,
	the algebra $\toba$ is generated by $y_i$, $i\in\I_4$, with defining relations
	\begin{align*}
	y_1^2&=\mu_1; & y_{13}&=0; &y_{14}&=\lambda_{14}; & y_{24}&=\lambda_{24}; & [[y_{(24)},&y_{3}]_c,y_{3}]_c=\lambda_{(24)};\\
	y_{2}^2&=\mu_2; & y_{4}^2&=\mu_4; & y_{332}&=0; & [y_{(13)},& y_{2}]_c=0; & [y_{334},&y_{34}]_c=\lambda_{34}.
	\end{align*}
	Here we apply Lemma \ref{lem:cut1} for $(i,j)=(1,2)$ and $\toba$ projects onto $\widetilde{\mathcal E}_{\bq'}$, $\bq'$ with two components, of types $\g(2,3)$ and $A_1$.

	 For $\xymatrix{ \overset{-1}{\underset{\ }{\circ}}\ar  @{-}[r]^{\zeta}  &
		\overset{-1}{\underset{\ }{\circ}} \ar  @{-}[r]^{\ztu}  & \overset{-\ztu}{\underset{\
			}{\circ}}
		\ar  @{-}[r]^{\ztu}  & \overset{-1}{\underset{\ }{\circ}}}$,
	the algebra $\toba$ is generated by $y_i$, $i\in\I_4$, with defining relations
	\begin{align*}
	y_{13}&=0;  & y_{2}^2&=\mu_2; &y_{14}&=\lambda_{14}; & y_{24}&=\lambda_{24}; & [y_{(13)},&y_2]_c=0;\\
	&  & y_{4}^2&=\mu_4; & y_{3332}&=\lambda_{3332}; & y_{3334}&=\lambda_{3334}; & y_{1}^2&=\mu_1;
	\end{align*}
	\vsp
	\begin{align*}
	[y_2,y_{334}]_c+q_{34}[y_{(24)},y_3]_c+(\zeta^2-\zeta)q_{23}y_3y_{(24)}=0.
	\end{align*}
	Here we apply Lemma \ref{lem:cut1} for $(i,j)=(1,2)$ and $\toba$ projects onto $\widetilde{\mathcal E}_{\bq'}$, $\bq'$ with two components, of types super $\g(2,3)$ and $A_1$.

 For $\xymatrix{ \overset{\zeta}{\underset{\ }{\circ}}\ar  @{-}[r]^{\ztu}  &
		\overset{-1}{\underset{\ }{\circ}} \ar  @{-}[r]^{\zeta}  & \overset{\zeta}{\underset{\
			}{\circ}}
		\ar  @{-}[r]^{\ztu}  & \overset{-1}{\underset{\ }{\circ}}}$,
the algebra $\toba$ is generated by $y_i$, $i\in\I_4$, with defining relations
\begin{align*}
y_{2}^2&=\mu_2; & y_{13}&=0; & y_{112}&=0; & [y_{332},&y_{32}]_c=\lambda_{32}; \\
y_{4}^2&=\mu_4; & y_{14}&=0; & y_{334}&=0; & [[y_{(24)},&y_3]_c,y_3]_c=\lambda_{(24)};\\
& & y_{24}&=\lambda_{24}; & [y_{(13)},&y_2]_c=0; & [y_{3321},&y_{32}]_c=0.
\end{align*}
Here we apply Lemma \ref{lem:cut1} for $(i,j)=(1,2)$ and $\toba$ projects onto $\widetilde{\mathcal E}_{\bq'}$, $\bq'$ with two components, of types super $\g(2,3)$ and $A_1$.
	
	 For $\xymatrix{ \overset{\zeta}{\underset{\ }{\circ}}\ar  @{-}[r]^{\ztu}  &
		\overset{-1}{\underset{\ }{\circ}} \ar  @{-}[r]^{\zeta}  & \overset{-1}{\underset{\
			}{\circ}}
		\ar  @{-}[r]^{\zeta}  & \overset{-1}{\underset{\ }{\circ}}}$,
	the algebra $\toba$ is generated by $y_i$, $i\in\I_4$, with defining relations
	\begin{align*}
	y_{13}&=0;  & y_{3}^2&=\mu_3; & [y_{(13)},&y_{2}]_c=0; & y_{24}&=\lambda_{24}; & [[y_{23},&y_{(24)}]_c,y_{3}]_c=\lambda_{234};\\
	y_{14}&=0;  & y_{4}^2&=\mu_4; & y_{112}&=0; & y_{2}^2&=\mu_2; & [[y_{43},&y_{432}]_c,y_{3}]_c=\lambda_{432}.
	\end{align*}
	Here we apply Lemma \ref{lem:cut1} for $(i,j)=(1,2)$ and $\toba$ projects onto $\widetilde{\mathcal E}_{\bq'}$, $\bq'$ with two components, of types super $\g(2,3)$ and $A_1$.

For $\xymatrix@R-6pt{  &    \overset{-1}{\circ} \ar  @{-}[d]_{\ztu}\ar  @{-}[dr]^{\ztu} & \\  \overset{\zeta}{\underset{\ }{\circ}} \ar  @{-}[r]^{\ztu}  & \overset{\zeta}{\underset{\ }{\circ}} \ar  @{-}[r]^{\ztu}  & \overset{\zeta}{\underset{\ }{\circ}}}$,
	the algebra $\toba$ is generated by $y_i$, $i\in\I_4$, with defining relations
	\begin{align*}
	y_{13}&=0; &y_{14}&=0; & y_{221}&=\lambda_{221}; & y_{223}&=\lambda_{223}; & y_{224}&=0;\\
	y_{332}&=\lambda_{332}; & y_{334}&=0;  & y_{112}&=\lambda_{112}; & y_4^2&=\mu_4;
	\end{align*}
	\vsp
	\begin{align*}
	y_{(24)}-\zeta q_{34}[y_{24},y_3]_c-(1-\ztu)q_{23}y_3y_{24}=0.
	\end{align*}
	Here we apply Lemma \ref{lem:cut2} for $(i,j,k)=(4,2,3)$ and $\toba$ projects onto $\widetilde{\mathcal E}_{\bq'}$, $\bq'$ with two components of types $A_3$ and $A_1$.

\subsubsection{\ Type $\g(2, 6)$, $\zeta \in \G'_3$}\label{subsec:type-g(2,6)}
The Weyl groupoid has four objects.

 For $\xymatrix{\overset{\zeta}{\underset{\ }{\circ}}\ar  @{-}[r]^{\ztu}  &
	\overset{\zeta}{\underset{\ }{\circ}} \ar  @{-}[r]^{\ztu}  & \overset{-1}{\underset{\
		}{\circ}}
	\ar  @{-}[r]^{\ztu}  & \overset{\zeta}{\underset{\ }{\circ}} \ar  @{-}[r]^{\ztu}  &
	\overset{\zeta}{\underset{\ }{\circ}}}$
the algebra $\toba$ is generated by $y_i$, $i\in\I_5$, with defining relations
\begin{align*}
& & y_{112}&=\lambda_{112}; & y_{223}&=0; & [[[y_{(14)},y_3]_c , y_2]_c , y_3]_c &=0; \\
& & y_{3}^2&=\mu_3; & y_{554}&=\lambda_{554}; &
[[[y_{5432},y_3]_c, y_4]_c , y_3]_c &=0; \\
y_{221}&=\lambda_{221}; & y_{443}&=0; &
y_{445}& =\lambda_{225}; & y_{ij}=\lambda_{ij}, \, i<j, \, & \widetilde{q}_{ij}=1.
\end{align*}
Here we apply Lemma \ref{lem:cut1} for $(i,j)=(2,3)$ and $\toba$ projects onto $\widetilde{\mathcal E}_{\bq'}$, $\bq'$ with two components of types $A_3$ and $A_2$.

 For $\xymatrix@R-8pt{  &    \overset{-1}{\circ} \ar  @{-}[d]_{\zeta}\ar  @{-}[dr]^{\zeta} &
	\\
	\overset{\zeta}{\underset{\ }{\circ}} \ar  @{-}[r]^{\ztu}  & \overset{-1}{\underset{\
		}{\circ}} \ar  @{-}[r]^{\zeta}  & \overset{-1}{\underset{\ }{\circ}}\ar  @{-}[r]^{\ztu}  &
	\overset{\zeta}{\underset{\ }{\circ}}}$
the algebra $\toba$ is generated by $y_i$, $i\in\I_5$, with defining relations
\begin{align*}
y_{112}&=0; & y_{3}^2&=\mu_3; & y_{ij}&=\lambda_{ij}, \, i<j, \, \widetilde{q}_{ij}=1;  & [y_{125},& y_2]_c=0;\\
y_{443}&=0; & y_{5}^2&=\mu_5; & [y_{(13)},& y_2]_c=0; & [y_{(24)},& y_3]_c=0;
\end{align*}
\vsp
\begin{align*}
y_{2}^2&=\mu_2; & [y_{435},& y_3]_c=0; & y_{235}-q_{35}\ztu [y_{25},y_3]_c-q_{23}(1-\zeta)y_3y_{25}=0.
\end{align*}
This algebra is nonzero by \cite[Lemma 5.16]{AAGMV}.

 For $\xymatrix@R-10pt{ & &    \overset{\zeta}{\circ} \ar  @{-}[d]^{\ztu} & \\
	\overset{\zeta}{\underset{\ }{\circ}} \ar  @{-}[r]^{\ztu} & \overset{\zeta}{\underset{\
		}{\circ}} \ar  @{-}[r]^{\ztu}  & \overset{-1}{\underset{\ }{\circ}} \ar  @{-}[r]^{\zeta}
	& \overset{-1}{\underset{\ }{\circ}}}$
the algebra $\toba$ is generated by $y_i$, $i\in\I_5$, with defining relations
\begin{align*}
y_{223}&=0; & y_{112}&=\lambda_{112}; & y_{221}&=\lambda_{221}; & y_{3}^2&=\mu_3; \quad
[y_{435}, y_3]_c=0;\\
 [y_{(24)},& y_3]_c=0; & y_{4}^2&=\mu_4; & y_{553}&=0; & y_{ij}&=\lambda_{ij}, \, i<j, \, \widetilde{q}_{ij}=1.
\end{align*}
Here we apply Lemma \ref{lem:cut1} for $(i,j)=(2,3)$ and $\toba$ projects onto $\widetilde{\mathcal E}_{\bq'}$, $\bq'$ with two components of types $A_3$ and $A_2$.

 For $\xymatrix@R-10pt{  &  &  \overset{\zeta}{\circ} \ar  @{-}[d]^{\ztu} & \\
	\overset{\zeta}{\underset{\ }{\circ}} \ar  @{-}[r]^{\ztu} & \overset{\zeta}{\underset{\
		}{\circ}} \ar  @{-}[r]^{\ztu}  & \overset{\zeta}{\underset{\ }{\circ}} \ar  @{-}[r]^{\ztu}
	& \overset{-1}{\underset{\ }{\circ}}}$
the algebra $\toba$ is generated by $y_i$, $i\in\I_5$, with defining relations
\begin{align*}
y_{4}^2&=\mu_4; \qquad y_{553}=\lambda_{553}; & y_{112}&=\lambda_{112}; & y_{332}&=\lambda_{332}; & y_{221}&=\lambda_{221};\\
y_{ij}&=\lambda_{ij}, \, i<j, \, \widetilde{q}_{ij}=1;  & y_{335}&=\lambda_{335}; & y_{223}&= \lambda_{223}; & y_{334}&=0.
\end{align*}
Here we apply Lemma \ref{lem:cut1} for $(i,j)=(3,4)$ and $\toba$ projects onto $\widetilde{\mathcal E}_{\bq'}$, $\bq'$ with two components of types $A_4$ and $A_1$.

\subsubsection{\ Type $\el(5;3)$, $\zeta \in \G'_3$}\label{subsec:type-el(5;3)}
The Weyl groupoid has fifteen objects.

 For $\xymatrix{\overset{\zeta}{\underset{\ }{\circ}}\ar  @{-}[r]^{\ztu}  &
	\overset{\zeta}{\underset{\ }{\circ}} \ar  @{-}[r]^{\ztu}  & \overset{\zeta}{\underset{\
		}{\circ}}
	\ar  @{-}[r]^{\ztu}  & \overset{\ztu}{\underset{\ }{\circ}} \ar  @{-}[r]^{\zeta}  &
	\overset{-1}{\underset{\ }{\circ}}}$
the algebra $\toba$ is generated by $y_i$, $i\in\I_5$, with defining relations
\begin{align*}
y_{221}&=\lambda_{221}; & y_{5}^2&=\mu_5; & y_{334}&=0; & [[y_{(35)},&y_4]_c,y_4]_c=0;\\
& & y_{445}&=0; & y_{112}&=\lambda_{112}; & y_{ij}&=\lambda_{ij}, \, i<j, \, \widetilde{q}_{ij}=1; \\
& & y_{223}&=\lambda_{223}; & y_{332}&=\lambda_{332}; & [y_{4432},&y_{43}]_c=\lambda_{432}.
\end{align*}
Here we apply Lemma \ref{lem:cut1} for $(i,j)=(4,5)$ and $\toba$ projects onto $\widetilde{\mathcal E}_{\bq'}$, $\bq'$ with two components of types $B_4$ and $A_1$.

 For $\xymatrix{\overset{\zeta}{\underset{\ }{\circ}}\ar  @{-}[r]^{\ztu}  &
	\overset{\zeta}{\underset{\ }{\circ}} \ar  @{-}[r]^{\ztu}  & \overset{\zeta}{\underset{\
		}{\circ}}
	\ar  @{-}[r]^{\ztu}  & \overset{-1}{\underset{\ }{\circ}} \ar  @{-}[r]^{\ztu}  &
	\overset{-1}{\underset{\ }{\circ}}}$
the algebra $\toba$ is generated by $y_i$, $i\in\I_5$, with defining relations
\begin{align*}
& & y_{4}^2&=\mu_4; & y_{112}&=\lambda_{112}; & y_{221}&=\lambda_{221}; &  &[[y_{54},y_{543}]_c,y_4]_c=0;\\
y_{334}&=0; & y_{5}^2&=\mu_5; & y_{223}&=\lambda_{223}; & y_{332}&=\lambda_{332}; & &y_{ij}=\lambda_{ij}, \, i<j, \, \widetilde{q}_{ij}=1.
\end{align*}
Here we apply Lemma \ref{lem:cut1} for $(i,j)=(4,5)$ and $\toba$ projects onto $\widetilde{\mathcal E}_{\bq'}$, $\bq'$ with two components of types $A_4$ and $A_1$.

 For $\xymatrix{\overset{-1}{\underset{\ }{\circ}} \ar  @{-}[r]^{\zeta} &
	\overset{-1}{\underset{\ }{\circ}} \ar  @{-}[r]^{\ztu}  & \overset{-1}{\underset{\
		}{\circ}} \ar  @{-}[r]^{\ztu}  & \overset{\zeta}{\underset{\ }{\circ}} \ar  @{-}[r]^{\ztu}
	& \overset{\zeta}{\underset{\ }{\circ}}}$
the algebra $\toba$ is generated by $y_i$, $i\in\I_5$, with defining relations
\begin{align*}
y_{443}&=0;      &  y_{445}&=\lambda_{445}; & [[y_{23},&y_{(24)}]_c,y_3]_c=0; & y_{1}^2&=\mu_1;\quad y_{2}^2=\mu_2\\
y_{3}^2&=\mu_3;  & y_{554}&=\lambda_{554};  & [y_{(13)},&y_2]_c=\nu_2;         & y_{ij}&=\lambda_{ij}, \, i<j, \, \widetilde{q}_{ij}=1.
\end{align*}
Here we apply Lemma \ref{lem:cut1} for $(i,j)=(3,4)$ and $\toba$ projects onto $\widetilde{\mathcal E}_{\bq'}$, $\bq'$ with two components of types $A_3$ and $A_2$.

 For $\xymatrix{\overset{\zeta}{\underset{\ }{\circ}} \ar  @{-}[r]^{\ztu} &
	\overset{-1}{\underset{\ }{\circ}} \ar  @{-}[r]^{\zeta}  & \overset{\ztu}{\underset{\
		}{\circ}} \ar  @{-}[r]^{\ztu}  & \overset{\zeta}{\underset{\ }{\circ}} \ar  @{-}[r]^{\ztu}
	& \overset{\zeta}{\underset{\ }{\circ}}}$
the algebra $\toba$ is generated by $y_i$, $i\in\I_5$, with defining relations $y_{ij}=\lambda_{ij}, \, i<j, \, \widetilde{q}_{ij}=1;$
\begin{align*}
y_{2}^2&=\mu_2; & y_{112}&=0; & y_{445}&=\lambda_{445}; & [[y_{(24)},&y_{3}]_c,y_3]_c=0; \\
y_{332}&=0; & y_{443}&=0; & y_{554}&=\lambda_{554}; & [y_{3345},&y_{34}]_c=\lambda_{345};& [y_{(13)},&y_2]_c=\nu_2.
\end{align*}
If $\nu_2=0$, then we apply Lemma \ref{lem:cut1} for $(i,j)=(2,3)$ and $\toba$ projects onto $ \widetilde{\mathcal E}_{\bq'}$, $\bq'$ with two components of types $B_3$ and $A_2$.
If $\nu_2\neq 0$, then $\lambda_{345}=0$, so we apply Lemma \ref{lem:cut1} for $(i,j)=(3,4)$ and $\toba$ projects onto $\widetilde{\mathcal E}_{\bq'}$, $\bq'$ with two components of types $A_3$ and $A_2$.

 For $\xymatrix@C-4pt{\overset{-1}{\underset{\ }{\circ}} \ar  @{-}[r]^{\ztu} &
	\overset{\zeta}{\underset{\ }{\circ}} \ar  @{-}[r]^{\ztu}  & \overset{-1}{\underset{\
		}{\circ}} \ar  @{-}[r]^{\ztu}  & \overset{\zeta}{\underset{\ }{\circ}} \ar  @{-}[r]^{\ztu}
	& \overset{\zeta}{\underset{\ }{\circ}}}$
the algebra $\toba$ is generated by $y_i$, $i\in\I_5$, with defining relations
$y_{ij}=\lambda_{ij}, \, i<j, \, \widetilde{q}_{ij}=1;$
\begin{align*}
& & y_{3}^2&=\mu_3; & y_{445}&=\lambda_{445}; & y_{221}&=0; & [[[y_{5432},&y_{3}]_c,y_4]_c,y_3]_c=0;
\\
y_{223}&=0; & y_{1}^2&=\mu_1; & y_{554}&=\lambda_{554}; & y_{443}&=0; & [[[y_{(14)},&y_{3}]_c,y_2]_c,y_3]_c=\lambda_{(14)}.
\end{align*}
If $\lambda_{(14)}=0$, then we apply Lemma \ref{lem:cut1} for $(i,j)=(3,4)$ and $\toba$ projects onto $ \widetilde{\mathcal E}_{\bq'}$, $\bq'$ with two components of types $A_3$, $A_2$. If $\lambda_{(14)}\neq 0$ but all the other $\lambda$'s are zero, then we apply Lemma \ref{lem:cut1} for $(i,j)=(4,5)$ and $\toba$ projects onto $\widetilde{\mathcal E}_{\bq'}$, $\bq'$ with two components of types $B_4$, $A_1$. The general case follows by Lemma \ref{lem:nonzero-k<l}.

 For $\xymatrix@R-8pt{  &  &  \overset{-1}{\circ} \ar  @{-}[dl]_{\zeta}\ar  @{-}[d]^{\zeta} &
	\\
	\overset{-1}{\underset{\ }{\circ}} \ar  @{-}[r]^{\ztu}  & \overset{-1}{\underset{\
		}{\circ}} \ar  @{-}[r]^{\zeta}  & \overset{-1}{\underset{\ }{\circ}} \ar  @{-}[r]^{\ztu}
	& \overset{\zeta}{\underset{\ }{\circ}}}$
the algebra $\toba$ is generated by $y_i$, $i\in\I_5$, with defining relations
$y_{ij}=\lambda_{ij}, \, i<j, \, \widetilde{q}_{ij}=1;$
\begin{align*}
y_{2}^2&=\mu_2; & [y_{125},&y_{2}]_c=\nu_2'; & y_{3}^2&=\mu_3; & [y_{(24)},&y_{3}]_c=0;\\
[y_{(13)},&y_{2}]_c=\nu_2; & y_{5}^2&=\mu_5; &y_{1}^2&=\mu_1; & [y_{435},&y_{3}]_c=0;
\end{align*}
\vsp
\begin{align*}
 y_{443}&=0; & y_{235}-q_{35}\ztu[y_{25},y_3]_c-q_{23}(1-\zeta)y_3y_{25}=0.
\end{align*}
Here we apply Lemma \ref{lem:cut1} for $(i,j)=(3,4)$ and $\toba$ projects onto $\widetilde{\mathcal E}_{\bq'}$, $\bq'$ with two components of types super $CD$ and $A_1$.

 For $\xymatrix@R-8pt{  &  & & \overset{-1}{\circ} \ar  @{-}[d]^{\zeta}\ar  @{-}[dl]_{\zeta}
	\\
	\overset{\zeta}{\underset{\ }{\circ}} \ar  @{-}[r]^{\ztu}  & \overset{\zeta}{\underset{\
		}{\circ}} \ar  @{-}[r]^{\ztu}  & \overset{-1}{\underset{\ }{\circ}} \ar  @{-}[r]^{\zeta}
	& \overset{\ztu}{\underset{\ }{\circ}}}$
the algebra $\toba$ is generated by $y_i$, $i\in\I_5$, with defining relations
$y_{ij}=\lambda_{ij}, \, i<j, \, \widetilde{q}_{ij}=1;$
\begin{align*}
y_{3}^2&=\mu_3; & [y_{(24)},&y_{3}]_c=\nu_3; & y_{445}&=0; & y_{112}&=\lambda_{112};\\
y_{221}&=\lambda_{221}; & y_{223}&=0; & y_{443}&=0; & [y_{235},&y_{3}]_c=0;
\end{align*}
\vsp
\begin{align*}
y_{5}^2&=\mu_5; & y_{(35)}-q_{45}\ztu[y_{35},y_4]_c-q_{34}(1-\zeta)y_4y_{35}=0.
\end{align*}
Here we apply Lemma \ref{lem:cut2} for $(i,j,k)=(5,3,4)$ and $\toba$ projects onto $ \widetilde{\mathcal E}_{\bq'}$, $\bq'$ with two components of types super $A_4$ and $A_1$.

\label{page:g26-example}
 For $\xymatrix@R-8pt{  &  &  \overset{-1}{\circ} \ar  @{-}[dl]_{\zeta}\ar  @{-}[d]^{\zeta} &
	\\
	\overset{-1}{\underset{\ }{\circ}} \ar  @{-}[r]^{\zeta}  & \overset{\ztu}{\underset{\
		}{\circ}} \ar  @{-}[r]^{\zeta}  & \overset{-1}{\underset{\ }{\circ}} \ar  @{-}[r]^{\ztu}
	& \overset{\zeta}{\underset{\ }{\circ}}}$
the algebra $\toba$ is generated by $y_i$, $i\in\I_5$, with defining relations
$y_{ij}=\lambda_{ij}, \, i<j, \, \widetilde{q}_{ij}=1;$
\begin{align*}
y_{1}^2&=\mu_1; & [y_{(24)},&y_{3}]_c=\nu_2; & y_{3}^2&=\mu_3; & y_{221}&=0;\\
y_{223}&=0; & y_{5}^2&=\mu_5; & y_{443}&=0; & [y_{435},&y_{3}]_c=0;
\end{align*}
\vsp
\begin{align*}
 y_{225}&=0; & y_{235}-q_{35}\ztu[y_{25},y_3]_c-q_{23}(1-\zeta)y_3y_{25}=0.
\end{align*}
Here we apply Lemma \ref{lem:cut2} for $(i,j,k)=(5,2,3)$ and $\toba$ projects onto $ \widetilde{\mathcal E}_{\bq'}$, $\bq'$ with two components of types super $A_4$ and $A_1$.

 For $\xymatrix@R-8pt{  &  &  \overset{\zeta}{\circ} \ar  @{-}[d]^{\ztu} & \\
	\overset{\zeta}{\underset{\ }{\circ}} \ar  @{-}[r]^{\ztu}  & \overset{-1}{\underset{\
		}{\circ}} \ar  @{-}[r]^{\zeta}  & \overset{\ztu}{\underset{\ }{\circ}} \ar
	@{-}[r]^{\zeta}  & \overset{-1}{\underset{\ }{\circ}}}$
the algebra $\toba$ is generated by $y_i$, $i\in\I_5$, with defining relations
$y_{ij}=\lambda_{ij}, \, i<j, \, \widetilde{q}_{ij}=1;$
\begin{align*}
y_{112}&=0; & y_{2}^2&=\mu_2; & [y_{(13)},&y_2]_c=\nu_2; &
y_{4}^2&=\mu_4; & [[y_{235},&y_{3}]_c,y_3]_c=0;\\
& & y_{553}&=0; & y_{332}&=0; & y_{334}&=0; & [[y_{435},&y_{3}]_c,y_3]_c=0.
\end{align*}
Here we apply Lemma \ref{lem:cut1} for $(i,j)=(3,5)$ and $\toba$ projects onto $\widetilde{\mathcal E}_{\bq'}$, $\bq'$ with two components of types super $A_4$ and $A_1$.

 For $\xymatrix@R-8pt{  &  &  \overset{\zeta}{\circ} \ar  @{-}[d]^{\ztu} & \\
	\overset{-1}{\underset{\ }{\circ}} \ar  @{-}[r]^{\zeta}  & \overset{-1}{\underset{\
		}{\circ}} \ar  @{-}[r]^{\ztu}  & \overset{-1}{\underset{\ }{\circ}} \ar  @{-}[r]^{\zeta}
	& \overset{-1}{\underset{\ }{\circ}}}$
the algebra $\toba$ is generated by $y_i$, $i\in\I_5$, with defining relations $y_{ij}=\lambda_{ij}, \, i<j, \, \widetilde{q}_{ij}=1;$
\begin{align*}
& & y_{2}^2&=\mu_2; & [[y_{23},&y_{235}]_c,y_3]_c=0; &
y_{3}^2&=\mu_3; & [y_{(13)},&y_2]_c=\nu_2;\\
y_{4}^2&=\mu_4; & y_{553}&=0; & [y_{(24)},&y_3]_c=\nu_3; & y_{1}^2&=\mu_1; & [y_{435},&y_{3}]_c=0.
\end{align*}
Here we apply Lemma \ref{lem:cut1} for $(i,j)=(3,5)$ and $\toba$ projects onto $\widetilde{\mathcal E}_{\bq'}$, $\bq'$ with two components of types super $A_4$ and $A_1$.

 For $\xymatrix@R-8pt{  &  &  \overset{\zeta}{\circ} \ar  @{-}[d]^{\ztu} & \\
	\overset{-1}{\underset{\ }{\circ}} \ar  @{-}[r]^{\ztu}  & \overset{\zeta}{\underset{\
		}{\circ}} \ar  @{-}[r]^{\ztu}  & \overset{-1}{\underset{\ }{\circ}} \ar  @{-}[r]^{\zeta}
	& \overset{-1}{\underset{\ }{\circ}}}$
the algebra $\toba$ is generated by $y_i$, $i\in\I_5$, with defining relations $y_{ij}=\lambda_{ij}, \, i<j, \, \widetilde{q}_{ij}=1;$
\begin{align*}
y_{1}^2&=\mu_1; & y_{3}^2&=\mu_3; & y_{223}&=0; & [[[y_{1235},&y_{3}]_c,y_2]_c,y_3]_c=\lambda_{1235};\\
y_{221}&=0; & [y_{(24)},&y_3]_c=0; &
y_{4}^2&=\mu_4; & y_{553}&=0; \quad [y_{435},y_{3}]_c=0.
\end{align*}
Here we apply Lemma \ref{lem:cut1} for $(i,j)=(3,4)$ and $\toba$ projects onto $\widetilde{\mathcal E}_{\bq'}$, $\bq'$ with two components of types super $CD$ and $A_1$.

 For $\xymatrix@R-8pt{  &  &  \overset{\zeta}{\circ} \ar  @{-}[d]^{\ztu} & \\
	\overset{-1}{\underset{\ }{\circ}} \ar  @{-}[r]^{\zeta}  & \overset{-1}{\underset{\
		}{\circ}} \ar  @{-}[r]^{\ztu}  & \overset{\zeta}{\underset{\ }{\circ}} \ar  @{-}[r]^{\ztu}
	& \overset{-1}{\underset{\ }{\circ}}}$
the algebra $\toba$ is generated by $y_i$, $i\in\I_5$, with defining relations
$y_{ij}=\lambda_{ij}, \, i<j, \, \widetilde{q}_{ij}=1;$
\begin{align*}
y_{1}^2&=\mu_1; & [y_{(13)},&y_2]_c=0; & y_{2}^2&=\mu_2; & y_{335}&=\lambda_{335}; \\
y_{332}&=0; & y_{4}^2&=\mu_4; & y_{334}&=0; &  y_{553}&=\lambda_{553}.
\end{align*}
Here we apply Lemma \ref{lem:cut1} for $(i,j)=(3,4)$ and $\toba$ projects onto $\widetilde{\mathcal E}_{\bq'}$, $\bq'$ with two components of types super $A_4$ and $A_1$.

 For $\xymatrix@R-8pt{  &  &  \overset{\zeta}{\circ} \ar  @{-}[d]^{\ztu} & \\
	\overset{\zeta}{\underset{\ }{\circ}} \ar  @{-}[r]^{\ztu}  & \overset{-1}{\underset{\
		}{\circ}} \ar  @{-}[r]^{\zeta}  & \overset{-1}{\underset{\ }{\circ}} \ar  @{-}[r]^{\ztu}
	& \overset{-1}{\underset{\ }{\circ}}}$
the algebra $\toba$ is generated by $y_i$, $i\in\I_5$, with defining relations
$y_{ij}=\lambda_{ij}, \, i<j, \, \widetilde{q}_{ij}=1;$
\begin{align*}
& & y_{2}^2&=\mu_2; & [[y_{43},&y_{435}]_c,y_3]_c=0; &
y_{3}^2&=\mu_3; & [y_{(24)},&y_3]_c=\nu_3;\\
y_{112}&=0; & y_{4}^2&=\mu_4; & [y_{(13)},&y_2]_c=0; & y_{553}&=0; & [y_{235},&y_3]_c=0.
\end{align*}
Here we apply Lemma \ref{lem:cut1} for $(i,j)=(3,5)$ and $\toba$ projects onto $\widetilde{\mathcal E}_{\bq'}$, $\bq'$ with two components of types super $A_4$ and $A_1$.

 For $\xymatrix@R-8pt{  &  &  \overset{\zeta}{\circ} \ar  @{-}[d]^{\ztu} & \\
	\overset{-1}{\underset{\ }{\circ}} \ar  @{-}[r]^{\ztu}  & \overset{\zeta}{\underset{\
		}{\circ}} \ar  @{-}[r]^{\ztu}  & \overset{\zeta}{\underset{\ }{\circ}} \ar  @{-}[r]^{\ztu}
	& \overset{-1}{\underset{\ }{\circ}}}$
the algebra $\toba$ is generated by $y_i$, $i\in\I_5$, with defining relations
$y_{ij}=\lambda_{ij}, \, i<j, \, \widetilde{q}_{ij}=1;$
\begin{align*}
y_{1}^2&=\mu_1; & y_{221}&=0; & y_{553}&=\lambda_{553}; & y_{332}&=\lambda_{332}; \\
y_{4}^2&=\mu_4; & y_{223}&=\lambda_{223};  & y_{334}&=0; & y_{335}&=\lambda_{335}.
\end{align*}
Here we apply Lemma \ref{lem:cut1} for $(i,j)=(3,4)$ and $\toba$ projects onto $\widetilde{\mathcal E}_{\bq'}$, $\bq'$ with two components of types super $A_4$ and $A_1$.

 For $\xymatrix@R-8pt{  &   \overset{\zeta}{\circ} \ar  @{-}[d]^{\ztu} & & \\
	\overset{\ztu}{\underset{\ }{\circ}} \ar  @{-}[r]^{\zeta}  & \overset{-1}{\underset{\
		}{\circ}} \ar  @{-}[r]^{\ztu}  & \overset{\zeta}{\underset{\ }{\circ}} \ar  @{-}[r]^{\ztu}
	& \overset{\zeta}{\underset{\ }{\circ}}}$
the algebra $\toba$ is generated by $y_i$, $i\in\I_5$, with defining relations
$y_{ij}=\lambda_{ij}, \, i<j, \, \widetilde{q}_{ij}=1;$
\begin{align*}
y_{552}&=0; & [y_{(13)},&y_2]_c=\nu_2; & y_{443}&=\lambda_{443}; & [[[y_{4325},&y_{2}]_c,y_3]_c,y_2]_c=0;
\\
y_{112}&=0; & [y_{125},&y_2]_c=\nu_2'; & y_{334}&=\lambda_{334}; & y_{332}&=0;
\quad  y_{2}^2=\mu_2.
\end{align*}
If $\nu_2'=0$, then we apply Lemma \ref{lem:cut1} for $(i,j)=(3,4)$ and $\toba$ projects onto $ \widetilde{\mathcal E}_{\bq'}$, $\bq'$ with two components of types $A_4$, $A_1$. If $\nu_2'\neq 0$ but all the other $\lambda$'s are zero, then we apply Lemma \ref{lem:cut1} for $(i,j)=(2,3)$ and $\toba$ projects onto $\widetilde{\mathcal E}_{\bq'}$, $\bq'$ with two components of types super $A_3$, $A_2$. The general case follows by Lemma \ref{lem:nonzero-k<l}.

\subsubsection{\ Type $\g(8,3)$, $\zeta \in \G'_3$}\label{subsec:type-g(8,3)}
The Weyl groupoid has 21 objects:

 For $\xymatrix@C-4pt{\overset{-1}{\underset{\ }{\circ}}\ar  @{-}[r]^{\zeta}  &
	\overset{\ztu}{\underset{\ }{\circ}} \ar  @{-}[r]^{\zeta}  & \overset{\ztu}{\underset{\
		}{\circ}}\ar  @{-}[r]^{\ztu}
	& \overset{\zeta}{\underset{\ }{\circ}} \ar  @{-}[r]^{\ztu}  & \overset{\zeta}{\underset{\
		}{\circ}}}$, the algebra $\toba$ is generated by $y_i$, $i\in\I_5$, with defining relations
\begin{align*}
y_{221}&=0; & y_{223}&=\lambda_{223}; & y_{332}&=\lambda_{332}; & &y_{ij}=\lambda_{ij}, \, i<j, \, \widetilde{q}_{ij}=1;\\
& & y_{443}&=0; & y_{445}&=\lambda_{445}; & &[y_{3345},y_{34}]_c=\lambda_{345};\\
& & y_{554}&=\lambda_{554}; & y_{1}^2&=\mu_1; & &[[y_{(24)},y_{3}]_c,y_3]_c=\nu_3;.
\end{align*}
Here we apply Lemma \ref{lem:cut1} for $(i,j)=(1,2)$ and $\toba$ projects onto $\widetilde{\mathcal E}_{\bq'}$, $\bq'$ with two components of types $F_4$ and $A_1$.

 For $\xymatrix@C-4pt{\overset{-1}{\underset{\ }{\circ}}\ar  @{-}[r]^{\ztu}  &
	\overset{\zeta}{\underset{\ }{\circ}} \ar  @{-}[r]^{\ztu}  & \overset{\zeta}{\underset{\
		}{\circ}}
	\ar  @{-}[r]^{\ztu}  & \overset{\ztu}{\underset{\ }{\circ}} \ar  @{-}[r]^{\zeta}  &
	\overset{-1}{\underset{\ }{\circ}}}$, the algebra $\toba$ is generated by $y_i$, $i\in\I_5$, with defining relations
\begin{align*}
y_{221}&=0; & y_{223}&=\lambda_{223}; & y_{332}&=\lambda_{332}; & &y_{ij}=\lambda_{ij}, \, i<j, \, \widetilde{q}_{ij}=1;\\
& & y_{334}&=0; & y_{445}&=0; & &[y_{4432},y_{43}]_c=\lambda_{432};\\
& & y_{1}^2&=\mu_1; & y_{5}^2&=\mu_5; & &[[y_{(35)},y_{4}]_c,y_4]_c=0;.
\end{align*}
Here we apply Lemma \ref{lem:cut1} for $(i,j)=(4,5)$ and $\toba$ projects onto $\widetilde{\mathcal E}_{\bq'}$, $\bq'$ with two components of types super $B_4$ and $A_1$.

 For $\xymatrix@C-4pt{\overset{-1}{\underset{\ }{\circ}}\ar  @{-}[r]^{\zeta}  &
	\overset{-1}{\underset{\ }{\circ}} \ar  @{-}[r]^{\ztu}  & \overset{\zeta}{\underset{\
		}{\circ}}
	\ar  @{-}[r]^{\ztu}  & \overset{\ztu}{\underset{\ }{\circ}} \ar  @{-}[r]^{\zeta}  &
	\overset{-1}{\underset{\ }{\circ}}}$, the algebra $\toba$ is generated by $y_i$, $i\in\I_5$, with defining relations
\begin{align*}
[y_{(13)},&y_2]_c=0; & y_{332}&=0; & y_{334}&=0; & &y_{ij}=\lambda_{ij}, \, i<j, \, \widetilde{q}_{ij}=1;\\
& & y_{445}&=0; & y_{1}^2&=\mu_1; & &[y_{4432},y_{43}]_c=0;\\
& & y_{2}^2&=\mu_2; & y_{5}^2&=\mu_5; & &[[y_{(35)},y_{4}]_c,y_4]_c=0.
\end{align*}
Here we apply Lemma \ref{lem:cut1} for $(i,j)=(4,5)$ and $\toba$ projects onto $\widetilde{\mathcal E}_{\bq'}$, $\bq'$ with two components of types super $B_4$ and $A_1$.

 For $\xymatrix@C-4pt{\overset{-1}{\underset{\ }{\circ}}\ar  @{-}[r]^{\ztu}  &
	\overset{\zeta}{\underset{\ }{\circ}} \ar  @{-}[r]^{\ztu}  & \overset{\zeta}{\underset{\
		}{\circ}}
	\ar  @{-}[r]^{\ztu}  & \overset{-1}{\underset{\ }{\circ}} \ar  @{-}[r]^{\ztu}  &
	\overset{-1}{\underset{\ }{\circ}}}$ the algebra $\toba$ is generated by $y_i$, $i\in\I_5$, with defining relations
\begin{align*}
& & y_{221}&=0; & y_{223}&=\lambda_{223}; & y_{334}&=0; & &y_{ij}=\lambda_{ij}, \, i<j, \, \widetilde{q}_{ij}=1;\\
y_{4}^2&=\mu_4; & y_{5}^2&=\mu_5; & y_{332}&=\lambda_{332}; & y_{1}^2&=\mu_1; & &[[y_{54},y_{543}]_c,y_4]_c=0.
\end{align*}
Here we apply Lemma \ref{lem:cut1} for $(i,j)=(3,4)$ and $\toba$ projects onto $\widetilde{\mathcal E}_{\bq'}$, $\bq'$ with two components of types super $A_3$ and $A_2$.

 For $\xymatrix@C-4pt{\overset{\zeta}{\underset{\ }{\circ}}\ar  @{-}[r]^{\ztu}  &
	\overset{-1}{\underset{\ }{\circ}} \ar  @{-}[r]^{\zeta}  & \overset{-1}{\underset{\
		}{\circ}}
	\ar  @{-}[r]^{\ztu}  & \overset{\ztu}{\underset{\ }{\circ}} \ar  @{-}[r]^{\zeta}  &
	\overset{-1}{\underset{\ }{\circ}}}$ the algebra $\toba$ is generated by $y_i$, $i\in\I_5$, with defining relations
\begin{align*}
y_{112}&=0; & [y_{443},&y_{43}]_c=\lambda_{43}; & y_{445}&=0; & & y_{ij}=\lambda_{ij}, \, i<j, \, \widetilde{q}_{ij}=1;\\
& & [y_{(13)},&y_2]_c=0; & y_{2}^2&=\mu_2; & &[y_{4432},y_{43}]_c=0;\\
y_{5}^2&=\mu_5; & [y_{(24)},&y_3]_c=0; & y_{3}^2&=\mu_3; & &[[y_{(35)},y_4]_c,y_4]_c=\lambda_{(35)}.
\end{align*}
Here we apply Lemma \ref{lem:cut1} for $(i,j)=(2,3)$ and $\toba$ projects onto $\widetilde{\mathcal E}_{\bq'}$, $\bq'$ with two components of types super $A_2$ and $\g(2,3)$.

 For $\xymatrix@C-4pt{\overset{-1}{\underset{\ }{\circ}}\ar  @{-}[r]^{\zeta}  &
	\overset{-1}{\underset{\ }{\circ}} \ar  @{-}[r]^{\ztu}  & \overset{\zeta}{\underset{\
		}{\circ}}
	\ar  @{-}[r]^{\ztu}  & \overset{-1}{\underset{\ }{\circ}} \ar  @{-}[r]^{\ztu}  &
	\overset{-1}{\underset{\ }{\circ}}}$ the algebra $\toba$ is generated by $y_i$, $i\in\I_5$, with defining relations
\begin{align*}
y_{5}^2&=\mu_5; & y_{332}&=0; & y_{334}&=0; & y_{1}^2&=\mu_1; & & y_{ij}=\lambda_{ij}, \, i<j, \, \widetilde{q}_{ij}=1;\\
& & [y_{(13)},&y_2]_c=0; & y_{2}^2&=\mu_2; & y_{4}^2&=\mu_4; & &[[y_{54},y_{543}]_c,y_4]_c=0.
\end{align*}
Here we apply Lemma \ref{lem:cut1} for $(i,j)=(3,4)$ and $\toba$ projects onto $\widetilde{\mathcal E}_{\bq'}$, $\bq'$ with two components of types super $A_2$ and $A_3$.

 For $\xymatrix@C-4pt{\overset{\zeta}{\underset{\ }{\circ}}\ar  @{-}[r]^{\ztu}  &
	\overset{\zeta}{\underset{\ }{\circ}} \ar  @{-}[r]^{\ztu}  & \overset{-1}{\underset{\
		}{\circ}}
	\ar  @{-}[r]^{\zeta}  & \overset{-\zeta}{\underset{\ }{\circ}} \ar  @{-}[r]^{\zeta}  &
	\overset{-1}{\underset{\ }{\circ}}}$ the algebra $\toba$ is generated by $y_i$, $i\in\I_5$, with defining relations
\begin{align*}
y_{3}^2&=\mu_3; & y_{221}&=\lambda_{221}; & y_{4443}&=\lambda_{4443}; & & y_{ij}=\lambda_{ij}, \, i<j, \, \widetilde{q}_{ij}=1;\\
y_{5}^2&=\mu_5; & y_{112}&=\lambda_{112}; & y_{4445}&=\lambda_{4445}; & &[y_{(24)},y_{3}]_c=0;
\end{align*}
\vsp
\begin{align*}
y_{223}&=0; & [y_3,y_{445}]_c+q_{45}[y_{(35)},y_4]_c+(\zeta-\ztu)q_{34}y_4y_{(35)}=0.
\end{align*}
Here we apply Lemma \ref{lem:cut1} for $(i,j)=(2,3)$ and $\toba$ projects onto $\widetilde{\mathcal E}_{\bq'}$, $\bq'$ with two components of types super $A_2$ and $\g(2,3)$.

 For $\xymatrix@C-4pt{\overset{\zeta}{\underset{\ }{\circ}}\ar  @{-}[r]^{\ztu}  &
	\overset{\zeta}{\underset{\ }{\circ}} \ar  @{-}[r]^{\ztu}  & \overset{-1}{\underset{\
		}{\circ}}\ar  @{-}[r]^{\zeta}
	& \overset{\ztu}{\underset{\ }{\circ}} \ar  @{-}[r]^{\ztu}  & \overset{-1}{\underset{\
		}{\circ}}}$ the algebra $\toba$ is generated by $y_i$, $i\in\I_5$, with defining relations
\begin{align*}
y_{3}^2&=\mu_3; & [y_{(24)},&y_3]_c=\nu_3; & y_{112}&=\lambda_{112}; & y_{223}&=0; \\
y_{5}^2&=\mu_5; & [y_{445},&y_{45}]_c=\lambda_{45}; & y_{221}&=\lambda_{221}; & y_{445}&=0;
\end{align*}
\vsp 
\begin{align*}
 y_{ij}&=\lambda_{ij}, \, i<j, \, \widetilde{q}_{ij}=1; & [[y_{(35)},y_{4}]_c,y_4]_c&=\lambda_{(35)}.
\end{align*}
If $\nu_3=0$, then we apply Lemma \ref{lem:cut1} for $(i,j)=(2,3)$ and $\toba$ projects onto $ \widetilde{\mathcal E}_{\bq'}$, $\bq'$ with two components of types super $A_2$ and $\g(2,3)$. If $\nu_3\neq 0$, then $\lambda_{(35)}= \lambda_{45}=0$, so we apply Lemma \ref{lem:cut1} for $(i,j)=(4,5)$ and $\toba$ projects onto $\widetilde{\mathcal E}_{\bq'}$, $\bq'$ with two components of types super $A_4$ and $A_1$.

 For $\xymatrix@C-4pt{\overset{\zeta}{\underset{\ }{\circ}}\ar  @{-}[r]^{\ztu}  &
	\overset{-1}{\underset{\ }{\circ}} \ar  @{-}[r]^{\zeta}  & \overset{-1}{\underset{\
		}{\circ}}\ar  @{-}[r]^{\ztu}
	& \overset{-1}{\underset{\ }{\circ}} \ar  @{-}[r]^{\ztu}  & \overset{-1}{\underset{\
		}{\circ}}}$ the algebra $\toba$ is generated by $y_i$, $i\in\I_5$, with defining relations
\begin{align*}
y_{4}^2&=\mu_4; & y_{5}^2&=\mu_5; & [y_{(24)},&y_3]_c=\nu_3; & & y_{ij}=\lambda_{ij}, \, i<j, \, \widetilde{q}_{ij}=1;\\
&& y_{2}^2&=\mu_2; &  [y_{(13)},&y_2]_c=0; &   &[[y_{34},y_{(35)}]_c,y_4]_c=\lambda_{345};\\
&& y_{3}^2&=\mu_3; & y_{112}&=0; & &[[y_{54},y_{543}]_c,y_4]_c=\lambda_{543}.
\end{align*}
Here we apply Lemma \ref{lem:cut1} for $(i,j)=(1,2)$ and $\toba$ projects onto $\widetilde{\mathcal E}_{\bq'}$, $\bq'$ with two components of types $A_1$ and $\g(4,3)$.

 For $\xymatrix@C-4pt{\overset{\ztu}{\underset{\ }{\circ}}\ar  @{-}[r]^{\zeta}  &
	\overset{-1}{\underset{\ }{\circ}} \ar  @{-}[r]^{\ztu}  & \overset{-1}{\underset{\
		}{\circ}}\ar  @{-}[r]^{\ztu}
	& \overset{\zeta}{\underset{\ }{\circ}} \ar  @{-}[r]^{\ztu}  & \overset{\zeta}{\underset{\
		}{\circ}}}$ the algebra $\toba$ is generated by $y_i$, $i\in\I_5$, with defining relations
\begin{align*}
[y_{(13)},&y_2]_c=0; & y_{2}^2&=\mu_2; & y_{112}&=0; & y_{445}&=\lambda_{445}; & & y_{ij}=\lambda_{ij}, \, i<j, \, \widetilde{q}_{ij}=1;\\
& & y_{3}^2&=\mu_3; &  y_{443}&=0; & y_{554}&=\lambda_{554}; & &[[y_{23},y_{(24)}]_c,y_3]_c=0.
\end{align*}
Here we apply Lemma \ref{lem:cut1} for $(i,j)=(2,3)$ and $\toba$ projects onto $\widetilde{\mathcal E}_{\bq'}$, $\bq'$ with two components of types super $A_2$ and $A_3$.

 For $\xymatrix@C-4pt{\overset{-1}{\underset{\ }{\circ}}\ar  @{-}[r]^{\ztu}  &
	\overset{-1}{\underset{\ }{\circ}} \ar  @{-}[r]^{\zeta}  & \overset{\ztu}{\underset{\
		}{\circ}}\ar  @{-}[r]^{\ztu}
	& \overset{\zeta}{\underset{\ }{\circ}} \ar  @{-}[r]^{\ztu}  & \overset{\zeta}{\underset{\
		}{\circ}}}$ the algebra $\toba$ is generated by $y_i$, $i\in\I_5$, with defining relations
\begin{align*}
y_{1}^2&=\mu_1; & [y_{(13)},&y_2]_c=0; & y_{332}&=0; & y_{445}&=\lambda_{445};\\
y_{2}^2&=\mu_2; & [y_{3345},&y_{34}]_c=\lambda_{(35)}'; & y_{443}&=0; & y_{554}&=\lambda_{554};
\end{align*}
\vsp 
\begin{align*}
 y_{ij}&=\lambda_{ij}, \, i<j, \, \widetilde{q}_{ij}=1; & [[y_{(35)},y_{4}]_c,y_4]_c&=\lambda_{(35)}.
\end{align*}
Here we apply Lemma \ref{lem:cut1} for $(i,j)=(2,3)$ and $\toba$ projects onto $\widetilde{\mathcal E}_{\bq'}$, $\bq'$ with two components of types super $A_2$ and $C_3$.

 For $\xymatrix@R-8pt{  &   & \overset{-1}{\circ} \ar  @{-}[d]^{\zeta}\ar@{-}[dr]^{\zeta}  & \\
	\overset{-1}{\underset{\ }{\circ}} \ar  @{-}[r]^{\ztu}  & \overset{\zeta}{\underset{\
		}{\circ}} \ar  @{-}[r]^{\ztu}  & \overset{-1}{\underset{\ }{\circ}} \ar  @{-}[r]^{\zeta}
	& \overset{\ztu}{\underset{\ }{\circ}}}$ the algebra $\toba$ is generated by $y_i$, $i\in\I_5$, with defining relations
\begin{align*}
y_1^2&=\mu_1; & [y_{(24)},&y_3]_c=\nu_3; & y_{221}&=0; & y_{223}&=0; & & y_{ij}=\lambda_{ij}, \, i<j, \, \widetilde{q}_{ij}=1;\\
y_3^2&=\mu_3; & y_5^2&=\mu_5; & y_{443}&=0; & y_{445}&=0; & &[y_{235},y_3]_c=0;
\end{align*}
\vsp
\begin{align*}
y_{(35)}-q_{45}\ztu[y_{35},y_4]_c-q_{34}(1-\zeta)y_4y_{35}=0.
\end{align*}
Here we apply Lemma \ref{lem:cut2} for $(i,j,k)=(4,3,5)$ and $\toba$ projects onto $ \widetilde{\mathcal E}_{\bq'}$, $\bq'$ with two components of types super $A_4$ and $A_1$.

 For $\xymatrix@R-8pt{  &   & \overset{-1}{\circ} \ar  @{-}[d]^{\zeta} \ar@{-}[dr]^{\zeta} &
	\\
	\overset{-1}{\underset{\ }{\circ}} \ar  @{-}[r]^{\zeta}  & \overset{-1}{\underset{\
		}{\circ}} \ar  @{-}[r]^{\ztu}  & \overset{-1}{\underset{\ }{\circ}} \ar  @{-}[r]^{\zeta}
	& \overset{\ztu}{\underset{\ }{\circ}}}$ the algebra $\toba$ is generated by $y_i$, $i\in\I_5$, with defining relations
\begin{align*}
[y_{(13)},&y_2]_c=\nu_2; & y_{443}&=0; & y_1^2&=\mu_1; & & y_{ij}=\lambda_{ij}, \, i<j, \, \widetilde{q}_{ij}=1;\\
[y_{(24)},&y_3]_c=0; & y_{445}&=0; & y_2^2&=\mu_2; & &[y_{235},y_3]_c=\nu_3';
\end{align*}
\vsp
\begin{align*}
y_3^2&=\mu_3; & y_5^2&=\mu_5; & y_{(35)}-q_{45}\ztu[y_{35},y_4]_c-q_{34}(1-\zeta)y_4y_{35}=0.
\end{align*}
Here we apply Lemma \ref{lem:cut2} for $(i,j,k)=(4,3,5)$ and $\toba$ projects onto $ \widetilde{\mathcal E}_{\bq'}$, $\bq'$ with two components of types super $A_4$ and $A_1$.

 For $\xymatrix@R-8pt{  &    \overset{-1}{\circ} \ar  @{-}[d]^{\zeta} \ar  @{-}[dr]^{\zeta}
	&  & \\
	\overset{\ztu}{\underset{\ }{\circ}} \ar  @{-}[r]^{\zeta}  & \overset{\ztu}{\underset{\
		}{\circ}} \ar  @{-}[r]^{\zeta}  & \overset{-1}{\underset{\ }{\circ}} \ar  @{-}[r]^{\ztu}
	& \overset{\zeta}{\underset{\ }{\circ}}}$ the algebra $\toba$ is generated by $y_i$, $i\in\I_5$, with defining relations
\begin{align*}
[y_{(24)},&y_3]_c=\nu_3; & y_{112}&=\lambda_{112}; & y_{223}&=0; & & y_{ij}=\lambda_{ij}, \, i<j, \, \widetilde{q}_{ij}=1;\\
y_{443}&=0; & y_{221}&=\lambda_{221}; & y_{225}&=0; & &[y_{435},y_3]_c=0;
\end{align*}
\vsp
\begin{align*}
y_3^2&=\mu_3; & y_5^2&=\mu_5; & y_{235}-q_{35}\ztu[y_{25},y_3]_c-q_{23}(1-\zeta)y_3y_{25}=0.
\end{align*}
Here we apply Lemma \ref{lem:cut2} for $(i,j,k)=(5,2,3)$ and $\toba$ projects onto $ \widetilde{\mathcal E}_{\bq'}$, $\bq'$ with two components of types super $A_4$ and $A_1$.

 For $\xymatrix@R-8pt{  &   & \overset{-1}{\circ} \ar  @{-}[d]^{\zeta} \ar@{-}[dr]^{\zeta} &
	\\
	\overset{\zeta}{\underset{\ }{\circ}} \ar  @{-}[r]^{\ztu}  & \overset{-1}{\underset{\
		}{\circ}} \ar  @{-}[r]^{\zeta}  &
	\overset{\ztu}{\underset{\ }{\circ}} \ar  @{-}[r]^{\zeta}  & \overset{\ztu}{\underset{\
		}{\circ}}}$ the algebra $\toba$ is generated by $y_i$, $i\in\I_5$, with defining relations
\begin{align*}
[y_{(13)},&y_2]_c=\nu_2; & y_{112}&=0; & y_{332}&=0; & & y_{ij}=\lambda_{ij}, \, i<j, \, \widetilde{q}_{ij}=1;\\
y_{445}&=0; & y_{334}&=\lambda_{334}; & y_{335}&=0; & &y_{443}=\lambda_{443};
\end{align*}
\vsp
\begin{align*}
y_2^2&=\mu_2; & y_5^2&=\mu_5; &
y_{(35)}-q_{45}\ztu[y_{35},y_4]_c-q_{34}(1-\zeta)y_4y_{35}=0.
\end{align*}
Here we apply Lemma \ref{lem:cut2} for $(i,j,k)=(5,3,4)$ and $\toba$ projects onto $ \widetilde{\mathcal E}_{\bq'}$, $\bq'$ with two components of types super $A_4$ and $A_1$.

 For $\xymatrix@R-8pt{  &   & \overset{\zeta}{\circ} \ar  @{-}[d]^{\ztu}  & \\
	\overset{-1}{\underset{\ }{\circ}} \ar  @{-}[r]^{\zeta}  & \overset{\ztu}{\underset{\
		}{\circ}} \ar  @{-}[r]^{\zeta}  & \overset{-1}{\underset{\ }{\circ}} \ar  @{-}[r]^{\ztu}
	& \overset{-1}{\underset{\ }{\circ}}}$ the algebra $\toba$ is generated by $y_i$, $i\in\I_5$, with defining relations
\begin{align*}
y_{3}^2&=\mu_3; & [y_{(24)},&y_3]_c=0; & y_{221}&=0; & y_{223}&=0; & & y_{ij}=\lambda_{ij}, \, i<j, \, \widetilde{q}_{ij}=1;\\
y_{4}^2&=\mu_4; & [y_{235},&y_3]_c=\nu_3'; &  y_{553}&=0; & y_1^2&=\mu_1; & &[[y_{43},y_{435}]_c,y_3]_c=0.
\end{align*}
Here we apply Lemma \ref{lem:cut1} for $(i,j)=(3,5)$ and $\toba$ projects onto $\widetilde{\mathcal E}_{\bq'}$, $\bq'$ with two components of types super $A_4$ and $A_1$.

 For $\xymatrix@R-8pt{  &   & \overset{\zeta}{\circ} \ar  @{-}[d]^{\ztu}  & \\
	\overset{\ztu}{\underset{\ }{\circ}} \ar  @{-}[r]^{\zeta}  & \overset{-1}{\underset{\
		}{\circ}} \ar  @{-}[r]^{\ztu}  & \overset{\zeta}{\underset{\ }{\circ}} \ar  @{-}[r]^{\ztu}
	& \overset{-1}{\underset{\ }{\circ}}}$ the algebra $\toba$ is generated by $y_i$, $i\in\I_5$, with defining relations
\begin{align*}
y_{112}&=0; & y_{2}^2&=\mu_2; & y_{332}&=0; & y_{335}&=\lambda_{335}; & & y_{ij}=\lambda_{ij}, \, i<j, \, \widetilde{q}_{ij}=1;\\
& & y_{4}^2&=\mu_4; & y_{334}&=0; & y_{553}&=\lambda_{553}; & &[y_{(13)},y_2]_c=\nu_2.
\end{align*}
Here we apply Lemma \ref{lem:cut1} for $(i,j)=(3,4)$ and $\toba$ projects onto $\widetilde{\mathcal E}_{\bq'}$, $\bq'$ with two components of types super $A_4$ and $A_1$.

 For $\xymatrix@R-8pt{  &   & \overset{\zeta}{\circ} \ar  @{-}[d]^{\ztu}  & \\
	\overset{-1}{\underset{\ }{\circ}} \ar  @{-}[r]^{\ztu}  & \overset{-1}{\underset{\
		}{\circ}} \ar  @{-}[r]^{\zeta}  & \overset{\ztu}{\underset{\ }{\circ}} \ar
	@{-}[r]^{\zeta}  & \overset{-1}{\underset{\ }{\circ}}}$ the algebra $\toba$ is generated by $y_i$, $i\in\I_5$, with defining relations
\begin{align*}
y_1^2&=\mu_1; & [[y_{235},&y_3]_c,y_3]_c=0; & [y_{(13)},&y_2]_c=0; &  y_{ij}&=\lambda_{ij}, \, i<j, \, \widetilde{q}_{ij}=1;\\
y_{2}^2&=\mu_2; & [[y_{435},&y_3]_c,y_3]_c=0; &  y_4^2&=\mu_4; & y_{332}&=y_{334}=y_{553}=0.
\end{align*}
This algebra is nonzero by \cite[Lemma 5.16]{AAGMV}.

 For $\xymatrix@R-8pt{  &   & \overset{\zeta}{\circ} \ar  @{-}[d]^{\ztu}  & \\
	\overset{-1}{\underset{\ }{\circ}} \ar  @{-}[r]^{\zeta}  & \overset{\ztu}{\underset{\
		}{\circ}} \ar  @{-}[r]^{\zeta}  & \overset{\ztu}{\underset{\ }{\circ}} \ar
	@{-}[r]^{\zeta}  & \overset{-1}{\underset{\ }{\circ}}}$
the algebra $\toba$ is generated by $y_i$, $i\in\I_5$, with defining relations
\begin{align*}
y_{1}^2&=\mu_1; & y_{221}&=y_{334}=y_{553}=0; & y_{223}&=\lambda_{223}; & & y_{ij}=\lambda_{ij}, \, i<j, \, \widetilde{q}_{ij}=1;\\
y_{4}^2&=\mu_4; & [[y_{435},&y_3]_c,y_3]_c=0; & y_{332}&=\lambda_{332}; &  &[[y_{235},y_3]_c,y_3]_c=\lambda_{235}.
\end{align*}
Here we apply Lemma \ref{lem:cut1} for $(i,j)=(3,4)$ and $\toba$ projects onto $\widetilde{\mathcal E}_{\bq'}$, $\bq'$ with two components of types super $CD$ and $A_1$.

 For $\xymatrix@R-8pt{  &   & \overset{\zeta}{\circ} \ar  @{-}[d]^{\ztu}  & \\
	\overset{\ztu}{\underset{\ }{\circ}} \ar  @{-}[r]^{\zeta}  & \overset{-1}{\underset{\
		}{\circ}} \ar  @{-}[r]^{\ztu}  & \overset{-1}{\underset{\ }{\circ}} \ar  @{-}[r]^{\zeta}
	& \overset{-1}{\underset{\ }{\circ}}}$ the algebra $\toba$ is generated by $y_i$, $i\in\I_5$, with defining relations
\begin{align*}
y_{112}&=0; & [y_{(24)},&y_3]_c=\nu_3; & y_{4}^2&=\mu_4; &  & y_{ij}=\lambda_{ij}, \, i<j, \, \widetilde{q}_{ij}=1;\\
&&     [y_{(13)},&y_2]_c=0; &  y_{3}^2&=\mu_3; &   &[[y_{23},y_{235}]_c,y_3]_c=0;\\
&& y_{553}&=0; & y_{2}^2&=\mu_2; & &[y_{435},y_3]_c=0.
\end{align*}
Herewe apply Lemma \ref{lem:cut1} for $(i,j)=(3,5)$ and $\toba$ projects onto $\widetilde{\mathcal E}_{\bq'}$, $\bq'$ with two components of types super $A_4$ and $A_1$.

 For $\xymatrix@R-8pt{  &   & \overset{\zeta}{\circ} \ar  @{-}[d]^{\ztu}  & \\
	\overset{-1}{\underset{\ }{\circ}} \ar  @{-}[r]^{\ztu}  & \overset{-1}{\underset{\
		}{\circ}} \ar  @{-}[r]^{\zeta}  & \overset{-1}{\underset{\ }{\circ}} \ar  @{-}[r]^{\ztu}
	& \overset{-1}{\underset{\ }{\circ}}}$ the algebra $\toba$ is generated by $y_i$, $i\in\I_5$, with defining relations
\begin{align*}
y_{1}^2&=\mu_1; & [y_{(24)},&y_3]_c=\nu_3; & y_{553}&=0; & & y_{ij}=\lambda_{ij}, \, i<j, \, \widetilde{q}_{ij}=1;\\
&&     [y_{(13)},&y_2]_c=\nu_2; &  y_{4}^2&=\mu_4; &   &[[y_{43},y_{435}]_c,y_3]_c=0;\\
&& y_{3}^2&=\mu_3; & y_{2}^2&=\mu_2; & &[y_{235},y_{3}]_c=0.
\end{align*}
Here we apply Lemma \ref{lem:cut1} for $(i,j)=(3,5)$ and $\toba$ projects onto $\widetilde{\mathcal E}_{\bq'}$, $\bq'$ with two components of types super $A_4$ and $A_1$.

\subsubsection{Type $\g(4,6)$, $\zeta \in \G'_3$}\label{subsec:type-g(4,6)}
The Weyl groupoid has seven objects.

 For $\xymatrix@C-4pt{\overset{\zeta}{\underset{\ }{\circ}}\ar  @{-}[r]^{\ztu}  &
\overset{\zeta}{\underset{\ }{\circ}} \ar  @{-}[r]^{\ztu}  & \overset{\zeta}{\underset{\
}{\circ}}
\ar  @{-}[r]^{\ztu}  & \overset{-1}{\underset{\ }{\circ}} \ar  @{-}[r]^{\ztu}  &
\overset{\zeta}{\underset{\ }{\circ}}  \ar  @{-}[r]^{\ztu}  & \overset{\zeta}{\underset{\
}{\circ}}}$ the algebra $\toba$ is generated by $y_i$, $i\in\I_6$, with defining relations
\begin{align*}
y_{112}&=\lambda_{112}; & y_{221}&=\lambda_{221}; & y_{223}&=\lambda_{223}; & & y_{ij}=\lambda_{ij}, \, i<j, \, \widetilde{q}_{ij}=1;\\
y_{665}&=\lambda_{665}; & y_{332}&=\lambda_{332}; & y_{334}&=0; & &[[[y_{(25)},y_{4}]_c,y_3]_c,y_4]_c=0;\\
y_{4}^2&=\mu_4; & y_{554}&=0; & y_{556}&=\lambda_{556}; & &[[[y_{6543},y_{4}]_c,y_5]_c,y_4]_c=0.
\end{align*}
Here we apply Lemma \ref{lem:cut1} for $(i,j)=(3,4)$ and $\toba$ projects onto $\widetilde{\mathcal E}_{\bq'}$, $\bq'$ with two components of types $A_3$.

 For $\xymatrix@C-5pt@R-8pt{  &   & & \overset{\zeta}{\circ} \ar  @{-}[d]^{\ztu}  & \\
\overset{\zeta}{\underset{\ }{\circ}} \ar  @{-}[r]^{\ztu}  & \overset{\zeta}{\underset{\
}{\circ}} \ar  @{-}[r]^{\ztu}  & \overset{\zeta}{\underset{\ }{\circ}} \ar  @{-}[r]^{\ztu}
 & \overset{-1}{\underset{\ }{\circ}}
\ar  @{-}[r]^{\zeta}  & \overset{-1}{\underset{\ }{\circ}}}$ the algebra $\toba$ is generated by $y_i$, $i\in\I_6$, with defining relations
\begin{align*}
y_{4}^2&=\mu_4; & y_{221}&=\lambda_{221}; & y_{223}&=\lambda_{223}; & & y_{ij}=\lambda_{ij}, \, i<j, \, \widetilde{q}_{ij}=1;\\
y_{5}^2&=\mu_5; & y_{334}&=0; & y_{664}&=0; & &[y_{546},y_{4}]_c=0;\\
[y_{(35)},y_{4}]_c&=0; & y_{112}&=\lambda_{112}; & y_{332}&=\lambda_{332}; & &[[[y_{2346},y_4]_c,y_3]_c,y_4]_c=0.
\end{align*}
Here we apply Lemma \ref{lem:cut1} for $(i,j)=(4,5)$ and $\toba$ projects onto $\widetilde{\mathcal E}_{\bq'}$, $\bq'$ with two components of types $CD$ and $A_1$.

 For $\xymatrix@C-5pt@R-8pt{  &   & \overset{\zeta}{\circ} \ar  @{-}[d]^{\ztu}  & &\\
\overset{\zeta}{\underset{\ }{\circ}} \ar  @{-}[r]^{\ztu}  & \overset{-1}{\underset{\
}{\circ}} \ar  @{-}[r]^{\zeta}  & \overset{-1}{\underset{\ }{\circ}} \ar  @{-}[r]^{\ztu}
& \overset{\zeta}{\underset{\ }{\circ}}
\ar  @{-}[r]^{\ztu}  & \overset{\zeta}{\underset{\ }{\circ}}}$ the algebra $\toba$ is generated by $y_i$, $i\in\I_6$, with defining relations
\begin{align*}
[y_{(13)},&y_{2}]_c=0; & y_{445}&=\lambda_{445}; & y_{443}&=0; & & y_{ij}=\lambda_{ij}, \, i<j, \, \widetilde{q}_{ij}=1;\\
y_{112}&=0; & y_{554}&=\lambda_{554}; & y_{663}&=0; & &[y_{236},y_{3}]_c=0;\\
[y_{(24)},&y_{3}]_c=0; & y_{2}^2&=\mu_2; & y_{3}^2&=\mu_3; & &[[[y_{5436},y_3]_c,y_4]_c,y_3]_c=0.
\end{align*}
Here we apply Lemma \ref{lem:cut1} for $(i,j)=(3,4)$ and $\toba$ projects onto $\widetilde{\mathcal E}_{\bq'}$, $\bq'$ with two components of types $A_4$ and $A_2$.

 For $\xymatrix@C-5pt@R-8pt{  &   & & \overset{\zeta}{\circ} \ar  @{-}[d]^{\ztu}  & \\
\overset{\zeta}{\underset{\ }{\circ}} \ar  @{-}[r]^{\ztu}  & \overset{\zeta}{\underset{\
}{\circ}} \ar  @{-}[r]^{\ztu}  & \overset{\zeta}{\underset{\ }{\circ}} \ar  @{-}[r]^{\ztu}
 & \overset{\zeta}{\underset{\ }{\circ}}
\ar  @{-}[r]^{\ztu}  & \overset{-1}{\underset{\ }{\circ}}}$ the algebra $\toba$ is generated by $y_i$, $i\in\I_6$, with defining relations
\begin{align*}
y_{112}&=\lambda_{112}; & y_{221}&=\lambda_{221}; & y_{223}&=\lambda_{223}; & & y_{ij}=\lambda_{ij}, \, i<j, \, \widetilde{q}_{ij}=1;\\
y_{332}&=\lambda_{332}; & y_{334}&=\lambda_{334}; & y_{443}&=\lambda_{443}; & & y_{445}=0;\\
y_{446}&=\lambda_{446}; & y_{664}&=\lambda_{664}; & y_{5}^2&=\mu_5.
\end{align*}
Here we apply Lemma \ref{lem:cut1} for $(i,j)=(4,5)$ and $\toba$ projects onto $\widetilde{\mathcal E}_{\bq'}$, $\bq'$ with two components of types $A_5$ and $A_1$.

 For $\xymatrix@C-5pt@R-8pt{  &   & \overset{\zeta}{\circ} \ar  @{-}[d]^{\ztu}  & &\\
\overset{-1}{\underset{\ }{\circ}} \ar  @{-}[r]^{\zeta}  & \overset{-1}{\underset{\
}{\circ}} \ar  @{-}[r]^{\ztu}  & \overset{\zeta}{\underset{\ }{\circ}} \ar  @{-}[r]^{\ztu}
 & \overset{\zeta}{\underset{\ }{\circ}}
\ar  @{-}[r]^{\ztu}  & \overset{\zeta}{\underset{\ }{\circ}}}$ the algebra $\toba$ is generated by $y_i$, $i\in\I_6$, with defining relations
\begin{align*}
y_{332}&=0; & y_{334}&=\lambda_{334}; & y_{336}&=\lambda_{336}; & & y_{ij}=\lambda_{ij}, \, i<j, \, \widetilde{q}_{ij}=1;\\
y_{663}&=\lambda_{663}; & y_{443}&=\lambda_{443}; & y_{445}&=\lambda_{443}; & &[y_{(13)},y_2]_c=0;\\
y_{554}&=\lambda_{554}; & y_{1}^2&=\mu_1; & y_{2}^2&=\mu_2.
\end{align*}
Here we apply Lemma \ref{lem:cut1} for $(i,j)=(2,3)$ and $\toba$ projects onto $\widetilde{\mathcal E}_{\bq'}$, $\bq'$ with two components of types $A_2$ and $A_4$.

 For $\xymatrix@C-5pt@R-8pt{  &    & \overset{\zeta}{\circ} \ar  @{-}[d]^{\ztu} &  & \\
\overset{-1}{\underset{\ }{\circ}} \ar  @{-}[r]^{\ztu}  & \overset{\zeta}{\underset{\
}{\circ}} \ar  @{-}[r]^{\ztu}  & \overset{\zeta}{\underset{\ }{\circ}} \ar  @{-}[r]^{\ztu}
 & \overset{\zeta}{\underset{\ }{\circ}}
\ar  @{-}[r]^{\ztu}  & \overset{-1}{\underset{\ }{\circ}}}$ the algebra $\toba$ is generated by $y_i$, $i\in\I_6$, with defining relations
\begin{align*}
y_{221}&=0; & y_{223}&=\lambda_{223}; & y_{332}&=\lambda_{332}; & & y_{ij}=\lambda_{ij}, \, i<j, \, \widetilde{q}_{ij}=1;\\
y_{334}&=\lambda_{334}; & y_{443}&=\lambda_{443}; & y_{445}&=0; & & y_{336}=\lambda_{336};\\
y_{663}&=\lambda_{663}; & y_{1}^2&=\mu_1; & y_{5}^2&=\mu_5.
\end{align*}
Here we apply Lemma \ref{lem:cut1} for $(i,j)=(4,5)$ and $\toba$ projects onto $\widetilde{\mathcal E}_{\bq'}$, $\bq'$ with two components of types $CD$ and $A_1$.

 For $\xymatrix@C-5pt@R-8pt{  &   & \overset{-1}{\circ} \ar  @{-}[d]^{\zeta} \ar
@{-}[dr]^{\zeta}  & &\\
\overset{\zeta}{\underset{\ }{\circ}} \ar  @{-}[r]^{\ztu}  & \overset{\zeta}{\underset{\
}{\circ}} \ar  @{-}[r]^{\ztu}  & \overset{-1}{\underset{\ }{\circ}} \ar  @{-}[r]^{\zeta}
& \overset{-1}{\underset{\ }{\circ}}
\ar  @{-}[r]^{\ztu}  & \overset{\zeta}{\underset{\ }{\circ}}}$
the algebra $\toba$ is generated by $y_i$, $i\in\I_6$, with defining relations
\begin{align*}
y_{112}&=\lambda_{112}; & y_{6}^2&=\mu_6; & y_{223}&=0; & & y_{ij}=\lambda_{ij}, \, i<j, \, \widetilde{q}_{ij}=1;\\
y_{554}&=0; & y_{4}^2&=\mu_4; & [y_{(24)},&y_3]_c=0; &  &[y_{236},y_3]_c=0;\\
y_{221}&=\lambda_{221}; & y_{3}^2&=\mu_3; & [y_{(35)},&y_4]_c=0; & &[y_{546},y_4]_c=0;
\end{align*}
\vsp
\begin{align*}
y_{346}-q_{46}\ztu [y_{36},y_4]_c-q_{34}(1-\zeta)y_4y_{36}=0.
\end{align*}
Here we apply Lemma \ref{lem:cut2} for $(i,j,k)=(6,3,4)$ and $\toba$ projects onto $ \widetilde{\mathcal E}_{\bq'}$, $\bq'$ with two components of types $A_5$ and $A_1$.

\subsubsection{Type $\g(6,6)$, $\zeta \in \G'_3$}\label{subsec:type-g(6,6)}
The Weyl groupoid has 21 objects.

 For $\xymatrix@C-6pt{\overset{\zeta}{\underset{\ }{\circ}}\ar  @{-}[r]^{\ztu}  &
\overset{\zeta}{\underset{\ }{\circ}} \ar  @{-}[r]^{\ztu}  & \overset{\zeta}{\underset{\
}{\circ}}
\ar  @{-}[r]^{\ztu}  & \overset{\zeta}{\underset{\ }{\circ}} \ar  @{-}[r]^{\ztu}  &
\overset{\ztu}{\underset{\ }{\circ}}  \ar  @{-}[r]^{\zeta}  & \overset{-1}{\underset{\
}{\circ}}}$ the algebra $\toba$ is generated by $y_i$, $i\in\I_6$, with defining relations
\begin{align*}
y_{112}&=\lambda_{112}; & y_{221}&=\lambda_{221}; & y_{223}&=\lambda_{223}; & & y_{ij}=\lambda_{ij}, \, i<j, \, \widetilde{q}_{ij}=1;\\
y_{556}&=0 & y_{332}&=\lambda_{332}; & y_{334}&=\lambda_{334}; &  &[[y_{(46)},y_5]_c,y_5]_c=0;\\
y_{6}^2&=\mu_6; & y_{443}&=\lambda_{443}; & y_{445}&=\lambda_{445}; & &[y_{5543},y_{54}]_c=\lambda_{543}.
\end{align*}
Here we apply Lemma \ref{lem:cut1} for $(i,j)=(5,6)$ and $\toba$ projects onto $\widetilde{\mathcal E}_{\bq'}$, $\bq'$ with two components of types $B_5$ and $A_1$.

 For $\xymatrix@C-6pt{\overset{\zeta}{\underset{\ }{\circ}}\ar  @{-}[r]^{\ztu}  &
\overset{\zeta}{\underset{\ }{\circ}} \ar  @{-}[r]^{\ztu}  & \overset{\zeta}{\underset{\
}{\circ}}
\ar  @{-}[r]^{\ztu}    & \overset{\zeta}{\underset{\ }{\circ}}  \ar  @{-}[r]^{\ztu} &
\overset{-1}{\underset{\ }{\circ}} \ar  @{-}[r]^{\ztu}  & \overset{-1}{\underset{\
}{\circ}}}$ the algebra $\toba$ is generated by $y_i$, $i\in\I_6$, with defining relations
\begin{align*}
y_{5}^2&=\mu_5; & y_{112}&=\lambda_{112}; & y_{221}&=\lambda_{221}; & & y_{ij}=\lambda_{ij}, \, i<j, \, \widetilde{q}_{ij}=1;\\
y_{6}^2&=\mu_6; & y_{223}&=\lambda_{223}; & y_{332}&=\lambda_{332}; &  &[[y_{65},y_{654}]_c,y_5]_c=0;\\
& & y_{334}&=\lambda_{334}; & y_{443}&=\lambda_{443}; & &y_{445}=0.
\end{align*}
Here we apply Lemma \ref{lem:cut1} for $(i,j)=(5,6)$ and $\toba$ projects onto $\widetilde{\mathcal E}_{\bq'}$, $\bq'$ with two components of types $A_5$ and $A_1$.

 For $\xymatrix@C-7pt{\overset{\zeta}{\underset{\ }{\circ}} \ar  @{-}[r]^{\ztu}  &
\overset{-1}{\underset{\ }{\circ}}\ar  @{-}[r]^{\zeta}  & \overset{-1}{\underset{\
}{\circ}}\ar  @{-}[r]^{\ztu}  &
\overset{-1}{\underset{\ }{\circ}}\ar  @{-}[r]^{\ztu}  & \overset{\zeta}{\underset{\
}{\circ}}
\ar  @{-}[r]^{\ztu}  & \overset{\zeta}{\underset{\ }{\circ}}}$ the algebra $\toba$ is generated by $y_i$, $i\in\I_6$, with defining relations
\begin{align*}
y_{112}&=0; & y_{2}^2&=\mu_2; & y_{554}&=0; & & y_{ij}=\lambda_{ij}, \, i<j, \, \widetilde{q}_{ij}=1;\\
y_{556}&=\lambda_{556}; & y_{3}^2&=\mu_3; & y_{665}&=\lambda_{665}; & &[[y_{34},y_{(35)}]_c,y_4]_c=0;\\
& & y_{4}^2&=\mu_4; & [y_{(13)},&y_2]_c=0; & &[y_{(24)},y_3]_c=\nu_3.
\end{align*}
Here we apply Lemma \ref{lem:cut1} for $(i,j)=(4,5)$ and $\toba$ projects onto $\widetilde{\mathcal E}_{\bq'}$, $\bq'$ with two components of types $A_4$ and $A_2$.

 For $\xymatrix@C-7pt{\overset{-1}{\underset{\ }{\circ}} \ar  @{-}[r]^{\zeta}  &
\overset{-1}{\underset{\ }{\circ}}
\ar  @{-}[r]^{\ztu}  & \overset{\zeta}{\underset{\ }{\circ}} \ar  @{-}[r]^{\ztu}  &
\overset{-1}{\underset{\ }{\circ}} \ar  @{-}[r]^{\ztu}  & \overset{\zeta}{\underset{\
}{\circ}}\ar  @{-}[r]^{\ztu}  & \overset{\zeta}{\underset{\ }{\circ}}}$ the algebra $\toba$ is generated by $y_i$, $i\in\I_6$, with defining relations
\begin{align*}
y_1^2&=\mu_1; & y_{334}&=0; & y_{556}&=\lambda_{556}; & & y_{ij}=\lambda_{ij}, \, i<j, \, \widetilde{q}_{ij}=1;\\
y_2^2&=\mu_2; & y_{554}&=0; & y_{665}&=\lambda_{665}; &  &[[[y_{(25)},y_4]_c,y_3]_c,y_4]_c=\lambda_{_{(25)}};\\
y_4^2&=\mu_4; & y_{332}&=0; & [y_{(13)},&y_2]_c=0; & &[[[y_{6543},y_4]_c,y_5]_c,y_4]_c=0.
\end{align*}
Here we apply Lemma \ref{lem:cut1} for $(i,j)=(1,2)$ and $\toba$ projects onto $\widetilde{\mathcal E}_{\bq'}$, $\bq'$ with two components of types $\el(5,3)$ and $A_1$.

 For $\xymatrix@C-6pt{\overset{-1}{\underset{\ }{\circ}}\ar  @{-}[r]^{\ztu}  &
\overset{\zeta}{\underset{\ }{\circ}} \ar  @{-}[r]^{\ztu}  & \overset{\zeta}{\underset{\
}{\circ}}
\ar  @{-}[r]^{\ztu}  & \overset{-1}{\underset{\ }{\circ}} \ar  @{-}[r]^{\ztu}  &
\overset{\zeta}{\underset{\ }{\circ}}\ar  @{-}[r]^{\ztu}  & \overset{\zeta}{\underset{\
}{\circ}}}$ the algebra $\toba$ is generated by $y_i$, $i\in\I_6$, with defining relations
\begin{align*}
y_{221}&=0; & y_{223}&=\lambda_{223}; & y_{332}&=\lambda_{332}; & & y_{ij}=\lambda_{ij}, \, i<j, \, \widetilde{q}_{ij}=1;\\
y_{1}^2&=\mu_1; & y_{334}&=0; & y_{554}&=0; & &[[[y_{(25)},y_4]_c,y_3]_c,y_4]_c=0;\\
y_{4}^2&=\mu_4; & y_{556}&=\lambda_{556}; & y_{665}&=\lambda_{665}; & &[[[y_{6543},y_4]_c,y_5]_c,y_4]_c=0.
\end{align*}
Here we apply Lemma \ref{lem:cut1} for $(i,j)=(4,5)$ and $\toba$ projects onto $\widetilde{\mathcal E}_{\bq'}$, $\bq'$ with two components of types $A_4$ and $A_2$.

 For $\xymatrix@C-7pt{\overset{\zeta}{\underset{\ }{\circ}}\ar  @{-}[r]^{\ztu}  &
\overset{\zeta}{\underset{\ }{\circ}} \ar  @{-}[r]^{\ztu}  & \overset{-1}{\underset{\
}{\circ}}
\ar  @{-}[r]^{\zeta}  & \overset{\ztu}{\underset{\ }{\circ}} \ar  @{-}[r]^{\ztu}  &
\overset{\zeta}{\underset{\ }{\circ}}\ar  @{-}[r]^{\ztu}  & \overset{\zeta}{\underset{\
}{\circ}}}$ the algebra $\toba$ is generated by $y_i$, $i\in\I_6$, with defining relations
\begin{align*}
y_{3}^2&=\mu_3; & y_{112}&=\lambda_{112}; & y_{221}&=\lambda_{221}; & & y_{ij}=\lambda_{ij}, \, i<j, \, \widetilde{q}_{ij}=1;\\
y_{665}&=\lambda_{665}; & y_{223}&=0; & y_{443}&=0; &  &[[y_{(35)},y_4]_c,y_4]_c=0;\\
[y_{(24)},y_3]_c&=\nu_2; & y_{554}&=0; & y_{556}&=\lambda_{665}; & &[y_{4456},y_{45}]_c=\lambda_{456}.
\end{align*}
If $\nu_2=0$, then we apply Lemma \ref{lem:cut1} for $(i,j)=(3,4)$ and $\toba$ projects onto $ \widetilde{\mathcal E}_{\bq'}$, $\bq'$ with two components of types $A_3$ and $B_3$. If $\nu_2\neq 0$, then $\lambda_{456}=0$, so we apply Lemma \ref{lem:cut1} for $(i,j)=(4,5)$ and $\toba$ projects onto $ \widetilde{\mathcal E}_{\bq'}$, $\bq'$ with two components of types $A_4$ and $A_2$.

 For $\xymatrix@R-8pt{  &  & & \overset{-1}{\circ} \ar  @{-}[d]^{\zeta}\ar@{-}[dr]^{\zeta}  & \\
\overset{\zeta}{\underset{\ }{\circ}} \ar  @{-}[r]^{\ztu} & \overset{\zeta}{\underset{\
}{\circ}} \ar  @{-}[r]^{\ztu} & \overset{\zeta}{\underset{\ }{\circ}} \ar  @{-}[r]^{\ztu}
& \overset{-1}{\underset{\ }{\circ}} \ar  @{-}[r]^{\zeta}  & \overset{\ztu}{\underset{\
}{\circ}}}$ the algebra $\toba$ is generated by $y_i$, $i\in\I_6$, with defining relations
\begin{align*}
y_{112}&=\lambda_{112}; & y_{221}&=\lambda_{221}; & y_{223}&=\lambda_{223}; & & y_{ij}=\lambda_{ij}, \, i<j, \, \widetilde{q}_{ij}=1;\\
y_{4}^2&=\mu_4; & y_{332}&=\lambda_{332}; & y_{334}&=0; & &[y_{346},y_4]_c=0;\\
y_{5}^2&=\mu_5; & y_{554}&=0; & y_{556}&=0; & &[y_{(35)},y_4]_c=\nu_4;
\end{align*}
\vsp
\begin{align*}
y_{(46)}-q_{56}\ztu [y_{46},y_5]_c-q_{45}(1-\zeta)y_5y_{46}=0.
\end{align*}
Here we apply Lemma \ref{lem:cut2} for $(i,j,k)=(6,4,5)$ and $\toba$ projects onto $ \widetilde{\mathcal E}_{\bq'}$, $\bq'$ with two components of types $A_5$ and $A_1$.

 For $\xymatrix@C-3pt@R-8pt{  &   & \overset{-1}{\circ} \ar  @{-}[d]^{\zeta} \ar
@{-}[dr]^{\zeta}  & \\
\overset{-1}{\underset{\ }{\circ}} \ar  @{-}[r]^{\ztu}  & \overset{\zeta}{\underset{\
}{\circ}} \ar  @{-}[r]^{\ztu}
& \overset{-1}{\underset{\ }{\circ}} \ar  @{-}[r]^{\zeta}
& \overset{-1}{\underset{\ }{\circ}}\ar  @{-}[r]^{\ztu}  & \overset{\zeta}{\underset{\
}{\circ}}}$ the algebra $\toba$ is generated by $y_i$, $i\in\I_6$, with defining relations
\begin{align*}
y_{221}&=0; & y_{223}&=0; & y_{554}&=0; & & y_{ij}=\lambda_{ij}, \, i<j, \, \widetilde{q}_{ij}=1;\\
y_{4}^2&=\mu_4; & y_{1}^2&=\mu_1; & [y_{(24)},&y_3]_c=0; & &[y_{236},y_3]_c=0;\\
y_{6}^2&=\mu_6; & y_{3}^2&=\mu_3; & [y_{(35)},&y_4]_c=0; & &[y_{546},y_4]_c=0;
\end{align*}
\vsp
\begin{align*}
y_{346}-q_{46}\ztu [y_{36},y_4]_c-q_{34}(1-\zeta)y_4y_{36}=0.
\end{align*}
Here we apply Lemma \ref{lem:cut2} for $(i,j,k)=(6,3,4)$ and $\toba$ projects onto $ \widetilde{\mathcal E}_{\bq'}$, $\bq'$ with two components of types $A_5$ and $A_1$.

 For $\xymatrix@C-3pt@R-8pt{  &   & \overset{-1}{\circ} \ar  @{-}[d]^{\zeta} \ar
@{-}[dr]^{\zeta}  & \\
\overset{-1}{\underset{\ }{\circ}} \ar  @{-}[r]^{\zeta}  & \overset{-1}{\underset{\
}{\circ}} \ar  @{-}[r]^{\ztu}
& \overset{-1}{\underset{\ }{\circ}} \ar  @{-}[r]^{\zeta}
& \overset{-1}{\underset{\ }{\circ}}\ar  @{-}[r]^{\ztu}  & \overset{\zeta}{\underset{\
}{\circ}}}$ the algebra $\toba$ is generated by $y_i$, $i\in\I_6$, with defining relations
\begin{align*}
y_{554}&=0; & y_{1}^2&=\mu_1; & [y_{(13)},&y_2]_c=\nu_2; & & y_{ij}=\lambda_{ij}, \, i<j, \, \widetilde{q}_{ij}=1;\\
y_{4}^2&=\mu_4; & y_{2}^2&=\mu_2; & [y_{(24)},&y_3]_c=\nu_3; &  &[y_{236},y_3]_c=\nu_3';\\
y_{6}^2&=\mu_6; & y_{3}^2&=\mu_3; & [y_{(35)},&y_4]_c=0; & &[y_{546},y_4]_c=0;
\end{align*}
\vsp
\begin{align*}
y_{346}-q_{46}\ztu [y_{36},y_4]_c-q_{34}(1-\zeta)y_4y_{36}=0.
\end{align*}
Here we apply Lemma \ref{lem:cut2} for $(i,j,k)=(3,4,6)$ and $\toba$ projects onto $\widetilde{\mathcal E}_{\bq'}$, $\bq'$ with two components of super type $A_3$.

 For $\xymatrix@C-3pt@R-8pt{  &   & \overset{-1}{\circ} \ar  @{-}[d]^{\zeta} \ar
@{-}[dr]^{\zeta}  & \\
\overset{\zeta}{\underset{\ }{\circ}} \ar  @{-}[r]^{\ztu}  & \overset{-1}{\underset{\
}{\circ}} \ar  @{-}[r]^{\zeta}
& \overset{\ztu}{\underset{\ }{\circ}} \ar  @{-}[r]^{\zeta}
& \overset{-1}{\underset{\ }{\circ}}\ar  @{-}[r]^{\ztu}  & \overset{\zeta}{\underset{\
}{\circ}}}$ the algebra $\toba$ is generated by $y_i$, $i\in\I_6$, with defining relations
\begin{align*}
y_{112}&=0; & y_{332}&=0; & y_{334}&=0; & & y_{ij}=\lambda_{ij}, \, i<j, \, \widetilde{q}_{ij}=1;\\
y_{4}^2&=\mu_4; & y_{336}&=0; & [y_{(13)},&y_2]_c=\nu_2; & &[y_{(35)},y_4]_c=\nu_4;\\
y_{6}^2&=\mu_6; & y_{554}&=0; & y_2^2&=\mu_2; & &[y_{546},y_4]_c=0;
\end{align*}
\vsp
\begin{align*}
y_{346}-q_{46}\ztu [y_{36},y_4]_c-q_{34}(1-\zeta)y_4y_{36}=0.
\end{align*}
Here we apply Lemma \ref{lem:cut2} for $(i,j,k)=(6,3,4)$ and $\toba$ projects onto $ \widetilde{\mathcal E}_{\bq'}$, $\bq'$ with two components of types $A_5$ and $A_1$.

 For $\xymatrix@C-3pt@R-8pt{  & &  & \overset{\zeta}{\circ} \ar  @{-}[d]^{\ztu}  & \\
\overset{\zeta}{\underset{\ }{\circ}} \ar  @{-}[r]^{\ztu} & \overset{\zeta}{\underset{\
}{\circ}} \ar  @{-}[r]^{\ztu}  & \overset{-1}{\underset{\ }{\circ}} \ar  @{-}[r]^{\zeta}
& \overset{-1}{\underset{\ }{\circ}} \ar  @{-}[r]^{\ztu}  & \overset{-1}{\underset{\
}{\circ}}}$ the algebra $\toba$ is generated by $y_i$, $i\in\I_6$, with defining relations
\begin{align*}
y_{112}&=\lambda_{112}; & y_{221}&=\lambda_{221}; & y_{223}&=0; & & y_{ij}=\lambda_{ij}, \, i<j, \, \widetilde{q}_{ij}=1;\\
y_{4}^2&=\mu_4; & y_{664}&=0; & [y_{(24)},&y_3]_c=0; & &[[y_{54},y_{546}]_c,y_4]_c=0;\\
y_{5}^2&=\mu_5; & y_{3}^2&=\mu_3; & [y_{(35)},&y_4]_c=\nu_4; & &[y_{346},y_4]_c=0.
\end{align*}
Here we apply Lemma \ref{lem:cut1} for $(i,j)=(4,6)$ and $\toba$ projects onto $\widetilde{\mathcal E}_{\bq'}$, $\bq'$ with two components of types $A_5$ and $A_1$.

 For $\xymatrix@R-8pt{  & &  & \overset{\zeta}{\circ} \ar  @{-}[d]^{\ztu}  & \\
\overset{\zeta}{\underset{\ }{\circ}} \ar  @{-}[r]^{\ztu} & \overset{-1}{\underset{\
}{\circ}} \ar  @{-}[r]^{\zeta}  & \overset{-1}{\underset{\ }{\circ}} \ar  @{-}[r]^{\ztu}
& \overset{\zeta}{\underset{\ }{\circ}} \ar  @{-}[r]^{\ztu}  & \overset{-1}{\underset{\
}{\circ}}}$ the algebra $\toba$ is generated by $y_i$, $i\in\I_6$, with defining relations
\begin{align*}
y_{3}^2&=\mu_3; & y_{112}&=0; & y_{443}&=0; & & y_{ij}=\lambda_{ij}, \, i<j, \, \widetilde{q}_{ij}=1;\\
y_{5}^2&=\mu_5; & y_{445}&=0; & y_{446}&=\lambda_{446}; & &[y_{(13)},y_2]_c=0;\\
& & y_{664}&=\lambda_{664}; & y_2^2&=\mu_2; & &[y_{(24)},y_3]_c=0.
\end{align*}
Here we apply Lemma \ref{lem:cut1} for $(i,j)=(3,4)$ and $\toba$ projects onto $\widetilde{\mathcal E}_{\bq'}$, $\bq'$ with two components of type $A_3$.

 For $\xymatrix@R-8pt{  & &  & \overset{\zeta}{\circ} \ar  @{-}[d]^{\ztu}  & \\
\overset{\zeta}{\underset{\ }{\circ}} \ar  @{-}[r]^{\ztu} & \overset{\zeta}{\underset{\
}{\circ}} \ar  @{-}[r]^{\ztu}  & \overset{-1}{\underset{\ }{\circ}} \ar  @{-}[r]^{\zeta}
& \overset{\ztu}{\underset{\ }{\circ}} \ar  @{-}[r]^{\zeta}  & \overset{-1}{\underset{\
}{\circ}}}$ the algebra $\toba$ is generated by $y_i$, $i\in\I_6$, with defining relations
\begin{align*}
y_{112}&=\lambda_{112}; & y_{221}&=\lambda_{221}; & y_{223}&=0; & & y_{ij}=\lambda_{ij}, \, i<j, \, \widetilde{q}_{ij}=1;\\
y_{3}^2&=\mu_3; & y_{443}&=0; & [y_{(24)},&y_3]_c=\nu_3; & &[[y_{346},y_4]_c,y_4]_c=0;\\
y_{5}^2&=\mu_5; & y_{445}&=0; & y_{664}&=0; & &[[y_{546},y_4]_c,y_4]_c=0.
\end{align*}
Here we apply Lemma \ref{lem:cut1} for $(i,j)=(4,6)$ and $\toba$ projects onto $\widetilde{\mathcal E}_{\bq'}$, $\bq'$ with two components of types $A_5$ and $A_1$.

 For $\xymatrix@R-8pt{  & &  & \overset{\zeta}{\circ} \ar  @{-}[d]^{\ztu}  & \\
\overset{-1}{\underset{\ }{\circ}} \ar  @{-}[r]^{\zeta} & \overset{-1}{\underset{\
}{\circ}} \ar  @{-}[r]^{\ztu}  & \overset{\zeta}{\underset{\ }{\circ}} \ar  @{-}[r]^{\ztu}
& \overset{\zeta}{\underset{\ }{\circ}} \ar  @{-}[r]^{\ztu}  & \overset{-1}{\underset{\
}{\circ}}}$ the algebra $\toba$ is generated by $y_i$, $i\in\I_6$, with defining relations
\begin{align*}
y_1^2&=\mu_1; & y_{332}&=0; & y_{334}&=\lambda_{334}; & & y_{ij}=\lambda_{ij}, \, i<j, \, \widetilde{q}_{ij}=1;\\
y_2^2&=\mu_2; & y_{443}&=\lambda_{443}; & y_{445}&=0; & &[y_{(13)},y_2]_c=0;\\
y_5^2&=\mu_5; & y_{446}&=\lambda_{446}; & y_{664}&=\lambda_{664}. & &
\end{align*}
Here we apply Lemma \ref{lem:cut1} for $(i,j)=(4,5)$ and $\toba$ projects onto $\widetilde{\mathcal E}_{\bq'}$, $\bq'$ with two components of types $A_5$ and $A_1$.

 For $\xymatrix@C-4pt@R-8pt{  & &  & \overset{\zeta}{\circ} \ar  @{-}[d]^{\ztu}  & \\
\overset{\zeta}{\underset{\ }{\circ}} \ar  @{-}[r]^{\ztu} & \overset{-1}{\underset{\
}{\circ}} \ar  @{-}[r]^{\zeta}  & \overset{-1}{\underset{\ }{\circ}} \ar  @{-}[r]^{\ztu}
& \overset{-1}{\underset{\ }{\circ}} \ar  @{-}[r]^{\zeta}  & \overset{-1}{\underset{\
}{\circ}}}$ the algebra $\toba$ is generated by $y_i$, $i\in\I_6$, with defining relations
\begin{align*}
y_{112}&=0; & y_{664}&=0; & [y_{(13)},&y_2]_c=0; & & y_{ij}=\lambda_{ij}, \, i<j, \, \widetilde{q}_{ij}=1;\\
y_{4}^2&=\mu_4; & y_{2}^2&=\mu_2; & [y_{(24)},&y_3]_c=\nu_3; & &[[y_{34},y_{346}]_c,y_4]_c=0;\\
y_{5}^2&=\mu_5; & y_{3}^2&=\mu_3; & [y_{(35)},&y_4]_c=\nu_4; & &[y_{546},y_4]_c=0.
\end{align*}
Here we apply Lemma \ref{lem:cut1} for $(i,j)=(4,6)$ and $\toba$ projects onto $\widetilde{\mathcal E}_{\bq'}$, $\bq'$ with two components of types $A_5$ and $A_1$.

 For $\xymatrix@R-8pt{  & &  & \overset{\zeta}{\circ} \ar  @{-}[d]^{\ztu}  & \\
\overset{-1}{\underset{\ }{\circ}} \ar  @{-}[r]^{\ztu} & \overset{\zeta}{\underset{\
}{\circ}} \ar  @{-}[r]^{\ztu}  & \overset{\zeta}{\underset{\ }{\circ}} \ar  @{-}[r]^{\ztu}
& \overset{\zeta}{\underset{\ }{\circ}} \ar  @{-}[r]^{\ztu}  & \overset{-1}{\underset{\
}{\circ}}}$ the algebra $\toba$ is generated by $y_i$, $i\in\I_6$, with defining relations
\begin{align*}
y_{1}^2&=\mu_1; & y_{221}&=0; & y_{223}&=\lambda_{223}; & & y_{ij}=\lambda_{ij}, \, i<j, \, \widetilde{q}_{ij}=1;\\
y_{5}^2&=\mu_5; & y_{332}&=\lambda_{332}; & y_{334}&=\lambda_{334}; &  &y_{443}=\lambda_{443}; \\
& & y_{445}&=0; & y_{446}&=\lambda_{446}; & &y_{664}=\lambda_{664}.
\end{align*}
Here we apply Lemma \ref{lem:cut1} for $(i,j)=(4,5)$ and $\toba$ projects onto $\widetilde{\mathcal E}_{\bq'}$, $\bq'$ with two components of types $A_5$ and $A_1$.

 For $\xymatrix@C-4pt@R-8pt{  & &  & \overset{\zeta}{\circ} \ar  @{-}[d]^{\ztu}  & \\
\overset{-1}{\underset{\ }{\circ}} \ar  @{-}[r]^{\zeta} & \overset{-1}{\underset{\
}{\circ}} \ar  @{-}[r]^{\ztu}  & \overset{\zeta}{\underset{\ }{\circ}} \ar  @{-}[r]^{\ztu}
& \overset{-1}{\underset{\ }{\circ}} \ar  @{-}[r]^{\zeta}  & \overset{-1}{\underset{\
}{\circ}}}$ the algebra $\toba$ is generated by $y_i$, $i\in\I_6$, with defining relations
\begin{align*}
y_{1}^2&=\mu_1; & y_{332}&=0; & y_{334}&=0; & & y_{ij}=\lambda_{ij}, \, i<j, \, \widetilde{q}_{ij}=1;\\
y_{2}^2&=\mu_2; & y_{664}&=0; & [y_{(13)},&y_2]_c=0; & &[y_{564},y_4]_c=0;\\
& & y_{4}^2&=\mu_4; & y_{5}^2&=\mu_5; & &[y_{(35)},y_4]_c=0.
\end{align*}
This algebra is nonzero by \cite[Lemma 5.16]{AAGMV}.

 For $\xymatrix@C-4pt@R-8pt{  &   & \overset{\zeta}{\circ} \ar  @{-}[d]^{\ztu}  & \\
\overset{\ztu}{\underset{\ }{\circ}} \ar  @{-}[r]^{\zeta}  & \overset{-1}{\underset{\
}{\circ}} \ar  @{-}[r]^{\ztu}  & \overset{\zeta}{\underset{\ }{\circ}} \ar  @{-}[r]^{\ztu}
& \overset{\zeta}{\underset{\ }{\circ}}\ar  @{-}[r]^{\ztu}  & \overset{\zeta}{\underset{\
}{\circ}}}$ the algebra $\toba$ is generated by $y_i$, $i\in\I_6$, with defining relations
\begin{align*}
y_{2}^2&=\mu_2; & y_{112}&=0; & y_{332}&=0; & & y_{ij}=\lambda_{ij}, \, i<j, \, \widetilde{q}_{ij}=1;\\
y_{554}&=\lambda_{554}; & y_{334}&=\lambda_{334}; & y_{336}&=\lambda_{336}; & &[y_{(13)},y_2]_c=\nu_2;\\
& & y_{663}&=\lambda_{663}; & y_{443}&=\lambda_{443}; & &y_{445}=\lambda_{445}.
\end{align*}
If $\nu_2=0$, then  we apply Lemma \ref{lem:cut1} for $(i,j)=(2,3)$ and $\toba$ projects onto $ \widetilde{\mathcal E}_{\bq'}$, $\bq'$ with two components of types $A_2$ and $A_4$. If $\nu_2\neq 0$, then $\lambda_{336}=\lambda_{663}=0$, so we apply Lemma \ref{lem:cut1} for $(i,j)=(3,6)$ and $\toba$ projects onto $\widetilde{\mathcal E}_{\bq'}$, $\bq'$ with two components of types $A_5$ and $A_1$.

 For $\xymatrix@C-3pt@R-8pt{  & &  & \overset{\zeta}{\circ} \ar  @{-}[d]^{\ztu}  & \\
\overset{-1}{\underset{\ }{\circ}} \ar  @{-}[r]^{\ztu} & \overset{\zeta}{\underset{\
}{\circ}} \ar  @{-}[r]^{\ztu}  & \overset{\zeta}{\underset{\ }{\circ}} \ar  @{-}[r]^{\ztu}
& \overset{-1}{\underset{\ }{\circ}} \ar  @{-}[r]^{\zeta}  & \overset{-1}{\underset{\
}{\circ}}}$ the algebra $\toba$ is generated by $y_i$, $i\in\I_6$, with defining relations
\begin{align*}
y_{1}^2&=\mu_1; & y_{223}&=\lambda_{223}; & y_{334}&=0; & & y_{ij}=\lambda_{ij}, \, i<j, \, \widetilde{q}_{ij}=1;\\
y_{4}^2&=\mu_4; & y_{332}&=\lambda_{332}; & [y_{(35)},&y_4]_c=0; & &[[[y_{2346},y_4]_c,y_3]_c,y_4]_c=0;\\
y_{5}^2&=\mu_5; & y_{664}&=0; & y_{221}&=0; & &[y_{546},y_4]_c=0.
\end{align*}
Here we apply Lemma \ref{lem:cut1} for $(i,j)=(3,4)$ and $\toba$ projects onto $\widetilde{\mathcal E}_{\bq'}$, $\bq'$ with two components of type $A_3$.

 For $\xymatrix@C-3pt@R-8pt{  &   & \overset{\zeta}{\circ} \ar  @{-}[d]^{\ztu}  & \\
\overset{-1}{\underset{\ }{\circ}} \ar  @{-}[r]^{\zeta}  & \overset{\ztu}{\underset{\
}{\circ}} \ar  @{-}[r]^{\zeta}
& \overset{-1}{\underset{\ }{\circ}} \ar  @{-}[r]^{\ztu}
& \overset{\zeta}{\underset{\ }{\circ}}\ar  @{-}[r]^{\ztu}  & \overset{\zeta}{\underset{\
}{\circ}}}$ the algebra $\toba$ is generated by $y_i$, $i\in\I_6$, with defining relations
\begin{align*}
y_{221}&=0; & y_{223}&=0; & y_{443}&=0; & & y_{ij}=\lambda_{ij}, \, i<j, \, \widetilde{q}_{ij}=1;\\
y_{1}^2&=\mu_1; & y_{445}&=\lambda_{445}; & [y_{(24)},&y_3]_c=0; & &[[[y_{5436},y_3]_c,y_4]_c,y_3]_c=0;\\
y_{3}^2&=\mu_3; & y_{554}&=\lambda_{554}; & y_{663}&=0; & &[y_{236},y_3]_c=\nu_3'.
\end{align*}
Here we apply Lemma \ref{lem:cut1} for $(i,j)=(3,4)$ and $\toba$ projects onto $\widetilde{\mathcal E}_{\bq'}$, $\bq'$ with two components of types $A_4$ and $A_2$.

 For $\xymatrix@C-3pt@R-8pt{  &   & \overset{\zeta}{\circ} \ar  @{-}[d]^{\ztu}   & \\
\overset{-1}{\underset{\ }{\circ}} \ar  @{-}[r]^{\ztu}  & \overset{-1}{\underset{\
}{\circ}} \ar  @{-}[r]^{\zeta}
& \overset{-1}{\underset{\ }{\circ}} \ar  @{-}[r]^{\ztu}
& \overset{\zeta}{\underset{\ }{\circ}}\ar  @{-}[r]^{\ztu}  & \overset{\zeta}{\underset{\
}{\circ}}}$ the algebra $\toba$ is generated by $y_i$, $i\in\I_6$, with defining relations
\begin{align*}
y_{1}^2&=\mu_1; & y_{445}&=\lambda_{445}; & y_{554}&=\lambda_{554}; & & y_{ij}=\lambda_{ij}, \, i<j, \, \widetilde{q}_{ij}=1;\\
y_{2}^2&=\mu_2; & y_{663}&=0; & [y_{(13)},&y_2]_c=\nu_2; & &[[[y_{5436},y_3]_c,y_4]_c,y_3]_c=0;\\
y_{3}^2&=\mu_3; & y_{443}&=0; & [y_{(24)},&y_3]_c=0; & &[y_{236},y_3]_c=0.
\end{align*}
Here we apply Lemma \ref{lem:cut1} for $(i,j)=(3,4)$ and $\toba$ projects onto $\widetilde{\mathcal E}_{\bq'}$, $\bq'$ with two components of types $A_4$ and $A_2$.

\subsubsection{Type $\g(8,6)$, $\zeta \in \G'_3$}\label{subsec:type-g(8,6)}
The Weyl groupoid has eight objects.

 For $\xymatrix{\overset{\zeta}{\underset{\ }{\circ}}\ar  @{-}[r]^{\ztu}  &
\overset{\zeta}{\underset{\ }{\circ}} \ar  @{-}[r]^{\ztu}  & \overset{\zeta}{\underset{\
}{\circ}}
\ar  @{-}[r]^{\ztu}  & \overset{\zeta}{\underset{\ }{\circ}}
\ar  @{-}[r]^{\ztu}  & \overset{-1}{\underset{\ }{\circ}}  \ar  @{-}[r]^{\ztu}  &
\overset{\zeta}{\underset{\ }{\circ}}  \ar  @{-}[r]^{\ztu}  & \overset{\zeta}{\underset{\
}{\circ}}
}$, $\toba$ is generated by $y_i$, $i\in\I_7$, with defining relations
\begin{align*}
y_{5}^2&=\mu_5; & y_{112}&=\lambda_{112}; & y_{221}&=\lambda_{221}; & y_{223}&=\lambda_{223}; \\
y_{665}&=0; & y_{667}&=\lambda_{667}; & y_{332}&=\lambda_{332}; & y_{334}&=\lambda_{334};\\
& & y_{443}&=\lambda_{443}; & y_{445}&=\lambda_{445}; & y_{776}&=\lambda_{776};
\end{align*}
\vsp
\begin{align*}
y_{ij}&=\lambda_{ij}, \, i<j, \, \widetilde{q}_{ij}=1; & [[[y_{(36)},y_5]_c,y_4]_c,y_5]_c&=[[[y_{7654},y_5]_c,y_6]_c,y_5]_c=0.
\end{align*}

Here we apply Lemma \ref{lem:cut1} for $(i,j)=(5,6)$ and $\toba$ projects onto $\widetilde{\mathcal E}_{\bq'}$, $\bq'$ with two components of types $A_5$ and $A_2$.

 For $\xymatrix@R-8pt{  &   & & \overset{-1}{\circ} \ar  @{-}[d]^{\zeta}\ar  @{-}[dr]^{\zeta}
& &\\
\overset{\zeta}{\underset{\ }{\circ}} \ar  @{-}[r]^{\ztu}  & \overset{\zeta}{\underset{\
}{\circ}} \ar  @{-}[r]^{\ztu}  & \overset{\zeta}{\underset{\ }{\circ}} \ar  @{-}[r]^{\ztu}
 & \overset{-1}{\underset{\ }{\circ}}
\ar  @{-}[r]^{\zeta}  & \overset{-1}{\underset{\ }{\circ}} \ar  @{-}[r]^{\ztu}  &
\overset{\zeta}{\underset{\ }{\circ}}}$ the algebra $\toba$ is generated by $y_i$, $i\in\I_7$, with defining relations
\begin{align*}
y_{4}^2&=\mu_4; & y_{112}&=\lambda_{112}; & y_{221}&=\lambda_{221}; & y_{ij}&=\lambda_{ij}, \, i<j, \, \widetilde{q}_{ij}=1;\\
y_{5}^2&=\mu_5; & y_{223}&=\lambda_{223}; & y_{332}&=\lambda_{332}; & [y_{(35)},&y_4]_c=[y_{347},y_4]_c=0;\\
y_{7}^2&=\mu_7; & y_{334}&=0; & y_{665}&=0; & [y_{(46)},&y_5]_c=[y_{657},y_5]_c=0
\end{align*}
\vsp
\begin{align*}
y_{457}-q_{57}\ztu [y_{47},y_5]_c-q_{45}(1-\zeta)y_5y_{47}=0.
\end{align*}
Here we apply Lemma \ref{lem:cut2} for $(i,j,k)=(7,4,5)$ and $\toba$ projects onto $ \widetilde{\mathcal E}_{\bq'}$, $\bq'$ with two components of types $A_6$ and $A_1$.

 For $\xymatrix@R-8pt{  &  & & & \overset{\zeta}{\circ} \ar  @{-}[d]^{\ztu}  & \\
\overset{\zeta}{\underset{\ }{\circ}} \ar  @{-}[r]^{\ztu}  & \overset{\zeta}{\underset{\
}{\circ}} \ar  @{-}[r]^{\ztu}  & \overset{\zeta}{\underset{\ }{\circ}} \ar  @{-}[r]^{\ztu}
 & \overset{-1}{\underset{\ }{\circ}}
\ar  @{-}[r]^{\zeta}  & \overset{-1}{\underset{\ }{\circ}} \ar  @{-}[r]^{\ztu}  &
\overset{-1}{\underset{\ }{\circ}}}$ the algebra $\toba$ is generated by $y_i$, $i\in\I_7$, with defining relations
\begin{align*}
y_{4}^2&=\mu_4; & y_{112}&=\lambda_{112}; & y_{221}&=\lambda_{221}; &  y_{ij}&=\lambda_{ij}, \, i<j, \, \widetilde{q}_{ij}=1;\\
y_{5}^2&=\mu_5; & y_{223}&=\lambda_{223}; & y_{332}&=\lambda_{332}; & [y_{(35)},&y_4]_c=[[y_{65},y_{657}]_c,y_5]_c=0;\\
y_{6}^2&=\mu_6; & y_{775}&=0; & [y_{(46)},&y_5]_c=\mu_5; & y_{334}&=[y_{457},y_5]_c=0.
\end{align*}
Here we apply Lemma \ref{lem:cut1} for $(i,j)=(5,7)$ and $\toba$ projects onto $\widetilde{\mathcal E}_{\bq'}$, $\bq'$ with two components of types $A_6$ and $A_1$.

 For $\xymatrix@R-8pt{  &  & & & \overset{\zeta}{\circ} \ar  @{-}[d]^{\ztu}  & \\
\overset{\zeta}{\underset{\ }{\circ}} \ar  @{-}[r]^{\ztu}  & \overset{\zeta}{\underset{\
}{\circ}} \ar  @{-}[r]^{\ztu}  & \overset{\zeta}{\underset{\ }{\circ}} \ar  @{-}[r]^{\ztu}
 & \overset{\zeta}{\underset{\ }{\circ}}
\ar  @{-}[r]^{\ztu}  & \overset{\zeta}{\underset{\ }{\circ}} \ar  @{-}[r]^{\ztu}  &
\overset{-1}{\underset{\ }{\circ}}}$ the algebra $\toba$ is generated by $y_i$, $i\in\I_7$, with defining relations
\begin{align*}
y_{112}&=\lambda_{112}; & y_{221}&=\lambda_{221}; & y_{223}&=\lambda_{223}; & & y_{ij}=\lambda_{ij}, \, i<j, \, \widetilde{q}_{ij}=1;\\
y_{557}&=\lambda_{557}; & y_{332}&=\lambda_{332}; & y_{334}&=\lambda_{334}; &  &y_{443}=\lambda_{443};\\
y_{775}&=\lambda_{775}; & y_{445}&=\lambda_{445}; & y_{554}&=\lambda_{554}; & &y_{556}=0; \qquad y_{6}^2=\mu_2.
\end{align*}
Herewe apply Lemma \ref{lem:cut1} for $(i,j)=(5,6)$ and $\toba$ projects onto $\widetilde{\mathcal E}_{\bq'}$, $\bq'$ with two components of types $A_6$ and $A_1$.

 For $\xymatrix@R-8pt{  &   & & \overset{\zeta}{\circ} \ar  @{-}[d]^{\ztu}  & &\\
\overset{\zeta}{\underset{\ }{\circ}} \ar  @{-}[r]^{\ztu}  & \overset{\zeta}{\underset{\
}{\circ}} \ar  @{-}[r]^{\ztu}  & \overset{-1}{\underset{\ }{\circ}} \ar  @{-}[r]^{\zeta}
& \overset{-1}{\underset{\ }{\circ}}
\ar  @{-}[r]^{\ztu}  & \overset{\zeta}{\underset{\ }{\circ}} \ar  @{-}[r]^{\ztu}  &
\overset{\zeta}{\underset{\ }{\circ}}}$ the algebra $\toba$ is generated by $y_i$, $i\in\I_7$, with defining relations
\begin{align*}
y_{112}&=\lambda_{112}; & y_{221}&=\lambda_{221}; & [y_{347},&y_4]_c=0; & & y_{ij}=\lambda_{ij}, \, i<j, \, \widetilde{q}_{ij}=1;\\
y_{3}^2&=\mu_3; & y_{556}&=\lambda_{556}; & [y_{(24)},&y_3]_c=0; &  &[[[y_{6547},y_4]_c,y_5]_c,y_4]_c=0;\\
y_{4}^2&=\mu_4; & y_{665}&=\lambda_{665}; & [y_{(35)},&y_4]_c=0; & &y_{223}=y_{774}=y_{554}=0.
\end{align*}
Here we apply Lemma \ref{lem:cut1} for $(i,j)=(4,7)$ and $\toba$ projects onto $\widetilde{\mathcal E}_{\bq'}$, $\bq'$ with two components of types $A_6$ and $A_1$.

For $\xymatrix@R-8pt{  &   & & \overset{\zeta}{\circ} \ar  @{-}[d]^{\ztu}  & & \\
\overset{\zeta}{\underset{\ }{\circ}} \ar  @{-}[r]^{\ztu}  & \overset{-1}{\underset{\
}{\circ}} \ar  @{-}[r]^{\zeta}  & \overset{-1}{\underset{\ }{\circ}} \ar  @{-}[r]^{\ztu}
& \overset{\zeta}{\underset{\ }{\circ}}
\ar  @{-}[r]^{\ztu}  & \overset{\zeta}{\underset{\ }{\circ}} \ar  @{-}[r]^{\ztu}  &
\overset{\zeta}{\underset{\ }{\circ}}}$ the algebra $\toba$ is generated by $y_i$, $i\in\I_7$, with defining relations
\begin{align*}
y_{112}&=y_{443}=0; & y_{665}&=\lambda_{665}; & y_{445}&=\lambda_{445}; & & y_{ij}=\lambda_{ij}, \, i<j, \, \widetilde{q}_{ij}=1;\\
y_{2}^2&=\mu_2; & y_{447}&=\lambda_{447}; & y_{774}&=\lambda_{774}; & &[y_{(13)},y_2]_c=0;\\
y_{3}^2&=\mu_3; & y_{554}&=\lambda_{554}; & y_{556}&=\lambda_{556}; & &[y_{(24)},y_3]_c=0.
\end{align*}
Here we apply Lemma \ref{lem:cut1} for $(i,j)=(3,4)$ and $\toba$ projects onto $\widetilde{\mathcal E}_{\bq'}$, $\bq'$ with two components of types $A_3$ and $A_4$.

 For $\xymatrix@R-8pt{  &   & & \overset{\zeta}{\circ} \ar  @{-}[d]^{\ztu}    & &\\
\overset{-1}{\underset{\ }{\circ}} \ar  @{-}[r]^{\zeta}  & \overset{-1}{\underset{\
}{\circ}} \ar  @{-}[r]^{\ztu}  & \overset{\zeta}{\underset{\ }{\circ}} \ar  @{-}[r]^{\ztu}
 & \overset{\zeta}{\underset{\ }{\circ}}
\ar  @{-}[r]^{\ztu}  & \overset{\zeta}{\underset{\ }{\circ}} \ar  @{-}[r]^{\ztu}  &
\overset{\zeta}{\underset{\ }{\circ}}}$ the algebra $\toba$ is generated by $y_i$, $i\in\I_7$, with defining relations
\begin{align*}
y_{1}^2&=\mu_1; & y_{334}&=\lambda_{334}; & y_{443}&=\lambda_{443}; & & y_{ij}=\lambda_{ij}, \, i<j, \, \widetilde{q}_{ij}=1;\\
y_{2}^2&=\mu_2; & y_{445}&=\lambda_{445}; & y_{447}&=\lambda_{447}; & &[y_{(13)},y_2]_c=0;\\
y_{665}&=\lambda_{665} & y_{774}&=\lambda_{774}; & y_{554}&=\lambda_{554}; & & y_{556}=\lambda_{556}; \qquad y_{332}=0.
\end{align*}
Here we apply Lemma \ref{lem:cut1} for $(i,j)=(2,3)$ and $\toba$ projects onto $\widetilde{\mathcal E}_{\bq'}$, $\bq'$ with two components of types $A_2$ and $D_5$.

 For $\xymatrix@R-8pt{  &   & & \overset{\zeta}{\circ} \ar  @{-}[d]^{\ztu}    & &\\
\overset{-1}{\underset{\ }{\circ}} \ar  @{-}[r]^{\ztu}  & \overset{\zeta}{\underset{\
}{\circ}} \ar  @{-}[r]^{\ztu}  & \overset{\zeta}{\underset{\ }{\circ}} \ar  @{-}[r]^{\ztu}
 & \overset{\zeta}{\underset{\ }{\circ}}
\ar  @{-}[r]^{\ztu}  & \overset{\zeta}{\underset{\ }{\circ}} \ar  @{-}[r]^{\ztu}  &
\overset{\zeta}{\underset{\ }{\circ}}}$ the algebra $\toba$ is generated by $y_i$, $i\in\I_7$, with defining relations
\begin{align*}
y_{221}&=0; & y_{223}&=\lambda_{223}; & y_{332}&=\lambda_{332}; & & y_{ij}=\lambda_{ij}, \, i<j, \, \widetilde{q}_{ij}=1;\\
y_{556}&=\lambda_{556}; & y_{334}&=\lambda_{334}; & y_{443}&=\lambda_{443}; & &y_{445}=\lambda_{445};\\
y_{665}&=\lambda_{665}; & y_{447}&=\lambda_{447}; & y_{774}&=\lambda_{774}; & &y_{554}=\lambda_{554}; \qquad y_{1}^2=\mu_1.
\end{align*}
Here we apply Lemma \ref{lem:cut1} for $(i,j)=(1,2)$ and $\toba$ projects onto $\widetilde{\mathcal E}_{\bq'}$, $\bq'$ with two components of types $A_1$ and $E_6$.

\subsection{Super modular type, characteristic
5}\label{sec:by-diagram-super-modular-char5}
\subsubsection{Type $\brj(2; 5)$, $\zeta \in \G'_5$}\label{subsec:type-brj(2;5)}
The Weyl groupoid has two objects.

For $\xymatrix{\overset{\zeta}{\underset{\ }{\circ}} \ar  @{-}[r]^{\zeta^2} &
\overset{-1}{\underset{\ }{\circ}}}$, $\toba$ is generated by $y_1, y_2$ with defining relations
\begin{align*}
[y_{1112},y_{112}]_c&=\lambda_1;  &
 y_2^2&=\lambda_2;  & y_{11112}&=0; &  [[[y_{112},y_{12}]_c,y_{12}]_c,y_{12}]_c&=0,
\end{align*}
with $\lambda_1\neq0$ only if $\bq=\left(\begin{smallmatrix}
\zeta & 1 \\\zeta^2 &-1
\end{smallmatrix}\right)$. See \texttt{brj25.log}, resp. \texttt{brj252.log}, for the deformation of the first, resp. last, relation.
If $\lambda_1=0$, then this algebra is nonzero by \cite[Lemma 5.16]{AAGMV}. If $\lambda_1,\lambda_2\neq0$, then this algebra is nonzero by
 \texttt{brj25b.log}.

For $\xymatrix{\overset{-\zeta^3}{\underset{\ }{\circ}} \ar  @{-}[r]^{\zeta^3} &
\overset{-1}{\underset{\ }{\circ}}}$, $\toba$ is generated by $y_1, y_2$ with defining relations
\begin{align*}
 y_2^2&=\lambda_2;  & y_{111112}&=0; & [y_1,&[y_{112},y_{12}]_c]_c+q_{12}y_{112}^2=0.
\end{align*}
This algebra is nonzero by \cite[Lemma 5.16]{AAGMV}.

\subsubsection{Type $\el(5;5)$, $\zeta \in \G'_{5}$}\label{subsec:type-el(5;5)}
The Weyl groupoid has seven objects.

 For $\xymatrix@C-4pt{\overset{\zeta^2}{\underset{\ }{\circ}} \ar  @{-}[r]^{\ztu^{\, 2}} &
\overset{\zeta^{2}}{\underset{\ }{\circ}} \ar  @{-}[r]^{\ztu^{\, 2}}  &
\overset{-1}{\underset{\ }{\circ}} \ar  @{-}[r]^{\ztu}  & \overset{\zeta}{\underset{\
}{\circ}} \ar  @{-}[r]^{\ztu^{\, 2}}  & \overset{\zeta^2}{\underset{\ }{\circ}}}$
the algebra $\toba$ is generated by $y_i$, $i\in\I_5$, with defining relations
\begin{align*}
y_{112}&=0; & y_{221}&=0; & y_{223}&=0; & & y_{ij}=\lambda_{ij}, \, i<j, \, \widetilde{q}_{ij}=1;\\
y_{554}&=\lambda_{554}; & y_{443}&=0; & y_{4445}&=\lambda_{4445}; &  &[[[y_{(14)},y_3]_c,y_2]_c,y_3]_c=0;
\end{align*}
\vsp
\begin{align*}
y_{3}^2&=\mu_3; & [[y_{5432},y_4]_c,y_3]_c-&q_{43}(\zeta^2-\zeta)[[y_{5432},y_3]_c,y_4]_c=0.
\end{align*}
Here we apply Lemma \ref{lem:cut1} for $(i,j)=(3,4)$ and $\toba$ projects onto $\widetilde{\mathcal E}_{\bq'}$, $\bq'$ with two components of types $A_3$ and $B_2$.

 For $\xymatrix@C-4pt{\overset{\zeta^2}{\underset{\ }{\circ}} \ar  @{-}[r]^{\ztu^{\, 2}} &
\overset{\zeta^{2}}{\underset{\ }{\circ}} \ar  @{-}[r]^{\ztu^{\, 2}}  &
\overset{\zeta^{2}}{\underset{\ }{\circ}} \ar  @{-}[r]^{\ztu^{\, 2}}  &
\overset{-1}{\underset{\ }{\circ}} \ar  @{-}[r]^{\zeta}  & \overset{-1}{\underset{\
}{\circ}}}$
the algebra $\toba$ is generated by $y_i$, $i\in\I_5$, with defining relations
\begin{align*}
y_{112}&=0; & y_{4}^2&=\mu_4; & y_{221}&=0; & y_{223}&=0; & & y_{ij}=\lambda_{ij}, \, i<j, \, \widetilde{q}_{ij}=1;\\
& & y_{5}^2&=\mu_5; & y_{332}&=0; & y_{334}&=0; &  &[[y_{54},y_{543}]_c,y_4]_c=0.
\end{align*}
This algebra is nonzero by \cite[Lemma 5.16]{AAGMV}.

 For $\xymatrix@C-4pt{\overset{\zeta^2}{\underset{\ }{\circ}} \ar  @{-}[r]^{\ztu^{\, 2}} &
\overset{\zeta^{2}}{\underset{\ }{\circ}} \ar  @{-}[r]^{\ztu^{\, 2}}  &
\overset{\zeta^{2}}{\underset{\ }{\circ}} \ar  @{-}[r]^{\ztu^{\, 2}}  &
\overset{\zeta}{\underset{\ }{\circ}} \ar  @{-}[r]^{\ztu}  & \overset{-1}{\underset{\
}{\circ}}}$ the algebra $\toba$ is generated by $y_i$, $i\in\I_5$, with defining relations
\begin{align*}
y_{5}^2&=\mu_5; & y_{112}&=0; & y_{221}&=0; & y_{223}&=0; & & y_{ij}=\lambda_{ij}, \, i<j, \, \widetilde{q}_{ij}=1;\\
& & y_{445}&=0; & y_{332}&=0; & y_{334}&=\lambda_{334}; & &y_{4443}=\lambda_{4443}.
\end{align*}
Here we apply Lemma \ref{lem:cut1} for $(i,j)=(4,5)$ and $\toba$ projects onto $\widetilde{\mathcal E}_{\bq'}$, $\bq'$ with two components of types $A_1$ and $B_4$.

 For $\xymatrix@R-8pt{  &   \overset{\zeta^2}{\circ} \ar  @{-}[d]^{\ztu^{\, 2}} & & \\
\overset{-1}{\underset{\ }{\circ}} \ar  @{-}[r]^{\zeta^2}  & \overset{-1}{\underset{\
}{\circ}} \ar  @{-}[r]^{\ztu^{\, 2}}  & \overset{\zeta^2}{\underset{\ }{\circ}} \ar
@{-}[r]^{\ztu^{\, 2}}  & \overset{\zeta^2}{\underset{\ }{\circ}}}$
the algebra $\toba$ is generated by $y_i$, $i\in\I_5$, with defining relations
\begin{align*}
y_{332}&=0; & y_{1}^2&=\mu_1; & y_{334}&=0; & y_{443}&=0; & & y_{ij}=\lambda_{ij}, \, i<j, \, \widetilde{q}_{ij}=1;\\
& & y_{2}^2&=\mu_2; & y_{552}&=0; & [y_{(13)},&y_2]_c=0; & &[y_{125},y_2]_c=0.
\end{align*}
This algebra is nonzero by \cite[Lemma 5.16]{AAGMV}.

For $\xymatrix@R-8pt{  &   \overset{\zeta^2}{\circ} \ar  @{-}[d]^{\ztu^{\, 2}} & & \\
\overset{-1}{\underset{\ }{\circ}} \ar  @{-}[r]^{\ztu^{\, 2}}  &
\overset{\zeta^2}{\underset{\ }{\circ}} \ar  @{-}[r]^{\ztu^{\, 2}}  &
\overset{\zeta^2}{\underset{\ }{\circ}} \ar  @{-}[r]^{\ztu^{\, 2}}  &
\overset{\zeta^2}{\underset{\ }{\circ}}}$
the algebra $\toba$ is generated by $y_i$, $i\in\I_5$, with defining relations
\begin{align*}
y_{1}^2&=\mu_1; & y_{ij}&=\lambda_{ij}, \, i<j, \, \widetilde{q}_{ij}=1; & y_{iij}&=0,  \, \widetilde{q}_{ij}\neq 1, \, i\neq 1.
\end{align*}
This algebra is nonzero by \cite[Lemma 5.16]{AAGMV}.

 For $\xymatrix@R-8pt{  &   \overset{-1}{\circ} \ar  @{-}[d]_{\zeta^{2}} \ar  @{-}[dr]^{\zeta}
& & \\
\overset{\zeta^2}{\underset{\ }{\circ}} \ar  @{-}[r]^{\ztu^{\, 2}}  &
\overset{-1}{\underset{\ }{\circ}} \ar  @{-}[r]^{\zeta^{2}}  & \overset{-1}{\underset{\
}{\circ}} \ar  @{-}[r]^{\ztu^{\, 2}}  & \overset{\zeta^2}{\underset{\ }{\circ}}}$
the algebra $\toba$ is generated by $y_i$, $i\in\I_5$, with defining relations
\begin{align*}
y_{112}&=0; & y_{443}&=0; & [y_{(13)},&y_2]_c=0; & & y_{ij}=\lambda_{ij}, \, i<j, \, \widetilde{q}_{ij}=1;\\
[y_{125},&y_2]_c=0; & y_{2}^2&=\mu_2; & [y_{(24)},&y_3]_c=0; & &[[y_{53},y_{534}]_c,y_3]_c=0;
\end{align*}
\vsp
\begin{align*}
y_{3}^2&=\mu_3; & y_{5}^2&=\mu_5; & y_{235}-\frac{q_{35}}{\zeta^2+\zeta}[y_{25},y_3]_c-q_{23}(1-\zeta)y_3y_{25}=0.
\end{align*}
This algebra is nonzero by \cite[Lemma 5.16]{AAGMV}.

 For $\xymatrix@R-8pt{& &  \overset{\zeta}{\circ} \ar  @{-}[d]_{\ztu} \ar  @{-}[dr]^{\ztu} &
\\
\overset{\zeta^{2}}{\underset{\ }{\circ}}  \ar  @{-}[r]^{\ztu^{\, 2}} &
\overset{\zeta^{2}}{\underset{\ }{\circ}} \ar  @{-}[r]^{\ztu^{\, 2}}  &
\overset{-1}{\underset{\ }{\circ}} \ar  @{-}[r]^{\zeta^{2}}  & \overset{-1}{\underset{\
}{\circ}} }$
the algebra $\toba$ is generated by $y_i$, $i\in\I_5$, with defining relations
\begin{align*}
y_{112}&=0; & y_{221}&=0; & y_{223}&=0; & & y_{ij}=\lambda_{ij}, \, i<j, \, \widetilde{q}_{ij}=1;\\
y_{554}&=0; & y_{553}&=0; & [y_{(24)},&y_3]_c=0; &  &[[[y_{1235},y_3]_c,y_2]_c,y_3]_c=0;
\end{align*}
\vsp
\begin{align*}
y_{3}^2&=\mu_3; & y_4^2&=\mu_4; & y_{(35)}-q_{45}\zeta[y_{35},y_4]_c-q_{34}(1-\ztu)y_4y_{35}=0.
\end{align*}
This algebra is nonzero by \cite[Lemma 5.16]{AAGMV}.

\subsection{Unidentified type}\label{sec:by-diagram-Unidentified}

\subsubsection{Type $\bgl(4,\alpha)$, $q\neq \pm 1$}\label{subsec:type-bgl(4,alpha)}
The Weyl groupoid has five objects.
Two of them are obtained with $-q^{-1}$
instead of $q$; hence we study the other three.
We set $N=\ord q$, $M=\ord -q^{-1}$.

For $\xymatrix@C-4pt{\overset{q}{\underset{\ }{\circ}}\ar  @{-}[r]^{q ^{-1}}  &
\overset{q}{\underset{\ }{\circ}}
\ar  @{-}[r]^{q^{-1}}  & \overset{-1}{\underset{\ }{\circ}}
\ar  @{-}[r]^{-q}  & \overset{-q^{-1}}{\underset{\ }{\circ}}}$, the algebra $\toba$ is generated by $y_i$, $i\in\I_4$, with defining relations
\begin{align*}
 y_3^2&=\lambda_3; & y_{13}&=0; & y_{14}&=0; & y_{24}&=0; \\
y_{112}&=\lambda_1; & y_{221}&=\lambda_2; & y_{223}&=\lambda_4; &  y_{443}&=\lambda_5.
\end{align*}
Here $\lambda_1\lambda_2\neq 0$ or $\lambda_4\neq 0$ only if $N=3$. Also, $\lambda_1\lambda_4=\lambda_2\lambda_4=0$. Similarly, $\lambda_5\neq 0$ only if $M=4$. In any case we can project onto a case of smaller rank and hence $\toba\neq0$.

For $\xymatrix@C-4pt{\overset{q}{\underset{\ }{\circ}}\ar  @{-}[r]^{q ^{-1}}  &
\overset{-1}{\underset{\ }{\circ}}
\ar  @{-}[r]^{-1}  & \overset{-1}{\underset{\ }{\circ}}
\ar  @{-}[r]^{-q}  & \overset{-q^{-1}}{\underset{\ }{\circ}}}$, the algebra $\toba$ is generated by $y_i$, $i\in\I_4$, with defining relations
\begin{align*}
 y_2^2&=\lambda_2; & y_{13}&=0; & y_{14}&=0; & y_{24}&=0; \\
y_{3}^2&=\lambda_3; & y_{23}^2&=\lambda_4; & y_{112}&=\lambda_1; &  y_{443}&=\lambda_5;
\end{align*}
\vsp
\begin{align*}
[[y_{(14)},y_2 ]_c,y_3]_c-q_{23}\frac{1+q}{1-q}[[ y_{(14)},y_3]_c,y_2]_c=0.
\end{align*}
Here, $\lambda_1\lambda_5=0$  and $\lambda_1$ or $\lambda_5$ are nonzero only if $N=4$. Thus, we can always project onto a case of smaller rank and hence $\toba\neq 0$. See \texttt{bgla.log}, resp. \texttt{bgla1.log}, \texttt{bgla2.log} for the deformation of the longest relation when one, resp. two, resp. three, of the relations $y_2^2,y_3^2,y_{23}^2$ are deformed.

For $\xymatrix@C-5pt{ \overset{q}{\underset{\ }{\circ}} \ar  @{-}[r]^{q
^{-1}}  & \overset{-1}{\underset{\ }{\circ}} \ar  @{-}[r]^{-1}
\ar@/^2pc/  @{-}[rr]^{q}
 & \overset{-1}{\underset{\ }{\circ}}
\ar  @{-}[r]^{-q^{-1}} & \overset{-1}{\circ}}$, the algebra $\toba$ is generated by $y_i$, $i\in\I_4$, with defining relations
\begin{align*}
 y_{112}&=\lambda_1; & y_{13}&=0; & y_{14}&=\lambda_5; & [y_{(13)},&y_2]_c=0; \\
 y_{2}^2&=\lambda_2; & y_{24}^2&=\lambda_6; & y_{3}^2&=\lambda_3; &  y_{4}^2&=\lambda_4;
\end{align*}
\vsp
\begin{align*}
y_{(24)}+\frac{(1+q)q_{43}}{2}[y_{24},y_3]_c-q_{23}(1+q^{-1})y_3y_{24}=0.
\end{align*}
Here $\lambda_1\neq 0$ only if $N=4$. Hence, if $\lambda_3=0$ we can project onto a case of smaller rank. If $\lambda_3\neq 0$, then $\lambda_4=\lambda_5=\lambda_6=0$ and we can project onto a case of smaller rank. Assume $\lambda_1=0$. If $\lambda_5\neq0$, then we can project. Otherwise, $\lambda_2=\lambda_4=0$ and we have that either $\lambda_6=0$ or $\lambda_3=0$.
In any case, we can project. Hence $\toba\neq0$.

\subsubsection{Type $\ufo(1)$, $\zeta\in\G_4$}\label{subsec:type-ufo(1)}
The Weyl groupoid has six objects.

For $\Dchainfive{\zeta }{\ztu }{\zeta }{\ztu }{-1}{-1}{-1}{\ztu }{\zeta }$ the algebra $\toba$ is generated by
$y_i$, $i\in\I_5$, with defining relations
\begin{align*}
y_3^2&=\lambda_3; & y_{112}&=0; & y_{221}&=0; & y_{ij}&=0, \, i<j, \, \widetilde{q}_{ij}=1;\\
y_4^2&=\lambda_4; & y_{223}&=\lambda_1; & y_{554}&=\lambda_2; &  [[[y_{(14)},&y_3]_c,y_2]_c,y_3]_c=0;
\end{align*}
\vsp
\begin{align*}
y_{34}^2&=\lambda_5; & [[y_{(25)},y_3]_c,y_4]_c-q_{34}\zeta[[y_{(25)},y_4]_c,y_3]_c=0.
\end{align*}
Here $\toba\neq 0$ as we can project onto a case of smaller rank by making $y_1=0$.

For $\Dchainfive{\zeta }{\ztu }{\zeta }{\ztu }{\zeta }{\ztu }{-1}{\ztu }{\zeta }$ the algebra $\toba$ is generated by $y_i$, $i\in\I_5$, with defining relations
\begin{align*}
 y_{112}&=0; & y_{221}&=0;  & y_{4}^2&=\lambda_1 & y_{ij}&=0, \, i<j, \, \widetilde{q}_{ij}=1;\\
y_{223}&=0; & y_{332}&=0; &  y_{334}&=\lambda_2; & y_{554}&=\lambda_3.
 \end{align*}
In this case $\toba\neq 0$ as we can project onto a case of smaller rank.

For $\xymatrix@R-8pt{  &   \overset{-1}{\circ} \ar  @{-}[d]_{\zeta} \ar  @{-}[dr]^{-1} & & \\
\overset{\zeta}{\underset{\ }{\circ}} \ar  @{-}[r]^{\ztu}  & \overset{-1}{\underset{\
}{\circ}} \ar  @{-}[r]^{\zeta}  & \overset{-1}{\underset{\ }{\circ}} \ar  @{-}[r]^{\ztu}
 & \overset{\zeta}{\underset{\ }{\circ}}}$ the algebra $\toba$ is generated by $y_i$, $i\in\I_5$, with defining relations
\begin{align*}
y_{2}^2&=\lambda_2; & y_3^2&=\lambda_3;& y_5^2&=\lambda_5; & [y_{(13)},&y_2]_c=0;  & y_{ij}&=0, \, i<j, \, \widetilde{q}_{ij}=1;\\
&& y_{112}&=\lambda_1; & y_{443}&=\lambda_4; &  [y_{125},&y_2]_c=0; &  [y_{(24)},& y_3]_c=0;
 \end{align*}
\vsp
\begin{align*}
y_{35}^2&=\lambda_6; & y_{235}+q_{35}(1+\zeta)[y_{25},y_3]_c-2q_{23}y_3y_{25}=0.
\end{align*}
If $\lambda_1=0$, then we can project. If $\lambda_1\lambda_3\neq 0$, then $\lambda_4=0$, and we can make $y_4=0$.
If $\lambda_1\lambda_4\neq 0$, then $\lambda_5=0$, and we can make $y_5=0$. Hence $\toba\neq 0$.

For $\xymatrix@R-8pt{& &  \overset{-1}{\circ} \ar  @{-}[d]_{-1} \ar  @{-}[dr]^{\zeta} &  \\
\overset{\zeta}{\underset{\ }{\circ}}  \ar  @{-}[r]^{\ztu} &
\overset{\zeta}{\underset{\ }{\circ}} \ar  @{-}[r]^{\ztu}  & \overset{-1}{\underset{\
}{\circ}} \ar  @{-}[r]^{\zeta}  & \overset{-1}{\underset{\ }{\circ}} }$ the algebra $\toba$ is generated by $y_i$, $i\in\I_5$, with defining relations
\begin{align*}
 y_{112}&=0; & y_{221}&=0; & y_{223}&=\lambda_1; & y_{ij}&=0, && i<j, \widetilde{q}_{ij}=1;\\
y_5^2&=\lambda_2; & y_{35}^2&=\lambda_5; & y_{3}^2&=\lambda_3; & y_4^2&=\lambda_4; &&  [[[y_{1235},y_3]_c,y_2]_c,y_3]_c=0;
 \end{align*}
\vsp
\begin{align*}
 [y_{(24)},&y_3]_c=0; & y_{(35)}+\frac{q_{45}(1+\zeta)}{2}[y_{35},y_4]_c-q_{34}(1-\zeta)y_4y_{35}=0.
\end{align*}
Hence $\toba\neq 0$ as we can apply Lemma \ref{lem:cut2} for $(i,j,k)=(4,3,5)$.

For $\xymatrix@R-8pt{  &   \overset{\zeta}{\circ} \ar  @{-}[d]^{\ztu} & & \\
\overset{-1}{\underset{\ }{\circ}} \ar  @{-}[r]^{\zeta}  & \overset{-1}{\underset{\
}{\circ}} \ar  @{-}[r]^{\ztu}  & \overset{\zeta}{\underset{\ }{\circ}} \ar
@{-}[r]^{\ztu}  & \overset{\zeta}{\underset{\ }{\circ}}}$ the algebra $\toba$ is generated by $y_i$, $i\in\I_5$, with defining relations
\begin{align*}
y_1^2&=\lambda_1; & y_2^2&=\lambda_2 & [y_{(13)},&y_2]_c=0; & y_{ij}&=0, &&\quad i<j, \widetilde{q}_{ij}=1;\\
 y_{332}&=\lambda_3; & y_{334}&=0; & y_{443}&=0; & y_{552}&=\lambda_4; &&  [y_{125},y_2]_c=0.
\end{align*}
We can apply Lemma \ref{lem:cut3} for $i=4$, hence $\toba\neq 0$.

For $\xymatrix@R-8pt{  &   \overset{\zeta}{\circ} \ar  @{-}[d]^{\ztu} & & \\
\overset{-1}{\underset{\ }{\circ}} \ar  @{-}[r]^{\ztu}  & \overset{\zeta}{\underset{\
}{\circ}} \ar  @{-}[r]^{\ztu}  & \overset{\zeta}{\underset{\ }{\circ}} \ar
@{-}[r]^{\ztu}  & \overset{\zeta}{\underset{\ }{\circ}}}$ the algebra $\toba$ is generated by $y_i$, $i\in\I_5$, with defining relations
\begin{align*}
 y_{221}&=\lambda_1; & y_{223}&=0; & y_{225}&=0; & y_{ij}&=0, & &  i<j, \widetilde{q}_{ij}=1;\\
y_1^2&=\lambda_2; & y_{552}&=0; & y_{332}&=0; & y_{334}&=0;&  & y_{443}=0.
\end{align*}
We can apply Lemma \ref{lem:cut3} for $i=3$, hence $\toba\neq 0$.

\subsubsection{Type $\ufo(2)$, $\zeta \in \G'_4$}\label{subsec:type-ufo(2)}
The Weyl groupoid has seven objects.

For $\xymatrix@C-4pt{\overset{\zeta}{\underset{\ }{\circ}}\ar  @{-}[r]^{\ztu}  &
\overset{\zeta}{\underset{\ }{\circ}} \ar  @{-}[r]^{\ztu}  & \overset{\zeta}{\underset{\
}{\circ}}
\ar  @{-}[r]^{\ztu}  & \overset{-1}{\underset{\ }{\circ}} \ar  @{-}[r]^{-1}  &
\overset{-1}{\underset{\ }{\circ}}  \ar  @{-}[r]^{\ztu}  & \overset{\zeta}{\underset{\
}{\circ}}}$ the algebra $\toba$ is generated by $y_i$, $i\in\I_6$, with defining relations
\begin{align*}
y_{112}&=y_{221}=0; & y_{334}&=\lambda_1; & y_4^2&=\lambda_3; & y_{ij}&=0, \, i<j, \, \widetilde{q}_{ij}=1;\\
y_{332}&=y_{223}=0; & y_{665}&=\lambda_2; & y_5^2&=\lambda_4; & [[[y_{(25)},&y_4]_c,y_3]_c,y_4]_c=0;
\end{align*}
\vsp
\begin{align*}
y_{45}^2&=\lambda_5; & [[y_{(36)},y_4]_c,y_5]_c-&q_{45}\zeta[[y_{(36)},y_5]_c,y_4]_c=0.
\end{align*}
We can apply Lemma \ref{lem:cut3} for $i=1$, hence $\toba\neq 0$.

For $\xymatrix@C-4pt{\overset{\zeta}{\underset{\ }{\circ}}\ar  @{-}[r]^{\ztu}  &
\overset{\zeta}{\underset{\ }{\circ}} \ar  @{-}[r]^{\ztu}  & \overset{\zeta}{\underset{\
}{\circ}}
\ar  @{-}[r]^{\ztu}  & \overset{\zeta}{\underset{\ }{\circ}} \ar  @{-}[r]^{\ztu}  &
\overset{-1}{\underset{\ }{\circ}}  \ar  @{-}[r]^{\ztu}  & \overset{\zeta}{\underset{\
}{\circ}}}$ the algebra $\toba$ is generated by $y_i$, $i\in\I_6$, with defining relations
\begin{align*}
y_{445}&=\lambda_1; & y_{112}&=0; & y_{221}&=0; & y_{223}&=0; & y_{ij}&=0, \, i<j, \, \widetilde{q}_{ij}=1;\\
y_{665}&=\lambda_2; & y_{332}&=0; &  y_{334}&=0; & y_{443}&=0; & y_{5}^2&=\lambda_3
\end{align*}
We can apply Lemma \ref{lem:cut3} for $i=2$, hence $\toba\neq 0$.

For $\xymatrix@C-5pt@R-8pt{  &   &  \overset{\zeta}{\circ} \ar  @{-}[d]^{\ztu} & & \\
\overset{\zeta}{\underset{\ }{\circ}} \ar  @{-}[r]^{\ztu}  & \overset{-1}{\underset{\
}{\circ}} \ar  @{-}[r]^{\zeta}  & \overset{-1}{\underset{\ }{\circ}} \ar  @{-}[r]^{\ztu}
& \overset{\zeta}{\underset{\ }{\circ}} \ar  @{-}[r]^{\ztu}  & \overset{\zeta}{\underset{\
}{\circ}}}$ the algebra $\toba$ is generated by $y_i$, $i\in\I_6$, with defining relations
\begin{align*}
y_{112}&=\lambda_1; & y_{443}&=\lambda_2; & y_2^2&=\lambda_4; & y_{ij}&=0, \, i<j, \, \widetilde{q}_{ij}=1; \quad y_{554}=0\\
y_{445}&=0; & y_{663}&=\lambda_3; & y_3^2&=\lambda_5; & [y_{(13)},&y_2]_c= [y_{(24)},y_3]_c= [y_{236},y_3]_c=0.
\end{align*}
Hence $\toba\neq 0$ since we can apply Lemma \ref{lem:cut3} for $i=4$.

For $\xymatrix@C-5pt@R-8pt{  &   & \overset{\zeta}{\circ} \ar  @{-}[d]^{\ztu}  & &\\
\overset{-1}{\underset{\ }{\circ}} \ar  @{-}[r]^{\zeta}  & \overset{-1}{\underset{\
}{\circ}} \ar  @{-}[r]^{\ztu}  & \overset{\zeta}{\underset{\ }{\circ}} \ar  @{-}[r]^{\ztu}
 & \overset{\zeta}{\underset{\ }{\circ}} \ar  @{-}[r]^{\ztu}  &
\overset{\zeta}{\underset{\ }{\circ}}}$ the algebra $\toba$ is generated by $y_i$, $i\in\I_6$, with defining relations
\begin{align*}
y_1^2&=\lambda_1; & y_{332}&=\lambda_2; & y_{336}&=y_{663}=0; & y_{ij}&=0, \, i<j, \, \widetilde{q}_{ij}=1;\\
y_2^2&=\lambda_3; & [y_{(13)},&y_2]_c=0; & y_{334}&=y_{443}=0; & y_{445}&=y_{554}=0.
\end{align*}
We can apply Lemma \ref{lem:cut3} for $i=4$, hence $\toba\neq 0$.

For $\xymatrix@C-5pt@R-8pt{  &   & & \overset{-1}{\circ} \ar  @{-}[d]^{-1} \ar
@{-}[dr]^{\zeta}  & \\
\overset{\zeta}{\underset{\ }{\circ}} \ar  @{-}[r]^{\ztu}  & \overset{\zeta}{\underset{\
}{\circ}} \ar  @{-}[r]^{\ztu}  & \overset{\zeta}{\underset{\ }{\circ}} \ar  @{-}[r]^{\ztu}
 & \overset{-1}{\underset{\ }{\circ}} \ar  @{-}[r]^{\zeta}  & \overset{-1}{\underset{\
}{\circ}}}$ the algebra $\toba$ is generated by $y_i$, $i\in\I_6$, with defining relations
\begin{align*}
y_4^2&=\lambda_1; & y_5^2&=\lambda_2; &  [y_{(35)},&y_4]_c=0; & y_{223}&=0; & &y_{ij}=0, \, i<j, \, \widetilde{q}_{ij}=1;\\
y_6^2&=\lambda_3; & y_{46}^2&=\lambda_4; & y_{332}&=0; & y_{334}&=\lambda_5; &  &[[[y_{2346},y_4]_c,y_3]_c,y_4]_c=0;
\end{align*}
\vsp
\begin{align*}
y_{112}&=y_{221}=0; & y_{(46)}+\frac{q_{56}(1+\zeta)}{2}[y_{46},y_5]_c-q_{45}(1-\zeta)y_5y_{46}=0.
\end{align*}
We can apply Lemma \ref{lem:cut3} for $i=2$, hence $\toba\neq 0$.

For $\xymatrix@C-5pt@R-8pt{  &   & \overset{\zeta}{\circ} \ar  @{-}[d]^{\ztu}  & &\\
\overset{-1}{\underset{\ }{\circ}} \ar  @{-}[r]^{\ztu}  & \overset{\zeta}{\underset{\
}{\circ}} \ar  @{-}[r]^{\ztu}  & \overset{\zeta}{\underset{\ }{\circ}} \ar  @{-}[r]^{\ztu}
 & \overset{\zeta}{\underset{\ }{\circ}} \ar  @{-}[r]^{\ztu}  &
\overset{\zeta}{\underset{\ }{\circ}}}$ the algebra $\toba$ is generated by $y_i$, $i\in\I_6$, with defining relations
\begin{align*}
y_{223}&=y_{332}=0; & y_{663}&=y_{336}=0; & y_{ij}&=0, \, i<j, \, \widetilde{q}_{ij}=1;\\
y_{334}&=y_{443}=0; & y_{445}&=y_{554}=0; & y_1^2&=\lambda_1; \quad y_{221}=\lambda_2.
\end{align*}
We can apply Lemma \ref{lem:cut3} for $i=4$, hence $\toba\neq 0$.

For $\xymatrix@C-5pt@R-8pt{  &   & \overset{-1}{\circ} \ar  @{-}[d]^{\zeta} \ar
@{-}[dr]^{-1}  & &\\
\overset{\zeta}{\underset{\ }{\circ}} \ar  @{-}[r]^{\ztu}  & \overset{\zeta}{\underset{\
}{\circ}} \ar  @{-}[r]^{\ztu}  & \overset{-1}{\underset{\ }{\circ}} \ar  @{-}[r]^{\zeta}
& \overset{-1}{\underset{\ }{\circ}}
\ar  @{-}[r]^{\ztu}  & \overset{\zeta}{\underset{\ }{\circ}}}$ the algebra $\toba$ is generated by $y_i$, $i\in\I_6$, with defining relations
\begin{align*}
y_{3}^2&=\lambda_1; & y_4^2&=\lambda_2; &  y_{223}&=\lambda_3; & [y_{(24)},&y_3]_c=0; & y_{ij}&=0, \, i<j, \, \widetilde{q}_{ij}=1;\\
y_6^2&=\lambda_4; & y_{46}^2&=\lambda_5; & y_{554}&=\lambda_6; & [y_{(35)},&y_4]_c=0; & y_{112}&=y_{221}=0;
\end{align*}
\vsp
\begin{align*}
[y_{236},&y_3]_c=0; & y_{346}+q_{46}(1+\zeta)&[y_{36},y_4]_c-2q_{34}y_4y_{36}=0.
\end{align*}
We can apply Lemma \ref{lem:cut3} for $i=2$, hence $\toba\neq 0$.

\subsubsection{Type $\ufo(3)$, $\zeta \in \G'_3$}\label{subsec:type-ufo(3)}
The Weyl groupoid has five objects.

For $\Dchainthree{-1}{\ztu}{\zeta }{-\ztu}{-\zeta }$, the algebra $\toba$ is generated by $y_i$, $i\in\I_3$, with defining relations
\begin{align*}
 y_1^2&=\lambda_1; & y_2^3&=\lambda_2; & y_{13}&=0; & y_{221}&=0; &  y_{332}&=\lambda_3,
\end{align*}
with $\lambda_1\lambda_3=\lambda_2\lambda_3=0$. If $\lambda_1=0$ then we apply Lemma \ref{lem:cut3} for $i=1$, otherwise $\lambda_3=0$ and we apply the same Lemma but for $i=3$.

For $\Dchainthree{-1}{\zeta }{-1}{-\ztu}{-\zeta }$, the algebra $\toba$ is generated by $y_i$, $i\in\I_3$, with defining relations
\begin{align*}
 y_{13}&=0; & y_{332}&=0; &  y_1^2&=\lambda_1; & y_2^2&=\lambda_2.
\end{align*}
This algebra is nonzero by \cite[Lemma 5.16]{AAGMV}.

For $\Dchainthree{\zeta }{-1}{-1}{-\ztu}{-\zeta }$, the algebra $\toba$ is generated by $y_i$, $i\in\I_3$, with defining relations
\begin{align*}
y_1^3&=\lambda_1; & y_2^2&=\lambda_2; & [y_{112},y_{12}]_c&=\lambda_3;  &  [[y_{12},&y_{(13)}]_c,y_2]_c= y_{13}=y_{332}=0.
\end{align*}
Here we can apply Lemma \ref{lem:cut3} for $i=3$.

For $\Dchainthree{\zeta }{-\ztu}{-\zeta }{-\ztu}{-\zeta }$, the algebra $\toba$ is generated by $y_i$, $i\in\I_3$, with defining relations
\begin{align*}
 y_{13}&=0; & y_{221}&=\lambda_2; &  y_1^3&=\lambda_1; & y_{223}&=0; & y_{332}&=0.
\end{align*}
This algebra is nonzero by \cite[Lemma 5.16]{AAGMV}.

For $\Dtriangle{\zeta }{-1}{-1}{\ztu}{-\zeta }{-1}$, the algebra $\toba$ is generated by $y_i$, $i\in\I_3$, with defining relations
\begin{align*}
 y_{113}&=0; & y_1^3&=\lambda_1; &  [y_{112},y_{12}]_c&=\lambda_3; & y_2^2&=\lambda_2; &  y_3^2&=\lambda_4
\end{align*}
\vsp
\begin{align*}
y_{(13)}-\frac{q_{23}\zeta}{1-\ztu}[y_{13},y_2]_c+q_{12}\ztu y_2y_{13}=0,
\end{align*}
with $\lambda_2\lambda_4=0$. Hence the algebra $\toba$ projects onto a rank-2 nonzero algebra by setting $y_2=0$ or $y_3=0$ according to whether $\lambda_4=0$ or $\lambda_2=0$.

\subsubsection{Type $\ufo(4)$, $\zeta \in \G'_3$}\label{subsec:type-ufo(4)}
The Weyl groupoid has nine objects.

For $\Dchainthree{-1}{-1}{-1}{\zeta }{-1}$, the algebra $\toba$ is generated by $y_i$, $i\in\I_3$, with defining relations
\begin{align*}
 y_1^2&=\lambda_1; & y_2^2&=\lambda_2; & y_3^2&=\lambda_3; &  y_{23}^3&=\lambda_6;\\
& & y_{12}^2&=\lambda_4; & y_{13}&=\lambda_5; & [y_{123},[y_{123},y_{23}]_c]_c&=\lambda_7
\end{align*}
with $\lambda_2\lambda_3=\lambda_7\lambda_5=0$. The deformation of $y_{23}^3$ follows from \texttt{ufo4a1.log} and of the longest relation by \texttt{ufo4a2.log}.

Assume $\lambda_7=0$. Then we can project onto a case of smaller rank unless $\lambda_1\lambda_2\lambda_4\lambda_6\neq 0$ or $\lambda_1\lambda_3\lambda_6\neq 0$, which gives
$\bq=
\left(\begin{smallmatrix}
-1& \pm 1 & \pm 1\\ \mp1 & -1 & -\zeta \\ \pm 1 &-1 &-1
\end{smallmatrix}\right)
$ or $\bq=
\left(\begin{smallmatrix}
-1& \pm 1 & \pm 1\\ \mp1 & -1 & -1 \\ \pm 1 &-\zeta &-1
\end{smallmatrix}\right)
$. In any case, $\toba\neq 0$, see \texttt{ufo4a.log}.

Now $\lambda_7\neq 0$ only if
$\bq=
\left(\begin{smallmatrix}
-1& a & a^{-1}\\ -a^{-1} & -1 & c\\ a & c^{-1}\zeta &-1
\end{smallmatrix}\right)$, with $c^3=-a^2$. In any case, $\toba\neq 0$, see \texttt{ufo4a7.log}. 

For $\Dchainthree{-1}{-1}{-1}{-\ztu}{\ztu}$, the algebra $\toba$ is generated by $y_i$, $i\in\I_3$, with defining relations
\begin{align*}
 y_1^2&=\lambda_1; & y_2^2&=\lambda_2; & y_3^3&=\lambda_3 & [y_{332},y_{32}]_c&=\lambda_4;  \\
 [y_{3321},y_{321}]_c&=\lambda_5; & y_{13}&=0; &  y_{12}^2&=\lambda_6; & [[y_{32},y_{321}]_c,y_2]_c&=0.
\end{align*}
We can project onto a case of rank two unless $\bq=
\left(\begin{smallmatrix}
-1& \pm 1 & 1 \\ \mp 1 & -1 & \zeta^{-1}\\ 1 & -1 &\zeta^{-1}
\end{smallmatrix}\right)
$, when all scalars can be nonzero. In this case, $\toba\neq0$ by \texttt{ufo4b.log}.

For $\Dchainthree{-1}{-1}{\zeta }{\ztu}{-1}$, the algebra $\toba$ is generated by $y_i$, $i\in\I_3$, with defining relations
\begin{align*}
 y_1^2&=\lambda_1; & y_2^3&=\lambda_2; & y_3^2&=\lambda_3; & [y_{221}, y_{21}]_c&=\lambda_5; \\
y_{223}&=0; & y_{13}&=\lambda_4;  & [y_{12},y_{(13)}]_c&=0,
\end{align*}
with $\lambda_1\lambda_2=0=\lambda_2\lambda_3=\lambda_2\lambda_5=\lambda_2\lambda_4$. Hence, if $\lambda_2\neq0$, then $\toba$ projects onto $\k\langle y_2| y_2^3=\lambda_2\rangle$. If $\lambda_1\neq0$, then $\lambda_5=0$ and thus $\toba$ projects onto the linking of two rank-one algebras, by making $y_2=0$.
A similar situation holds when $\lambda_3\neq 0$ and $\lambda_4=0$. On the other hand, $\lambda_3\lambda_4\neq 0$ only when
$
\bq=
\left(\begin{smallmatrix}
-1&\mp1&-1\\\pm1&\zeta&\pm 1\\\mp1&\pm\zeta^{2}&-1
\end{smallmatrix}\right)
$ and $\toba\neq0$ by \texttt{ufo4c.log}.

For $\Dchainthree{-1}{\zeta }{-1}{-\zeta }{-1}$, the algebra $\toba$ is generated by $y_i$, $i\in\I_3$, with defining relations
\begin{align*}
 y_1^2&=\lambda_1; & y_2^2&=\lambda_2; & [[y_{32},y_{321}]_c,y_2]_c&=\lambda_5; &
y_3^2&=\lambda_3; &  y_{13}&=\lambda_4.
\end{align*}
Here $\lambda_5\neq0$ only if $\bq=\left(
\begin{smallmatrix}
-1& -\texttt{c}^2\zeta & -\texttt{c}^{-3}\\ -\texttt{c}^{-2}& -1 & \texttt{c}\\ -\texttt{c}^3 & -\texttt{c}^{-1}\zeta & -1
\end{smallmatrix}
\right)$, some $\texttt{c}\in\k^\times$. Otherwise $\toba\neq 0$ by by applying Lemma \ref{lem:cut1} for $(i,j)=(2,3)$. When $\lambda_5\neq0$ is the only nonzero scalar then $\toba\neq0$ by \texttt{ufo4d.log}. Then $\toba\neq 0$ by Lemma \ref{lem:nonzero-k<l}.

For $\Dchainthree{-1}{\zeta }{-\zeta }{-\ztu}{-1}$, the algebra $\toba$ is generated by $y_i$, $i\in\I_3$, with defining relations
\begin{align*}
 y_1^2&=\lambda_1; &  y_{2221}&=\lambda_2; & y_3^2&=\lambda_3; &  y_{13}&=\lambda_4 &  y_{223}&=0.
\end{align*}
Here $\lambda_2\neq 0$ only if $\bq=\left(
\begin{smallmatrix}
-1& -\zeta & *\\ -1& -\zeta & *\\ *&* & -1
\end{smallmatrix}
\right)$. If $\lambda_2=0$, $\toba\neq 0$, by making $y_2=0$.
A similar argument shows that $\toba\neq 0$ when $\lambda_2\neq0$: $\toba$ arises as the linking of a rank-two case (generated by $y_1,y_2$) and a rank-one case.

For $\Dchainthree{-1}{\ztu}{\ztu}{-\ztu}{-1}$, the algebra $\toba$ is generated by $y_i$, $i\in\I_3$, with defining relations
\begin{align*}
 y_1^2&=\lambda_1; & y_2^3&=\lambda_2; & y_3^2&=\lambda_3; &  y_{13}&=\lambda_4; & [y_{223}, y_{23}]_c&=\lambda_5;
\end{align*}
\vsp
\begin{align*}
[y_1,y_{223}]_c+q_{23}[y_{123},y_2]_c -q_{12} y_2y_{123}=0,
\end{align*}
with $\lambda_2\lambda_4=0$. We may project onto a case of smaller rank if $\lambda_4=\lambda_1=0$ or $\lambda_4=\lambda_3=\lambda_5=0$, so $\toba\neq 0$ in these cases. The same holds if $\lambda_2=\lambda_5=0$. Now, if  $\lambda_1\lambda_4\neq0$, then $\lambda_5=0$; hence $\toba\neq0$.

If $\lambda_2\lambda_5\neq0$ or $\lambda_2\lambda_3\neq0$ and besides $\lambda_1\neq0$, then $\bq=\left(
\begin{smallmatrix}
-1& \zeta^2 & \pm1\\ 1& \zeta^2 & 1\\ \pm1 & \zeta^2 & -1
\end{smallmatrix}
\right)$.
In any case $\toba\neq0$: on the one hand this holds if $\lambda_5\neq 0$ is the unique nonzero scalar, by making $y_1=0$, and on the other when $\lambda_4=\lambda_5=0$ by \cite[Lemma 5.16]{AAGMV}; hence $\toba\neq0$ using Lemma \ref{lem:nonzero-k<l}.

For $\Dchainthree{-1}{\ztu}{\zeta }{-\zeta }{-1}$, the algebra $\toba$ is generated by $y_i$, $i\in\I_3$, with defining relations
\begin{align*}
 y_1^2&=\lambda_1; & y_2^3&=\lambda_2; & y_3^2&=\lambda_3;  \\
y_{221}&=0; & y_{13}&=\lambda_4; & [y_{223},y_{23}]_c&=\lambda_5.
\end{align*}
Here $\toba\neq 0$ by making $y_{12}=0$.

For $\Dtriangle{-1}{-1}{\zeta }{-1}{\ztu}{-\zeta }$, the algebra $\toba$ is generated by $y_i$, $i\in\I_3$, with defining relations
\begin{align*}
 y_1^2&=\lambda_1; & y_2^3&=\lambda_2; & y_3^2&=\lambda_3; &  [y_{221},y_{21}]_c&=\lambda_4;
& y_{13}^2&=\lambda_5; &  y_{223}&=0;
\end{align*}
\vsp
\begin{align*}
y_{(13)}+\frac{q_{23}(1-\zeta)}{2}[y_{13},y_2]_c-q_{12}(1-\ztu)y_2y_{13}=0.
\end{align*}
We can project onto a case of smaller rank and thus show that $\toba$ is nonzero whenever $\bq\neq\left(
\begin{smallmatrix}
-1& \zeta & \pm1\\ -1& \zeta & 1\\ \mp1 & \zeta^2 & -1
\end{smallmatrix}
\right)$. Still, $\toba\neq0$, by \texttt{ufo4h.log}.

For $\Dtriangle{-\zeta }{\zeta }{-1}{-\ztu}{\ztu}{-\ztu}$, the algebra $\toba$ is generated by $y_i$, $i\in\I_3$, with defining relations
\begin{align*}
y_{113}&=\lambda_1; & y_2^2&=\lambda_2; & y_3^3&=\lambda_3 & y_{112}&=0; &  y_{332}&=0;
\end{align*}
\vsp
\begin{align*}
y_{(13)}+q_{23}(1-\ztu)[y_{13},y_2]_c-q_{12}(1-\ztu)y_2y_{13}=0.
\end{align*}
Here $\lambda_1\neq 0$ only if $\bq=\left(
\begin{smallmatrix}
-\zeta & * & \zeta\\ *& -1 & *\\ -\zeta&* & \zeta
\end{smallmatrix}
\right)$ and $\lambda_1\lambda_2=0$. Hence $\toba\neq0$, by projecting onto a rank 2 case.

\subsubsection{\ Type $\ufo(5)$, $\zeta \in \G'_3$}\label{subsec:type-ufo(5)}
The Weyl groupoid has six objects.

For $\Dchainfour{-\zeta }{-\ztu}{-\zeta }{-\ztu}{-\zeta }{-\ztu}{\zeta }$, the algebra $\toba$ is generated by $y_i$, $i\in\I_4$, with defining relations
\begin{align*}
 y_{112}&=0; & y_{221}&=0; & y_{4}^3&=\lambda_1;\\
 y_{223}&=0; &  y_{332}&=0; & y_{334}&=\lambda_2; & y_{ij}&=0, \, i<j, \, \widetilde{q}_{ij}=1.
\end{align*}
We can apply Lemma \ref{lem:cut3} for $i=2$, hence $\toba\neq 0$.

For $\Dchainfour{-\zeta }{-\ztu}{-\zeta }{-\ztu}{-1}{-1}{\zeta }$, the algebra $\toba$ is generated by $y_i$, $i\in\I_4$, with defining relations
\begin{align*}
 y_{112}&=0; & y_{221}&=0; & [y_{443},& y_{43}]_c=0;  & y_{ij}&=0, \, i<j, \, \widetilde{q}_{ij}=1;\\
y_{223}&=0; &  y_3^2&=\lambda_1; & y_{4}^3&=\lambda_2; & [[y_{43},&y_{432}]_c,y_3]_c=0.
\end{align*}
We can apply Lemma \ref{lem:cut3} for $i=1$, hence $\toba\neq 0$.

For $\Dchainfour{-\zeta }{-\ztu}{\zeta }{\ztu}{-1}{-\ztu}{-\zeta }$, the algebra $\toba$ is generated by $y_i$, $i\in\I_4$, with defining relations
\begin{align*}
 y_{112}&=\lambda_1; & y_{443}&=0; & y_{ij}&=0, \, i<j, \, \widetilde{q}_{ij}=1; &
y_2^3&=\lambda_2; & y_{3}^2&=\lambda_3.
\end{align*}
We can apply Lemma \ref{lem:cut3} for $i=4$, hence $\toba\neq 0$.

For $\Dchainfour{-\zeta }{-\ztu}{-\zeta }{-\ztu}{-1}{-\ztu}{-\zeta }$, the algebra $\toba$ is generated by $y_i$, $i\in\I_4$, with defining relations
\begin{align*}
 y_{112}&=y_{221}=y_{223}=y_{443}=0; & y_{ij}&=0, \, i<j, \, \widetilde{q}_{ij}=1; &
y_{3}^2&=\lambda_1.
\end{align*}
We can apply Lemma \ref{lem:cut3} for $i=4$, hence $\toba\neq 0$.

For $\Drightofway{-\zeta }{-\ztu}{-1}{-\zeta }{\ztu}{-1}{-1}{\zeta }$, the algebra $\toba$ is generated by $y_i$, $i\in\I_4$, with defining relations
\begin{align*}
[y_{334},&y_{34}]_c=\lambda_1; & y_{112}&=0; & y_{2}^2&=\lambda_2; & y_{ij}&=0, \, i<j, \, \widetilde{q}_{ij}=1;\\
[y_{124},& y_2]_c=0; & y_{332}&=0; & y_3^3&=\lambda_3; & y_{4}^2&=\lambda_4;
\end{align*}
\vsp
\begin{align*}
y_{(24)}-2q_{34}\zeta[y_{24},y_3]_c-2q_{23}y_3y_{24}=0.
\end{align*}
We can apply Lemma \ref{lem:cut3} for $i=1$, hence $\toba\neq 0$.

For $\Drightofway{-\zeta }{-\ztu}{-1}{-\zeta }{\zeta }{-1}{-\zeta }{-1}$, the algebra $\toba$ is generated by $y_i$, $i\in\I_4$, with defining relations
\begin{align*}
y_{2}^2&=\lambda_1; & y_{3}^2&=\lambda_2;  & y_4^2&=\lambda_3; & y_{ij}&=0, \, i<j, \, \widetilde{q}_{ij}=1;
\end{align*}
\vsp
\begin{align*}
y_{112}&=0; &  [y_{124},&y_2]_c=0; & y_{(24)}+q_{34}\ztu [y_{24},y_3]_c+q_{23}\ztu y_3y_{24}=0.
\end{align*}
We can apply Lemma \ref{lem:cut3} for $i=1$, hence $\toba\neq 0$.

\subsubsection{\ Type $\ufo(6)$, $\zeta \in \G'_4$}\label{subsec:type-ufo(6)}
The Weyl groupoid has eight objects.

For $\Dchainfour{\ztu }{\zeta }{-1}{\ztu }{\zeta }{\zeta }{\ztu }$ the algebra $\toba$ is generated by $y_i$, $i\in\I_4$, with defining relations
\begin{align*}
 y_{112}&=\lambda_1;  & y_{2}^2&=\lambda_2;  &   y_3^4&=\lambda_7 & y_{332}&=\lambda_3;   \\
  y_{13}&=\lambda_5 &   y_{443}&=0; & y_{14}&=y_{24}=0;  & [y_{(13)},& y_2]_c=\lambda_4;
\end{align*}
\vsp
 \begin{align*}
    [[[y_{(24)}, y_3]_c,y_3]_c,y_3]_c=\lambda_6.
\end{align*}
Here $\lambda_1\lambda_3=\lambda_1\lambda_4=\lambda_1\lambda_5=0=\lambda_1\lambda_6$. 
If $\lambda_{1}\neq0$ then we apply Lemma \ref{lem:cut3} for $i=4$. Hence we assume $\lambda_1=0$. If $\lambda_4=0$, then we apply Lemma \ref{lem:cut1} for $(i,j)=(1,2)$. If $\lambda_4$ is the unique non-zero scalar, then we apply Lemma \ref{lem:cut3} for $i=4$. The general case follows by Lemma \ref{lem:nonzero-k<l}. See \texttt{ufo6a.log} for the deformation of the last relation.

For $\Dchainfour{-1}{\ztu }{-1}{\zeta }{-1}{\zeta }{\ztu }$ the algebra $\toba$ is generated by $y_i$, $i\in\I_4$, with defining relations
\begin{align*}
y_1^2&=\lambda_1; &   y_2^2&=\lambda_2; & y_3^2&=\lambda_3; & [y_{(13)},&y_2]_c=\lambda_4; \\
y_{443}&=\lambda_5; & y_{13}&=\lambda_6 & y_{14}&=y_{24}=0; & [[y_{23},&[y_{23},y_{(24)}]_c]_c,y_3]_c=0.
 \end{align*}
If $\lambda_5=0$, then $\toba\neq0$ by making $y_4=0$. Similarly if $\lambda_5$ is the unique nonzero scalar. The general case now follows from Lemma \ref{lem:nonzero-k<l}.

For $\Dchainfour{-1}{\zeta }{\ztu }{\zeta }{-1}{\zeta }{\ztu }$
the algebra $\toba$ is generated by $y_i$, $i\in\I_4$, with defining relations
\begin{align*}
y_1^2&=\lambda_1;  & y_{221}&=\lambda_2; & y_3^2&=\lambda_3; & y_{223}&=\lambda_4; \\
 y_{443}&=\lambda_5; & y_{13}&=\lambda_6;  &  y_{14}&=y_{24}=0.
 \end{align*}
If $\lambda_5\lambda_4\neq 0$, then $\lambda_1=\lambda_2=\lambda_6=0$ and hence we can apply Lemma \ref{lem:cut3} for $i=1$. If $\lambda_5\neq 0$ and $\lambda_4=0$, then $\toba\neq 0$ by Lemma \ref{lem:cut1} for $(i,j)=(2,3)$. If $\lambda_5=0$, then $\toba\neq 0$ by Lemma \ref{lem:cut3} for $i=4$.

For $\Dchainfour{-1}{\ztu }{\zeta }{-1}{-1}{\zeta }{\ztu }$ the algebra $\toba$ is generated by $y_i$, $i\in\I_4$, with defining relations
\begin{align*}
 y_1^2&=\lambda_1; & y_{221}&=\lambda_2 & y_3^2&=\lambda_3; & y_{13}&=\lambda_4; \\
 y_{443}&=\lambda_5; &    y_{2223}&=0; & y_{14}&=0 & y_{24}&=0;
 \end{align*}
 \vsp
\begin{align*}
[y_2, & [y_{(24)},y_3]_c]_c = \frac{q_{23}q_{43}}{1+\zeta}[y_{23},y_{(24)}]_c+(\zeta-1)q_{23}q_{24} y_{(24)}y_{23}.
\end{align*}
If $\lambda_5=0$, then $\toba\neq0$ by by Lemma \ref{lem:cut3} for $i=4$. Similarly if $\lambda_5$ is the unique nonzero scalar. The general case now follows from Lemma \ref{lem:nonzero-k<l}.

For $\Dchainfour{-1}{\zeta }{-1}{-1}{-1}{\zeta }{\ztu }$ the algebra $\toba$ is generated by $y_i$, $i\in\I_4$, with defining relations
\begin{align*}
y_1^2&=\lambda_1; & y_2^2&=\lambda_2;  &  y_3^2&=\lambda_3; &  y_{443}&=\lambda_4; \\
y_{23}^2&=\lambda_5; & y_{13}&=\lambda_6 & y_{14}&=y_{24}=0; &  [[y_{12},&y_{(13)}]_c, y_2]_c=\lambda_7.
\end{align*}
If $\lambda_4=0$, then $\toba\neq0$ by by Lemma \ref{lem:cut3} for $i=4$. Similarly if $\lambda_4$ is the unique nonzero scalar. The general case now follows from Lemma \ref{lem:nonzero-k<l}.

For $\Drightofway{\ztu }{\zeta }{-1}{\ztu }{-1}{\zeta }{\ztu }{-1}$ the algebra $\toba$ is generated by $y_i$, $i\in\I_4$, with defining relations
\begin{align*}
 y_{112}&=\lambda_1; & y_2^2&=\lambda_2; & y_3^2&=\lambda_3; & y_{443}&=\lambda_4; & y_{442}&=\lambda_5; \\
y_{23}^2&=\lambda_6;   &  y_{14}&=\lambda_7; & [y_{124},&y_2]_c=\lambda_8; & y_{13}&=0;
\end{align*}
\vsp
\begin{align*}
y_{(24)}-q_{34}\zeta [y_{24},y_3]_c-q_{23}(1+\zeta)y_3y_{24}=0.
\end{align*}
We analyze the cases in which we cannot project onto a case of smaller rank, notice that, for instance, $\lambda_1\lambda_5\neq 0$ or $\lambda_1\lambda_5\neq 0$, imply that $\lambda_3=\lambda_4=\lambda_6=0$ and we can make $y_3=0$. Similarly if $\lambda_7\neq 0$ or $\lambda_8\neq 0$.

Finally, if $\lambda_1\lambda_4\neq 0$, then
$\bq=\left(\begin{smallmatrix}
\ztu & -1 & 1 & \pm1\\
\ztu & -1 & -1 & \ztu \\
1 & 1 & -1 & \zeta \\
\pm1  & 1 & -1 & \zeta \\
\end{smallmatrix}\right)$ and $\lambda_5=\lambda_7=\lambda_8=0$; we could also have $\lambda_2\lambda_3\lambda_6\neq 0$. In any case, $\toba\neq 0$ by \texttt{ufo6f.g}.

For $\Drightofway{-1}{\ztu }{\zeta }{\ztu }{-1}{-1}{\ztu }{-1}$ the algebra $\toba$ is generated by $y_i$, $i\in\I_4$, with defining relations
\begin{align*}
 y_1^2&=\lambda_1; & y_{13}&=\lambda_2 & y_3^2&=\lambda_3; & y_4^2&=\lambda_4; \\
 y_{221}&=\lambda_5; & y_{224}&=\lambda_6 & y_{14}&=\lambda_7; & y_{2223}&=0;
\end{align*}
\vsp
\begin{align*}
y_{(24)}-q_{34}\zeta [y_{24},y_3]_c-q_{23}(1+\zeta)y_3y_{24}=0.
\end{align*}
Here $\toba\neq0$ as we can apply Lemma \ref{lem:cut2} for $(i,j,k)=(3,2,4)$.

For $\Drightofway{-1}{\zeta }{-1}{\ztu }{-1}{-1}{\ztu }{-1}$ the algebra $\toba$ is generated by $y_i$, $i\in\I_4$, with defining relations
\begin{align*}
 y_1^2&=\lambda_1; & y_2^2&=\lambda_2; & y_3^2&=\lambda_3; &  y_4^2&=\lambda_4; & [y_{124},& y_2]_c=\lambda_5;  \\
& &  y_{14}&=\lambda_6; & y_{23}^2&=\lambda_7;  &  y_{13}&=\lambda_8;  & [[y_{12},&y_{(13)}]_c,y_2]_c=\lambda_9;
\end{align*}
\vsp
\begin{align*}
y_{(24)}-q_{34}\zeta [y_{24},y_3]_c-q_{23}(1+\zeta)y_3y_{24}=0.
\end{align*}
If $\lambda_5=0$, then $\toba\neq 0$ since we can we can apply Lemma \ref{lem:cut2} for $(i,j,k)=(4,2,3)$. Similarly if $\lambda_4$ is the unique nonzero scalar. The general case now follows from Lemma \ref{lem:nonzero-k<l}.

\subsubsection{Type $\ufo(7)$, $\zeta \in \G'_{12}$}\label{subsec:type-ufo(7)}
The Weyl groupoid has five objects but two of them are obtained with $\zeta^5$
instead of $\zeta$. Hence we focus in the remaining three cases.

For $\Dchaintwo{-\overline{\zeta}^{\,2}}{-\zeta ^3}{-\zeta ^2}$,
$\toba$ is generated by $y_1, y_2$ with defining relations
\begin{align*}
y_1^3&=\lambda_1; & y_2^3&=\lambda_2; &  [y_1,y_{122}]_c -\frac{\zeta^{10}(1+\zeta)q_{12}}{1+\zeta^3}y_{12}^2&=0,
\end{align*}
with $\lambda_1\lambda_2=0$. This is nonzero by \cite[Lemma 5.16]{AAGMV}.

For $\Dchaintwo{-\overline{\zeta}^{\,2}}{\overline{\zeta}}{-1}$,
$\toba$ is generated by $y_1, y_2$ with defining relations
\begin{align*}
y_1^3&=\lambda_1; & y_2^2&=\lambda_2; & [[y_{112},y_{12}]_c,y_{12}]_c&=0,
\end{align*}
with $\lambda_1\lambda_2=0$. The last relation is undeformed in any case, see
\texttt{ufo(7).log}. This algebra is nonzero by \cite[Lemma 5.16]{AAGMV}.

For $\Dchaintwo{-\zeta ^3}{\zeta}{-1}$, $\toba$ is generated by $y_1, y_2$ with defining relations
\begin{align*}
y_1^4&=\lambda_1; & y_2^2&=\lambda_2; &  [ y_{112},y_{12}]_c&=0,
\end{align*}
with $\lambda_1\lambda_2=0$. This is nonzero by \cite[Lemma 5.16]{AAGMV}.

\subsubsection{Type $\ufo(8)$, $\zeta \in \G'_{12}$}\label{subsec:type-ufo(8)}
The Weyl groupoid has three objects.

For $\Dchaintwo{-\zeta ^2}{\zeta }{-\zeta ^2}$, $\toba$ is generated by $y_1, y_2$ with defining relations
\begin{align*}
y_1^3&=\lambda_1; & y_2^3&=\lambda_2; & [y_1,y_{122}]_c-(1+\zeta+\zeta^2)\zeta^4q_{12}y_{12}^2&=0,
\end{align*}
with $\lambda_1\lambda_2=0$. This algebra is nonzero by \cite[Lemma 5.16]{AAGMV}.

For $\Dchaintwo{-\zeta ^2}{\zeta ^3}{-1}$, $\toba$ is generated by $y_1, y_2$ with defining relations
\begin{align*}
y_1^3&=\lambda_1; & y_2^2&=\lambda_2; & [[y_{112},y_{12}]_c,y_{12}]_c&=0.
\end{align*}
with $\lambda_1\lambda_2=0$. The last relation is undeformed in any case, see
\texttt{ufo8.log}. This algebra is nonzero by \cite[Lemma 5.16]{AAGMV}.

For $\Dchaintwo{-\ztu}{-\zeta ^3}{-1}$, $\toba$ is generated by $y_1, y_2$ with defining relations
\begin{align*}
y_{11112}&=0; & y_2^2&=\lambda_2; &   [y_{112},y_{12}]_c&=0.
\end{align*}
This algebra is nonzero by \cite[Lemma 5.16]{AAGMV}.

\subsubsection{Type $\ufo(9)$, $\zeta \in \G'_{24}$}\label{subsec:type-ufo(9)}
The Weyl groupoid has four objects.

For $\Dchaintwo{\zeta ^6}{-\ztu}{\ \ -\ztu^{\, 4}}$, $\toba$ is generated by $y_1, y_2$ with defining relations
\begin{align*}
y_1^4&=\lambda_1; & y_2^3&=\lambda_2; & [y_1,y_{122}]_c - \frac{1+\zeta^7}{1+\zeta}\zeta^{10}q_{12}y_{12}^2&=0,
\end{align*}
with $\lambda_1\lambda_2=0$. This algebra is nonzero by \cite[Lemma 5.16]{AAGMV}.

For $\Dchaintwo{\zeta ^6}{\zeta }{\ztu}$, $\toba$ is generated by $y_1, y_2$ with defining relations
\begin{align*}
y_1^4&=\lambda_1; \qquad  y_{221}=0; \\
[[y_{112}&,y_{12}]_c,y_{12}]_c= -\lambda_1q_{21}^{-1}(2\zeta+\zeta^8+2\zeta^{17}-\zeta^{22}) y_2^3 .
\end{align*}
See \texttt{ufo9b.log} and \texttt{ufo9b2.log}. This algebra is nonzero, as it projects onto $\k\langle y_1|y_1^4=\lambda_1 \rangle$.

For $\Dchaintwo{-\ztu^{\, 4}}{\zeta ^5}{-1}$, $\toba$ is generated by $y_1, y_2$ with defining relations
\begin{align*}
y_1^3&=\lambda_1; & y_2^2&=\lambda_2;  & [y_{112},[[y_{112},y_{12}]_c,y_{12}]_c]_c -\alpha [y_{112},y_{12}]_c^2 &=0,
\end{align*}
for $\alpha=\frac{1+\zeta+\zeta^6+2\zeta^7+\zeta^{17}}{1+\zeta^4+\zeta^6+\zeta^{11}} \zeta^{9}q_{12}$. Here, $\lambda_1\lambda_2=0$.
The last relation remains undeformed, see
\texttt{ufo9c.log}. This algebra is nonzero by \cite[Lemma 5.16]{AAGMV}.

For $\Dchaintwo{\zeta }{\ztu^{\, 5}}{-1}$, $\toba$ is generated by $y_1, y_2$ with defining relations
\begin{align*}
y_2^2&=\lambda_1; & y_{1111112}&=0; & [y_{112},y_{12}]_c&=0
\end{align*}
This algebra is nonzero by \cite[Lemma 5.16]{AAGMV}.
The last relation remains undeformed by direct computation.

\subsubsection{Type $\ufo(10)$, $\zeta\in \G'_{20}$}\label{subsec:type-ufo(10)}
The Weyl groupoid has four objects but two of them are obtained with $-\zeta$
instead of $\zeta$. Hence we focus in just two cases.

For $\Dchaintwo{\zeta }{\ztu^{\, 3}}{-1}$, $\toba$ is generated by $y_1, y_2$ with defining relations
\begin{align*}
y_2^2&=\lambda_1;  &  y_{11112}&=0;  & [[[y_{112},y_{12}]_c,y_{12}]_c,y_{12}]_c&=0.
\end{align*}
This algebra is nonzero by \cite[Lemma 5.16]{AAGMV}.  The last relation remains undeformed by
\texttt{ufo10a.log}.

For $\Dchaintwo{-\ztu^{\, 2}}{\zeta ^3}{-1}$, $\toba$ is generated by $y_1, y_2$ with defining relations
\begin{align*}
y_1^5&=\lambda_1; & y_2^2&=\lambda_2 &
[y_1,[y_{112},y_{12}]_c]_c+\frac{1-\zeta^{17}}{1-\zeta^2}q_{12}y_{112}^2&=0,
\end{align*}
with $\lambda_1\lambda_2=0$. The last relation remains undeformed, see  \texttt{ufo10c.log}, and thus this algebra is nonzero by \cite[Lemma 5.16]{AAGMV}.

\subsubsection{Type $\ufo(11)$, $\zeta\in \G'_{15}$}\label{subsec:type-ufo(11)}
The Weyl groupoid has four objects.

For $\Dchaintwo{-\zeta }{-\ztu^{\, 3}}{\zeta ^5}$, $\toba$ is generated by $y_1, y_2$ with defining relations
\begin{align*}
y_2^3&=\lambda_2; & y_{11112}&=0; & [[y_{112},y_{12}]_c,y_{12}]_c &=0;
&
[y_1,y_{122}]_c +\alpha y_{12}^2&=0,
\end{align*}
for $\alpha=\frac{1+\zeta^{13}}{1+\zeta^{12}} \zeta^{10}q_{12}$. This algebra is nonzero by \cite[Lemma 5.16]{AAGMV}.
 The third relation  is undeformed by \texttt{ufo11a1.log}.
  We recall that the fourth one is primitive in $T(V)$.

For $\Dchaintwo{\zeta ^3}{-\zeta ^4}{-\ztu^{\,4}}$, $\toba$ is generated by $y_1, y_2$ with defining relations
\begin{align*}
[y_1,[y_{112},y_{12}]_c]_c &=\frac{1-\zeta^2}{1+\zeta^7}\zeta^9q_{12} y_{112}^2; \qquad y_1^5=\lambda_1; \qquad  y_{221}=0; \\
[[[y_{112},y_{12}]_c,y_{12}]_c,y_{12}]_c&= \lambda_1q_{21}(\zeta+\zeta^4+3\zeta^8+4\zeta^{13}-3\zeta^{14})y_2^4,
\end{align*}
See \texttt{ufo11b.log} for the rule of deformation of the bottom quantum Serre relation.
This algebra is nonzero as it projects onto $\k\langle y_1| y_1^5=\lambda_1 \rangle$.

For $\Dchaintwo{\zeta ^5}{-\ztu^{\, 2}}{-1}$, $\toba$ is generated by $y_1, y_2$ with defining relations
\begin{align*}
 y_1^3&=\lambda_1; & y_2^2&=\lambda_2; & [[[y_{112},y_{12}]_c,y_{12}]_c,y_{12}]_c& =0,
\end{align*}
with $\lambda_1\lambda_2=0$. The third relation is undeformed by \texttt{ufo11c.log} and thus
 this algebra is nonzero by  \cite[Lemma 5.16]{AAGMV}.

For $\Dchaintwo{\zeta ^3}{-\zeta ^2}{-1}$, $\toba$ is generated by $y_1, y_2$ with defining relations
\begin{align*}
 y_1^5&=\lambda_1; & y_2^2&=\lambda_2; &  [y_{112},y_{12}]_c&=0; &  [[y_{1112},y_{112}]_c,y_{112}]_c&=0,
\end{align*}
$\lambda_1\lambda_2=0$. The third relation is undeformed by direct computation.
As for the fourth, see  \texttt{ufo11d.log}.
Then this algebra is nonzero \cite[Lemma 5.16]{AAGMV}.

\subsubsection{Type $\ufo(12)$, $\zeta\in \G'_{7}$}\label{subsec:type-ufo(12)}
The Weyl groupoid has two objects.

For $\Dchaintwo{-\zeta }{-\ztu^{\, 3}}{-1}$, $\toba$ is generated by $y_1, y_2$ with defining relations
\begin{align*}
 y_{11112}&=0; & y_2^2&=\lambda_2; &[y_{112},[[y_{112},y_{12}]_c,y_{12}]_c]_c
 -\alpha[y_{112},y_{12}]_c^2& =0,
\end{align*}
for $\alpha=-q_{12}\frac{-\zeta+3\zeta^2+3\zeta^3-\zeta^4+3\zeta^6}{-2\zeta+2\zeta^3-\zeta^5+\zeta^6}$.
By \texttt{ufo12a.log}, the last relation is undeformed.
The algebra is nonzero as it projects onto $\k\langle y_2^2|y_2^2=\lambda_2\rangle$.

For $\Dchaintwo{-\ztu^{\, 2}}{-\zeta ^3}{-1}$, $\toba$ is generated by $y_1, y_2$ with defining relations
\begin{align*}
 y_{1111112}&=0; & y_2^2&=\lambda_2; &
[y_1,[y_{112},y_{12}]_c]_c +q_{12}\frac{1+\zeta^4}{1-\zeta^2}y_{112}^2 &=0.
\end{align*}
The last relation is undeformed by direct computation. The algebra is nonzero as it projects onto $\k\langle y_2^2|y_2^2=\lambda_2\rangle$.


\begin{thebibliography}{AAG}
\bibitem[AA1]{AA}  {\sc N. Andruskiewitsch}, {\sc I. Angiono}.
\emph{Finite dimensional Nichols algebras of diagonal type}. To appear in Bull.~Math.~Sci., \texttt{arXiv:1707.08387}.

\bibitem[AA2]{AA-new}  \bysame.
\emph{On Nichols algebras over basic Hopf algebras}, \texttt{arXiv:1802.00316}.


\bibitem[A+]{AAGMV} {\sc N. Andruskiewitsch}, {\sc I. Angiono}, {\sc A. Garc\'ia Iglesias}, {\sc A. Masuoka},
{\sc C. Vay}. \emph{Lifting via cocycle
deformation}. J. Pure Appl. Alg. {\bf 218} (4), 684--703 (2014).

\bibitem[AAG]{AAG}  {\sc N. Andruskiewitsch}, {\sc I. Angiono}, {\sc A. Garc\'ia Iglesias}.
\emph{ Liftings of Nichols algebras of diagonal type I. Cartan type A}. Int. Math. Res. Not. IMRN {\bf 2017} (9), 2793--2884 (2017).

\bibitem[AG]{AG}  {\sc N. Andruskiewitsch}, {\sc A. Garc\'ia Iglesias}.
\emph{Twisting Hopf algebras from cocycle deformations}. Ann. Univ. Ferrara {\bf 63} (2), 221--247 (2017).

\bibitem[AS1]{AS-lift-meth}  {\sc N. Andruskiewitsch}, {\sc H.-J. Schneider}.
\emph{Finite quantum groups over abelian groups of prime exponent}, Ann. Sci. Ec. Norm. Super. {\bf 35}, 1--26, (2002).

\bibitem[AS2]{AS}  \bysame
\emph{Isomorphism classes and automorphisms
of finite Hopf algebras of type $A_n$}. Proceedings of the XVIth Latin American Algebra Colloquium (Spanish), 201--226, Bibl. Rev. Mat. Iberoamericana, Rev. Mat. Iberoamericana, Madrid, (2007).


\bibitem[AS3]{AS-annals} \bysame
\emph{On the classification of finite-dimensional pointed Hopf
algebras}. Ann. Math. \textbf{171}, 375--417, (2010).

\bibitem[AV]{AV}  {\sc N. Andruskiewitsch}, {\sc C. Vay}.
\emph{Finite dimensional
Hopf algebras over the dual group algebra of the symmetric group in three letters},
Comm. Alg. {\bf 39}, 4507--4517, (2011).

\bibitem[A1]{A-nichols} {\sc I. Angiono}.
{\em On Nichols algebras of diagonal type}. J. Reine Angew. Math.  \textbf{683} (2013), 189--251.

\bibitem[A2]{A-distinguished} \bysame \emph{Distinguished Pre-Nichols algebras}, Transf. Groups {\bf 21}, 1--33, (2016).

\bibitem[AG]{AGI} {\sc I. Angiono}, {\sc A. Garc\'ia Iglesias}. 
\emph{Pointed Hopf algebras with standard braiding are generated in degreeone}, Contemp. Math. {\bf 537} (2011).


\bibitem[B]{B - diamond} {\sc G. Bergman}.
\emph{The diamond lemma for ring theory}. Adv. Math. \textbf{29} (1978), 178--218.

\bibitem[CG]{CG}  {\sc A. M. Cohen}, {\sc D. A. H. Gijsbers}. \emph{GBNP
0.9.5 (Non-commutative Gr\"obner bases)}, \texttt{http://www.win.tue.nl/\~\,amc}.

\bibitem[JG]{JG} {\sc J.M. Jury Giraldi}, {\sc A. Garc\'ia Iglesias}.
 \emph{Liftings of Nichols algebras of diagonal type III. Cartan type $G_2$}, J. Algebra {\bf 478}, 506--568 (2017).

\bibitem[GAP]{GAP} {\sc The GAP Group}, \emph{GAP --- Groups, Algorithms and
Programming}. Version \textbf{4.4.12}, (2008), \texttt{http://www.gap-system.org}.


\bibitem[Gu]{Gu} {\sc G\"unther,R.}, {\it Crossed products for pointed Hopf algebras}.
Comm.  Algebra, \textbf{27}, 4389--4410,  (1999).

\bibitem[H]{H-classif} {\sc I. Heckenberger}. {\em Classification of arithmetic root systems}. Adv. Math. \textbf{220} (2009), 59--124.

\bibitem[Mo]{Mo} {\sc S. Montgomery}. \emph{Hopf algebras and their action on rings}.
CBMS Lecture Notes 82, American Math Society, Providence, RI, (1993).



\bibitem[S]{S} {\sc Schauenburg, P.}, {\it Hopf bi-Galois extensions}. Comm. Alg.
\textbf{24}, 3797--3825, (1996).
\end{thebibliography}
\end{document}